\documentclass{memo-l}

\usepackage{amsmath,amssymb,mathrsfs,tikz,multirow,hyperref}

\newtheorem{theorem}{Theorem}[chapter]
\newtheorem{lemma}[theorem]{Lemma}
\newtheorem{conjecture}[theorem]{Conjecture}
\newtheorem{proposition}[theorem]{Proposition}
\newtheorem{corollary}[theorem]{Corollary}

\theoremstyle{definition}
\newtheorem{definition}[theorem]{Definition}

\theoremstyle{remark}

\numberwithin{section}{chapter}
\numberwithin{equation}{chapter}
\numberwithin{table}{chapter}
\newcommand{\PSL}{\mathrm{PSL}}
\newcommand{\PSU}{\mathrm{PSU}}
\newcommand{\SU}{\mathrm{SU}}
\newcommand{\PGL}{\mathrm{PGL}}
\newcommand{\GL}{\mathrm{GL}}
\newcommand{\SL}{\mathrm{SL}}
\newcommand{\PSp}{\mathrm{PSp}}
\newcommand{\Sp}{\mathrm{Sp}}
\newcommand{\Spin}{\mathrm{Spin}}
\newcommand{\POmega}{\mathrm{P}\Omega}
\newcommand{\cf}{\mathrm{cf}}
\newcommand{\Irr}{\mathrm{Irr}}

\newcommand{\slf}{\mathfrak{sl}}
\newcommand{\topp}{\mathrm{top}}
\newcommand{\Ann}{\mathrm{Ann}}
\newcommand{\Aut}{\mathrm{Aut}}
\newcommand{\F}{\mathbb{F}}
\newcommand{\Out}{\mathrm{Out}}
\newcommand{\gen}[1]{\langle#1\rangle}
\newcommand{\Hom}{\mathrm{Hom}}
\newcommand{\soc}{\mathrm{soc}}
\newcommand{\rad}{\mathrm{rad}}
\newcommand{\Ext}{\mathrm{Ext}}
\newcommand{\Alt}{\mathrm{Alt}}
\newcommand{\Sym}{\mathrm{Sym}}
\newcommand{\N}{\mathbb{N}}
\newcommand{\Q}{\mathbb{Q}}
\newcommand{\I}{\mathrm{i}}

\newcommand{\bG}{\mathbf{G}}
\newcommand{\bH}{\mathbf{H}}
\newcommand{\bX}{\mathbf{X}}
\newcommand{\bY}{\mathbf{Y}}
\newcommand{\bT}{\mathbf{T}}
\newcommand{\bL}{\mathbf{L}}
\newcommand{\bP}{\mathbf{P}}
\newcommand{\bC}{\mathbf{C}}

\begin{document}

\frontmatter

\title{Maximal \texorpdfstring{$\mathrm{PSL}_2$}{PSL2} Subgroups of Exceptional Groups of Lie Type}

\author{David A. Craven\\Mem. Amer. Math. Soc., to appear}
\address{School of Mathematics, University of Birmingham, Birmingham, B15 2TT, United Kingdom}
\email{d.a.craven@bham.ac.uk}
\thanks{The author is a Royal Society University Research Fellow, and gratefully acknowledges the financial support of the Society.}

%    \date is required; it is the date received by the editor.
\date{December 14, 2017}

%    Recognition of the 2010 edition of the Mathematics Subject
%    Classification requires a version of amsbook.cls from July 2009
%    or later.  If "2010" is not recognized, please upgrade.


%\dedicatory{Dedication text (use \\[2pt] for line break if necessary)}

\begin{abstract}
We study embeddings of $\mathrm{PSL}_2(p^a)$ into exceptional groups $G(p^b)$ for $G=F_4,E_6,{}^2\!E_6,E_7$, and $p$ a prime with $a,b$ positive integers. With a few possible exceptions, we prove that any almost simple group with socle $\mathrm{PSL}_2(p^a)$, that is maximal inside an almost simple exceptional group of Lie type $F_4$, $E_6$, ${}^2\!E_6$ and $E_7$, is the fixed points under the Frobenius map of a corresponding maximal closed subgroup of type $A_1$ inside the algebraic group.

Together with a recent result of Burness and Testerman for $p$ the Coxeter number plus one, this proves that all maximal subgroups with socle $\mathrm{PSL}_2(p^a)$ inside these finite almost simple groups are known, with three possible exceptions ($p^a=7,8,25$ for $E_7$).

In the three remaining cases we provide considerable information about a potential maximal subgroup.
\end{abstract}

\maketitle

\tableofcontents

%    Include unnumbered chapters (preface, acknowledgments, etc.) here.
%\include{}

\mainmatter
%    Include main chapters here.
\chapter{Introduction}

Classifying the maximal subgroups of a finite group is one of the most fundamental problems in the field of finite group theory. Aschbacher and Scott \cite{aschbacherscott1985} reduced the problem for all finite groups to understanding the $1$-cohomology groups $H^1(G,M)$ for all simple modules $M$ for all finite almost simple groups $G$, and classifying all maximal subgroups of almost simple groups.

This paper is a contribution towards the latter, ambitious goal. For alternating and classical groups there is in some sense no complete answer, since the dimensions of the classical groups (and degrees of the alternating groups) tend to infinity, although there is substantial work in this direction. However, for sporadic and exceptional groups there is a possibility of a complete answer being known.

For sporadic groups, a complete answer is known for all groups but the Monster, and here we concentrate on exceptional groups of Lie type. There is a classification of maximal subgroups for exceptional groups $G=G(q)$ for $G$ not of type $F_4$, $E_6$, ${}^2\!E_6$, $E_7$ and $E_8$ already (see \cite{wilsonrob} for example), and so we focus on the remaining cases. What is known in the literature so far is summarized in Chapter \ref{ch:maxsubgroups}, but broadly speaking, all maximal subgroups are known for these groups apart possibly from various almost simple maximal subgroups, and these are either a small list of simple groups that are not Lie type in defining characteristic, or if the potential maximal is Lie type in defining characteristic then what is left are groups of small rank and small field size, together with a large collection of possible subgroups $\PSL_2(p^a)$, the focus of this paper.

The following general theorem is a summary of our results, although we have much more detail about a putative maximal subgroup arising in (\ref{thmi:maind}).

\begin{theorem}\label{thm:generaltheorem} Let $p$ be a prime and $a,b\geq 1$ be integers. Let $G$ be the simple group $F_4(p^b)$, $E_6(p^b)$, ${}^2\!E_6(p^b)$ or $E_7(p^b)$, and let $H$ be a simple group $\PSL_2(p^a)$ contained in $G$. Write $\bar G$ for an almost simple group with socle $G$. If $N_{\bar G}(H)$ is an almost simple, maximal subgroup of $\bar G$, then one of the following holds:
\begin{enumerate}
\item\label{thmi:maina} $G=F_4(p^b)$ for some $p\geq 13$, $a=b$, and $N_{\bar G}(H)$ is unique up to conjugation;
\item\label{thmi:mainb} $G=E_7(p^b)$ for some $p\geq 17$, $a=b$, and $N_{\bar G}(H)$ is unique up to conjugation for $p=17$, and there are two classes for $p\geq 19$;
\item\label{thmi:mainc} $G=E_7(p^b)$ for some $p$, $a=7b$, and $N_{\bar G}(H)$ is unique up to conjugation;
\item\label{thmi:maind} $G=E_7(p^b)$, $p^a$ is one of $7,8,25$.
\end{enumerate}
In (\ref{thmi:maina}) and (\ref{thmi:mainb}), the subgroup $\PSL_2(p^a)$ arises from an $A_1$ subgroup of the algebraic group $F_4$ and $E_7$. In (\ref{thmi:mainc}), the subgroup $\PSL_2(p^a)$ arises from an $A_1^7$ subgroup of the algebraic group $E_7$. There are no known examples of maximal subgroups in (\ref{thmi:maind}).
\end{theorem}

We now give more information about what we prove for each group $G$.

\begin{theorem}\label{thm:f4} Let $p$ be a prime and $a,b\geq 1$ be integers. Let $G$ be an almost simple group with socle $F_4(p^b)$, and suppose that $H$ is a subgroup of $G$ with $F^*(H)=\PSL_2(p^a)$. If $H$ is maximal in $G$ then one of the following holds:
\begin{enumerate}
\item\label{thmi:f4a} $p^a=13$, $H=\PSL_2(13)$ and is a Serre embedding;
\item\label{thmi:f4b} $q=p^a$, $p\geq 13$, $F^*(H)=\PSL_2(q)$, and $H$ is the normalizer in $G$ of the fixed points $\bX^\sigma$ of an algebraic $A_1$ subgroup of the algebraic group $F_4$ under a Frobenius endomorphism $\sigma$.
\end{enumerate}
\end{theorem}

The definition of a Serre embedding is given formally in Definition \ref{defn:serreembedding}, but informally it is a copy of $\PSL_2(h+1)$ where $h$ is the Coxeter number of $G$ and this subgroup contains a regular unipotent element. (This subgroup is named after Serre as he constructed copies of $\PSL_2(h+1)$ (if $h+1$ is a prime) over all fields in \cite{serre1996}.) In recent work of Burness and Testerman \cite{burnesstesterman2017un}, Serre embeddings have been shown to come from algebraic $A_1$s, and so (\ref{thmi:f4a}) is a subcase of (\ref{thmi:f4b}) above. Thus Theorem \ref{thm:f4} implies Theorem \ref{thm:generaltheorem} for $G=F_4(p^b)$. Kay Magaard \cite{magaardphd} proved Theorem \ref{thm:f4} for $p\geq 5$ in his Ph.D. thesis, in addition proving that $b=1$ in (\ref{thmi:f4b}).

For $E_6$ we have a complete theorem, without relying on \cite{burnesstesterman2017un}, as we show that the Serre embedding lies in $F_4$.

\begin{theorem}\label{thm:e6} Let $p$ be a prime, $a,b\geq 1$ be integers, and let $G$ be an almost simple group with socle either $E_6(p^b)$ or ${}^2\!E_6(p^b)$. There does not exist a maximal subgroup $H$ of $G$ with $F^*(H)=\PSL_2(p^a)$.
\end{theorem}

Almost all of this theorem was obtained by Aschbacher \cite{aschbacherE6Vun} using geometric techniques, where only the case $p^a=p^b=11$ and $H$ contains a semiregular unipotent element, from class $E_6(a_1)$, is left open; here we prove the whole result again, using representation theory, and remove this final case using the Lie algebra structure of the adjoint module $L(E_6)$. Of course, Theorem \ref{thm:e6} implies Theorem \ref{thm:generaltheorem} for $G$ isomorphic to $E_6(p^b)$ and ${}^2\!E_6(p^b)$.

For $E_7$, here we have some potential exceptions. The difficult cases are the Serre embedding $p^a=19$ and $p^a=7,8,25$.

\begin{theorem}\label{thm:e7} Let $p$ be a prime and $a,b\geq 1$ be integers. Let $G$ be an almost simple group with socle $E_7(p^b)$, and suppose that $H$ is a maximal subgroup of $G$ with $F^*(H)=\PSL_2(p^a)$. One of the following holds:
\begin{enumerate}
\item\label{thmi:e7a} $p^a=7$, $p^a=8$ or $p^a=25$;
\item\label{thmi:e7b} $p^a=19$, $H=\PSL_2(p^a)$ and is a Serre embedding;
\item\label{thmi:e7c} $p^a=p^b$, $p\geq 17$, and $H$ is the normalizer in $G$ of the fixed points $\bX^\sigma$ of an algebraic $A_1$ subgroup of the algebraic group $E_7$. (There is one class for $p=17$, and two for $p\geq 19$.)
\item\label{thmi:e7d} $p^a=p^{7b}$ and $H$ is the normalizer in $G$ of the fixed points $\bX^\sigma$ of an algebraic $A_1^7$ subgroup of the algebraic group $E_7$.
\end{enumerate}
\end{theorem}

Again, Burness and Testerman have showed that (\ref{thmi:e7b}) is a subcase of (\ref{thmi:e7c}), and (\ref{thmi:e7d}) is the fixed points of a maximal-rank subgroup given in \cite[Table 5.1]{liebecksaxlseitz1992}. In the case (\ref{thmi:e7a}) where $p^a=8$, we can give the precise module structure of $H$ on the minimal module for $E_7$. For $p^a=7$, there are unresolved cases of potential copies of $\PSL_2(7)$ where the preimage of the subgroup in the simply connected version of $E_7$ is both $2\times \PSL_2(7)$ and $\SL_2(7)$. In both cases the module structures on both minimal and adjoint modules can be given precisely, but it seems difficult to progress further using these techniques. In the case of $p^a=25$, this is a copy of $\SL_2(25)$ inside the simply connected version of $E_7$ with the centre of the subgroup being equal to the centre of $E_7$, and we have complete information about the module structures on both the minimal and adjoint modules. It exists inside the $A_1A_1$ maximal subgroup of $E_7$, and if it is unique up to $\bG$-conjugacy then it is strongly imprimitive.\footnote{Uniqueness has since been proved, and will appear in forthcoming work. This removes this possibility from Theorems \ref{thm:generaltheorem} and \ref{thm:e7}.} Theorem \ref{thm:e7} implies Theorem \ref{thm:generaltheorem} for $G$ isomorphic to $E_7(p^b)$.

We do not deal with maximal subgroups of $E_8$ here, and only consider it for certain lemmas, which will be useful in a later treatment of this case. For exceptional groups other than $E_8$, the minimal module has dimension much smaller than the dimension of the group (as an algebraic group) and we can use representation theory to analyse this module. We can still do things with the Lie algebra for $E_8$, as we did in \cite{craven2015un2}, but these are postponed to avoid making this work even longer.

\medskip

The strategy for the proofs of these theorems is given in Chapter \ref{ch:strategy}, and relies heavily on computer calculations in three ways:
\begin{enumerate}
\item The first is to compute the traces of semisimple elements of large order on various modules for exceptional groups. Tables of these traces are available for elements of small order, but we need them for very large orders, sometimes in the hundreds. For this we can use the program that Alastair Litterick produced in his Ph.D.\ thesis \cite[Chapter 7]{litterick}, or construct the normalizer of a torus explicitly in Magma and take the conjugacy classes, then compute their eigenvalues. (Litterick has produced a much faster algorithm for computing traces of elements on $p$-restricted modules, but we do not need this for our cases.)
\item The second is to do large linear algebra problems. To find all sets of composition factors that could arise as the composition factors of the restriction of a $kG$-module to a subgroup $H$ involves checking many possible combinations against the large lists of traces of semisimple elements. This is done to reduce the possible module structures for the subgroup on the minimal and adjoint modules, and was also used in \cite{litterick}.
\item The third is to construct explicit modules for finite groups, and show that certain module structures cannot exist. This would be possible by hand, at least in some cases, but incredibly complicated and prone to mistakes. In each case, a clear recipe is given for how to reproduce the module we construct to ease verifiability. The Magma commands \texttt{ProjectiveCover(M)} and \texttt{Ext(A,B)} compute projective modules and the space of extensions between two modules. If \texttt{M} is a module and \texttt{N} is a simple module, then one can construct the maximal extension of \texttt{N} by \texttt{M} with the code
\begin{verbatim}E,rho:=Ext(M,N);
Mnew:=MaximalExtension(M,N,E,rho);
\end{verbatim}
This can be used to easily verify statements made in the paper about the structures of certain modules.
\end{enumerate}
With these three uses of a computer in mind, the rest of the argument is done by hand, in Chapters \ref{ch:f4} to \ref{ch:e7oddsl}.

The structure of this article is as follows: in the next chapter we give notation and some preliminary results. In Chapter \ref{ch:maxsubgroups} we give information about maximal subgroups of finite and algebraic exceptional groups. Chapter \ref{ch:maxsubspace} gives results about how to prove results about almost simple subgroups of algebraic groups, given information about their simple socles. Chapter \ref{ch:semisimple} proves results about `blueprints', finite subgroups of an algebraic group $\bG$ that stabilize exactly the same subspaces of a (non-trivial) rational $k\bG$-module as some infinite subgroup. We obtain upper bounds on the orders of semisimple elements that are not blueprints for the minimal module for $F_4$, $E_6$ and $E_7$. In Chapter \ref{ch:unipssemis} after, we give lots of information about unipotent and semisimple elements of exceptional groups, together with information about $\slf_2$-subalgebras of exceptional Lie algebras. Chapter \ref{ch:sl2modules} gives information about modules for $\SL_2(p^a)$, and the chapter after gives some constructions of $\PSL_2$s inside $E_6$ in characteristic $3$.

We then launch into the proof proper, with Chapter \ref{ch:strategy} giving an outline of the strategy of the proof, Chapters \ref{ch:f4} and \ref{ch:e6} proving the results for $F_4$ and $E_6$, and then the three chapters after doing $E_7$ first in characteristic $2$, and then $E_7$ in odd characteristic, split into two chapters according as the embedding into the simply connected group is $2\times \PSL_2(p^a)$ or $\SL_2(p^a)$.

The first appendix gives some widely known information about the composition factors of the reductive and parabolic maximal subgroups of $F_4$, $E_6$ and $E_7$ on the minimal and adjoint modules, information that is well known but given here for ease of reference. The second gives the traces of semisimple elements of small order on the minimal and adjoint modules for the algebraic groups $F_4$, $E_6$ and $E_7$.

\bigskip

\textbf{Acknowledgement}: The author would like to thank the referee for many detailed and helpful comments that have greatly improved the exposition of the manuscript.

\chapter{Notation and Preliminaries}
\label{ch:notation}

In this chapter we give the notation that we need, for both groups and modules, and give a few preliminary results.

Throughout this paper, $p\geq 2$ is a prime number, $q$ is a power of $p$, and $G=G(q)$ denotes an exceptional finite group of Lie type defined over $\F_q$. More specifically, let $\bG$ be a simple, simply connected algebraic group of exceptional type over the algebraic closure $k$ of $\F_p$, equipped with a Frobenius endomorphism $\sigma$, and set $G=\bG^\sigma$. The precise types of $G$ that we are interested in are those exceptional groups whose maximal subgroups are not yet known, i.e., $F_4(q)$, $E_6(q)$, ${}^2\!E_6(q)$, $E_7(q)$ and $E_8(q)$, although we do not do much in the case of $E_8(q)$, and often will exclude it from consideration.

Notice that we consider the simply connected versions of $\bG$ and $G$, so $E_7(q)$ possesses a centre when $p$ is odd, and $E_6(q)$ does for $3\mid (q-1)$. We want the simply connected versions in order to work with the minimal module and the adjoint module simultaneously. Where this is particularly important we will remind the reader, for example when considering $\PSL_2(p^a)$ embedded in the simple group of type $E_7$, where we can embed either $\SL_2(p^a)$ in $E_7$ with the centres coinciding or $2\times \PSL_2(p^a)$ into $E_7$ with the centres coinciding, representing the two possible preimages of a copy of $\PSL_2(p^a)$ in the simple group. If $G$ possesses a graph automorphism of order $2$, denote this by $\tau$; we will remind the reader of this notation when we use it.

We let $\bar G$ be an almost simple group with socle $G/Z(G)$, which embeds into $\Aut(G)$. The maximal subgroups $M$ of $\bar G$ split into three categories: $M\cap (G/Z(G))$ is a maximal subgroup of $G/Z(G)$, $M\cap (G/Z(G))$ is not a maximal subgroup of $G/Z(G)$, and $(G/Z(G))\leq M$. The third collection is easily computed, and the first can be deduced from a list of maximal subgroups of $G$ by taking normalizers. However, the second, called \emph{novelty} maximal subgroups, cannot easily be seen from the maximal subgroups of $G$. They arise in the following manner: let $H$ be a subgroup that is not maximal in a simple group $X$, but $H$ is normalized by a group of automorphisms $A$ of $X$ while every proper subgroup of $X$ properly containing $H$ is not normalized by it. In this case, $H.A$ is a maximal subgroup of $X.A$. However, it is of course very difficult to understand these if one is simply given a list of maximal subgroups of $X$, so we will prove more than simply that a given subgroup is not maximal in the simple group, but that it is contained in stabilizers of various subspaces of a given module, enough that we can see that it cannot form a novelty maximal subgroup.

Let $L(\lambda)$ denote the irreducible highest weight module of weight $\lambda$. The notation for the weight lattice is `standard', consistent with the main references in this work and can be found in for example \cite[Chapter VI]{bourbakilie2}. The modules that we normally consider are the two smallest non-trivial ones. Write $M(\bG)$ for one of the minimal modules for $\bG$, namely $L(\lambda_4)$ for $F_4$, either $L(\lambda_1)$ or $L(\lambda_6)$ for $E_6$ and ${}^2\!E_6$, $L(\lambda_7)$ for $E_7$ and not defined for $E_8$. We write $L(\bG)$ for the Lie algebra or adjoint module, which is $L(\lambda_1)$, $L(\lambda_2)$, $L(\lambda_1)$ and $L(\lambda_1)$ respectively. If $L(\bG)$ has a trivial composition factor so is not irreducible, which occurs in $E_7$ in characteristic $2$ and $E_6$ in characteristic $3$, let $L(\bG)^\circ$ denote the non-trivial composition factor, and in other cases let $L(\bG)^\circ=L(\bG)$. These two modules have the following dimensions:
\begin{center}
\begin{tabular}{ccc}
\hline Group & $\dim(M(\bG))$ & $\dim(L(\bG)^\circ)$
\\\hline $F_4$ & $26-\delta_{p,3}$ & $52$
\\ $E_6$ & $27$ & $78-\delta_{p,3}$
\\ $E_7$ & $56$ & $133-\delta_{p,2}$
\\ $E_8$ & $-$ & $248$
\\ \hline
\end{tabular}
\end{center}
In characteristic $2$, $L(F_4)$ has factors $L(\lambda_4)=M(F_4)$ and $L(\lambda_1)=M(F_4)^\tau$, where $\tau$ denotes the graph automorphism of $G$ (which does not extend `nicely' to a morphism of $\bG$, see the definition of $\Aut^+(\bG)$ in Chapter \ref{ch:maxsubgroups}), so in this case we can consider $L(\lambda_1)$ and $L(\lambda_4)$ or $L(F_4)$. In all other cases, $L(\bG)^\circ$ is irreducible.

\bigskip

We now introduce some notation for modules. All modules will be finite dimensional and are defined over $k$. If $H$ is a group, let $\Irr(H)$ denote the set of irreducible modules over the field, which is always $k$. We also write $k$ for the trivial module for any group over the field $k$, although we will also denote it by `$1$'. As usual write `$\oplus$' and `$\otimes$' for the direct sum and tensor product of two modules. Let $\Lambda^i$ and $S^i$ denote the exterior and symmetric powers.  Write $\soc^i(M)$ for the $i$th socle layer and $\rad^i(M)$ for the $i$th radical layer of $M$. Write $\topp(M)$ for the top of $M$, i.e., $M/\rad(M)$, and $\cf(M)$ for the composition factors of $M$ as a multiset. Let $H^1(H,M)$ denote the $1$-cohomology group of $M$, and in general $\Ext^1(M,M')$ denote the group of extensions with submodule $M'$ and quotient $M$. The projective cover of a module $M$ will be denoted by $P(M)$. Let $M^*$ denote the dual of $M$.

We write $M\downarrow_H$ for the restriction of $M$ to $H$. Let $H$ and $G$ be groups, let $M$ be a $kG$-module, and let $\mathcal N$ be a set of simple $kH$-modules. Write $N$ for the direct sum of all members of $\mathcal N$. The set $\mathcal N$ is \emph{conspicuous} for $M$ if, for every $p$-regular element $x$  (i.e., order not divisible by $p$) of $H$, there is a $p$-regular element $y$ of $G$ such that the eigenvalues of $x$ on $M$ coincide with the eigenvalues of $y$ on $N$. Informally, this means that $N$ could be the composition factors of $M\downarrow_H$, if $H$ were a subgroup of $G$. Normally we use this by specifying $H$ to be a subgroup of $G$, but whose composition factors on $M$ are unknown. We might not want to check all $p$-regular elements of $H$, so we will say that $\mathcal N$ is conspicuous for elements of particular orders when we do not check all elements.

If $u$ is an element of a group $H$ of $p$-power order, and $M$ is a $kH$-module, then $u$ acts on $M$ as a sum of Jordan blocks of various dimensions. If the dimensions are, say, $5$, $5$, and $1$, we write that the structure is $5^2,1$. This is in keeping with the notation from \cite{lawther1995}, which is our main reference for the actions of unipotent classes on $M(\bG)$ and $L(\bG)$. The element $u$ \emph{acts projectively} on $M$ if all blocks of $u$ on $M$ have size $o(u)$, where $o(u)$ denotes the order of $u$. This is equivalent to the restriction $M\downarrow_{\gen u}$ being a projective $k\gen{u}$-module.

We will often have to talk about the structures of modules, as in their socle layers. If $M$ is a module with socle $A$ and second socle $B$ then we can write
\[ \begin{array}{c}B\\A\end{array}\]
for this structure; however, this is often too space-consuming when we have many socle layers, and so we also write $B/A$ for this module. Generalizing this, we delineate between socle layers by `$/$', so that $A/B,C/D,E$ is a module with socle $D\oplus E$, second socle $B\oplus C$, and third socle $A$. Because of the potential for confusion with the quotient, we will always mention if we mean a quotient.

We also introduce the concepts of radical and residual. If $I$ is a subset of $\Irr(H)$, then the \emph{$I$-radical} of $M$ is the largest submodule of $M$ whose composition factors lie in $I$, and the \emph{$I$-residual} of $M$ is the smallest submodule such that every composition factor of the quotient lies in $I$. Write $I'$ for $\Irr(H)\setminus I$.

One lemma that we occasionally use, that can be quite powerful, relates the minimal and adjoint modules for exceptional groups. We place it here because there seems no more appropriate place.

\begin{lemma}\label{lem:relatingvminlg} Let $\bG$ be one of $F_4$, $E_6$ and $E_7$.
\begin{enumerate}
\item\label{lemi:vmina} Let $\bG=F_4$. If $p=3$ then $L(F_4)$ is a submodule of $\Lambda^2(M(F_4))$. If $p\geq 5$ then $L(F_4)$ is a summand of $\Lambda^2(M(F_4))$.
\item\label{lemi:vminb} Let $\bG=E_6$. If $p=2$ then $L(E_6)$ is a submodule of $M(E_6)\otimes M(E_6)^*$. If $p=3$ then the socle of 
$M(E_6)\otimes M(E_6)^*$ is $1$-dimensional, and quotienting out by this, $L(E_6)^\circ$ is a submodule. If $p\geq 5$ then $L(E_6)$ is a summand of $M(E_6)\otimes M(E_6)^*$.
\item\label{lemi:vminc} Let $\bG=E_7$. If $p=2$ then the socle of 
$\Lambda^2(M(E_7))$ is $1$-dimensional, and quotienting out by this, $L(E_7)^\circ$ is a submodule. If $p=3$ then $L(E_7)$ is a submodule of $S^2(M(E_7))$. If $p\geq 5$ then $L(E_7)$ is a summand of $S^2(M(E_7))$.
\end{enumerate}
\end{lemma}
\begin{proof} In characteristic $0$, the highest weight modules in $\Lambda^2(M(F_4))$ are $L(\lambda_1)=L(F_4)$ and $L(\lambda_3)$, of dimension $273$. If $p\geq 5$ these two modules remain irreducible, so the exterior square is the sum of these, but $L(\lambda_3)$ has dimension $196$ when $p=3$, so that the exterior square of $M(F_4)$ -- which has dimension $300$ as $M(F_4)$ has dimension $25$ -- has three composition factors: two copies of $L(\lambda_1)$ and one of $L(\lambda_3)$. As the exterior square is self-dual, either this module has the form $52/196/52$ or there is a $52$- or $196$-dimensional summand. However, the regular unipotent element of $F_4$ acts on $M(F_4)$ with blocks $15,9,1$ (see Table \ref{t:unipotentF4}) and thus on the exterior square with blocks $27,21^2,18^6,15^2,9^{10},3$, all of which are divisible by $3$, so there is no simple summand of the exterior square. This proves (\ref{lemi:vmina}).

For characteristic $0$, $M(E_6)\otimes M(E_6)^*=L(\lambda_1)\otimes L(\lambda_6)$ has composition factors the highest weight modules $L(0)$, $L(\lambda_2)=L(E_6)$ and $L(\lambda_1+\lambda_6)$ (of dimension $650$). These again remain irreducible for $p\geq 5$, so we obtain the result. In characteristic $2$, $L(\lambda_1+\lambda_6)$ has dimension $572$ and a copy of $L(\lambda_2)$ is the other composition factor of the Weyl module $W(\lambda_1+\lambda_6)$. Thus we get $L(0)\oplus (L(\lambda_2)/L(\lambda_1+\lambda_6)/L(\lambda_2))$: the trivial breaks off as a summand as the module is self-dual and has a trivial submodule, and the easiest way to check that the second summand is indecomposable is to prove it for $E_6(2)$ with a computer, where we indeed obtain 
\[1\oplus (78/572/78).\]
In characteristic $3$, the composition factors become $L(0)^3,L(\lambda_2)^2,L(\lambda_1+\lambda_6)$. Again, the easiest way to deduce the structure is to test it for $E_6(3)$ using a computer, and note that it has form
\[ 1/77/1,572/77/1,\]
hence must for any group of type $E_6$ because it must be a refinement of this structure and remain self-dual.

For $E_7$, in characteristic $0$, $S^2(M(E_7))$ has composition factors the highest weight modules $L(\lambda_1)=L(E_7)$ and $L(2\lambda_7)$ of dimension $1463$. These both remain irreducible for $p\geq 5$, and so $S^2(M(E_7))$ is the sum of these. For $p=3$, $L(2\lambda_7)$ has dimension $1330$, and the structure is $L(\lambda_1)/L(2\lambda_7)/L(\lambda_1)$. Again, this is most easily checked with a computer for $E_7(3)$. For $p=2$ we take the exterior square, which has composition factors of dimension $1,1,132,132,1274$, and as with $E_6$, the easiest way to check the structure is with a computer for $E_7(2)$, where we get
\[ L(0)/L(\lambda_1)/L(\lambda_6)/L(\lambda_1)/L(0).\]
This completes the proof.
\end{proof}

In many cases we want to prove that a module has a particular composition factor as a submodule or quotient, often the trivial module, which we denote by $k$ or $1$. Thus we need a method of proving that a particular composition factor is always a submodule or quotient in any module with those factors. This is the idea of pressure.

Suppose that $H$ is a finite group such that $O^p(H)=H$, and such that for all simple modules $M$ over a field $k$, $H^1(H,M)=H^1(H,M^*)$. The \emph{pressure} of a module $V$ for $H$ is the quantity
\[ \sum_{M\in\cf(V)} (\dim H^1(H,M)-\delta_{M,k}),\]
where $\delta$ is the Kronecker delta. Results on pressure have occurred in the literature before, with the most general so-far being \cite[Lemma 1.8]{craven2015un2}. Another generalization of this allows us to understand the situation of forcing a module from a collection $\mathcal M$ of simple modules to be a submodule of a given module $V$. If $\mathcal M$ is a collection of simple modules for a group $H$, with $\Ext^1(M,M')=0$ for all $M,M'\in\mathcal M$, and such that $\Ext^1(A,M)=\Ext^1(M,A)$ for all simple modules $A$ and $M$ with $M\in\mathcal M$, then the \emph{$\mathcal M$-pressure} of a module $V$ is the quantity
\[ \sum_{M'\in \cf(V)}\sum_{M\in \mathcal M}(\dim\Ext^1(M,M')-\delta_{M,M'}).\]
The result \cite[Lemma 1.8]{craven2015un2} directly generalizes to $\mathcal M$-pressure, with the exact same proof, so we simply state the result.

\begin{lemma}\label{lem:pressure} Suppose that $H$ is a finite group, and let $\mathcal M$ be a set of simple modules for $H$ such that $\Ext^1(M,M')=0$ for all $M,M'\in \mathcal M$, and $\Ext^1(M,A)=\Ext^1(A,M)$ for all $M\in \mathcal M$ and all simple modules $A$. Let $V$ be a module for $H$ of $\mathcal M$-pressure $n$.
\begin{enumerate}
\item If $n<0$ then $\Hom(M,V)\neq 0$ for some $M\in \mathcal M$, i.e., $V$ has a simple submodule isomorphic to some $M\in\mathcal M$. If $n=0$ then either $\Hom(M,V)\neq 0$ or $\Hom(V,M)\neq 0$, i.e., $V$ has either a simple submodule or quotient isomorphic to some member of $\mathcal M$.

\item More generally, if a composition factor of $V$ has $\mathcal M$-pressure greater than $n$, then either $\Hom(M,V)\neq 0$ or $\Hom(V,M)\neq 0$ for some $M\in\mathcal M$.

\item If $\Hom(M,V)=\Hom(V,M)=0$ for all $M\in\mathcal M$, then any subquotient $W$ of $V$ has $\mathcal M$-pressure between $-n$ and $n$.
\end{enumerate}
\end{lemma}

The concept of pressure can be used to prove that either $M(\bG)$ or $L(\bG)$ possesses a trivial submodule or quotient when restricted to $H$. We therefore would like to know whether that is enough in some circumstances to conclude that $H$ is contained within a positive-dimensional subgroup of $\bG$. The next result is \cite[Lemma 1.4]{craven2015un2}.

\begin{lemma}\label{lem:smallsubspaces} Let $\bG$ be one of $F_4$, $E_6$, ${}^2\!E_6$, $E_7$ or $E_8$. Let $H\leq \bG$. If one of the following holds, then $H$ is contained in a positive-dimensional subgroup of $\bG$:
\begin{enumerate}
\item $H$ stabilizes a $1$-space or hyperplane of $M(\bG)$ or $L(\bG)$;
\item $\bG=F_4$, $E_6$, ${}^2\!E_6$ or $E_7$, and $H$ stabilizes a $2$-space or a space of codimension $2$ in $M(\bG)$;
\item $\bG=E_6$ or ${}^2\!E_6$, and $H$ stabilizes a $3$-space or a $24$-space of $M(E_6)$.
\end{enumerate}
\end{lemma}

In Propositions \ref{prop:fixlineonLG}, \ref{prop:fixlineonMG} and \ref{prop:fix2spaceonMG}, we extend the statements about stabilizing a line of $L(\bG)^\circ$ or $M(\bG)$, or a $2$-space of $M(\bG)$, to include stability under outer automorphisms of the finite group $G$. This allows such statements to be used to deduce the results in the introduction about maximal subgroups.

We end with giving the line stabilizers for the minimal modules for the finite groups $E_6(q)$ and $E_7(q)$. These have appeared in the literature before, and we take these from \cite[Lemmas 5.4 and 4.3]{liebecksaxl1987}.

\begin{lemma}\label{lem:e6stabs} Let $G=E_6(q)$. There are three orbits of lines of the action of $G$ on $M(E_6)$, with line stabilizers as follows:
\begin{enumerate}
\item $F_4(q)$ acting on $M(E_6)$ as $L(\lambda_4)\oplus L(0)$ ($L(0)/L(\lambda_4)/L(0)$ in characteristic $3$);
\item a $D_5$-parabolic subgroup; $q^{16}D_5(q).(q-1)$, acting uniserially as \[L(\lambda_1)/L(\lambda_4)/L(0);\]
\item a subgroup $q^{16}.B_4(q).(q-1)$ acting indecomposably as \[L(0)/L(\lambda_4)/L(0),L(\lambda_1).\]
\end{enumerate}
\end{lemma}

\begin{lemma}\label{lem:e7stabs} Let $G=E_7(q)$. There are five orbits of lines of the action of $G$ on $M(E_7)$, with line stabilizers as follows:
\begin{enumerate}
\item $E_6(q).2$ (the graph automorphism) acting semisimply with composition factors of dimensions $54,1,1$;
\item ${}^2\!E_6(q).2$ (the graph automorphism) acting semisimply with composition factors of dimensions $54,1,1$;
\item an $E_6$-parabolic subgroup $q^{27}.E_6(q).(q-1)$ acting uniserially as \[L(0)/L(\lambda_1)/L(\lambda_6)/L(0);\]
\item a subgroup $q^{1+32}.B_5(q).(q-1)$ acting indecomposably as \[L(0),L(\lambda_1)/L(\lambda_5)/L(0),L(\lambda_1);\]
\item a subgroup $q^{26}.F_4(q).(q-1)$ acting indecomposably as \[L(0),L(0)/L(\lambda_4)/L(\lambda_4)/L(0),L(0).\]
\end{enumerate}
\end{lemma}

\chapter{Maximal Subgroups}
\label{ch:maxsubgroups}
This chapter summarizes what is known about the maximal subgroups of the finite groups $\bar G$, and also the algebraic group $\bG$, about which complete information on positive-dimensional maximal subgroups is known.

The maximal subgroups of positive dimension in $\bG$ are given in \cite{liebeckseitz2004}, and given $\bG$ we denote this collection by $\mathscr X$; write $\mathscr X^\sigma$ for the fixed points $X=\bX^\sigma$ for $\bX$ a $\sigma$-stable member of $\mathscr X$. If $Z(G)\neq 1$, we also write $\mathscr X^\sigma$ for their images modulo the centre of $G$. If $\bar G$ is almost simple, the set $\mathscr X^\sigma$ shall be taken to mean the normalizers in $\bar G$ of the elements of $\mathscr X^\sigma$ for $F^*(\bar G)$. Because each definition applies to a different group ($G$, $G/Z(G)$ and $\bar G$) no confusion should arise.

For technical reasons, we do not include in $\mathscr X^\sigma$ the fixed points of $G$ under a field, graph, or field-graph automorphism of prime order (so, for example, ${}^2\!E_6(p^2)$ and $E_6(p)$ inside $E_6(p^2)$). Such subgroups are said to \emph{have the same type} as $\bG$. 

While the maximal subgroups of $\bG$ are known, the maximal subgroups of $G$ and $\bar G$ are of course not. We start with a broad characterization of the maximal subgroups of $\bar G$, given in \cite{borovik1989} and \cite[Theorem 2]{liebeckseitz1990}.

\begin{theorem}\label{thm:classificationMaximal}
Let $M$ be a maximal subgroup of $\bar G$ not containing $F^*(\bar G)$. One of the following holds:
\begin{enumerate}
\item\label{thmi:maxa} $M$ is a member of $\mathscr X^\sigma$ or $M$ has the same type as $\bG$;
\item\label{thmi:maxb} $M$ is the normalizer of an elementary abelian $r$-group for some $r\neq p$ (an \emph{exotic $r$-local subgroup});
\item\label{thmi:maxc} $F^*(M)=\Alt(5)\times \Alt(6)$ and $G=E_8$ with $p>5$;
\item\label{thmi:maxd} $M$ is an almost simple group whose socle is not a group of Lie type in characteristic $p$;
\item\label{thmi:maxe} $M$ is an almost simple group whose socle is a group of Lie type in characteristic $p$, that does not appear in (\ref{thmi:maxa}).
\end{enumerate}
\end{theorem}

The subgroups in (\ref{thmi:maxa}) are known and are the fixed points of those in \cite{liebeckseitz2004}, together with the (normalizers of) fixed points under field, graph, and field-graph automorphisms of $G$; the subgroups in (\ref{thmi:maxb}) are known and given in \cite{clss1992}; the subgroup (\ref{thmi:maxc}) was discovered by Borovik \cite{borovik1989} and is unique up to conjugacy. The potential subgroups in (\ref{thmi:maxd}) have been steadily reduced over the last two decades. Here the list is fairly short and given in \cite{liebeckseitz1999}, but note that a fair number of these have been eliminated in a variety of papers, too numerous to list here, but we mention the papers \cite{litterickmemoir} and \cite{craven2015un2} for all Lie type groups, and with $F_4$ and $E_6$ having almost all possibilities for $M$ removed by Magaard and Aschbacher in \cite{magaardphd} and \cite{aschbacherE6Vun} respectively. The author has also made progress on eliminating still more of this list and proving uniqueness of various maximal subgroups, with details to appear elsewhere.

It is probable that (\ref{thmi:maxe}) is empty. To express this another way, we introduce a definition here that will be used throughout the text. As in \cite[Section 2]{liebeckmartinshalev2005}, define $\Aut^+(\bG)$ to be the group generated by inner, diagonal, graph and $p$-power field automorphisms of $\bG$. As seen in \cite[Section 12.5]{cartersimple} if $G=\bG^\sigma$ then any automorphism of $G$ extends to an element of $\Aut^+(\bG)$, but if  $p=2$ and $\bG=B_2,F_4$, and $p=3$ and $\bG=G_2$ we must be more careful. In these cases there are `exceptional' graph morphisms, which can be included in $\Aut^+(\bG)$ (see \cite[Sections 12.3 and 12.4]{cartersimple}). However, because the graph morphism powers to a field automorphism ($\Out(F_4(2^b))$ is cyclic) we can only add a single graph morphism to $\Aut^+(\bG)$ at a time. Since the almost simple group $\bar G$ can only induce one graph morphism on $G$, this restriction is merely formal, and does not affect our conclusions. Unlike \cite{liebeckmartinshalev2005}, we will include these in $\Aut^+(\bG)$: the only case where this will make a difference for us is in $F_4$ for $p=2$. Of course this means that $\Aut^+(\bG)$ depends on the specific Frobenius morphism $\sigma$.

\begin{definition} A subgroup $H$ of $\bG$ that is not contained in any member of $\mathscr X$, and is not of the same type as $\bG$, is called \emph{Lie primitive}, and otherwise \emph{Lie imprimitive}. If $\sigma$ is a Frobenius endomorphism on $\bG$ and $H$ is contained in $G=\bG^\sigma$, then $H$ is called \emph{strongly imprimitive} if $H$ is contained in a $\sigma$-stable, $N_{\Aut^+(\bG)}(H)$-stable member of $\mathscr X$.
\end{definition}

The condition on $\sigma$-stability is just that $H$ is contained in a member of $\mathscr X^\sigma$. It is not necessary as it is subsumed under the second condition, but we include it for emphasis. The condition on automorphisms is needed so that $N_{\bar G}(H)$ is contained in a member of $\mathscr X^\sigma$. Hence the statement that (\ref{thmi:maxe}) is empty is equivalent to the following.

\begin{conjecture}\label{conj:morphismextension} Let $H$ be a simple subgroup of $\bG$, and that $H$ is a group of Lie type in characteristic $p$, not of the same type as $\bG$. Then $H$ is strongly imprimitive.
\end{conjecture}

Such a statement is true for classical algebraic groups \cite[Theorem 4]{seitz1988}, and for exceptional groups and $p>113$ in \cite{seitztesterman1990}, in which the above conjecture was suggested. It also follows for the group $G_2$ from the classification of maximal subgroups of $G_2(q)$, given by Cooperstein \cite{cooperstein1981} and Kleidman \cite{kleidman1988}. We will need the fact that Conjecture \ref{conj:morphismextension} is true for classical groups and $G_2$ later on, so we mention it now.

\begin{proposition}\label{prop:morphismextensiontrueforclassical} If $\bG$ is reductive, and is a product of classical algebraic groups and copies of $G_2$, and $H\cong \PSL_2(p^a)$ is a $\bG$-irreducible subgroup, then $H$ is contained in a $\bG$-irreducible $A_1$ subgroup of $\bG$.
\end{proposition}

Later work reduced the possibilities for (\ref{thmi:maxe}) still further: The rank of $M$ can only be at most half the rank of $\bG$ by \cite{lst1996,liebeckseitz2005}. Furthermore, for those groups we have the following possibilities by \cite{liebeckseitz1998,lawther2014}:
\begin{enumerate}
\item $M(p^a)$ has semisimple rank at most half that of $G$, $p^a\leq 9$, and $M(p^a)$ is not one of $\PSL_2(p^a)$, ${}^2\!B_2(p^a)$ and ${}^2\!G_2(p^a)$;
\item $\PSL_3(16)$ and $\PSU_3(16)$;
\item $\PSL_2(p^a)$, ${}^2\!B_2(p^a)$ and ${}^2G_2(p^a)$ for $p^a\leq \gcd(2,p-1)\cdot t(\bG)$, where
\[ t(G_2)=12,\quad t(F_4)=68,\quad t(E_6)=124,\quad t(E_7)=388,\quad t(E_8)=1312.\]
\end{enumerate}

If $H$ is a subgroup of $\bG$ that is a group of Lie type in characteristic $p$, we therefore call $H$ \emph{large rank}, \emph{medium rank}, and \emph{small rank}, according as the semisimple rank of $H$ is more than half that of $\bG$, between $2$ and half that of $\bG$, and $1$ respectively. Thus there are no large-rank members of (\ref{thmi:maxe}), the medium-rank members of (\ref{thmi:maxe}) are only over small fields (although the Suzuki and small Ree groups are also medium rank) and the small-rank subgroups are $\PSL_2(p^a)$, where $p^a$ can be quite large compared with medium-rank groups. Theorem \ref{thm:generaltheorem} states that Conjecture \ref{conj:morphismextension} is true for $H\cong \PSL_2(p^a)$ and $\bG$ of types $F_4$, $E_6$ and $E_7$, but with potential counterexamples for $\bG$ of type $E_7$ and $p^a=7,8,25$.

\medskip

The constant $t(\bG)$ is linked to the eigenvalues of semisimple elements of $\bG$ on $L(\bG)$, but to state it uses the following definition. It appears in work of Liebeck and Seitz \cite{liebeckseitz1998}, but is not specifically defined there.

\begin{definition}\label{defn:blueprints} Let $\bG$ be an infinite group and let $V$ be a module for $\bG$. A finite subgroup $H$ of $\bG$ is a \emph{blueprint} for $V$ if there exists an infinite subgroup $\bX$ of $\bG$ such that $\bX$ and $H$ stabilize the same subspaces of $V$. An element $x$ is a blueprint for $V$ if $\gen x$ is.
\end{definition}

If $H$ is a blueprint for $V$ a simple $k\bG$-module, then one of two cases hold: $H$ is irreducible on $V$, or the intersection of the stabilizers of all $H$-stable subspaces of $\bG$ is a proper, positive-dimensional subgroup $\bX$ of $\bG$ containing $H$. If $V$ is either $M(\bG)$ or $L(\bG)^\circ$ then the subgroups of $\bG$ acting irreducibly on $V$ are classified (but not up to conjugacy) in \cite{liebeckseitz2004a}.

\medskip

If one wants to prove that a finite subgroup $H$ of $\bG$ is strongly imprimitive, one may directly construct an infinite subgroup containing it. An alternative is to use the following statement, which easily follows from the definitions.

\begin{lemma}\label{lem:containstrongimp} Let $H$ and $K$ be finite subgroups of $\bG$, with $H\leq K$. If $K$ is Lie imprimitive then so is $H$. If $K$ is both strongly imprimitive and $N_{\Aut^+(\bG)}(H)$-stable, then $H$ is strongly imprimitive.
\end{lemma}

Fix $\bG$. Up to isomorphism, there are only finitely many candidates for simple subgroups of $\bG$ that are not strongly imprimitive: subgroups of the exotic $r$-local subgroups, the medium- and small-rank subgroups above, and the almost simple groups that are not Lie type in characteristic $p$ from \cite{liebeckseitz1999}. Write $\mathscr P$ for this set, and partially order it by inclusion.

Let $H$ be a subgroup of $\bG$ that is a maximal member of $\mathscr P$. Suppose that $H$ is not strongly imprimitive. Suppose also that $K$ is another finite subgroup of $\bG$ containing $H$, and that $K$ is $N_{\Aut^+(\bG)}(H)$-stable. By Lemma \ref{lem:containstrongimp}, either $K\leq N_{\bG}(H)$ or $K$ has the same type as $\bG$, because $K$ cannot be strongly imprimitive, but does not lie in $\mathscr P$.

If we construct a subgroup $K=\gen{H,x}$ so that $K$ stabilizes a proper, non-zero subspace $W$ of either $M(\bG)$ or $L(\bG)^\circ$, then $K$ cannot have the same type as $\bG$. Since $\Aut(H)$ is known, one may find the orders of all of its elements, and if $o(x)$ is not one of them, $K\not\leq N_{\bG}(H)$. If $K$ stabilizes every subspace in the $N_{\Aut^+(\bG)}(H)$-orbit of $W$ then $K$ is $N_{\Aut^+(\bG)}(H)$-stable by Proposition \ref{prop:intersectionstabilizers} below. We have therefore derived a contradiction, and so $H$ must be strongly imprimitive.

\begin{proposition}\label{prop:maximalnotinP} Let $H$ be a maximal member of $\mathscr P$, and suppose that $K$ is a finite subgroup of $\bG$ containing $H$. If $K\not\leq N_{\bG}(H)$, and $K$ stabilizes a non-zero, proper subspace of an irreducible $k\bG$-module, then $H$ is Lie imprimitive. If, moreover, $K$ is $N_{\Aut^+(\bG)}(H)$-stable, then $H$ is strongly imprimitive.
\end{proposition}

In order to use this argument, we need to know the maximal members of $\mathscr P$. Those members of $\mathscr P$ that are not Lie type in characteristic $p$ rarely contain $\PSL_2(p^a)$, so we wish to exclude some medium-rank groups from $\mathscr P$ in order to force more $A_1$ subgroups to be maximal in $\mathscr P$.

We will prove this here, because this is the most natural place for it. However, it uses techniques and results from Chapters \ref{ch:maxsubspace} to \ref{ch:unipssemis}, and it would be better to read the proof after those chapters.

\begin{proposition}\label{prop:eliminatemediumrank} Let $\bG$ be one of $F_4$, $E_6$ and $E_7$.
\begin{enumerate}
\item\label{propi:medra} For $p=5,7$, any copy of $H=\PSp_4(p)$ in $\bG$, or $H=\Sp_4(p)$ in $\bG=E_7$ with $Z(H)=Z(\bG)$, is a blueprint for $M(\bG)\oplus M(\bG)^*$.

\item\label{propi:medrb} For $p=5,7$, any copy of $H=\PSL_4(p)$ or $\PSU_4(p)$ in $\bG$, or $H=2\cdot \PSL_4(p)$ or $H=2\cdot \PSU_4(p)$ in $\bG=E_7$ with $Z(H)=Z(\bG)$, is a blueprint for $M(\bG)\oplus M(\bG)^*$.

\item\label{propi:medrc} Let $p$ be an odd prime and $a\geq 1$. Any copy of $H=\PSp_6(p^a)$ in $\bG$, or $H=\Sp_6(p^a)$ in $\bG=E_7$ with $Z(H)=Z(\bG)$, is a blueprint for $M(\bG)\oplus M(\bG)^*$.

\item\label{propi:medrd} Let $p$ be an odd prime and $a\geq 1$. Any copy of $H=\Omega_7(p^a)$ in $\bG$, or $H=\Spin_7(p^a)$ in $\bG=E_7$ with $Z(H)=Z(\bG)$, is a blueprint for $M(\bG)\oplus M(\bG)^*$.
\end{enumerate}
Therefore if $p=5,7$ then $\PSL_4(p)$, $\PSU_4(p)$ and $\PSp_4(p)$ do not appear in $\mathscr P$, and for $p^a=3,5,7,9$, $\PSp_6(p^a)$ and $\Omega_7(p^a)$ do not appear in $\mathscr P$.
\end{proposition}
\begin{proof} We first prove the conclusion. None of the groups in this result appears in \cite[Tables 1.2 and 1.3]{liebeckseitz2004a}, so $H$ cannot act irreducibly on $M(\bG)$. Thus since $H$ is a blueprint for $M(\bG)$, $H$ is strongly imprimitive by Proposition \ref{prop:blueprintissi}.

We now prove the individual statements. Note that, unless $M(\bG)\not\cong M(\bG)^*$, taking the direct sum is irrelevant for whether a subgroup is a blueprint. We will therefore prove the result for $M(E_7)$, and then use the embeddings of $F_4$ and $E_6$ into $E_7$ to descend to those groups.

We prove (\ref{propi:medra}) and (\ref{propi:medrb}) for $p=5$ first. For (\ref{propi:medra}), we compute the conspicuous sets of composition factors for $M(E_7)\downarrow_H$. The simple modules for $\PSp_4(5)$ of dimension at most $56$ are $1$, $5$, $10$, $13$, $30$, $35_1$, $35_2$ and $55$. The only conspicuous set of composition factors is $10^2,5^6,1^6$, and since these composition factors have no extensions with each other, $M(E_7)\downarrow_H$ is semisimple. Let $u$ denote an element of order $5$ with the largest centralizer. This acts on $M(E_7)$ with Jordan blocks $3^2,2^{16},1^{18}$, which is the generic class $2A_1$ (see Definition \ref{defn:genericunipotent}) by \cite[Table 7]{lawther1995}. By Lemma \ref{lem:genericunipotent}, $u$, and therefore $H$, is a blueprint for $M(E_7)$.

If $H=\Sp_4(5)\leq E_7(k)$ with centres coinciding, then the involutions in $H$ act on faithful modules with trace $0$, which is not allowed since the trace of an involution in $E_7$ is $\pm 8$ (see Appendix \ref{app:traces}). This proves (\ref{propi:medra}).

As $\SL_4(5)$ and $\SU_4(5)$ contain $\Sp_4(5)$, and the centres of $\SL_4(5)$ and $\SU_4(5)$ contain the centre of $\Sp_4(5)$, we have that $\PSp_4(5)\leq \PSL_4(5),\PSU_4(5)$, and therefore (\ref{propi:medrb}) holds as subgroups that contain blueprints are themselves blueprints.

For $p=7$ the exact same proof holds, except that the dimensions of the simple modules for $\PSp_4(7)$ are now $1$, $5$, $10$, $14$, $25$, $35_1$, $35_2$ and $54$.

\medskip

We now prove (\ref{propi:medrc}). For $p^a\neq 3,5$, the largest semisimple element of odd order has order $171$ and $365$ respectively. These are greater than $75$, so these elements are blueprints for $M(E_7)$ by Theorem \ref{thm:blueprintsminimal}. Hence $H$ is a blueprint for $M(E_7)$ for $p^a=7,9$.

If $H=\PSp_6(3)$ then there are only three simple modules of dimension at most $56$, with dimensions $1$, $14$ and $21$. The traces of elements of orders $5$ and $7$ are enough to prove that $H$ does not embed in $\bG=E_7$, and hence not in $F_4$ or $E_6$. If $H=\Sp_6(3)$ then the appropriate simple modules have dimensions $6$, $14$ and $50$, and traces of elements of order $5$ are enough to prove that the only conspicuous set of composition factors for $M(E_7)\downarrow_H$ is $14,6^7$. There are no extensions between composition factors, so this is semisimple. The action on $M(E_7)$ of $u$ of order $3$ with the largest centralizer in $H$ is $2^{12},1^{32}$. This is the generic class $A_1$ by \cite[Table 7]{lawther1995}. We conclude from Lemma \ref{lem:genericunipotent} that $H$ is a blueprint for $M(E_7)$, and therefore so is any subgroup containing $H$, as needed.

For $p=5$ all of the same statements hold, except we only need traces of elements of order $3$ to prove that $\PSp_6(5)$ does not embed, and for $\Sp_6(5)$ elements of orders $2$ and $3$ suffice.

Finally, we consider (\ref{propi:medrd}). Since the semisimple elements have the same orders in $\Omega_7(p^a)$ as $\PSp_6(p^a)$, we again need only consider $p^a=3,5$. For $p^a=3$, the simple modules for $H=\Omega_7(3)$ of dimension at most $56$ are $1$, $7$, $27$ and $35$. The traces of elements of orders $2$ and $4$ are enough to find the unique conspicuous set of composition factors, $21^2,7^2$. and since there are no extensions between these modules $M(E_7)\downarrow_H$ is semisimple. An element $u$ of order $3$ with maximal centralizer size acts on this module with blocks $3^2,2^{16},1^{18}$, which is the generic class $2A_1$ by \cite[Table 7]{lawther1995}. Hence $u$ and $H$ are blueprints for $M(E_7)$ by Lemma \ref{lem:genericunipotent}. In the other case of $H=\Spin_7(3)$, a non-central involution in $H$ has trace $0$ on all faithful modules, and so since an involution in $E_7$ has trace $\pm 8$ from Appendix \ref{app:traces}, we cannot get this case.

The exact same proof works for $p=5$ except we use traces of elements of orders $2$ and $3$ to eliminate all but one set of composition factors.
\end{proof}

Now that these medium-rank subgroups have been removed from $\mathscr P$, we also need to remove $\PSL_2(p^a)$ for some large $p^a$. This is performed by Theorem \ref{thm:goodblueprint}: let $v(\bG)$ be defined by
\[ v(\bG)=\begin{cases}
18 & \bG=F_4,\; E_6,
\\75 & \bG=E_7.\end{cases}\]
If $p^a>v(\bG)\cdot \gcd(2,p-1)$ then $H=\PSL_2(p^a)$ does not lie in $\mathscr P$.

\begin{proposition}\label{prop:possiblemaximals}Let $H=\PSL_2(p^a)$ for $p^a\leq v(\bG)\cdot \gcd(2,p-1)$, and suppose that $H$ lies in $\mathscr P$. Then $H$ is a maximal member of $\mathscr P$ unless $p^a$ is one of the following:
\begin{itemize}
\item $G=F_4$, $p^a=4,5,7,8,9,11$;
\item $G=E_6$, $p^a=4,5,7,8,9,11,16,25$ \item $G=E_7$, $p^a=4,5,7,8,9,11,16,25,27,64,81$.
\end{itemize}
For $p^a=25$, if $\bG=E_6$ then $H$ is only contained in ${}^2\!F_4(2)'$ from $\mathscr P$, and if $\bG=E_7$ then $H$ is only contained in $Ru$ from $\mathscr P$. For $p^a=27$, $H$ is only contained in a Ree group ${}^2\!G_2(27)$ from $\mathscr P$.
\end{proposition}
\begin{proof} The possible members of $\mathscr P$ are: \begin{enumerate}
\item\label{li:mempa} $\PSL_2(p^b)$ for $p^b\leq v(\bG)\cdot \gcd(2,p-1)$;
\item\label{li:mempb} medium-rank subgroups not eliminated in Proposition \ref{prop:eliminatemediumrank};
\item\label{li:mempc} simple subgroups of exotic $r$-local subgroups;
\item\label{li:mempd} simple groups not of Lie type in characteristic $p$, given in \cite[Section 10]{liebeckseitz1999}.
\end{enumerate}
Suppose that $H$ is contained in another member $K$ of $\mathscr P$. If $K$ comes from (\ref{li:mempa}), then $K=\PSL_2(p^b)$ for $b\geq 2a$, and hence $p^{2a}\geq v(\bG)\cdot \gcd(2,p-1)$. Thus $p^a=4,5$ for $\bG=F_4,E_6$, and $p^a=4,5,7,8,11$ for $\bG=E_7$. Each of these appears in the proposition, so we may assume $K$ does not occur in (\ref{li:mempa}).

\medskip

Suppose that $K$ occurs in (\ref{li:mempb}). The groups $\PSL_3(p^b)$, $\PSU_3(p^b)$ and $G_2(p^b)$ only contain $\PSL_2(p^b)$, whence $p^a\leq 9$ if $K$ is one of these groups. The Suzuki groups cannot contain $H$ and so we can exclude these groups. The Ree group ${}^2\!G_2(3^{2n+1})$ contains $\PSL_2(3^{2n+1})$, so if $K$ is one of these groups then $H$ may only be $\PSL_2(27)$. Furthermore, ${}^2\!G_2(3^{2n+1})$ contains a real semisimple element of order $3^{2n+1}+3^n+1$. By Theorem \ref{thm:blueprintsminimal}, this means that $K$ lies in $\mathscr P$ for $\bG=F_4$ unless $n=1$, and for $\bG=E_7$ unless $n\leq 3$. For $\bG=E_6$, if $K$ is not irreducible on $M(E_6)$ then $n\leq 1$, and if $K$ is irreducible on $M(E_6)$ then $K$ is strongly imprimitive by \cite{liebeckseitz2004a}.

If $p=2$ then the only case not on our list is $\PSL_2(32)$ for $\bG=E_7$, and this cannot be contained in any medium-rank subgroup, so assume that $p$ is odd.

If $p=3$ then we need to consider $\PSL_4(3^a)$ and $\PSU_4(3^a)$ for $a=1,2$, as $\PSp_4(3^a)$ is contained in $\PSL_4(3^a)$. These contain $\PSL_2(3^a)$ and $\PSL_2(3^{2a})$, so we obtain $3^a=3,9,81$. If this is less than $2\cdot v(\bG)$ then it appears on our list.

Thus $p=5,7$, but now there are no more medium-rank subgroups to consider, so this deals with those $K$ in (\ref{li:mempb}).

\medskip

Suppose that $K$ occurs in (\ref{li:mempc}). The exotic $r$-local subgroups for $\bG\neq E_8$ have composition factors either cyclic groups, or $\SL_3(2)$ ($G_2$ and above), $\SL_3(3)$ ($F_4$ and above) and $\SL_3(5)$ ($E_6$ and above); the first two are minimal simple groups anyway, and the third contains only $\SL_2(4)$, so we get $p^a=4$ for $E_6,E_7$. 

\medskip

Suppose that $K$ occurs in (\ref{li:mempd}). If $K$ is alternating then $K$ is either $\Alt(6)$ or $\Alt(7)$ by \cite[Theorem 1]{craven2015un2}. These contain $\PSL_2(p^a)$ for $p^a=4,5,7,9$, and these appear on all three lists. If $K=\PSL_2(r^b)$ then $H$ must be $\PSL_2(4),\PSL_2(5)$, so this case is easily considered.

If $K$ is Lie type but not $\PSL_2$ then very few may occur, and even fewer that contain a copy of $\PSL_2(p^a)$ that we have not already seen. We use \cite{atlas} to find subgroups of these simple groups. Table \ref{t:simplecontainspsl2} lists the $p^a$ such that $K=K(r^b)$ contains $\PSL_2(p^a)$, $r\neq p$, and $K$ exists in $\bG$ one of $F_4$, $E_6$ and $E_7$ in characteristic $p$.
\begin{table}
\begin{center}
\begin{tabular}{cccc}
\hline Group & Prime powers & Group & Prime powers
\\ \hline $\PSL_3(3)$& - & $\POmega_7(3)$ ($p=2$) & $4,8$
\\ $\PSL_3(4)$ & $5,7,9$ & $G_2(3)$ ($p=2$) & $8$
\\ $\PSU_3(3)$& $7$ & ${}^2\!B_2(32)$ ($p=5$) & $5$
\\ $\PSU_3(8)$& - & $M_{11}$ & $4,5,9,11$
\\ $\PSU_4(2)=\PSp_4(3)$&$5$ & $M_{12}$ & $4,5,9,11$
\\ $\PSp_6(2)$&$5,7,9$ & $J_2$ & $4,5,7,9$
\\ $\Omega_8^+(2)$&$5,7,9$ & $M_{22}$ ($p=2,5$) & $4,5$
\\ ${}^3\!D_4(2)$& $7$ & $J_1$ ($p=11$) & $11$
\\ ${}^2\!F_4(2)'$& $5,9,25$ & $Ru$ ($p=5$) & $5,25$
\\ $\PSL_4(3)$ ($p=2$) & $4$ & $HS$ ($p=5$) & $5$
\\ $\PSU_4(3)$ ($p=2$) & $4$ & $Fi_{22}$ ($p=2$) & $4,8$
\\ \hline
\end{tabular}
\end{center}
\caption{The prime powers $p^a$ for which $\PSL_2(p^a)$ is a cross-characteristic subgroup of various simple groups involved in $E_7(k)$} \label{t:simplecontainspsl2}
\end{table}
From \cite{liebeckseitz1999}, we see that $Ru$ only embeds in $E_7$, and ${}^2\!F_4(2)'$ in $E_6$ and $E_7$. The result is complete, except for why ${}^2\!F_4(2)'$ does not lie in $\mathscr P$ for $\bG=E_7$. We see from \cite[Table 6.231]{litterickmemoir} that ${}^2\!F_4(2)'$ only embeds in $E_7$ with two trivial composition factors on $M(E_7)$, hence it has pressure $0$ and therefore stabilizes a line on $M(E_7)$ by Lemma \ref{lem:pressure}. We now apply Proposition \ref{prop:fixlineonMG} to see that $K$ is strongly imprimitive, and therefore ${}^2\!F_4(2)'$ does not lie in $\mathscr P$ for $\bG=E_7$.
\end{proof}

Propositions \ref{prop:maximalnotinP} and \ref{prop:possiblemaximals} together mean that we may prove that $\PSL_2(p^a)$ is strongly imprimitive, at least for certain $p^a$, just by exhibiting a larger subgroup that stabilizes a proper subspace of an irreducible module for $\bG$. This will come in useful for $\bG=E_7$ in particular.

\chapter{Maximal Subgroups and Subspace Stabilizers}
\label{ch:maxsubspace}

The purpose of this chapter is to prove some sufficiency criteria for a subgroup to be strongly imprimitive, essentially making upgrading `Lie imprimitive' to `strongly imprimitive' a formal process in many cases.

Lemma \ref{lem:smallsubspaces} stated that if a subgroup $H$ stabilizes a line or hyperplane on $M(\bG)$ or $L(\bG)^\circ$ then $H$ is Lie imprimitive. Propositions \ref{prop:fixlineonLG} and \ref{prop:fixlineonMG} generalize this to show that $H$ is strongly imprimitive, but we have to replace `stabilize' by `centralize', a condition that is of course equivalent if $H$ is a simple group. We also show the same result for certain $H$ stabilizing a $2$-space on $M(\bG)$ in Proposition \ref{prop:fix2spaceonMG}. If a subgroup $H$ is a blueprint for a special class of module, usually including $M(\bG)$ or $L(\bG)^\circ$, then we again can use this to prove strong imprimitivity, see Proposition \ref{prop:blueprintissi}. The general criterion is Proposition \ref{prop:intersectionstabilizers}, which is reproduced from \cite[Section 4.2]{litterickmemoir}, itself generalizing work of \cite[Proposition 1.12]{liebeckseitz1998}.

These results will be used in Chapters \ref{ch:f4} to \ref{ch:e7oddsl} to prove strong imprimitivity given information about subgroups $H$ of $\bG$.

\medskip

We recall some definitions: a subgroup $H$ of $\bG$ is \emph{$\bG$-completely reducible} if, whenever $H$ is contained in a parabolic subgroup of $\bG$, it is contained in the corresponding Levi subgroup; it is \emph{$\bG$-irreducible} if $H$ is not contained in any parabolic subgroup of $\bG$; as we have seen before, it is \emph{Lie primitive} if it is not contained in any proper positive-dimensional subgroup of $\bG$.

The first result that we need is by Liebeck, Martin and Shalev \cite[Proposition 2.2 and Remark 2.4]{liebeckmartinshalev2005}, and this is translated into our notation.

\begin{proposition}\label{prop:paraisgood} Let $H$ be a finite subgroup of $\bG$ such that $H$ is not $\bG$-completely reducible and $H\leq \bG^\sigma$. There exists a proper parabolic subgroup $\bP$ of $\bG$ such that $H$ is contained in $\bP$ and $\bP$ is both $\sigma$-stable and $N_{\Aut^+(\bG)}(H)$-stable. In particular, $H$ is strongly imprimitive.
\end{proposition}

If $R$ is an elementary abelian $r$-group such that $N_{\bG}(R)$ is an exotic $r$-local subgroup then of course $R$ is contained in the normalizer of a torus but $N_{\bG}(R)$ is Lie imprimitive, so that $R$ is Lie imprimitive but not strongly imprimitive. Thus one cannot prove that every Lie imprimitive subgroup is strongly imprimitive. Of course, $R$ stabilizes $1$-spaces and $2$-spaces on all non-trivial simple modules, in particular the minimal and adjoint modules, so stabilizing such 
subspaces cannot be enough to conclude that a subgroup is strongly imprimitive.

\begin{proposition}\label{prop:Gcrisgood} Suppose that $H$ is a finite subgroup of $\bG$. If $H$ is $\bG$-reducible then $H$ is strongly imprimitive.
\end{proposition}
\begin{proof} Suppose that $H$ is not strongly imprimitive, and that $H$ is $\bG$-reducible. By Proposition \ref{prop:paraisgood}, $H$ is $\bG$-completely reducible, so $H$ is contained in a proper ($\sigma$-stable) Levi subgroup $\bL$. Since $C_{\bG}(\bL)$ contains a torus $\bT$ of $\bG$, $\bY=C_{\bG}(H)^\circ\neq 1$. Note that $\bY$ is $N_{\Aut^+(\bG)}(H)$-stable, since $H$ is, and therefore $H\cdot \bY$ is a $N_{\Aut^+(\bG)}(H)$-stable, positive-dimensional subgroup containing $H$. Since $H$ is normal in this subgroup, it is proper, and therefore $H$ is strongly imprimitive.
\end{proof}

Thus if $H$ is a subgroup of $\bG^\sigma$ such that $N_{\bar G}(H)$ is almost simple, $H$ is $\bG$-irreducible.

The next result comes entirely from \cite[Chapter 4]{litterickmemoir}, particularly \cite[Corollary 4.5]{litterickmemoir}, and is a direct generalization of the methods of \cite[Proposition 1.12]{liebeckseitz1998}. If $W$ is a subspace of a $\bG$-module $V$, write $\bG_W$ for the stabilizer of $W$ in $\bG$, and if $\mathcal W$ is a collection of subspaces of $V$, write $\bG_{\mathcal W}$ for the intersection of the $\bG_W$ for $W\in \mathcal W$. Note that this also applies for $\bG=F_4$ and $p=2$, where the `graph morphism' interchanges $\lambda_i$ and $\lambda_{5-i}$. 

\begin{proposition}\label{prop:intersectionstabilizers} Let $\bG$ be a simple algebraic group over an algebraically closed field, and let $V$ be a semisimple module such that the highest weights of the composition factors are stable under graph morphisms of $\bG$. Let $H$ be a finite subgroup of $\bG$.

If $\phi\in N_{\Aut^+(\bG)}(H)$ then $\phi$ permutes the $H$-invariant subspaces of $V$. Hence if $\mathcal W$ is an orbit of $H$-invariant subspaces of $V$, then $H\leq \bG_{\mathcal W}$ and $\bG_{\mathcal W}$ is $\phi$-stable.

If $\mathcal W$ is one of
\begin{enumerate}
\item all $H$-invariant subspaces of $V$,
\item all $H$-invariant subspaces of $V$ whose dimension lies in some subset $I\subseteq \N$,
\item all simple $H$-invariant subspaces of $V$,
\item all simple $H$-invariant subspaces of $V$ whose dimension lies in some subset $I\subseteq \N$,
\end{enumerate}
then $\mathcal W$ is a union of orbits of $H$-invariant subspaces of $V$. In particular, $\bG_{\mathcal W}$ is a subgroup of $\bG$ containing $H$ and stable under every element of $N_{\Aut^+(\bG)}(H)$.

Hence if $\bG_{\mathcal W}$ is positive dimensional for one of these sets $\mathcal W$ of subspaces, not stabilized by $\bG$ itself, then $H$ is strongly imprimitive.
\end{proposition}

This also gives us the following result, applying the statement to blueprints.

\begin{proposition}
\label{prop:blueprintissi} Suppose that $H$ is a finite subgroup of $\bG$, and let $V$ be a semisimple module that is $N_{\Aut^+(\bG)}(H)$-stable, i.e., $V$ is a sum of highest-weight modules whose weights form orbits under the action of the graph morphisms in $N_{\Aut^+(\bG)}(H)$. If $H$ is a blueprint for $V$, then either $H$ is strongly imprimitive or $H$ acts irreducibly on every composition factor of $V$.

In particular, if $H$ is either $\PSL_2(p^a)$ or $\SL_2(p^a)$ in $\bG=E_7$ with $Z(H)=Z(\bG)$, and $H$ is a blueprint on either $M(\bG)$ or $L(\bG)^\circ$, then $H$ is strongly imprimitive.
\end{proposition}
\begin{proof} The first part is a direct consequence of Proposition \ref{prop:intersectionstabilizers}. The second statement follows from the first if $H$ cannot act irreducibly on either $M(\bG)$ or $L(\bG)^\circ$. But this follows from \cite{liebeckseitz2004a}.
\end{proof}

We now prove results about stabilizing lines on $L(\bG)$ and $M(\bG)$, and $2$-spaces on $M(\bG)$.

\begin{proposition}\label{prop:fixlineonLG} If $H$ is a finite subgroup of $\bG$ such that $H$ centralizes a line on $L(\bG)^\circ$, then $H$ is strongly imprimitive.
\end{proposition}
\begin{proof} Suppose that the fixed-point subspace of $H$ on $L(\bG)$ is $W\neq 0$, and we consider $\bG$ to be of adjoint type, so that $L(\bG)$ cannot have a trivial submodule \cite[Lemma 2.1.1]{liebeckseitz2004}. By \cite[(1.3)]{seitz1991}, if $W$ contains a nilpotent element then $H$ is contained in a proper parabolic subgroup, hence by Proposition \ref{prop:Gcrisgood} $H$ is strongly imprimitive.

Thus $W$ consists of semisimple elements. Indeed, since $W$ is the space of fixed points of $H$ on $L(\bG)$, it is a subalgebra of $L(\bG)$. This means that $W$ is contained in a maximal torus of $L(\bG)$, hence (as in \cite[(1.3)]{seitz1991}) $\bC=C_{\bG}(W)$ contains a maximal torus of $\bG$. Since $H\leq \bC$, this means that $\bC$ is a maximal-rank subgroup of $\bG$, and is $\sigma$-stable since $W$ is. To show invariance under $N_{\bar G}(H)$, we now apply Proposition \ref{prop:intersectionstabilizers}, since $\bC$ is the intersection of the centralizers of all $H$-centralized $1$-spaces of $L(\bG)$. This shows that $\bC$ is $N_{\Aut^+(\bG)}(H)$-stable, so $H$ is strongly imprimitive.
\end{proof}

We now do the same thing with $M(\bG)$ in place of $L(\bG)$. We add a small condition on $H$ that is certainly satisfied for simple groups. This can probably be removed, at the expense of making the proof more complicated.

\begin{proposition}\label{prop:fixlineonMG} If $H$ is a finite subgroup of $\bG$ such that $H$ has no subgroup of index $2$, and $H$ centralizes a line on $M(\bG)$, then $H$ is strongly imprimitive.
\end{proposition}
\begin{proof} We may assume that $H$ is $\bG$-irreducible by Proposition \ref{prop:Gcrisgood}. Suppose that the fixed-point subspace $M(\bG)^H$ of $H$ on $M(\bG)$ is non-zero. If $\bG=F_4$ and $p=2$ then as $M(F_4)$ is a submodule of $L(F_4)$, we get the result in this case by Proposition \ref{prop:fixlineonLG}. If $p$ is odd, then $H$ is contained in a line stabilizer of $M(F_4)$, which from \cite[Lemma 2.2(iii)]{liebeckseitz2005} are contained in either a maximal parabolic or $B_4$. The first case cannot hold by assumption, so $H$ is contained in $B_4$.

If $M(\bG)^H$ has dimension $1$ or $2$ then the centralizer of $M(\bG)^H$ is positive dimensional, and hence we may apply Proposition \ref{prop:intersectionstabilizers} to obtain the result, as in the previous proposition. Thus $M(\bG)^H$ has dimension at least $3$, and so centralizes a $2$-space on the sum of the natural and spin modules for $B_4$. Note that $H$ cannot lie in a parabolic subgroup of $B_4$ since it does not lie in a parabolic subgroup of $F_4$, and hence $H$ cannot lie in a $B_3$ subgroup of $B_4$, which is the $2$-space centralizer on the natural module. On the other hand, since the spin module for $B_4$ appears in both $M(F_4)$ and $L(F_4)$, if $H$ centralizes a line on the spin module then $H$ centralizes a line on $L(F_4)$, and we apply Proposition \ref{prop:fixlineonLG}.

If $G=E_6$ then from Lemma \ref{lem:e6stabs} the line stabilizers are $F_4$ or subgroups contained in parabolic subgroups, so $H\leq F_4$. If $H$ centralizes a point on $M(F_4)$ then we apply the previous paragraph, and if $H$ does not then $H$ centralizes a unique line on $M(E_6)$, whence we apply Proposition \ref{prop:intersectionstabilizers} again.

If $G=E_7$ then from Lemma \ref{lem:e7stabs}, we may assume that $H$ only centralizes lines whose stabilizers are of the form $E_6.2$, as all others are contained in maximal parabolic subgroups and $H$ is $\bG$-irreducible. However, since $H$ has no subgroup of index $2$, this means that $H$ is contained in an $E_6$-Levi subgroup, but is $\bG$-irreducible. This completes the proof.
\end{proof}

We also will have cause to use $2$-space stabilizers on $M(\bG)$ for $\bG=F_4,E_6,E_7$. By dimension counting (see Lemma \ref{lem:smallsubspaces}) we see that the stabilizer of such a space is positive dimensional so $H$ is imprimitive, but we need to show that $H$ is strongly imprimitive. Because the proof of this goes case by case, we only prove the precise cases that we need, rather than produce a general method. Even with this restriction, the proof is long and technical.

\begin{proposition}\label{prop:fix2spaceonMG} If $H=\SL_2(p^a)$ is a subgroup of $\bG$ for $\bG=F_4,E_6,E_7$ and $p=2$, or $\bG=E_7$ and $p^a=5,7,9,25,27,49$, and $H$ stabilizes a $2$-space on $M(\bG)$, then $H$ is strongly imprimitive.
\end{proposition}
\begin{proof} If $p$ is odd then any group that normalizes $H$ centralizes $Z(H)$, so either lies in the centralizer of an involution -- hence $H$ is strongly imprimitive -- or $Z(H)=Z(\bG)$. Thus for $p$ odd we assume that $Z(H)=Z(\bG)$. In addition, if $H$ is $\bG$-reducible then $H$ is strongly imprimitive by Proposition \ref{prop:Gcrisgood}, so we may assume that $H$ is $\bG$-irreducible. By Lemma \ref{lem:smallsubspaces}, $H$ is not Lie primitive.

Let $\bX$ be a minimal $\bG$-irreducible subgroup of $\bG$ containing $H$, so that $H$ is Lie primitive in $\bX$. If $\bX$ is a product of groups then the projections along each must also be Lie primitive, else we could replace $\bX$ with a smaller subgroup. There are no Lie primitive subgroups $\SL_2(p^a)$ in classical groups or $G_2$ by Proposition \ref{prop:morphismextensiontrueforclassical}, so if $\bX$ has a classical or $G_2$ component then this is $A_1$. Thus $\bX$ can only have components $A_1$, $F_4$ and $E_6$. Furthermore, since there are no Lie primitive copies of $H$ in a product of $A_1$s, there is at most one $A_1$ component in $\bX$. Note that, of course, $\bX$ is not necessarily $\sigma$-stable, never mind $N_{\Aut^+(\bG)}(H)$-stable.

\medskip

If $\bG=F_4$ (and hence $p=2$) then $\bX$ is an $A_1$ subgroup. If $\bX\leq B_4$ then $\bX$ and hence $H$ stabilizes a line on $M(F_4)$, whence $H$ is strongly imprimitive by Proposition \ref{prop:fixlineonMG}. However, $\bX\leq B_4$ by \cite[Table 5]{thomas2016}, completing the proof for $F_4$.

\medskip

If $\bG=E_6$ (and again $p=2$), then instead of this analysis, which would also work, we use the tables of Cohen--Cooperstein \cite{cohencooperstein1988} to determine the possibilities for the $2$-space $W$ being stabilized, and see that in all cases $C_{\bG}(W)$ must have a unipotent radical, hence $H$ lies in a parabolic, hence is strongly imprimitive by Proposition \ref{prop:Gcrisgood}.

\medskip

For $\bG=E_7$ and all primes $p$, then as $\bX\neq \bG$, either $\bX=A_1$, or it involves $F_4$ or $E_6$. If $\bX$ involves $E_6$ then $\bX=E_6$ is a Levi subgroup, and if $\bX=F_4$ then $\bX\leq E_6$, so again $H$ is not $\bG$-irreducible. Thus $\bX$ is either $F_4A_1$ or $A_1$. 

If $\bX=F_4A_1$ then the projection of $H$ along $F_4$ must be Lie primitive in $F_4$, else $\bX$ is contained in a proper positive-dimensional subgroup of $F_4A_1$, and these are all classical, so $\bX=A_1$. In this case, the best way to proceed is to use Theorem \ref{thm:f4}, which we prove in Chapter \ref{ch:f4}, which states that there are no Lie primitive copies of $\PSL_2(p^a)$ in $F_4$. Using this result is valid and not circular reasoning, since we will not need to use this result for $E_7$ until Chapter \ref{ch:e7oddsl}, by which time Theorem \ref{thm:f4} will have been proved.

Thus for the rest of the proof, $\bX$ is an irreducible $A_1$ subgroup of $E_7$. If $p=2$ then $\bX\leq D_6A_1$ by \cite[Table 7]{thomas2016}, so since $D_6A_1$ stabilizes a line on $L(E_7)^\circ$, $H$ is strongly imprimitive by Proposition \ref{prop:fixlineonLG}. Thus $p$ is odd.

The information given in \cite[Tables 7 and 12]{thomas2016} is significant, but does not include the (non-trivial) unipotent class intersecting $\bX$, and we also need the module action of $H$ on $M(E_7)$ in certain cases.

The conditions on the embeddings of the $A_1$s in \cite{thomas2016} are important, as if those conditions are not satisfied then the subgroup is not $\bG$-irreducible. If $H$ is contained in a diagonal $A_1$ subgroup of a product of $A_1$s, then $H$ is contained in other diagonal $A_1$s, obtained by twisting the representations of the projections along the various $A_1$s. Every one of the diagonal $A_1$s containing $H$ must be $\bG$-irreducible.

The effect of this is, for each class of $A_1$ in the tables, we require that the conditions are satisfied for integers lying between $0$ and $a-1$, else we may replace one, say $i$, by $i-a$, and yield another $A_1$ subgroup containing $H$. The consequence is that $a$ needs to be large enough so that the conditions can still hold.

For example, if $p=3$ then from \cite[Table 7]{thomas2016} we see that $A_1$ subgroups with labels $3$, $10$, $11$ and $12$ exist. For subgroup $11$, the centre of it does not coincide with the centre of $\bG$, so this can be excluded for all primes. Subgroup $3$ requires $a\geq 3$ for the conditions to be satisfied, subgroup $10$ requires $a\geq 2$, and subgroup $12$ cannot occur for $a\leq 3$.

In what follows, we use the notation for simple modules for $\SL_2$ from Chapter \ref{ch:sl2modules}, because it is much easier to understand the dimensions of the composition factors than using the notation $L(\lambda)$. We will show that the simultaneous stabilizer of all $2$-dimensional submodules of $M(E_7)\downarrow_H$ is positive dimensional, and therefore apply Proposition \ref{prop:intersectionstabilizers}, as we did in Propositions \ref{prop:fixlineonLG} and \ref{prop:fixlineonMG}, to prove that $H$ is strongly imprimitive, completing the proof of the result. Note that if $H$ stabilizes a unique $2$-space on $M(E_7)$ then the stabilizer is positive dimensional by Lemma \ref{lem:smallsubspaces}, so we are done in this case.

\medskip

We first set $p=3$, and so $a=2,3$. The subgroup $\bX$ must be either subgroup $3$ or subgroup $10$ from \cite[Table 7]{thomas2017un}, both lying inside $A_1D_6$. For subgroup $3$ we require $a=3$, and may assume that the action of $H$ along $M(D_6)$ is $L(1)\otimes L(1)^{[1]}\otimes L(2)^{[2]}$ (where $[i]$ denotes $i$ field twists, so $L(1)^{[1]}=L(3)$), and the action along $A_1$ is $L(1)^{[i]}$ for $i=0,1,2$. Thus in the notation of \cite{thomas2016}, $r=0,1,2$, ${u,t}={0,1}$ (both possibilities are allowed), $s=2$, and in our notation from Chapter \ref{ch:sl2modules} the two modules are $12_{1,2,3}$ and $2_i$ for $i=1,2,3$.

Since $\bX$ lies in $A_1D_6$, $\bX$ acts on $M(E_7)$ as the sum of (the restrictions to $\bX$ of) $(M(A_1),M(D_6))$ and $(L(0),L(\lambda_5))$. The three possibilities for the action of $H$ on $(M(A_1),M(D_6))$ are
\[ 6_{2,3}\oplus 18_{2,1,3},\quad 6_{1,3}\oplus 18_{1,2,3},\quad 8/2_2,6_{2,1}/8.\]
The other summand of $M(E_7)\downarrow_H$ does not depend on $i$, but does depend on which of $u$ and $t$ is $0$ and $1$. From \cite[Table 12]{thomas2016} we can read off the factors of $M(E_7)\downarrow_{\bX}$ and hence of $M(E_7)\downarrow_H$. These are
\[ 18_{2,1,3},8,2_2^2,2_3\quad\text{and}\quad 18_{1,2,3},6_{3,1},2_1^2,2_2,2_3\]
according as $t=0$ and $t=1$ respectively. The (non-trivial) unipotent element in $H$ belongs to class $2A_2+A_1$ by \cite[4.10]{lawther2009}, which acts on $M(E_7)$ with Jordan blocks $3^{18},1^2$ by \cite[Table 7]{lawther1995}. In particular, since there are no blocks of size $2$, there are no $2$-dimensional summands of $M(E_7)\downarrow_H$, and since $2_2$ and $2_1$ respectively are the only composition factors to appear more than once, they must lie in the socle and top of $M(E_7)\downarrow_H$. In particular, $H$ stabilizes a unique $2$-space, so we are done.

For subgroup $10$, we can have either $p^a=9$ or $p^a=27$. If $p^a=9$, then the parameters of $\bX$ from \cite{thomas2016} are $s=0$, $t=1$, and $r,u=0,1$, up to field automorphism. The unipotent class of $E_7$ to which a non-trivial unipotent element of $H$ belongs is again $2A_2+A_1$, so again there are no $2$-dimensional summands of $M(E_7)\downarrow_H$. The restriction of $(M(A_1),M(D_6))$ to $H$  is the sum of the modules $2_r\otimes 9$ and $2_{r+1}\otimes 3_{u+1}$. The first of these is the projective cover of a $6$-dimensional module, and the second of these contains a $2$-dimensional submodule if and only if $r=u$, in which case it is of the form $2_{u+1}/2_{u+2}/2_{u+1}$. The restriction of $(L(0),L(\lambda_5))$ to $H$ has the form
\[(4/1,1,3_1,3_2/4)\otimes 2_{u+1}.\]
The submodule $1,1,3_1,3_2/4$ is uniquely defined, and its product with $2_{u+1}$ has a single $2_{u+2}$ in the socle and no $2_{u+1}$. The top of this tensor product is $4\otimes 2_{u+1}=2_{u+2}\oplus 6_{u+2,u+1}$, but the $2_{u+2}$ cannot become a submodule for then it would be a summand, and there are no $2$-dimensional summands. Thus the $2$-dimensional factors of the socle of $M(E_7)\downarrow_H$ are $2_{u+2}$ if $r\neq u$, and $2_{u+1}\oplus 2_{u+2}$ if $r=u$.

In the former case $H$ stabilizes a unique $2$-space on $M(E_7)$ so we are done, so assume the latter holds. Notice that the irreducible $A_1$ subgroup $\bX$ with the parameters above also stabilizes the $2_{u+1}$ and the $2_{u+2}$, so the simultaneous stabilizer of both contains $\bX$ and we are done.

For $p^a=27$, we need $s<t$, and so $s=0$, $t=1$ or $s=0$, $t=2$ up to field automorphism, with $0\leq u,r\leq 2$. Suppose first that $s=0$ and $t=1$; this is very similar to the case for $p^a=9$. The restriction of $(M(A_1),M(D_6))$ to $H$ contains the modules $2_r\otimes 3_s\otimes 3_t$ and $2_r\otimes 3_u$, and the only way to get a $2$-dimensional submodule is if $r=u$, in which case we get the summand $2_{u+1}/2_{u+2}/2_{u+1}$. For the contribution from $(L(0),L(\lambda_5))$, the structure is
\[ (4_{1,2}/1,3_2,4_{1,3}/4_{1,2})\otimes 2_{u+1}.\]
This tensor product has at most one copy of a $2$-dimensional module in the socle, and the subgroup $\bX$ stabilizes each of these potential $2$-dimensional submodules if they exist, and we are done.

If $s=0$ and $t=2$ instead then the same statements for $(M(A_1),M(D_6))$ hold. The structure of $(L(0),L(\lambda_5))$ also very similar, being
\[ (4_{1,3}/1,3_1,4_{2,3}/4_{1,3})\otimes 2_{u+1}.\]
The same statements about $2$-dimensional submodules hold as well, proving the result in this case.

\medskip

We now set $p=5$. For $p^a=5$, we see from \cite[Table 7]{thomas2016} that there are no $\bG$-irreducible subgroups $H$. For $p^a=25$, we can satisfy the conditions for subgroups $8$, $9$, $10$ and $19$ from Thomas's list. The first three lie in $A_1D_6$, the last in $A_1A_1$.

The first case is subgroup $8$: here $s\neq t$ and $u\leq v$. By field automorphism we may assume that $s=0$ and $t=1$, and $u,v=1,2$. The action of $H$ on the module $(M(A_1),M(D_6))$ is
\[ 2_{r+1}\otimes (5_1\oplus 3_2\oplus 2_{u+1}\otimes 2_{v+1}).\]
Any $2$-dimensional submodule of this is stabilized by $\bX$, so we consider the restriction of $(L(0),L(\lambda_5))$ to $\bX$ and $H$. This is
\[ 4_1\otimes 2_2\otimes(2_{u+1}\oplus 2_{v+1}),\]
and this yields $2$-dimensional submodules if $u=1$ or $v=1$, again stabilized by $\bX$ as well. Hence $H$ and $\bX$ stabilize the same $2$-dimensional subspaces of $M(E_7)$, so we are done.

For subgroup $9$, we may take $s=0$, $t=u=1$, and $r=1,2$. The composition factors of the restriction of $(L(0),L(\lambda_5))$ to $H$ are $10_{2,1}^2,6_{2,1}^2$, so we may concentrate on $(M(A_1),M(D_6))$. The action of both $H$ and $\bX$ on this is
\[ 2_{r+1}\otimes (5_1\oplus 3_1\oplus 2_2^{\otimes 2}).\]
Both the $H$-action and $\bX$-action have two $2$-dimensional composition factors in the socle for either $r=1$ or $r=2$, so again the intersection of the stabilizers of all $2$-dimensional submodules of $M(E_7)\downarrow_H$ is positive dimensional.

For subgroup $10$, here we may take $s=0$, $t=1$ and $r,u=1,2$. We have four composition factors of $M(E_7)\downarrow_{\bX}$, and the restrictions to $H$ of these are:
\begin{itemize}
\item $2_{r+1}\otimes 3_1\otimes 3_2$, which is $(2_{r+1}\oplus 4_{r+1})\otimes 3_{r+2}$;
\item $2_{r+1}\otimes 3_{u+1}$, which is $6_{u+1,r+1}$ if $u\neq r$ and $2_{u+1}\oplus 4_{u+1}$ if $u=r$;
\item $4_1\otimes 2_2\otimes 2_{u+1}$ and $4_2\otimes 2_1\otimes 2_{u+1}$, and these two are $(3_{u+1}\oplus 5_{u+1})\otimes 2_{u+2}$ and $(1\oplus 3_{u+1})\otimes 4_{u+2}$.
\end{itemize}
We see that there is a unique $2$-dimensional submodule if $u=r$, and none otherwise, and so we are done.

If the subgroup is $19$, then the action of $\SL_2(5)$ on $M(E_7)$ is
\[P(4)^{\oplus 4}\oplus (2/2,4)\oplus (2,4/2),\]
and one of the $2$-dimensional submodules is distinguishable as being contained in a module $2/2$. Thus either $H$ stabilizes a unique $2$-space, and we are done by Proposition \ref{prop:intersectionstabilizers}, or it stabilizes the same number of subspaces as the copy of $L=\SL_2(5)$ inside it, in which case it stabilizes the subspace that must lie in a singleton orbit, and again we are done.

\medskip

The final case is $p=7$. In this case, we simply consider all subgroups from Thomas's list that occur for $p=7$, and tabulate their restrictions to $H=\SL_2(7)$.
\begin{center}\begin{tabular}{ccc}
\hline Group Number & Unipotent class & Action on $M(E_7)$
\\ \hline $1$, $2$, $6$, $16$, $19$ & $E_7(a_5)$ & $P(6)^{\oplus 2}\oplus P(4)\oplus 6\oplus 4^{\oplus 2}$
\\ $3$, $8$, $9$, $10$ & $A_3+A_2+A_1$ & $6^{\oplus 3}\oplus 4^{\oplus 7}\oplus 2^{\oplus 5}$
\\ $7$ & $D_5(a_1)+A_1$ & $P(6)^{\oplus 3}\oplus 4\oplus 2^{\oplus 5}$
\\ $12$ & $A_2+3A_1$ & $4^{\oplus 7}\oplus 2^{\oplus 14}$
\\ $15$, $17$ & $A_6$ & $P(6)\oplus P(4)^{\oplus 2}\oplus P(2)$
\\ \hline
\end{tabular}\end{center}
We exclude those with class $E_7(a_5)$ as $H$ does not stabilize a $2$-space in this case. If the unipotent elements of $H$ come from class $A_6$ then $H$ stabilizes a unique $2$-space on $M(E_7)$, and we are done. From \cite[Table 7]{thomas2016} we see that subgroup 7 is $\bG$-reducible unless $a\geq 3$, so this case cannot occur. The remaining two unipotent classes, $A_3+A_2+A_1$ and $A_2+3A_1$, are generic for $M(E_7)$ in the sense of Definition \ref{defn:genericunipotent} below. Lemma \ref{lem:genericunipotent} states that subgroups containing generic unipotent elements for $M(E_7)$ are strongly imprimitive, and this completes the proof of the result.
\end{proof}

\chapter{Blueprint Theorems for Semisimple Elements}
\label{ch:semisimple}

In this chapter, we will consider analogues of the bounds given in \cite{liebeckseitz1998} for a semisimple element $x$ of an exceptional algebraic group $\bG$ to be a blueprint for $L(\bG)$, by producing bounds for $x$ to be a blueprint for $M(\bG)$, or $M(\bG)\oplus M(\bG)^\tau$ if a graph automorphism $\tau$ lies in $\Aut^+(\bG)$ and $M(\bG)$ is not $\tau$-stable, as in Proposition \ref{prop:intersectionstabilizers}. In particular, we can take $M(\bG)$ if $\bG=E_7$ or $p$ is odd and $\bG=F_4$, we can take $M(E_6)\oplus M(E_6)^*$ if $\bG=E_6$, and we take $M(F_4)\oplus M(F_4)^\tau$ if $\bG=F_4$ and $p=2$, where $\tau$ is a graph automorphism. (The last module in this list is $L(\lambda_1)\oplus L(\lambda_4)$, which has the same composition factors as $L(F_4)$.) Write $\overline{M(\bG)}$ for $M(\bG)$ if $M(\bG)$ is stable under graph morphisms of $\bG$, and $M(\bG)\oplus M(\bG)^\tau$ for the two cases given above.

Thus let $\bG$ be a simply connected, simple algebraic group of type $G_2$, $F_4$, $E_6$, or $E_7$. In \cite{liebeckseitz1998}, the constant $t(\bG)$ is introduced, and \cite{lawther2014} produces a set $T(\bG)$ of positive integers with $t(\bG)=\max(T(\bG))$. The set $T(\bG)$ is split into odd and even integers, and is defined to be
\[ T(\bG)=\begin{cases}
\{1,3,\dots,9\}\cup\{2,4\dots,12\}&\bG=G_2,
\\\{1,3,\dots,57\}\cup\{2,4\dots,68\}&\bG=F_4,
\\\{1,3,\dots,105\}\cup\{2,4\dots,120,124\}&\bG=E_6,
\\\{1,3,\dots,317\}\cup\{2,4\dots,364,370,372,388\}&\bG=E_7.
\end{cases}\]
We then have the following theorem.

\begin{theorem}[Liebeck--Seitz, Lawther]\label{thm:lls} Suppose that $x$ is a semisimple element of $\bG$, and that the order of $x$ does not lie in $T(\bG)$. There exists a positive-dimensional subgroup $\bX$ of $\bG$, containing $x$, such that $\bX$ and $x$ stabilize the same subspaces of $L(\bG)$, i.e., $x$ is a blueprint for $L(\bG)$.
\end{theorem}

We show in Section \ref{sec:appprelims} the easy result that if $V$ is any rational $k\bG$-module then there is an analogous finite set $\mathcal X_V$ of integers such that if $x$ is a semisimple element of order not in $\mathcal X_V$, then $x$ is a blueprint for $V$.

The section afterwards computes $\mathcal X_V$, or at least gives bounds on it, for $V=M(\bG)$, which are independent of the underlying characteristic of the group $\bG$. However, for applications to maximal subgroups, since $M(\bG)$ is not always stable under graph automorphisms of $\bG$, we replace $M(\bG)$ by $\overline{M(\bG)}$, the direct sum of $M(\bG)$ and its image under the graph automorphism (if there is a graph automorphism). This module satisfies the hypothesis of Proposition \ref{prop:blueprintissi}, so if $x$ is a blueprint for $\overline{M(\bG)}$ and $H$ is a finite subgroup of $\bG$ containing $x$, then either $H$ is strongly imprimitive in $\bG$ or $H$ acts irreducibly on both $M(\bG)$ and $M(\bG)^\tau$. (This is no extra condition if $\bG=E_6$, as clearly $H$ acts irreducibly on $M(E_6)^*$ whenever it acts irreducibly on $M(E_6)$, but it is an extra condition if $\bG=F_4$ in characteristic $2$.)

The final part of this chapter applies these results to produce Theorem \ref{thm:goodblueprint}, which gives new, smaller, bounds on $p^a$ such that $\PSL_2(p^a)$ lies in the set $\mathscr P$ from Chapter \ref{ch:maxsubgroups}. In almost all cases, the proof is a formal generalization of the proof given by Liebeck and Seitz in \cite{liebeckseitz1998} for the original bounds from Chapter \ref{ch:maxsubgroups}, that if $H=\PSL_2(p^a)$ and $p^a\geq t(\bG)\cdot\gcd(2,p-1)$ then $H$ is strongly imprimitive.

\section{Preliminary Results}
\label{sec:appprelims}
Let $V$ be a $k\bG$-module of dimension $d$, and fix a basis $e_1,\dots,e_d$ of $V$. Let $\bT$ denote a maximal torus in $\bG$, and assume that $\bT$ acts diagonally on $V$ with respect to the basis $e_1,\dots,e_d$. Let $x$ be a semisimple element of order $n$ in $\bT$, and let $\zeta$ denote a primitive $n$th root of unity. For each $i$, the $\zeta^i$-eigenspace of the action of $x$ on $V$ is a subspace of $V$ spanned by some subset of the $e_i$, since $x$ acts diagonally.

If $y$ is another element of $\bG$ that stabilizes every subspace of $V$ that $x$ stabilizes, then in particular each $e_i$ is an eigenvector for $y$, so $y$ acts diagonally on $V$, and therefore $y\in \bT$. Thus if $x$ is contained in some subgroup $X$ stabilizing the same subspaces of $V$ as $x$, then $X\leq \bT$. Thus it makes sense to focus our attention on subgroups of $\bT$.

\begin{proposition}
Let $V$ be a $k\bG$-module. There exists a finite set of integers $\mathcal X_V$ such that if $n\notin \mathcal X_V$ then for any semisimple element $x\in \bT$ of order $n$ there exists an infinite subgroup $\bY$ of $\bT$ such that $x$ and $\bY$ stabilize the same subspaces of $V$, and conversely, if $n\in \mathcal X_V$ then there exists a choice of $x$ of order $n$ such that there is no infinite subgroup $\bY$ stabilizing the same subspaces of $V$ as $x$.
\end{proposition}
\begin{proof} Let $e_1,\dots,e_d$ denote a basis of $V$ with respect to which $\bT$ acts diagonally. Let $\sim$ be the equivalence relation on $\bT$ given by $y\sim y'$ if and only if $y$ and $y'$ have the same eigenspaces in their actions on $V$. Note that this is equivalent to $y$ and $y'$ stabilizing the same subspaces of $V$.

Since there are only finitely many possible invariant subspaces (as they are spanned by subsets of the $e_i$) there are only finitely many options for the eigenspaces of $x\in \bT$. In this case there are only finitely many equivalence classes $A_1,\dots,A_s$ for $\sim$. Let $Y_i$ denote the subgroup generated by $A_i$; then this stabilizes the same subspaces as any element of $A_i$, and is equal to the intersection of the subspace stabilizers
\[ S_x=\bigcap \{\bG_W\mid W\leq V,\;W\text{ stabilized by $x$}\}\]
for any $x\in A_i$. Let $\mathcal X_V$ denote the set of orders of elements of those $Y_i$ that are finite.

If $n\not\in \mathcal X_V$ then any $x\in \bT$ of order $n$ must lie in an equivalence class $A_i$ whose corresponding subgroup $Y_i$ is infinite, so take $\bY=Y_i$. Now suppose that $n\in \mathcal X_V$, and let $x$ be an element of order $n$ such that $x$ lies in a finite $Y_i$. Note that a priori $x$ need not lie in $A_i$, but this does not matter for our proof. The subspaces that are stabilized by $y\in A_i$ are also stabilized by $x$, so $S_x\leq Y_i$ and $S_x$ is therefore also finite. This completes the proof.
\end{proof}

What we see from the proposition is that, although there are many subgroups of $\bT$, only finitely many of them appear as the stabilizers of the set of subspaces that are stabilized by an element of $\bT$. Theorem \ref{thm:lls} states that $\mathcal X_{L(\bG)}=T(\bG)$. We will use $M(\bG)$ rather than $L(\bG)$ in an attempt to obtain better bounds, at the expense of having to use more effort. However, much of the effort of millions of calculations in abelian groups is done via computer.

We introduce a few pieces of notation and some definitions to help our discussion later. We start with an omnibus definition, containing many of the basic ideas we need.

\begin{definition} Let $e_1,\dots,e_d$ be a basis for a module $V$ such that $\bT$ acts diagonally on the $e_i$. A \emph{block system} $B$ is a set partition of $\{1,\dots,d\}$, whose constituent sets are called \emph{blocks}, and the \emph{stabilizer} of $B$ is the subgroup of $\bT$ consisting all of elements $x$ such that, whenever $i$ and $j$ lie in the same block of $B$, $e_i$ and $e_j$ are eigenvectors for the action of $x$ with the same eigenvalue, i.e., those elements $x$ that act as a scalar on every subspace of $V$ spanned by the $e_i$ for $i$ in a block of $B$. The \emph{dimension} of a block system $B$, denoted $\dim(B)$, is the dimension of its stabilizer as an algebraic group.

If $B$ and $B'$ are block systems then $B'$ is a \emph{coarsening} of $B$ if any block of $B$ is a subset of a block of $B'$, in other words, if the blocks of $B'$ are obtained by amalgamating blocks of $B$.

If $x$ is an element of $\bT$, then the \emph{block system associated to $x$} is the block system where $i$ and $j$ lie in the same set if and only if $e_i$ and $e_j$ are eigenvectors with the same eigenvalue under the action of $x$.

If $B$ is a block system then the \emph{closure} of $B$ is the block system $B'$ such that the stabilizers of $B$ and $B'$ are the same and if $B''$ is any coarsening of $B$ such that $B$ and $B''$ have the same stabilizer, then $B''$ is a coarsening of $B'$. In other words, the closure of $B$ is the coarsest block system with the same stabilizer as $B$. 
\end{definition}

We see that $x$ is contained in an infinite subgroup $\bY$ of $\bT$ such that $x$ and $\bY$ stabilize the same subspaces of $V$ if and only if the block system associated to $x$ has positive dimension.

For a given module $V$, we therefore wish to construct all block systems with finite stabilizers and compute the exponents of such groups. (Since the stabilizers are abelian, we only need the exponents as a finite abelian group contains an element of order $m$ if and only if $m$ divides the exponent of the group.)

In order to construct finite stabilizers for the module $M(\bG)$, we need some representation for torus elements acting on the minimal module. While this is easy for classical groups, for exceptional groups it is not necessarily so easy to obtain a representation of a maximal torus acting on $M(\bG)$.

Our solution to this is to choose a maximal-rank subgroup $H$ of an exceptional algebraic group that is a product of classical groups, and then use the maximal torus from that, which we understand. One downside to this is that the subgroup that we have, for example $A_2A_2A_2$ inside $E_6$, is not really $\SL_3^{\times 3}$ but a central product of $\SL_3$s. This just means that we have to be a bit more careful with the elements; we will discuss this more in Section \ref{sec:detbounds}.

\section{Determination of the Bounds for the Minimal Module}
\label{sec:detbounds}

In this section we compute $\mathcal X_{M(\bG)}$ for $\bG$ an exceptional algebraic group (other than $E_8$); we do so for $\bG=G_2$ by hand. For the larger groups, the computations are too cumbersome to do by hand, and we use a computer. The files and outputs of these are available on the author's website, and the algorithm used will be described here. As an independent check of our algorithm we also compute $T(\bG)$, i.e., $\mathcal X_{L(\bG)}$ for $\bG=G_2,F_4$, and we obtain the same answer as at the start of the chapter.

\subsection{$G_2$}
We start with $G_2$, where we let $\bH$ be the maximal-rank $A_2$ subgroup. The representation of the $A_2$ subgroup on $M(G_2)$ has composition factors the highest weight modules $L(01)$, $L(10)$ and $L(00)$, i.e., the trivial, the natural and its dual. (If $p=2$ then the $L(00)$ does not occur, but this does not affect the rest of the proof.)

Let $x$ be an element of order $n$ in $A_2$, with eigenvalues $\zeta^a,\zeta^b,\zeta^{-a-b}$ on $L(01)$, where $\zeta$ is a primitive $n$th root of unity. The eigenvalues of $x$ on $M(G_2)$ are therefore
\[ \zeta^a,\zeta^b,\zeta^{-a-b},\zeta^{-a},\zeta^{-b},\zeta^{a+b},1.\]
(We label the basis elements $e_1$ to $e_7$ in the order above.) The block system associated to $x$ is, generically, the singleton partition of $\{1,\dots,7\}$, and of course has dimension $2$. By considering just the exponents of the eigenvalues we obtain the list
\[ a,b,(-a-b),-a,-b,(a+b),0,\]
and equalities between these yield systems of linear equations. (Note that these should be taken modulo $n$.) Up to the $\Sym(3)$ automorphism group acting on the torus of $\bH$, we may assume that in any non-trivial block system $B$ the first element does not lie in a singleton set, and still be able to swap $b$ and $(-a-b)$. We conclude that $a$ is equal to one of $b$, $-a$, $-b$ or $0$.

Suppose first that $a=0$. The eigenvalue exponents are therefore $0$, $b$, $-b$, $0$, $b$, $-b$ and $0$, but the dimension of the block system
\[ \left\{\vphantom{x^2}\{1,4,7\},\{2,5\},\{3,6\}\right\}\]
is still $1$, so we need to make more eigenvalues equal. This means that $b$ is equal to either $-b$ or $0$, yielding $o(x)=1,2$. Hence $1,2\in \mathcal X_{M(G_2)}$.

We may therefore assume that the $1$-eigenspace of $x$ is $1$-dimensional, i.e., none of $a,b,a+b$ is equal to $0$, so $\{7\}$ is a set in the block system. We will remove the $1$-eigenspace from our lists from now on to remind us that it has been eliminated.

We still have that $a$ is equal to another eigenvalue exponent, say $b$. In this case the eigenvalue exponents are
\[ a,a,-2a,-a,-a,2a,\]
and the corresponding block system $\{\{1,2\},\{4,5\},\{3\},\{6\},\{7\}\}$ has dimension $1$ again. Setting two of $\pm a$ and $\pm 2a$ equal to one another yields $\alpha a=0$ for some $\alpha=1,2,3,4$, so $o(x)=1,2,3,4$. Hence $\{1,2,3,4\}\subseteq \mathcal X_{M(G_2)}$.

We may therefore assume that no two of $a$, $b$ and $-a-b$ are equal to each other, by applying the $\Sym(3)$ automorphism group. This means that $a$ must be equal, via the automorphisms, to either $-a$ or $-b$.

If $a=-b$ then $\zeta^{a+b}=1$, which is not allowed since we already assume that the $1$-eigenspace is $1$-dimensional. Thus $a=-a$, so $\zeta^a=-1$ (as $\zeta^a\neq 1$), and therefore the eigenvalues are
\[ -1,\zeta^b,-\zeta^b,-1,\zeta^{-b},-\zeta^{-b},1.\]
In particular, $p$ is odd, since else the $1$-eigenspace is at least $3$-dimensional. Setting $\zeta^b$ equal to $-1$ gives $o(x)=2$, to $-\zeta^b$ is impossible as $p$ is odd, to $\zeta^{-b}$ gives $o(x)=2$, and to $-\zeta^{-b}$ gives $\zeta^b=\pm\mathrm{i}$, so that $o(x)=4$. We therefore have proved the following proposition, valid for all primes.

\begin{proposition}\label{prop:numbersforG2}
If $\bG=G_2$ and $V=M(G_2)=L(10)$ then $\mathcal X_V=\{1,2,3,4\}$.
\end{proposition}

Comparing this set with $T(G_2)=\{1,\dots,10\}\cup\{12\}$, we see it offers a significant reduction in the set, so we will move on to larger exceptional groups.

\subsection{$F_4$}

Let $\bG=F_4$ and let $V=M(F_4)$. Since every semisimple element $x$ of $F_4$ lies in $D_4$, $1$ is always an eigenvalue of any semisimple element. Since $D_4$ has $24$ dimensions of non-trivial composition factors on $M(F_4)$, the $1$-eigenspace of $x$ on $M(F_4)$ has dimension at least $2$ for $p\neq 3$, and at least $1$ for $p=3$, since then $\dim(M(F_4))=25$. As for $G_2$, the precise dimension of the $1$-eigenspace is irrelevant to the calculations, so $p=3$ will follow the general pattern, but with the $1$-eigenspace dimension being one smaller.

As with $G_2$, we choose a maximal-rank subgroup $\bH$ in which we can easily represent a torus. In this case we choose two different subgroups: $A_2\tilde A_2$ and $A_1^4$. In the finite version of these groups, for example the group of the form $3\cdot (\PSL_3(q)\times \PSL_3(q))\cdot 3$ in the first case, one cannot guarantee that a semisimple element $x$ in this subgroup lies in the $\SL_3(q)\circ \SL_3(q)$ subgroup, only $x^3$. However, for the algebraic group this is not an issue. (One way to see this is to note that $\bT$ has no subgroups of index $3$.)

Thus we may work in $\bH$ rather than $\bG$. Note that there is a kernel of the map $\SL_3\times \SL_3\to \bH$, of order $3$, which will mean that our calculations will need to be modified to become exact. We will solve this problem when we get to it.

We let $\bH=A_1^4$ first: the representation of this on $V$ is the sum of two trivials and all six possible ways of tensoring two natural modules and two trivial modules for the four $A_1$ factors.

(The nicest way to see this is to take the $A_1^7$ inside $E_7$, which acts on the natural module for $E_7$ as the sum of seven modules, each a tensor product of three naturals distributed according to the Fano plane \cite[Proof of Lemma 2.1]{liebeckseitz2004a}, and then centralize one of the summands. On the other hand, the easiest way to see this is to take the maximal-rank $A_1^4$ inside $D_4$, and note that $D_4$ acts on $V$ as two trivials and the sum of the three non-isomorphic simple $8$-dimensional modules.)

Letting $x\in \bH$ have order some integer $n$, and writing $\zeta$ for a primitive $n$th root of unity, if $x$ has eigenvalues $\zeta^{\pm a_i}$ on the natural module for the $i$th copy of $A_1$, then the exponents for the eigenvalues of $x$ are
\[ \{\pm a_i\pm a_j:1\leq i<j\leq 4\}\cup\{0\},\]
with $0$ appearing twice, although this is not important for considering coincidences of eigenvalues. The symmetry group we apply here is $\Sym(4)$ acting in the obvious way on the $A_1$ factors.

We use a computer to analyse this situation, finding all coarsenings of the block system of dimension $4$ consisting of singletons, taking their closures, working up to automorphism, and continuing until we only have block systems of dimension $0$, finding $1264$ distinct block systems, with the following dimensions:
\[ 4^2,\;\; 3^{11},\;\; 2^{113},\;\;1^{538},\;\;0^{600}.\]
The exponents of the six-hundred torsion subgroups (i.e., stabilizers of block systems of dimension $0$) are the set of even numbers $\{2,\dots,36\}$.

This means that $\mathcal X_V$ is a subset of $\{2,\dots,18\} \cup\{20,22,24,\dots,36\}$. To obtain the actual set, we need to deal with the kernel of the map from $\SL_2(k)^4$ to $\bH$, which is a little bit tricky. In practice, we construct the abelian group as generated by four elements, the $a_i$, subject to the relations imposed by making sums of the $a_i$ equal to other sums of the $a_i$.

Let $g_1,\dots,g_n$ be a basis for the abelian group, and write
\[ a_i=\sum \alpha_{i,j} g_j.\]
If $g_j$ has order $n_j$, let $\zeta_j$ be a primitive $n_j$th root of unity. The matrix corresponding to $g_j$ should be diagonal, with coefficients $\zeta_j^{\alpha_{i,j}}$. Now take the group generated by these matrices. This will be a finite subgroup of $\bT$, thought of as $4\times 4$ diagonal matrices. We quotient out by the scalar matrix $-1$, and this is the image of our abelian subgroup in $\bH$.

Doing this reduces the exponents of the abelian groups to the set $\{1,\dots,18\}$.

\medskip

We will confirm this by choosing a different subgroup and getting the same answer. For $\bH=A_2\tilde A_2$, the representation of $\bH$ on $V$ is as the sum of three modules: the tensor product of the two naturals, the tensor product of the two duals, and the trivial for the $A_2$ by the adjoint representation $L(11)$ for the $\tilde A_2$ \cite[Lemma 11.11]{liebeckseitzbook}. We can more easily write down the exponents of the eigenvalues in terms of six variables
\begin{equation*} \begin{split}
\{a_i+a_j:1\leq i\leq 3,\;4\leq j\leq 6\}&\cup \{a_i+a_j:1\leq i\leq 3,\;4\leq j\leq 6\}
\\ &\quad\cup \{a_i-a_j:4\leq i\neq j\leq 6\}\cup\{0\},
\end{split}\end{equation*}
where $a_1+a_2+a_3=0$ and $a_4+a_5+a_6=0$. This time we have $\Sym(3)\times \Sym(3)$ acting by permuting the $a_i$ for $\{1,2,3\}$ and $\{4,5,6\}$.

We again use a computer to analyse this in the same way as before, finding $9278$ distinct block systems up to automorphism, with the following dimensions:
\[ 4,\;\; 3^{17},\;\; 2^{255},\;\;1^{2123},\;\;0^{6882}.\]
The exponents of the nearly seven thousand torsion subgroups are all multiples of $3$ between $3$ and $54$, but again we must remove the subgroup of order $3$ that forms the kernel of the map from $\SL_3\times \SL_3$ to $\bH$.

A similar approach works: this time we have $6\times 6$ matrices (placing $\bT$ inside $\SL_3\times \SL_3\leq \SL_6$, i.e., a torus of rank $6$), and we obtain the abelian group as a group of diagonal matrices. The kernel is now not a scalar matrix, but a block scalar matrix with two blocks of size $3$, each with coefficient a (different) cube root of unity.

Quotienting out by this yields exponents all integers between $1$ and $18$, agreeing with the previous calculation. We have therefore proved the following result, twice.

\begin{proposition}\label{prop:numbersforF4}
If $\bG=F_4$ and $V=M(F_4)=L(\lambda_1)$ then $\mathcal X_V=\{1,\dots,18\}$.
\end{proposition}

\subsection{$E_6$}

Let $\bG=E_6$ be simply connected and let $V=M(E_6)$. If $x$ is a real semisimple element then $x$ lies inside $F_4$ by Proposition \ref{prop:f4real}, and the eigenspaces of $x$ on $V$ are the same as for the minimal module of $F_4$, except the $1$-eigenspace has dimension one greater than for $F_4$. Changing the multiplicity of $1$ as an eigenvalue (as long as it is not changed to $0$) does not affect the calculations of the previous section. Hence if $x$ is $E_6$ is real semisimple and has order at least $19$, then $x$ is a blueprint for $V$. The same holds for $\overline{V}=V\oplus V^*$: since $x$ is real, the eigenvalues of $x$ on $V$ and $V^*$ are the same.

We will also produce a result that works for non-real elements as well. As with $F_4$, we consider two maximal-rank subgroups $\bH$ to confirm our results, settling on $A_5A_1$ and $A_2A_2A_2$ because their tori are simpler to write down.

We start with $A_5A_1$. The action of this subgroup on $V$ is as two composition factors, $(L(\lambda_4),L(0))$ and $(L(\lambda_1),L(1))$, where $L(\lambda_1)$ is the natural module and $L(\lambda_4)$ the exterior square of its dual \cite[Table 11.3]{liebeckseitzbook}. As with $F_4$, we label the eigenvalue exponents for the $A_5$ factor by $a_1,\dots,a_6$ such that $\sum_{i=1}^6 a_i=0$, and the other $A_1$ as $\pm a_7$. The eigenvalue exponents of $A_5A_1$ on $V$ therefore become
\[ \{-(a_i+a_j):1\leq i<j\leq 6\}\cup\{a_i\pm a_7:1\leq i\leq 6\}.\]
Now we have $\Sym(6)$ acting by permuting the $a_i$ for $\{1,\dots,6\}$.

We again use a computer to analyse this, finding $33365$ distinct block systems up to automorphism, with the following dimensions:
\[ 6,\;\;5^7,\;\;4^{68},\;\;3^{630},\;\;2^{4154},\;\;1^{12488},\;\;0^{16017}.\]
When considering the exponents of the finite such groups, we need to consider the centre of $\bG$, which of course will appear in every one of these subgroups, as well as the kernel of the map from $\SL_6\times \SL_2$ to $\bH$. First, the exponents of the finite stabilizers are all multiples of $6$ from $6$ to $156$. Removing the kernel using the same method as above (it acts as the scalar $-1$ in this case) yields abelian groups with exponent all multiples of $3$ from $3$ to $78$. Hence $\mathcal X_V$ consists of all divisors of $3i$ for $1\leq i\leq 26$, i.e.,
\[ \{1,\dots,27\}\cup\{30,33,36,\dots,78\}.\]
However, we are also interested in elements of that do not power to a non-identity central element: checking this is also easy inside the 16017 groups, simply by quotienting out by the centre. Doing so yields groups of exponent between $1$ and $27$, so if $o(x)> 27$ and $\gen x\cap Z(\bG)=1$, then $x$ is a blueprint for $V$.

Now we turn to $\bH=A_2A_2A_2$. The action of this subgroup on $V$ is as three composition factors, $(10,01,00)$, $(00,10,01)$ and $(01,00,10)$, where $L(10)$ is the natural module and $L(01)$ is its dual \cite[Proposition 2.3]{liebeckseitz1996}. We label the eigenvalue exponents by $a_i$ for $i=1,\dots,9$, such that the sums $a_1+a_2+a_3$, $a_4+a_5+a_6$ and $a_7+a_8+a_9$ are all zero. The eigenvalue exponents of $x$ on $V$ are therefore
\begin{equation*} \begin{split} \{a_i-a_j:1\leq i\leq 3,\;4\leq j\leq 6\}&\cup\{a_i-a_j:4\leq i\leq 6,\;7\leq j\leq 9\}
\\ &\quad\cup \{a_i-a_j:7\leq i\leq 9,\;1\leq j\leq 3\}.
\end{split}\end{equation*}
In this case, the group $\Sym(3)\wr \Sym(3)\leq \Sym(9)$ acts by preserving the set partition $\{\{1,2,3\},\{4,5,6\},\{7,8,9\}\}$.

The computer now finds $26498$ distinct block systems, with the following dimensions:
\[ 6,\;\;5^4,\;\;4^{51},\;\;3^{565},\;\;2^{4002},\;\;1^{12162},\;\;0^{9713}.\]
The exponents of the torsion groups are all multiples of $3$ between $3$ and $81$. Again there is a kernel, and the centre. The kernel is generated by a scalar matrix that is a root of unity, and the centre is given by a block scalar matrix, with each block of size $3$ one of the three cube roots of unity. Quotienting out by the kernel and taking exponents again yields the set
\[ \{1,\dots,27\}\cup\{30,33,36,\dots,78\}.\]
Again, quotienting out by the centre and taking exponents yields the set $\{1,\dots,27\}$. Thus we again have two proofs of the following proposition. This particular case does not appear in Theorem \ref{thm:blueprintsminimal} later since we need to consider $\overline V$ rather than $V$.

\begin{proposition}\label{prop:numbersforE6}
If $\bG=E_6$ and $V=M(E_6)=L(\lambda_1)$ then
\[ \mathcal X_V=\{1,\dots,27\}\cup \{30,33,36,\dots,78\}.\]
Furthermore, if $x\in \bG$ is semisimple, and $|\gen x\cdot Z(\bG)/Z(\bG)|>27$, then $x$ is a blueprint for $V$.
\end{proposition}

Since $\overline V$ is a submodule of the restriction of the minimal module for $E_7$ to the $E_6$ Levi subgroup, we will be able to use our results in the following section to get bounds on $\mathcal X_{\overline V}$.

\subsection{$E_7$}

Let $\bG=E_7$ be simply connected, and let $V=M(E_7)$. Since $V$ has dimension $56$ and the torus has rank $7$, one expects the number of block systems to be much higher, and for the programs to take much longer to run, which is true. It also means that there are too many block systems to store them all efficiently, and so we have to alter our algorithm for computing these slightly.

We let $\bH$ be the maximal-rank $A_7$ subgroup. The representation of this on $V$ has composition factors $L(\lambda_2)$ and $L(\lambda_6)$, i.e., the exterior square of the natural and its dual \cite[Lemma 11.8]{liebeckseitzbook}. The exact form of the finite groups $H(q)$ is $4\cdot \PSL_8(q)\cdot 2$ if $q\equiv 1 \bmod 8$ \cite[Table 5.1]{liebecksaxlseitz1992}. By restricting our attention to odd-order elements $x$, we avoid questions about whether the element $x$ powers to the central involution of $\bG$ if the characteristic $p$ is odd, and also do not have to consider the kernel of the map from $\SL_8$ to $\bH$.

Let $x$ be an element of order $n$ in $\bH$, with eigenvalues on the natural module for $\bH$ being $\zeta^{a_i}$ for $i=1,\dots,8$ with $\sum a_i=0$, where $\zeta$ is a primitive $n$th root of unity. The eigenvalues of $x$ on $V$ are therefore $\{\zeta^{\pm(a_i+a_j)}\mid 1\leq i<j\leq 8\}$. We again use a computer to analyse this situation, although we stop when we reach the block systems of dimension $2$, yielding
\[ 7,\;\; 6^5,\;\;5^{47},\;\;4^{626},\;\; 3^{9781},\;\; 2^{116170}.\]
We obviously do not want to try to store the likely million block systems of dimension $1$, so for each block system of dimension $2$ we find all coarsenings, and repeat the process until we reach block systems of dimension $0$. This of course introduces computational repetition but reduces the space requirement.

Doing this produces the set of exponents of the finite subgroups of all multiples of $8$ up to $264$, and all $a\equiv 4\bmod 8$ up to $300$. In particular, the odd divisors of these numbers are all (odd) integers up to $75$, and this means that we obtain the following proposition.

\begin{proposition}\label{prop:numbersforE7}
If $\bG=E_7$ and $V=M(E_7)=L(\lambda_1)$, then the odd elements of $\mathcal X_V$ are $\{1,3,\dots,75\}$.
\end{proposition}

A much simpler calculation is to find $\mathcal X_V$ when $V$ is the restriction of $M(E_7)$ to the $A_4$ Levi subgroup, which we will need for the proof of Proposition \ref{prop:f4a4} later.

The restriction of $M(E_7)$ to $\bH=A_4$ has composition factors (with various multiplicities) $L(0000)$, $L(1000)$, $L(0100)$, $L(0010)$ and $L(0001)$. Since multiplicities are not important when computing the set $\mathcal X_V$, we may assume that $V$ is the sum of a single copy of each of these modules.

Thus in this case $x$ is an element of order $n$ in $\bH$, with eigenvalues on $M(A_4)$ being $\zeta^{a_i}$ for $i=1,\dots,5$ with $\sum a_i=0$, where $\zeta$ is a primitive $n$th root of unity. The eigenvalues of $x$ on $V$ are therefore
\[ \{1\}\cup \{\zeta^{\pm a_i}\mid 1\leq i\leq 5\}\cup \{\zeta^{\pm (a_i+a_j)}\mid 1\leq i<j\leq 5\}.\]
In this case, $\Sym(5)$ acts on the $a_i$. As usual, we use a computer to find $886$ distinct block systems, with the following dimensions:
\[ 4,\;\;3^8,\;\; 2^{61},\;\; 1^{305},\;\; 0^{511}.\]
The exponents of the $511$ finite abelian groups are easily found, and we obtain the following result.

\begin{proposition}\label{prop:numbersforA4} If $\bG=A_4$ and $V$ is the sum of the modules $L(1000)$, $L(0100)$, $L(0010)$, $L(0001)$, and $L(0000)$, then
\[ \mathcal X_V=\{1,\dots,28\}\cup\{30\}.\]
\end{proposition}

We amalgamate the results of the whole section.

\begin{theorem}\label{thm:blueprintsminimal} Suppose that $\bG$ is an algebraic group of type $F_4,E_6,E_7$ in characteristic $p$, and that $V=M(\bG)$, with $\overline V=M(\bG)\oplus M(\bG)^\tau$ where $\tau$ is the graph automorphism of $\bG$ if such an automorphism exists. Let $x$ be a semisimple element of $\bG$.
\begin{enumerate}
\item\label{thmi:bluea} If $\bG=F_4$ and $p$ is odd, then $x$ is a blueprint for $V$ if $o(x)>18$.
\item\label{thmi:blueb} If $\bG=F_4$ and $p=2$, then $x$ is a blueprint for $\overline V$ if $o(x)>57$.
\item\label{thmi:bluec} If $\bG=E_6$ and $x$ is real, then $x$ is a blueprint for $\overline V$ if $o(x)>18$.
\item\label{thmi:blued} If $\bG=E_6$ and $x$ is non-real, then $x$ is a blueprint for $\overline V$ if $o(x)>75$ and $o(x)$ is odd.
\item\label{thmi:bluee} If $\bG=E_7$, then $x$ is a blueprint for $V$ if $o(x)>75$ and $o(x)$ is odd.
\end{enumerate}
\end{theorem}
\begin{proof}(\ref{thmi:bluea}), (\ref{thmi:bluec}) and (\ref{thmi:bluee}) are proved directly in Section \ref{sec:detbounds}. For (\ref{thmi:blueb}), note that $\overline V$ has the same composition factors as $L(\bG)$, so $\mathcal X_{\overline V}=\mathcal X_{L(\bG)}=T(\bG)$. As $o(x)$ is odd (as $x$ is semisimple) we obtain the result.

For (\ref{thmi:blued}), if $o(x)>75$ and $o(x)$ is odd, then by placing $x$ inside $E_7$ via the $E_6$-Levi subgroup, we see that $x$ is a blueprint for the minimal module for $E_7$, which restricts to $E_6$ as $L(0)^{\oplus 2}\oplus \overline V$. We need to check that if $\bX$ is a positive-dimensional subgroup stabilizing the same subspaces of $\overline{V}$ as $x$ then $\bX\leq E_6$, and then we are done.

To see this, note that $\bX$ must act trivially on the $2$-dimensional subspace $L(0)^{\oplus 2}$ in the restriction of the minimal module to $E_6$. By Lemma \ref{lem:e7stabs}, we see that $\bX$ lies inside the subgroup $E_6.2$ where the $2$ is the graph automorphism. Since the graph automorphism acts non-trivially on this $2$-space and $\bX$ acts trivially on it, $\bX$ must lie inside $E_6$, as claimed.
\end{proof}

\section{Consequences for Maximal Subgroups}

In this short section we apply the results about semisimple elements being blueprints for $M(\bG)$ to obtain better bounds on when $H=\PSL_2(p^a)$ is strongly imprimitive in $\bG$ than $p^a>\gcd(2,p-1)\cdot t(\bG)$, which we saw in Chapter \ref{ch:maxsubgroups}.

The set of potential maximal subgroups that are irreducible on $M(\bG)$ is given in \cite{liebeckseitz2004a}, and for $\bG=F_4,E_6,E_7$ there are no irreducible subgroups $\PSL_2(p^a)$, but there is a copy of $\PSL_2(128).7$ acting irreducibly on $M(E_7)$. Since the $\PSL_2(128)$ subgroup still acts reducibly, this is no barrier to proving strong imprimitivity.

\begin{theorem}\label{thm:goodblueprint} Let $v(\bG)$ be given by
\[ v(\bG)=\begin{cases}
18 & \bG=F_4,\; E_6,
\\75 & \bG=E_7.\end{cases}\]
If $H$ is a subgroup of $\bG$ such that $H\cdot Z(\bG)/Z(\bG)\cong \PSL_2(p^a)$ for some $p^a>\gcd(2,p-1)\cdot v(\bG)$, then $H$ is strongly imprimitive.
\end{theorem}
\begin{proof} Let $H$ be as described in the result, and note that $H$ cannot act irreducibly on either $M(\bG)$ or $L(\bG)^\circ$ by \cite{liebeckseitz2004a}. Suppose first that $p$ is odd, and note that $H$ contains a real semisimple element $x$ of order $(p^a-1)/2$ and $y$ of order $(p^a+1)/2$.

The next three statements all follow from Theorem \ref{thm:blueprintsminimal}:
\begin{enumerate}
\item If $\bG=F_4$ then $y$, and hence $H$, are blueprints for $M(F_4)$. 
\item If $\bG=E_6$ then $y$, and hence $H$, are blueprints for $M(E_6)\oplus M(E_6)^*$.
\item If $\bG=E_7$ then one of $x$ and $y$ have odd order so it, and hence $H$, are blueprints for $M(E_7)$.
\end{enumerate}
In all three cases, $H$ is a blueprint for a module satisfying the hypothesis for Proposition \ref{prop:blueprintissi}, and so $H$ is strongly imprimitive.

Thus suppose $p=2$, and note that $H$ contains an element of order $2^a+1$. If $\bG=E_6$ then $H$ is a blueprint for $M(E_6)\oplus M(E_6)^*$, and if $\bG=E_7$ then $H$ is a blueprint for $E_7$, again by Theorem \ref{thm:blueprintsminimal}, and again strongly imprimitive by Proposition \ref{prop:blueprintissi}. Thus $\bG=F_4$. If $a\geq 6$ then $H$ contains an element of order $65$, whence $H$ is a blueprint for $M(F_4)\oplus M(F_4)^\tau$ by Theorem \ref{thm:blueprintsminimal} and so strongly imprimitive by Proposition \ref{prop:blueprintissi}. For $a=5$ then we do not know that $H$ is a blueprint for $M(F_4)\oplus M(F_4)^\tau$ (where $\tau$ is the graph automorphism), but we do know from Proposition \ref{prop:numbersforF4} that $H$ is a blueprint for $M(F_4)$.

Let $\bar G$ be an almost simple group with socle $G=F_4(2^b)$. If $\bar G$ does not induce a graph automorphism on $G$ then $M(F_4)$ is stable under $\Aut_{\bar G}(\bG)$, so we may apply Proposition \ref{prop:blueprintissi} to see that $H$ is strongly imprimitive. On the other hand, if $\bar G$ does induce a graph automorphism on $G$, let $G_0$ denote the (unique) subgroup of index $2$ in $\bar G$, so that $G_0$ only induces field (and inner) automorphisms on $G$.

If $N_{\bar G}(H)$ lies inside $G_0$ then we replace $\bar G$ by $G_0$ and get that $H$ is strongly imprimitive in $\bG$. Since $|\Out(H)|=5$, if $N_{\bar G}(H)\not\leq G_0$ then an element of $\bar G\setminus G_0$ centralizes $H$. But the centralizer of a graph morphism on $F_4(2^b)$ is simply a group ${}^2\!F_4(2^c)$ for some $c\leq b$. One may proceed by replacing the Frobenius automorphism $\sigma$ by one whose fixed points are ${}^2\!F_4(2^c)$, but the easier way is to simply use the list of maximal subgroups of the almost simple Ree groups, and particularly \cite[Proposition 2.7]{malle1991}. This shows that the normalizer of any copy of $\SL_2(32)$ lies in a subgroup ${}^2\!A_2$ (either simply connected or adjoint type). This completes the proof.
\end{proof}

\chapter{Unipotent and Semisimple Elements}
\label{ch:unipssemis}

This chapter collects together a variety of facts about unipotent and semisimple elements in groups of Lie type. We consider criteria for unipotent and semisimple elements to be blueprints. We then move on to considering modules for $\SL_2$, and how the weight spaces of the module and the eigenvalues of elements of $\SL_2$ interact, with the aim of finding blueprint elements and blueprint subgroups $\PSL_2(p^b)$.

\section{Actions of Unipotent Elements}

Let $\bG$ be a simple algebraic group in characteristic $p$. The Bala--Carter--Pommerening labelling system for the unipotent classes, as used in a slightly modified form (to deal with interpolation of extra classes in certain bad characteristics) in our main reference \cite{lawther1995} for unipotent classes of exceptional groups, gives us a way to discuss unipotent classes that is independent of the characteristic $p$ of $\bG$. We may therefore compare the action of a unipotent class on a fixed simple module for different primes.

As is well known, any matrix of order a power of a prime $p$ defined over a field of characteristic $p$ can be written in Jordan normal form, with the conjugacy class in the general linear group being determined by the sizes of the Jordan blocks. Thus, if $u$ is a unipotent element of $\bG$ then for every module for $\bG$ of dimension $n$ we can associate a partition of $n$, the sizes of the various Jordan blocks in the action of $u$ on the module. We use the notation for this, and unipotent classes, from \cite{lawther1995}, which determines the Jordan block structure of the action of all unipotent classes of exceptional groups on the minimal and adjoint modules.

The only cases we will need that are not covered in \cite{lawther1995} are when $L(\bG)\neq L(\bG)^\circ$ and $M(F_4)$ for $p=3$. The next lemma gives the actions of the unipotent classes on the $25$-dimensional simple module $M(F_4)$, on the $26$-dimensional Weyl module $25/1$, and on the $27$-dimensional minimal module $M(E_6)$ for $E_6$, which has structure $1/25/1$ when restricted to $F_4$.

\begin{lemma}\label{lem:unipotentF4}
Let $u$ be a unipotent element in $F_4(3^n)$. The Jordan blocks of the action of $u$ on the $25$-dimensional minimal module $M(F_4)$, together with the extension $25/1$ and the minimal module $M(E_6)$ for $E_6$ is one of those given in Table \ref{t:unipotentF4}.
\begin{table}
\begin{center}
\begin{tabular}{lccc}
\hline Class in $F_4$ & Action on $M(F_4)$ & Action on $25/1$ & Action on $M(E_6)=1/25/1$
\\\hline $A_1$ & $2^6,1^{13}$ & $2^6,1^{14}$ & $2^6,1^{15}$
\\ $\tilde A_1$ & $3,2^8,1^6$ & $3,2^8,1^7$ & $3,2^8,1^8$
\\ $A_1+\tilde A_1$ & $3^3,2^6,1^4$ & $3^3,2^6,1^5$ & $3^3,2^6,1^6$
\\ $A_2$ & $3^6,1^7$ & $3^6,1^8$ & $3^6,1^9$
\\ $A_2+\tilde A_1$ & $3^7,2^2$ & $3^7,2^2,1$ & $3^7,2^2,1^2$
\\ $\tilde A_2$, $\tilde A_2+A_1$ & $3^8,1$ & $3^8,2$ & $3^9$
\\ \hline $B_2$ & $5,4^4,1^4$ & $5,4^4,1^5$ & $5,4^4,1^6$
\\ $C_3(a_1)$ & $5^2,4^2,3,2^2$ & $5^2,4^2,3,2^2,1$ & $5^2,4^2,3,2^2,1^2$
\\ $F_4(a_3)$ & $5^3,3^3,1$ & $5^3,3^3,1^2$ & $5^3,3^3,1^3$
\\ $B_3$ & $7^3,1^4$ & $7^3,1^5$ & $7^3,1^6$
\\ $C_3$, $F_4(a_2)$ & $9,6^2,3,1$ & $9,6^2,3,2$ & $9,6^2,3^2$
\\ $F_4(a_1)$ & $9^2,7$ & $9^2,7,1$ & $9^2,7,1^2$
\\ \hline $F_4$ & $15,9,1$ & $15,9,2$ & $15,9,3$
\\ \hline
\end{tabular}
\end{center}
\label{t:unipotentF4}\caption{Actions of unipotent elements on $M(F_4)$ and its extensions for $F_4$ in characteristic $3$}
\end{table}
\end{lemma}
\begin{proof}
The Jordan blocks of the actions of the unipotent elements on the $26$-dimensional module are given in \cite[Table 3]{lawther1995}, and using a computer, a representative of each of the classes was constructed in $F_4(3)$. The Jordan blocks of their actions on the $25$-dimensional composition factor were then computed, and are as above. The classes on the $25/1$ are exactly those in \cite[Table 3]{lawther1995}, and the corresponding classes for $E_6$ are in \cite[Table 5]{lawther1995}.
\end{proof}

Using a computer and constructing classes manually is the method by which we prove the next two lemmas, which we include for completeness.

\begin{lemma}\label{lem:unipotentE6}
Let $u$ be a unipotent element in $E_6(3^n)$. The Jordan blocks of the action of $u$ on the $77$-dimensional Lie algebra module $L(E_6)^\circ$ are obtained from the action on $L(E_6)$ by removing a Jordan block of size 1, except in the cases listed in Table \ref{t:unipotentE6}.
\begin{table}
\begin{center}
\begin{tabular}{lcc}
\hline Class in $E_6$ & Action on $L(E_6)^\circ$ & Action on $L(E_6)$
\\\hline $2A_2$&$3^{23},1^8$&$3^{23},2,1^7$
\\$2A_2+A_1$&$3^{24},2^2,1$&$3^{24},2^3$
\\$A_5$&$9^3,8^2,6^4,3^2,1^4$&$9^3,8^2,6^4,3^2,2,1^3$
\\$E_6(a_3)$&$9^4,7,6^4,3^3,1$&$9^4,7,6^4,3^3,2$
\\$E_6(a_1)$&$9^8,5$&$9^8,6$
\\$E_6$&$19,15^2,9^3,1$&$19,15^2,9^3,2$
\\\hline
\end{tabular}
\end{center}
\label{t:unipotentE6}\caption{Actions of unipotent elements on $L(E_6)^\circ$ and $L(E_6)$ for $E_6$ in characteristic $3$, where one does not obtain the former from the latter by removing a trivial Jordan block}
\end{table}
\end{lemma}

\begin{lemma}\label{lem:unipotentE7} Let $u$ be a unipotent element in $E_7(2^n)$. The Jordan blocks of the action of $u$ on the $132$-dimensional Lie algebra module $L(E_7)^\circ$ are obtained from the action on $L(E_7)$ by removing a Jordan block of size 1, for every unipotent class.
\end{lemma}

We can see from the tables in \cite{lawther1995} that for every unipotent class there is a set of primes $P$ such that, for any prime $p\not\in P$ the partition describing the Jordan block structure on a fixed module $V$ is the same.

\begin{definition}\label{defn:genericunipotent} Let $\bG$ be an algebraic group and let $V$ be a highest weight module for $\bG$. Let $u$ be a unipotent element of $\bG$ belonging to a unipotent class that exists for all primes. If the Jordan block structure of $u$ on $V$ is the same as for cofinitely many primes, then $u$ is said to be \emph{generic} on $V$.
\end{definition}

Thus, informally, the non-generic classes are those where the prime is in the set $P$ described above, where the partition differs from the `usual' one.

The reason that generic unipotent classes are interesting is that we can find `nice' $A_1$ subgroups containing them, at least if the class has elements of order $p$. In \cite[Lemma 1.2]{craven2015un2} we show that such unipotent elements are blueprints for $M(\bG)$ and, indeed, $M(E_6)\oplus M(E_6)^*$ if $\bG=E_6$.

\begin{lemma} [{{\cite[Lemma 1.2]{craven2015un2}}}]
\label{lem:genericunipotent}
Suppose that $\bG=F_4,E_6,E_7$ with $p$ odd for $\bG=F_4$. Let $H$ be a finite subgroup of $\bG$ such that $H$ contains a non-trivial unipotent element whose action on a module $V$, one of $M(\bG)$, $M(E_6)\oplus M(E_6)^*$ if $\bG=E_6$, or $L(\bG)^\circ$ is generic. Then $u$ and $H$ are blueprints for $V$. In particular, $H$ is either strongly imprimitive or irreducible on $V$.
\end{lemma}
(The consequence follows from Proposition \ref{prop:blueprintissi}.)

Thus, if any subgroup $H$ of an exceptional algebraic group $\bG$ contains a unipotent element of order $p$ that is generic for either the minimal or adjoint module, then $H$ is either strongly imprimitive or $H$ acts irreducibly on this module. However, such subgroups are listed in \cite{liebeckseitz2004a}, so if $H$ does not appear on the list in \cite{liebeckseitz2004a} then $H$ is strongly imprimitive.

\medskip

For large primes, we will occasionally prove that $H$ stabilizes a unique $3$-dimensional submodule of $L(\bG)$, which must then be a subalgebra of the Lie algebra. If this $3$-dimensional submodule of $L(\bG)\downarrow_H$ is a summand then we may apply Proposition \ref{prop:sl2ifsplitoff} below and show that it is a copy of $\slf_2$, but if the $3$-dimensional submodule is not a summand then we cannot easily prove that it is an $\slf_2$, as it need not be simple. There is one case in particular where this occurs, which we refer to as a Serre embedding. These are embeddings of $\PSL_2(h+1)$ into an algebraic group in characteristic $h+1$, where $h$ is the Coxeter number of the group.

\begin{definition}\label{defn:serreembedding} Let $\bG$ be an exceptional algebraic group with Coxeter number $h$, and let $p=h+1$. A subgroup $H=\PSL_2(p)$ is a \emph{Serre embedding} if the following conditions hold:
\begin{enumerate}
\item on $L(\bG)$, $H$ stabilizes a unique $3$-dimensional subspace;
\item $H$ contains a regular unipotent element.
\end{enumerate}
\end{definition}

The $3$-dimensional subspace is in fact a subalgebra: the exterior square of a $3$-dimensional $kH$-submodule $W$ is $W^*\cong W$, and $\Hom_{kH}(\Lambda^2(W),L(E_8))$ is $1$-dimensional. From \cite[Lemma 1]{ryba2002}, $W$ is a Lie subalgebra of $L(E_8)$. Since this subspace need not (in fact, will not be) a summand of $L(E_8)\downarrow_H$, we cannot directly apply Proposition \ref{prop:sl2ifsplitoff} below.

\section{Blueprints and Element Orders}
\label{sec:blueprints}

The first result states that whether a semisimple class contains blueprints for a fixed Weyl module is independent of the characteristic of the underlying field. To prove this requires some of the theory of semisimple elements, for example \cite[Section 3.1.6]{litterickmemoir}. It also looks quite technical, but that is merely because we have to set up some bijection between characteristic $0$ and characteristic $p$, and also need to worry about the fact that the eigenvalues of a semisimple element on a Weyl module do not always determine the conjugacy class uniquely.

\begin{proposition}\label{prop:blueprintinvariant} Let $\bG_p$ and $\bG_0$ be simple, simply connected algebraic groups of the same type in characteristics $p$ and $0$ respectively. Let $V_p$ and $V_0$ denote the Weyl modules of a fixed weight for $\bG_p$ and $\bG_0$ respectively.

Let $a_1,\dots,a_d$ be integers. Let $n$ be a positive integer prime to $p$, and let $\theta_n$ and $\zeta_n$ be primitive $n$th roots of unity in $k$ and $\mathbb{C}$ respectively. Let $X_p$ and $X_0$ denote the set of elements of $\bG_p$ and $\bG_0$ with eigenvalues 
\[\{\theta_n^{a_i}\mid 1\leq i\leq d\}\quad \text{and}\quad \{\zeta_n^{a_i}\mid 1\leq i\leq d\}\]
on $V_p$ and $V_0$ respectively. The set $X_p$ consists solely of blueprints for $V_p$ if and only if the set $X_0$ consists solely of blueprints for $V_0$.
\end{proposition}
\begin{proof} Being a blueprint for a module is an invariant of the semisimple class, so we may assume that our elements lie in maximal tori $\bT_p$ and $\bT_0$ of $\bG_p$ and $\bG_0$ respectively. Choose bases of $V_p$ and $V_0$ so that $\bT_p$ and $\bT_0$ act diagonally. From the theory in \cite[Section 3.1.6]{litterickmemoir}, there is a bijection $f$ between $X_p\cap \bT_p$ and $X_0\cap \bT_0$, and a bijection between elements of $\bT_p$ of order $an$ for $a\geq 1$ powering to $x\in X_p\cap \bT_p$ and elements of $\bT_0$ of order $an$ powering to $f(x)$. Furthermore, this bijection preserves Brauer characters on $V_p$ and $V_0$, i.e., preserves eigenvalues under the assignment $\theta\mapsto \zeta$.

Suppose that $X_p$ consists of blueprints for $V_p$. Thus for each $x\in X_p\cap \bT_p$, there exist elements of arbitrarily large order in $\bT_p$ powering to $x$ and with the same number of distinct eigenvalues on $V_p$. Taking the corresponding elements of $\bT_0$, we find elements of arbitrarily large order $an$ for $a\geq 1$ in $\bT_0$ powering to $f(x)$. Hence each $f(x)$ is a blueprint for $V_0$.

The same argument would work as a converse if we could always choose $a$ to be prime to $p$. This is the case: if $x_0=f(x)$ lies in $\bT_0$ and the set of all elements of $\bT_0$ that stabilize all of the same subspaces of $V_0$ as $x_0$ forms a group of diagonal matrices, hence is a direct product of a torus and a finite abelian group by \cite[Theorem 16.2]{humphreyslag}. As a torus in characteristic $0$ contains a product of groups $\mathbb{\bar\Q}^\times$, and $\mathbb{\bar\Q}^\times$ contains elements of all orders, the result holds.
\end{proof}

We can also push being a blueprint for the minimal module of $F_4$ up into $E_6$ and $E_7$.

\begin{lemma}\label{lem:f4blueprintisblueprint} Let $\bG$ be $E_6$ or $E_7$, and let $x$ be a semisimple element of $\bG$ that lies in $F_4$. If $x$ is a blueprint for $M(F_4)$, then $x$ is a blueprint for $M(E_6)$, $M(E_6)\oplus M(E_6)^*$ and $M(E_7)$.
\end{lemma}
\begin{proof} The restrictions of the modules in question to $F_4$ have composition factors copies of $M(F_4)$ and trivial modules. Since $1$ is always an eigenvalue of any semisimple element of $F_4$ on $M(F_4)$, any element of $F_4$ with the same eigenspaces on $M(F_4)$ as $x$ also has the same eigenspaces on $M(E_6)$, $M(E_6)\oplus M(E_6)^*$ and $M(E_7)$. This completes the proof.
\end{proof}

The semisimple elements of $E_6$ that lie in $F_4$ are the real elements. Although the author has looked for a reference for this, beyond \cite[Theorem 3.1]{cohenwales1997}, which shows this only for order at most $7$, it does not appear to be easy to find explicitly. It follows almost immediately from results of Moody and Patera \cite{moodypatera1984}.

\begin{proposition}\label{prop:f4real} Let $\bG$ be the simply connected form of $E_6$. If $x$ is a real semisimple element of $\bG$ then $x$ is conjugate to an element of $F_4$.
\end{proposition}
\begin{proof} We use the results and notation of \cite[Section 4]{moodypatera1984}: if $x\in F_4$ is semisimple and its class corresponds to the Kac co-ordinates $(s_0,\dots,s_4)$, then the corresponding Kac co-ordinates in $E_6$ are $(s_0,s_1,s_2,s_3,s_2,s_1,s_4)$. By \cite[Proposition 5.3]{moodypatera1984}, a conjugacy class in $E_6$ is real if and only if its Kac co-ordinates have exactly this form, and hence all real semisimple classes of $E_6$ intersect $F_4$ non-trivially, as desired.
\end{proof}

If $\bG=E_7$ then $v(\bG)=75$ is fairly large, and in certain circumstances we can bring this down. Here is one such circumstance.

\begin{proposition}\label{prop:f4a4} Let $\bG$ be the simply connected form of $E_7$, and let $x$ be a semisimple element of $\bG$. If the $1$-eigenspace of $x$ on $M(E_7)$ has dimension at least $6$ then $x$ lies inside a conjugate of either an $F_4$ or $A_4$ subgroup. If in addition $o(x)>30$ then $x$ is a blueprint for $M(E_7)$.
\end{proposition}
\begin{proof} As in the proof of Proposition \ref{prop:numbersforE7}, we first place $x$ into the subgroup $A_7$. If $x$ centralizes a point on $\Lambda^2(M(A_7))$, then this requires one eigenvalue of $x$ on $M(A_7)$ to be the inverse of another. (In the notation of the proof of Proposition \ref{prop:numbersforE7}, $a_1=-a_2$.) This places $x$ into $A_1A_5$. However, by Lemma \ref{lem:e7stabs}, $x$ must lie in $E_6$ or $B_5$ (as if $x$ lies in a parabolic it lies in a Levi), whence $x$ lies in $E_6$. Since $x$ centralizes a $4$-space on $M(E_6)\oplus M(E_6)^*$, it must centralize a $2$-space $W$ on $M(E_6)$ (as it acts semisimply).

If $x$ lies in $F_4\leq E_6$ then we are done, so suppose not. From Lemma \ref{lem:e6stabs}, this means that $x$ cannot centralize an $F_4$-point, or a $B_4$-point, since $B_4\leq F_4$. In the notation of \cite{cohencooperstein1988}, this means that $W$ is a purely white $2$-space. By \cite[(P.3)]{cohencooperstein1988}, $W$ is unique up to action of $E_6$, and the stabilizer is an $A_1A_4$-parabolic of $E_6$, so the centralizer is an $A_4$-parabolic. Thus $x$ lies in $A_4$, as claimed.

If $x\in F_4$ and $o(x)>30$ then $x$ is a blueprint for $M(F_4)$ by Proposition \ref{prop:numbersforF4}, and hence a blueprint for $M(E_7)$ by Lemma \ref{lem:f4blueprintisblueprint}. On the other hand, if $x\in A_4$ then one uses Proposition \ref{prop:numbersforA4}. (The composition factors of $A_4$ on $M(E_7)$ are\[L(1000)^3,L(0100),L(0010),L(0001)^3,L(0000)^6\]
by \cite[Table 21]{thomas2016}, so the result may be used.) 
\end{proof}

Suppose we want to find the eigenvalues on $M(E_7)$ of semisimple elements of order $63$ inside $E_7$, which we will need to do when considering $\SL_2(64)$. There are too many to construct and store them all effectively, but we can take an element $x$ of order $21$ and consider all $3^7=2187$ elements $\hat x$ that cube to $x$ in a torus. Since we have the eigenvalues of all elements of order $21$, given a potential multiset of eigenvalues for an element $x$ of order $63$ in $E_7$, we take the eigenvalues of $x^3$, find all semisimple classes of elements of order $21$ with those eigenvalues, then consider all preimages of representatives of each of those classes. The eigenvalues of $x$ are valid for coming from $E_7$ if and only if one of those elements of order $21$ has a preimage with those values.

This idea to find the eigenvalues of elements of large composite order will be called the \emph{preimage trick} in the rest of this paper.

\section{Blueprints inside Subgroups of Type \texorpdfstring{$A_1$}{A1}}

We now prove that certain semisimple elements, and subgroups of the form $\SL_2(p^a)$ and $\PSL_2(p^a)$ of exceptional groups, are blueprints for a given module by examining the constituents of the restriction of the module to an $A_1$ subgroup containing the element or subgroup.

The first lemma deals with modules for the algebraic group $\SL_2$, and when the eigenspaces of semisimple elements match the weight spaces. Generalizing the idea of $p^a$-restricted, a highest weight module for $\SL_2$ is \emph{$n$-restricted} for some $n\in\N$ if its highest weight is $\lambda$ for $\lambda<n$.

\begin{lemma}\label{lem:restrictedsl2} Let $M$ be a module for $\SL_2$ with composition factors highest weight modules $L(\lambda_1),\dots,L(\lambda_r)$, arranged so that $\lambda_i\leq \lambda_{i+1}$. Let $\bT$ be a maximal torus of $\SL_2$, and let $x\in \bT$ be a semisimple element of order $n$. If $\lambda_r<n/2$ then the eigenvalues of $x$ on $M$ are the same as the weight spaces of $\bT$. In particular, $x$ is a blueprint for $M$.
\end{lemma}
\begin{proof} Since all maximal tori are conjugate, we may assume that $x$ is the matrix
\[ \begin{pmatrix}\zeta & 0\\0&\zeta^{-1}\end{pmatrix},\]
where $\zeta$ is a primitive $n$th root of $1$. The eigenvalues of $x$ on $L(1)$ are $\zeta,\zeta^{-1}$, and the eigenvalues of $x$ on $L(\lambda_i)$ are roots of unity $\zeta^{\pm j}$ for $0\leq j\leq \lambda_i$. If $\lambda_i<n/2$ for all $i$ then the eigenspaces of $x$ are simply the weight spaces of the $L(\lambda_i)$, and so $x$ and $\bT$ stabilize the same subspaces of $M$, thus $x$ is a blueprint for $M$.
\end{proof}

We will apply this lemma to $A_1$ subgroups of algebraic groups. We often will end up with composition factors that do not precisely satisfy the hypotheses of this lemma though: if one composition factor has slightly larger highest weight, then although the eigenspaces do not correspond to weight spaces, with some weight spaces being merged, these all take place within one composition factor of the module, and so the finite subgroup $A_1(q)$ of the $A_1$ is still a blueprint for the module in question, even if the element of order $n$ is not.

\begin{lemma}\label{lem:q/2restricted} Let $\bG$ be the simply connected form of an exceptional algebraic group, and let $\bX$ be a positive-dimensional subgroup of $\bG$ of type $A_1$. Let $x$ be a semisimple element of $\bX$ of order $n$. Let $M$ be a module for $\bG$.
\begin{enumerate}
\item\label{lemi:q2resa} If the composition factors of $\bX$ on $M$ are $n/2$-restricted then $x$ and a maximal torus $\bT$ containing $x$ stabilize the same subspaces of $M$, so that $x$ is a blueprint for $M$.
\item\label{lemi:q2resb} Suppose that the highest weights of $\bX$ on $M$ are $\lambda_1,\dots,\lambda_r$, with $\lambda_i\leq \lambda_{i+1}$, and let $H=A_1(q)$ be a finite subgroup of $\bX$ containing $x$. If $\lambda_{r-1}+\lambda_r<n$ then $H$ and $\bX$ stabilize the same subspaces of $M$, so that $H$ is a blueprint for $M$.
\end{enumerate}
\end{lemma}
\begin{proof} The first part follows immediately from Lemma \ref{lem:restrictedsl2}, so we concentrate on the second statement. Letting $\bT$ be a maximal torus of $\bX$ containing $x$, if $\lambda$ and $\mu$ are two weights of $\bT$ on $M$ that are equal when taken modulo $n$ (i.e., yield the same eigenvalue for the action of $x$), then $\lambda$ and $\mu$ differ by a multiple of $n$. By assumption on the $\lambda_i$, since $\lambda-\mu\geq n$, both $\lambda$ and $\mu$ must be weights for the composition factor $L(\lambda_r)$, since if $\lambda$ is a weight for one of the other $L(\lambda_i)$ then it lies between $-\lambda_i$ and $+\lambda_i$, and cannot differ by $n$ from any other weight for any other $\lambda_j$.

Let $N$ be any $kH$-submodule of $M$. If $N$ does not contain the factor $L(\lambda_r)$ then the eigenvectors of $x$ on $N$ all come from weight vectors for $\bT$, by the previous paragraph, and so $\bT$ stabilizes $N$. If $N$ contains $L(\lambda_r)$, then $\bT$ also stabilizes $N$ by taking duals: $\bT$ stabilizes the submodule $\Ann(N)$ of $M^*$, which is isomorphic to $(M/N)^*$. Thus $\bT$ stabilizes every $H$-submodule of $M$, and so $\gen{\bT,H}=\bX$ and $H$ stabilize the same subspaces of $M$, as claimed.
\end{proof}

\section{Traces of Modules for \texorpdfstring{$\PGL_2$}{PGL2}}
\label{sec:pgl}
For this section $p$ is odd. Here we produce a technical result about extending simple modules for $\PSL_2(p^a)$ to $\PGL_2(p^a)$, and the traces and eigenvalues of the elements on such an extension.

There are two extensions of any simple module from $\PSL_2(p^a)$ to $\PGL_2(p^a)$. We will give a way of telling these apart if the dimension of the simple module is odd, which is all that we need in what follows.

If $L(i)$ is a simple module for $\PSL_2(p^a)$ of odd dimension, then of the two extensions of $L(i)$ to $\PGL_2(p^a)$, for one all defining matrices have determinant $1$ and for the other half have determinant $-1$: to see this, note that all elements in $\PSL_2(p^a)$ act with determinant $1$, and the non-trivial $1$-dimensional representation acts like $-1$ for elements outside of $\PSL_2(p^a)$, so a given extension and its product with this $1$-dimensional representation give us the two cases. Write $L(i)^+$ for the module for $\PGL_2(p^a)$ for which all matrices have determinant $+1$, and $L(i)^-$ for the other extension. This notation will be used in the proof of the next lemma.

\begin{lemma}\label{lem:pglnotappearing} Let $p$ be an odd prime and $a\geq 1$ an integer. Let $M$ be a simple module for $H=\PGL_2(p^a)$ with Brauer character $\phi$, and let $g$ be an element of order $p^a\pm 1$ in $H$. Let $t$ be an involution in $\PSL_2(p^a)$, and let $h$ be the involution in $\gen g$.
\begin{enumerate}
\item There are two conjugacy classes of involutions in $H$. If $o(g)$ is twice an odd number then $t$ and $h$ are representatives of these two classes, and otherwise $t$ and $h$ are conjugate.
\item If $\dim(M)$ is even then $\phi(t)=\phi(h)=0$.
\item If $\dim(M)$ is odd, then the dimensions of the $(+1)$-eigenspace and $(-1)$-eigenspace of the action of $g$ differ by $1$.
\item If $\dim(M)$ is odd, then $\phi(t)=\pm\phi(h)$. Furthermore, if $+1$ is not an eigenvalue of $g$ on $M$, then $\phi(t)=\phi(h)$ if and only if $t$ and $h$ are conjugate. If $-1$ is not an eigenvalue of $g$ on $M$ then $\phi(t)=\phi(h)$.
\end{enumerate}
\end{lemma}
\begin{proof}
\begin{enumerate}
\item That $H$ has two classes of involutions is well known, and one is a class of complements, the other is in $\PSL_2(p^a)$. Thus the second statement follows easily.
\item We use Steinberg's tensor product theorem, lifting all modules to $\GL_2(p^a)$: $M$ has even dimension if and only if, as a tensor product, at least one of the factors has even dimension, and the Brauer character is $0$ for a given element if and only if one of the factors has Brauer character $0$ for the same element. Thus we need to check this for the symmetric powers of the natural module $S^i(M')$ for $0\leq i\leq p-1$ (where $M'$ is the natural module for $\GL_2$), where it is trivial to see that the trace of an involution is $0$ on even-dimensional modules and $\pm 1$ on modules of odd dimension.

\item Since $\dim(M)$ is odd and $M$ is self-dual, of course one of $\pm 1$ is an eigenvalue for the action of $g$ on $M$, and the dimensions of the $(+1)$- and $(-1)$-eigenspaces must differ by an odd number. It is an easy exercise to compute the eigenvalues of $g$ on the Steinberg module $L(p^a-1)^\pm$ for $\PGL_2(p^a)$, and we see that these are all distinct if $g$ has order $p^a+1$, and if $g$ has order $p^a-1$ then $\pm 1$ appears twice and $\mp1$ appears once.

For other modules, from the definition of the Steinberg module as a tensor product of twists of the $p$-restricted module $L(p-1)^\pm$, and the fact that the eigenvalues of $g$ on $L(i)^\pm$ all appear in $L(p-1)^\pm$, we see that the eigenvalues of $g$ on any simple module appear in the eigenvalues of $g$ on the Steinberg. Thus the result holds since the sum of the dimensions of the $(+1)$- and $(-1)$-eigenspaces must be odd.

\item Return to $\PGL_2(p^a)$, and suppose that $t$ and $h$ are not conjugate, so that $g$ has twice odd order. As we saw above, for a $p$-restricted module $L(i)^\pm$ for $0\leq i\leq p-1$ an even integer, the Brauer character values of $L(i)^+$ on $t$ and $h$ have the same sign, and the Brauer character values of $L(i)^-$ on $t$ and $h$ have opposite signs. Notice that $+1$ is an eigenvalue of $M$ if and only if there are an even number of minus-type modules in the tensor decomposition, and this happens if and only if $\phi(t)=\phi(h)$, as needed.
\end{enumerate}
\end{proof}

Using this, the following result is now clear.

\begin{corollary}\label{cor:sl2nopgl} Let $p$ be an odd prime and $a\geq 1$ an integer. Let $\bG$ be the simply connected form of $E_7$, and let $H$ be a copy of $\SL_2(p^a)$ in $\bG$ with $Z(\bG)=Z(H)$. Suppose that $g$ is an element of $\bG$ such that $o(g)=p^a\pm 1$ is twice an odd number, and $-1$ is not an eigenvalue for the action of $g$ on $L(E_7)$. Then the group $\bar H=\gen{H,g}$ does not satisfy $\bar H/Z(H)=\PGL_2(p^a)$.
\end{corollary}
\begin{proof} Suppose that $\bar H/Z(H)=\PGL_2(p^a)$. Let $t$ be an element of $H$ that is an involution in $H/Z(H)$, so $o(t)=4$. The trace of $t$ on $L(E_7)$ is $-7$ or $25$, depending on the class of $t$ in $\bG$. The involution $h$ in $\gen g$ has trace $5$ on $L(E_7)$, since it is an involution in $\bG$ rather than $\bG/Z(\bG)$. We now show that $h$ and $t$ must in fact have the same trace, which is a contradiction. By Lemma \ref{lem:pglnotappearing}, any even-dimensional composition factors of $L(E_7)\downarrow_H$ yield trace $0$ for both $t$ and $h$, and they have the same trace on odd-dimensional composition factors since $-1$ is not an eigenvalue of $g$ on $L(E_7)$. Thus the trace of $t$ and $h$ on $L(E_7)$ is the same, but this is a contradiction.
\end{proof}

This is how we will use Corollary \ref{cor:sl2nopgl}: suppose that $H\cong \SL_2(p^a)$ lies inside the simply connected form $\bG$ of $E_7$ with $Z(\bG)=Z(H)$. We find some element $x\in \bG$ such that $x$ and $H$ both stabilize some proper subspace $W$ of either $M(E_7)$ or $L(E_7)$, so that $\gen{H,x}$ is not of the same type as $\bG$. If $-1$ is not an eigenvalue of $x$ on $L(E_7)$ then $\gen{H,x}$, modulo $Z(\bG)$, is not $\PGL_2(p^a)$ either. If the composition factors of $L(E_7)\downarrow_H$ are not invariant under a field automorphism of $H$, then $x$ cannot induce a field automorphism on $H$, and therefore $\gen{H,x}$ is not almost simple modulo $Z(\bG)$ either. If, in addition, $H$ is a maximal member of $\mathscr P$, then we apply Proposition \ref{prop:maximalnotinP} to see that $H$ is Lie imprimitive. Furthermore, if we chose $W$ so that it is stable under $N_{\Aut^+(\bG)}(H)$, or our element $x$ stabilizes all submodules in the $N_{\Aut^+(\bG)}(H)$-orbit of $W$, then $H$ is also strongly imprimitive.

\section{The Graph Automorphism of \texorpdfstring{$F_4$}{F4}}

In this short section we describe how semisimple elements of odd order in $F_4$ react to the graph automorphism in characteristic $2$. Since the graph automorphism $\tau$ does not stabilize the minimal module $M(F_4)$, and $L(F_4)$ has composition factors $M(F_4)$ and $M(F_4)^\tau$, we can see the effect of the graph automorphism on semisimple classes by taking the eigenvalues of an element $x$ on $L(F_4)$ and removing those from $M(F_4)$.

Since the graph automorphism squares to a field automorphism, however, it is slightly more complicated to understand those classes that are left invariant under a graph automorphism, since we need to check whether the eigenvalues of $x^\tau$ and $x^i$ match for some $i$, rather than whether the eigenvalues of $x$ and $x^\tau$ match. This is still not difficult using a computer, however; we give two special cases, where a conjugacy class is stable under the graph automorphism (up to powers) and where the classes have integral traces.

\begin{lemma} Let $k$ be a field of characteristic $2$. Let $x$ be a semisimple element in $G=F_4(k)$ such that $x^\tau$ is conjugate to a power of $x$. If $x$ has order at most $9$, then a power of $x$ has trace on $M(F_4)$ given below.

\begin{center}
\begin{tabular}{cc}
\hline $o(x)$ & Possible traces on $M(F_4)$
\\\hline $3$ & $-1$
\\ $5$ & $1$
\\ $7$ & $4(\zeta_7+\zeta_7^{-1})+3(\zeta_7^2+\zeta_7^{-2})+5$, $-(\zeta_7+\zeta_7^{-1})$
\\ $9$ & $2-3(\zeta_9+\zeta_9^{-1})$ 
\\\hline
\end{tabular}
\end{center}
\end{lemma}

\begin{lemma}\label{lem:swappingintegralf4} Let $k$ be a field of characteristic $2$. Let $x$ be a semisimple element in $G=F_4(k)$ such that the trace of both $x$ and $x^\tau$ is an integer, and $x$ and $x^\tau$ are not conjugate. If $x$ has order at most $9$, then the traces of $x$ and $x^\tau$ are as below, where we give $x$ up to graph automorphism.

\begin{center}
\begin{tabular}{ccc}
\hline $o(x)$ & Trace of $x$ on $M(F_4)$ & Trace of $x^\tau$ on $M(F_4)$
\\\hline $3$ & $8$ & $-1$
\\ $5$ & None & None
\\ $7$ & $-2$ & $5$ 
\\ $9$ & $-1$ ($x^3$ has trace $-1$) & $2$ ($x^3$ has trace $-1$)
\\\hline
\end{tabular}
\end{center}
\end{lemma}

\section{Rank-1 Subalgebras of the Lie Algebra}
\label{sec:sl2subalgebra}

In this section we consider subalgebras of $L(\bG)$, specifically $\slf_2$-subalgebras. The stabilizers of $\slf_2$-subalgebras will be shown to be positive dimensional if $p$ is not too small. Thus if a subgroup $H$ of $\bG$ stabilizes a unique $3$-space on $L(\bG)$, and this is an $\slf_2$-subalgebra, then $H$ is Lie imprimitive, and must be strongly imprimitive via Proposition \ref{prop:intersectionstabilizers}, since the subspace is unique.

To begin with, we prove a proposition that gives us a criterion for a subgroup $H$ to stabilize an $\slf_2$-subalgebra in the first place. This proposition is a restatement of results of Ryba from \cite{ryba2002}, particularly Lemma 10 from that paper.

\begin{proposition}\label{prop:sl2ifsplitoff} Let $V$ be a $3$-dimensional subspace of $L(\bG)$, and let $H$ be a subgroup of $\bG$ such that $HZ(\bG)/Z(\bG)=\PSL_2(p^a)$ for some $p\geq 5$. If $V$ is $H$-stable and a complement for $V$ is also $H$-stable (i.e., $V$ is a summand of $L(\bG)\downarrow_H$), and $\Hom_{kH}(V,L(\bG))$ is $1$-dimensional (i.e., there are no other submodules of $L(\bG)\downarrow_H$ isomorphic to a quotient of $V$) then $V$ is a subalgebra of $L(\bG)$ isomorphic to $\slf_2$.
\end{proposition}
\begin{proof} As $\Hom_{kH}(V,L(\bG))$ is $1$-dimensional, $V$ is a non-trivial module, and therefore $V$ is isomorphic to the module $3_i$ for some $i$ for $H/Z(H)$. In particular, it is self-dual.

Suppose that $L(\bG)\downarrow_H$ has a unique submodule isomorphic to $V$, and that this is a summand, so that the quotient $L(\bG)\downarrow_H/V$ has no quotient isomorphic to $V^*\cong V$. By \cite[Lemma 6]{ryba2002}, we have that $V$ possesses a non-singular trace form, and then we apply Block's theorem \cite{block1962} to see that $V$ is a simple Lie algebra of type $\slf_2$.
\end{proof}

In order to use this proposition, we need to know something about the $\slf_2$-subalgebras of the Lie algebras of exceptional groups. The following is a theorem of Stewart and Thomas \cite[Theorem 1.2]{stewartthomas2018}, specialized to the case of $\bG=E_6,E_7,E_8$, for use in this paper and a later one on $\SL_2$-subgroups of $E_8$.

\begin{theorem}\label{thm:sl2subalgebra}
Let $\bG=E_6$ and $p\geq 7$, or $\bG=E_7,E_8$ and $p\geq 11$. The classes of $\slf_2$-subalgebras of $L(\bG)$ are in one-to-one correspondence with the nilpotent orbits of $L(\bG)$, with a bijection being realized by sending an $\slf_2$-subalgebra to the nilpotent orbit of largest dimension intersecting it non-trivially.
\end{theorem}

To prove \cite[Theorem 1.2]{stewartthomas2018}, representatives for $\slf_2$-triples $(e,h,f)$ for each class of $\slf_2$-subalgebra were constructed using a computer and the `Tools' referred to in \cite[Proof of Theorem 1.2]{stewartthomas2018}. David Stewart has kindly shared this computer file, and so it was possible to compute the Jordan block structures of the nilpotent elements $f$ in each triple. If $e^{[p]}=0$ then $e$ and $f$ lie in the same orbit, but if $e^{[p]}\neq 0$ then the orbits of $e$ and $f$ are different, $f^{[p]}=0$, and the precise orbit of $f$ for $p\geq 5$ is given in Tables \ref{t:sl2classesE6}, \ref{t:sl2classesE7} and \ref{t:sl2classesE8}. This yields the following corollary.

\begin{table}
\begin{center}
\begin{tabular}{cccc}
\hline Class & $p=5$ & $p=7$ & $p=11$
\\\hline$D_{4}$&$3A_{1}$&&
\\$A_{5}$&$A_{3}$&&
\\$D_{5}(a_{1})$&$A_{2}+A_{1}$&&
\\$E_{6}(a_{3})$&$A_{3}+A_{1}$&&
\\$D_{5}$&$A_{3}$&$A_{3}+A_{1}$&
\\$E_{6}(a_{1})$&$2A_{2}+A_{1}$&$A_{5}$&
\\$E_{6}$&$A_{2}+2A_{1}$&$2A_{2}+A_{1}$&$A_5$
\\ \hline
\end{tabular}
\end{center}
\caption{Second nilpotent class intersecting a non-restricted $\slf_2$-subalgebra of $L(E_6)$ for $p\geq 5$}
\label{t:sl2classesE6}
\end{table}

\begin{table}\begin{center}
\begin{tabular}{cccccc}
\hline Class & $p=5$ & $p=7$ & $p=11$ & $p=13$ & $p=17$
\\\hline $D_{4}$&$(3A_{1})'$&&&&
\\$(A_{5})''$&$A_{3}$&&&&
\\$D_{4}+A_{1}$&$4A_{1}$&&&&
\\$D_{5}(a_{1})$&$A_{2}+A_{1}$&&&&
\\$(A_{5})'$&$A_{3}$&&&&
\\$A_{5}+A_{1}$&$(A_{3}+A_{1})'$&&&&
\\$D_{5}(a_{1})+A_{1}$&$A_{2}+2A_{1}$&&&&
\\$D_{6}(a_{2})$&$D_{4}(a_{1})+A_{1}$&&&&
\\$E_{6}(a_{3})$&$(A_{3}+A_{1})'$&&&&
\\$D_{5}$&$A_{3}$&$(A_{3}+A_{1})'$&&&
\\$E_{7}(a_{5})$&$A_{3}+A_{2}$&&&&
\\$A_{6}$&$A_{2}+2A_{1}$&&&&
\\$D_{5}+A_{1}$&$(A_{3}+A_{1})'$&$A_{3}+2A_{1}$&&&
\\$D_{6}(a_{1})$&$A_{3}+2A_{1}$&$D_{4}(a_{1})+A_{1}$&&&
\\$E_{7}(a_{4})$&$D_{4}(a_{1})+A_{1}$&$A_{3}+A_{2}$&&&
\\$D_{6}$&$A_{3}$&$A_{3}+2A_{1}$&&&
\\$E_{6}(a_{1})$&$2A_{2}+A_{1}$&$(A_{5})'$&&&
\\$E_{6}$&$A_{2}+2A_{1}$&$2A_{2}+A_{1}$&$(A_{5})'$&&
\\$E_{7}(a_{3})$&$(A_{3}+A_{1})'$&$D_{4}(a_{1})+A_{1}$&&&
\\$E_{7}(a_{2})$&$2A_{2}+A_{1}$&$A_{3}+2A_{1}$&$D_{6}(a_{2})$&&
\\$E_{7}(a_{1})$&$A_{2}+2A_{1}$&$(A_{5})'$&$D_{6}$&$D_6(a_1)$&
\\$E_{7}$&$A_{4}+A_{2}$&$A_{6}$&$A_{5}+A_{1}$&$D_5+A_1$&$D_6$
\\ \hline
\end{tabular}
\end{center}
\caption{Second nilpotent class intersecting a non-restricted $\slf_2$-subalgebra of $L(E_7)$ for $p\geq 5$}
\label{t:sl2classesE7}
\end{table}

\begin{table}\begin{center}
\begin{tabular}{ccccc}
\hline Class & $p=5$ & $p=7$ & $p=11$ & $p=13$
\\\hline $D_{4}$&$3A_{1}$&&&
\\$A_{5}$&$A_{3}$&&&
\\$D_{4}+A_{2}$&$A_{2}+3A_{1}$&&&
\\$E_{6}(a_{3})$&$A_{3}+A_{1}$&&&
\\$D_{5}$&$A_{3}$&$A_{3}+A_{1}$&&
\\$A_{5}+A_{1}$&$A_{3}+A_{1}$&&&
\\$D_{5}(a_{1})+A_{2}$&$2A_{2}+A_{1}$&&&
\\$E_{6}(a_{3})+A_{1}$&$A_{3}+2A_{1}$&&&
\\$D_{5}+A_{1}$&$A_{3}+A_{1}$&$A_{3}+2A_{1}$&&
\\$E_{8}(a_{7})$&$A_4+A_3$&&&
\\$A_{6}+A_{1}$&$A_{2}+3A_{1}$&&&
\\$E_{6}(a_{1})$&$2A_{2}+A_{1}$&$A_{5}$&&
\\$D_{5}+A_{2}$&$A_{3}+A_{2}$&$A_{3}+A_{2}+A_{1}$&&
\\$E_{6}$&$A_{2}+2A_{1}$&$2A_{2}+A_{1}$&$A_{5}$&
\\$D_{7}(a_{2})$&$2A_{3}$&$A_{4}+A_{1}$&&
\\$A_{7}$&$2A_{2}+A_{1}$&$A_{5}$&&
\\$E_{6}(a_{1})+A_{1}$&$2A_{2}+2A_{1}$&$A_{5}+A_{1}$&&
\\$E_{7}(a_{3})$&$A_{3}+A_{1}$&$D_{4}(a_{1})+A_{1}$&&
\\$E_{8}(b_{6})$&$D_{4}(a_{1})+A_{1}$&$E_{7}(a_{5})$&&
\\$D_{7}(a_{1})$&$A_{3}+2A_{1}$&$A_{3}+A_{2}$&&
\\$E_{6}+A_{1}$&$A_{2}+3A_{1}$&$2A_{2}+2A_{1}$&$A_{5}+A_{1}$&
\\$E_{8}(a_{6})$&$2A_{3}$&$A_{4}+2A_{1}$&&
\\$D_{7}$&$A_{2}+3A_{1}$&$A_{5}$&$D_{5}+A_{1}$&
\\$E_{8}(b_{5})$&$2A_{2}+2A_{1}$&$D_{4}(a_{1})+A_{1}$&$E_{7}(a_{5})$&
\\$E_{7}(a_{1})$&$A_{2}+2A_{1}$&$A_{5}$&$D_{6}$&$D_{6}(a_{1})$
\\$E_{8}(a_{5})$&$2A_{2}+A_{1}$&$E_{6}(a_{3})+A_{1}$&$E_{7}(a_{4})$&
\\$E_{8}(b_{4})$&$A_{2}+3A_{1}$&$A_{5}+A_{1}$&$E_{7}(a_{3})$&$E_{7}(a_{4})$
\\$E_{8}(a_{4})$&$A_{4}+A_{3}$&$A_{5}$&$D_{6}(a_{2})$&$E_{7}(a_{3})$
\\$E_{8}(a_{3})$&$A_{4}+A_{2}+A_{1}$&$A_{6}+A_{1}$&$E_{6}(a_{3})+A_{1}$&$D_{6}(a_{1})$
\\$E_{8}(a_{2})$&$2A_{3}$&$A_{4}+A_{3}$&$A_{5}+A_{1}$&$D_{7}$
\\$E_{8}(a_{1})$&$2A_{2}+2A_{1}$&$A_{4}+A_{2}+A_{1}$&$D_{5}(a_{1})+A_{2}$&$A_{7}$
\\$E_{8}$&$A_{3}+A_{2}+A_{1}$&$A_{3}+A_{2}+A_{1}$&$A_{4}+A_{3}$&$A_{6}+A_{1}$
\\ \hline
\\ \hline Class & $p=17$ & $p=19$ & $p=23$ & $p=29$
\\ \hline $E_8(a_3)$ & $E_7(a_3)$&&&
\\ $E_8(a_2)$ & $E_7(a_1)$&$E_7(a_2)$&&
\\$E_{8}(a_{1})$&$D_{7}$&$E_{7}$&$E_{7}(a_{1})$&
\\$E_{8}$&$A_{7}$&$E_{6}+A_{1}$&$D_{7}$&$E_{7}$
\\ \hline
\end{tabular}
\end{center}
\caption{Second nilpotent class intersecting a non-restricted $\slf_2$-subalgebra of $L(E_8)$ for $p\geq 5$. (Missing classes, $D_4+A_1$, $D_5(a_1)$, $D_5(a_1)+A_1$, $D_6(a_2)$, $E_7(a_5)$, $A_6$, $D_6(a_1)$, $E_7(a_4)$, $D_6$, $E_7(a_2)$ and $E_7$, are exactly as in Table \ref{t:sl2classesE7})}
\label{t:sl2classesE8}
\end{table}

\begin{corollary}\label{cor:sl2coxeter-1} Let $\bG=E_6$ and $p\geq 7$, or $\bG=E_7,E_8$ and $p\geq 11$. Let $H\cong\PSL_2(p^a)$ be a subgroup of $\bG/Z(\bG)$. If $H$ stabilizes an $\slf_2$-subalgebra $\mathfrak h$ of $L(\bG)$, then $\mathfrak{h}$ is restricted and $H$ is Lie imprimitive. If, in addition, $H$ stabilizes a unique $3$-space on $L(\bG)$ then $H$ is strongly imprimitive.
\end{corollary}
\begin{proof} Since there is a unique conjugacy class of subgroups $\PSL_2(p^a)$ inside $\PSL_2$, we see that because the standard $\PSL_2(p^a)$ inside $\PSL_2$ swaps $e$ and $f$ in an $\slf_2$-triple $(e,h,f)$ for $\mathfrak{h}$, thus $H$ must swap the two nilpotent orbits of $\mathfrak{h}$. From the discussion above this result, if $\mathfrak{h}$ is not restricted then $e$ and $f$ lie in different nilpotent orbits of $L(\bG)$, so $\mathfrak{h}$ must be restricted.

Now we argue exactly as in \cite[Proof of Theorem 1.6]{stewartthomas2018}, using \cite[Proposition 4.1]{seitz2000} to find $\bX$ a good $A_1$ in $\bG$ that stabilizes $\mathfrak{h}$. Hence $\gen{H,\bX}$ is positive dimensional and stabilizes $\mathfrak h$, so $H$ is Lie imprimitive. If $H$ stabilizes a unique $3$-space on $L(\bG)$ then Proposition \ref{prop:intersectionstabilizers} shows that $\bX$ is $N_{\Aut^+(\bG)}(H)$-stable, and therefore $H$ is strongly imprimitive.
\end{proof}

\chapter{Modules for \texorpdfstring{$\SL_2$}{SL2}}
\label{ch:sl2modules}
The purpose of this chapter is to describe everything we need to know about the simple modules and extensions between them for the groups $\SL_2(p^a)$ for $p$ a prime and $a\geq 1$.

\section{Modules for \texorpdfstring{$\SL_2(2^a)$}{SL(2,2a)}}

We construct certain modules for $H=\SL_2(2^a)$ for some $a\leq 10$, and prove that various configurations of module do not exist. (The reason we choose $a\leq 10$ is so that these results may be used in work for $E_8$, for which $v(E_8)=t(E_8)=1312$.) The main motivation for this is to achieve better bounds on the number of occurrences of certain composition factors that are needed to prevent a particular simple module appearing in the socle of a given module $M$.

We begin with some notation. Let $u$ be an element of order $2$ in $H$. Denote by $1$ the trivial module. By $2_1$ we denote the natural module for $H$, and define $2_i$ by the equation
\[ 2_{i-1}^{\otimes 2}=1/2_i/1,\]
i.e., $2_i$ is the twist under the field automorphism of $2_{i-1}$. Given this, if $I$ is a subset of $\{1,\dots,a\}$, of cardinality $b$, we define, 
\[2^b_I=\bigotimes_{i\in I} 2_i,\]
for example, $4_{1,2}=2_1\otimes 2_2$; the modules $2^b_I$ for all $I\subseteq \{1,\dots,a\}$ furnish us with a complete set of irreducible modules for $H$, by Steinberg's tensor product theorem.

We first recall a result of Alperin \cite[Theorem 3]{alperin1979}, that determines the dimensions of $\Ext^1(A,B)$ for $A,B$ simple modules for $H$.

\begin{lemma}\label{lem:extforsl22a} Let $A$ and $B$ be simple $H$-modules, corresponding to the subsets $I$ and $J$ of $\{1,\dots,a\}$. The dimension of $\Ext^1(A,B)$ is always $0$, unless
\begin{enumerate}
\item $|I\cap J|+1=|I\cup J|<a$, and
\item if $i\in I\cup J$ and $i-1\notin I\cap J$, then $i-1\notin I\cup J$,
\end{enumerate}
and in this case the dimension is $1$.

In particular, if $\Ext^1(A,B)\neq 0$ then the dimension of $A$ is either half or double that of $B$.
\end{lemma}

Using Lemma \ref{lem:pressure}, we see that if a $kH$-module $M$ has no trivial submodule or quotient, then $M$ either has no trivial composition factors, or requires at least one more $2$-dimensional composition factor than trivial factor. We can do better than this in some circumstances.

If a module has pressure $1$, then we can still say something about the module. This is important for $F_4$ and $E_6$ because there are no involutions acting projectively on $M(F_4)$ or $M(E_6)$ (but there are involutions of $E_7$ acting projectively on $M(E_7)$). The next lemma is a special case of Lemma \ref{lem:char2pressuregeneral} below, but we provide a full proof in this simple case for the benefit of the reader.

\begin{lemma}\label{lem:char2pressure1} Let $M$ be a $kH$-module that has at least one trivial composition factor but no trivial submodules or quotients. If $M$ has pressure $1$, then an involution in $H$ acts projectively on $M$ if $\dim(M)$ is even and with a single Jordan block of size $1$ if $\dim(M)$ is odd.
\end{lemma}
\begin{proof} Note that, since $M$ has pressure $1$, it cannot have $2_i\oplus 2_j$ or $1^{\oplus 2}$ as a subquotient without stabilizing a line or hyperplane by Lemma \ref{lem:pressure}. We proceed by induction on $\dim(M)$, starting with the even-dimensional case. If $\dim(M)<6$ then $M$ cannot satisfy the hypotheses of the lemma, so our induction starts. We may assume that $\soc(M)=2_i$ for some $i$: first there are no composition factors of $\soc(M)$ of dimension greater than $2$ because the quotient by one would still satisfy the hypotheses of the lemma, and $2_i\oplus 2_j$ cannot be in the socle by the note above.

The quotient module $M/\soc(M)$ has pressure $0$, so $H$ must stabilize a line or hyperplane by Lemma \ref{lem:pressure}, but cannot stabilize a hyperplane by assumption, so $M/\soc(M)$ has a trivial submodule, and it must be unique by the note at the start of this proof. Quotient out by any possible factors of dimension at least $4$ in the socle of $M/\soc(M)$ to obtain a module $N$ of pressure $0$ and with $\soc(N)=1$. (If there is a $2_i$ in $\soc(N)$ then we find a submodule of pressure $2$, which is not allowed.)

The socle of the quotient module $N/\soc(N)$ must be $2_j$ for some $j$, since $2_j\oplus 2_l$ cannot be a subquotient and $1$ only has extensions with simple modules of dimension $2$. Now $N/\soc^2(N)$ again has pressure $0$, so by Lemma \ref{lem:pressure} has a trivial submodule as it cannot have a trivial quotient (it is a quotient of $M$), and we have constructed a submodule $1/2_j/1$ inside $N$. Letting $L$ be the quotient of $N$ by this submodule, we have removed $2_i,2_j,1^2$ from $M$, and possibly some other modules, and so an involution acts projectively on $L$ by induction, but it also acts projectively on the kernel of the map $N\to L$, namely $1/2_j/1$, and on the kernel of the map $M\to N$ since that has no trivial factors at all, so an involution acts projectively on all of $M$, as needed.

For odd-dimensional modules, we now simply find any submodule $N$ with a single trivial composition factor and such that $1$ is a quotient of $N$. The quotient module $M/N$ must have even dimension and has no trivial submodule as otherwise $M$ would have $1\oplus 1$ as a subquotient. Also, $N$ has pressure $0$ since otherwise $N$ with the $1$ removed from the top has pressure $2$, contradicting Lemma \ref{lem:pressure}. Hence $M/N$ has pressure $1$: thus an involution acts projectively on $M/N$ and with a single $1$ on $N$, as needed.
\end{proof}

We can generalize this result to modules of larger pressure, but for our proof of this we need a computation about modules for $\SL_2(2^a)$, which we have only checked in the range $2\leq a\leq 7$. It is certainly true for all $a$, but the author cannot see a purely theoretical proof. 

\begin{lemma}\label{lem:char2pressuregeneral} Let $2\leq a\leq 7$ and let $H=\SL_2(2^a)$. Let $M$ be a $kH$-module that has at least one trivial composition factor but no trivial submodule or quotient. If $M$ has pressure $n$ then an involution in $H$ acts on $M$ with at most $n$ Jordan blocks of size $1$.
\end{lemma}
\begin{proof} As with the previous lemma, we proceed by induction on $\dim(M)$, and note that if $\dim(M)\leq 4$ then $M$ cannot satisfy the hypotheses of the lemma. If $M$ is a minimal counterexample to the lemma, the socle and top of $M$ consist solely of $2$-dimensional composition factors. Notice that, by choice of minimal counterexample, there cannot exist a submodule $N$ such that $N$ has no trivial quotients and the quotient module $M/N$ has no trivial submodules, since otherwise one of $N$ and $M/N$ would also be a counterexample to the lemma, by a pressure argument using Lemma \ref{lem:pressure}.

Let $I$ denote the set of all simple $kH$-modules of dimension at least $4$. Let $N_1$ denote the $\{2_i\}$-radical of $M$ (so the socle), $N_2$ denote the preimage in $M$ of the $I$-radical of the quotient $M/N_1$, $N_3$ denote the preimage in $M$ of the $\{2_i\}$-radical of the quotient $M/N_2$, and let $N$ denote the preimage in $M$ of the $\{1\}$-radical of the quotient $M/N_3$.

We claim that all trivial composition factors of $N$ lie in the second socle layer. The proof of this claim uses a computer calculation in $P(2_1)$, which has been checked for $2\leq a\leq 7$. Build up a submodule $V_i$ of $P(2_i)$ by taking the socle, adding all copies of modules in $I$ on top of it, and then all $2$-dimensional factors on top of that, and finally all trivial factors on top of that. There is a single trivial composition factor of $V_i$, which lies in the second socle layer of $V_i$ (i.e., comes from a $1/2_i$ submodule). Applying this statement to the situation above, we see that $N$ is a submodule of a sum of $V_i$ for various $i$, and therefore all trivial composition factors of $N$ lie in the second socle layer of $N$.

Suppose that not all $2$-dimensional composition factors of $N$ lie in the socle, and let $\bar N$ denote a minimal submodule of $N$ subject to having a $2$-dimensional factor not in the socle and such that the quotient $N/\bar N$ has no trivial submodule. Such a module $\bar N$ has a single trivial composition factor and pressure $1$, hence the quotient $M/\bar N$ has pressure $n-1$ and no trivial submodule or quotient. Notice that therefore the action of $u$ on $M/\bar N$ has at most $n-1$ Jordan blocks of size $1$, and therefore $u$ acts on $M$ with at most $n$ Jordan blocks of size $1$, as needed.

Thus we may assume that we are in the following situation: all $2$-dimensional factors of $N$ lie in its socle, all trivial composition factors of $N$ lie in its second socle layer, all composition factors of $\soc(M/N)$ have dimension $2$ by Lemma \ref{lem:extforsl22a}.

This means that there is a $2$-dimensional submodule of $M/N$, so let $L$ denote its preimage in $M$, which therefore has a quotient $2_i/1$ for some $i$. If the quotient $M/L$ has no trivial submodules then $L$ is a submodule such that both $L$ and $M/L$ have no trivial submodules or quotients, contradicting the first paragraph of this proof. Therefore $M/L$ has a trivial submodule, and the preimage $L_1$ of this in $M$ must have a quotient $1/2_i/1$. Thus $M/L_1$ has no trivial submodules or quotients, and $u$ acts projectively on $L_1$. By induction $M/L_1$ satisfies the conclusion of the result, and since $u$ acts projectively on $L_1$ it has the same action on $(M/L_1)\oplus L_1$ as $M$, so $M$ satisfies the conclusion of the result.
\end{proof}

\begin{lemma}\label{lem:sl2832sfor2trivials} Let $a=3$. If $M$ is an even-dimensional module with $2n>0$ trivial composition factors and no trivial submodule or quotient, then it has at least $3n$ composition factors of dimension $2$.
\end{lemma}
\begin{proof} If $\dim(M)\leq 6$ then $M$ cannot satisfy the hypotheses of the lemma, so our induction starts. Note that if $M=M_1\oplus M_2$ with the $M_i$ both even-dimensional, then by induction $M$ satisfies the conclusion of the lemma: thus $M$ is either indecomposable or the sum of two odd-dimensional indecomposable modules.

The projective cover of $2_1$ is
\[ 2_1/1,4_{1,3}/2_1,2_2,2_3/1,1,4_{2,3}/2_1,2_2,2_3/1,4_{1,3}/2_1.\]
Remove any $4$-dimensional factors from the top and socle of $M$, so that $M$ is a submodule of a sum of copies of projectives $P(2_i)$. If $M$ has seven socle layers then it must have a submodule $P(2_i)$ for some $i$. This must be a summand, since projective submodules are always summands. Thus $M=P(2_i)$ and we are done. Hence $M$ has at most five socle layers. The number of $2$-dimensional factors in the first and third socle layers must be at least as many as the number of $1$s in the second layer, and there are at least as many $2$s in the third and fifth socle layers as $1$s in the fourth layer. We therefore must have that there are at least $3n$ composition factors of dimension $2$ in $M$, as claimed.
\end{proof}

Lemma \ref{lem:smallsubspaces} shows that if a subgroup $H$ of $\bG$ stabilizes a $1$- or $2$-space on $M(\bG)$ for $F_4$, $E_6$ and $E_7$ then $H$ is not Lie primitive. By Lemma \ref{lem:extforsl22a} we see that $2$-dimensional modules have non-split extensions only with modules of dimension $1$ and $4$, so we would like a similar result to the previous one, counting the number of $4$-dimensional factors in a module $M$ that has $2$-dimensional composition factors but no $2$-dimensional submodules or quotients. We start with the easier case, where there are no trivial composition factors in $M$ at all. Notice that we can use $\mathcal M$-pressure here as well, but we can do a bit better using the structure of modules for $\SL_2(2^a)$.

(We do not need to consider $a>6$ here as these lemmas do not appear to be of use for $E_8$: the stabilizers of $2$-spaces of $L(E_8)$ are not obviously positive dimensional.)

\begin{lemma}\label{lem:largest2and4} Let $H=\SL_2(2^a)$ for $4\leq a\leq 6$. The largest submodule of $P(4_{i,j})$ whose composition factors have dimension $2$ and $4$ is as follows: for $j=i\pm 1$, we have a $10$-dimensional module
\[ 4_{i-1,i+1}/2_{i+1}/4_{i,i+1};\]
for $a=4$ we have a $28$-dimensional module
\[ 4_{1,3},4_{1,3}/2_1,2_3/4_{1,4},4_{2,3}/2_1,2_3/4_{1,3};\]
for $j=i\pm2$ and $a>4$ we have a $32$-dimensional module
\[   4_{i,i+2},4_{i,i+2}/2_i,2_{i+2}/4_{i+1,i+2},4_{i,i+3},4_{i-1,i+2}/2_i,2_{i+2}/4_{i,i+2},\]
with $4_{i-1,i+2}$ as a quotient. In all other cases, we have the module
\[   4_{i,j},4_{i,j}/2_i,2_j/4_{i,j-1},4_{i,j+1},4_{i+1,j},4_{i-1,j}/2_i,2_j/4_{i,j},\]
with $4_{i,j-1}$ and $4_{i-1,j}$ as quotients.

Consequently, if $M$ is a module with no trivial or $8$-dimensional composition factors, with $c>0$ composition factors of dimension $2$, and no $2$-dimensional submodule or quotient, then $M$ has at least $c$ composition factors of dimension $4$.
\end{lemma}
\begin{proof} The statements for individual $a$ are verified by computer, so we concentrate on the consequence. If $M$ has no $8$-dimensional composition factors, then $M$ splits as the direct sum of two modules, one with composition factors of dimensions $1$, $2$ and $4$, (although there are no trivial factors in $M$) and one of dimensions $16$ and above, which we can ignore. Thus $M$ can be assumed to only have factors of dimensions $2$ and $4$. As $M$ has no $2$-dimensional submodules by hypothesis, it is a submodule of a sum of modules of the above form.

We cannot produce a module $4/2,2/4/2,2/4$ since the $4$s in the middle of the modules above do not have extensions with both $2$s by Lemma \ref{lem:extforsl22a}. Thus we have at least $4/2,2/4,4/2,2/4$, and so we need as many $4$s as $2$s.
\end{proof}

Of course, unlike the $2_i$, the $4_{i,j}$ are not all the same up to field automorphism. Thus for specific choices of composition factors, it is possible to achieve better bounds than the previous lemma.

The next lemma considers the case where we want to know how many $1$s and $2$s we can stack on top of a given simple module of dimension $4$. This lemma gives that answer, and hence how many $4$s one needs to `hide' all $1$ and $2$s inside the middle of the module.

\begin{lemma}\label{lem:largestsubmodofP4}
Let $H=\SL_2(2^a)$ for some $2\leq a\leq 6$. The largest submodule of $P(4_{1,2})$ whose composition factors modulo the socle have dimensions $1$ or $2$ is
\[ 2_2/1/2_3/1/2_2/4_{1,2},\]
and an involution acts projectively on this module.

For $a=4$ and $a\geq 5$ we have
\[ 2_2,2_4/1/2_1,2_3/4_{1,3},\quad\text{and}\quad 2_2/1/2_1,2_3/4_{1,3}\]
respectively. For $a=6$ and $i=4,5$ we have $1/2_1,2_i/4_{1,i}$.

In particular, if $M$ is a module for $H$ with no trivial or $2$-dimensional submodules or quotients, and it has $2n$ trivial composition factors for some $n>0$, then it has $n'\geq n+1$ factors of dimension $4$, and between $2n+1$ and $4(n'-1)$ composition factors of dimension $2$.
\end{lemma}
\begin{proof} The facts about the largest submodule of $P(4_{i,j})$ can easily be checked with a computer. For the conclusion, we proceed by induction, with the result holding vacuously if $\dim(M)< 14$, since no module can satisfy the hypotheses of the lemma. By removing submodules and quotients of dimension $8$ and above, we may assume that the socle and top of $M$ consist entirely of $4$-dimensional modules.

Let $M_1=\soc(M)$ and $M_2$ be the preimage in $M$ of the $\{4_{i,j}\}'$-radical of the quotient module $M/M_1$. There are no composition factors of dimension $4$ in the quotient $M_2/M_1$, and there are no extensions between simple modules of dimensions at most $2$ and at least $8$ by Lemma \ref{lem:extforsl22a}. Hence $M_2/M_1$ is the direct sum of its $\{1,2_i\}$-radical and $\{1,2_i,4_{i,j}\}'$-radical. Let $M_2'$ denote the preimage in $M$ of the $\{1,2_i\}$-radical of $M_2/M_1$, so that the quotient $M_2'/M_1$ only has composition factors of dimension $1$ and $2$, and $M_1$ only has composition factors of dimension $4$.

The module $M_2'$ is therefore a submodule of a sum of modules as in the first part of the lemma. Thus if $M_2$ has $2m$ trivial modules then $M_1$ has at least $m$ copies of $4$-dimensional modules to support the $2m$ trivials, and from the structure of the modules in the lemma the number of $2$-dimensionals is at most $4m$. Thus our result holds for $M_2'$, and hence for $M_2$ as the number of factors of dimensions $1$, $2$ and $4$ are the same in $M_2'$ and $M_2$.

Notice that the quotient $M/M_2$ also has no trivial or $2$-dimensional submodules or quotients, so satisfies the conclusion of the lemma by induction. Thus there are at least $n+1$ different $4$-dimensional factors in $M$ and at most $4n'-1$ factors of dimension $2$; there are at least $2n+1$ factors of dimension $2$ since $M$ must have positive pressure, by Lemma \ref{lem:pressure}.
\end{proof}

\section{Modules for \texorpdfstring{$\SL_2(3^a)$}{SL(2,3a)}}

In this section we describe the simple modules for $H=\SL_2(3^a)$ for $1\leq a\leq 7$, describe various extensions between some of the simple modules, and prove the existence or non-existence of various indecomposable modules.

Let $L=\SL_2(3)\leq H$. The simple modules for $L$ have dimension $1$, $2$ and $3$, with only the $2$-dimensional being faithful. Therefore, the non-trivial simple modules for $H$ are tensor products of modules of dimension $2$ and $3$, with a module of dimension $2^m3^n$ being faithful if and only if $m$ is odd.

Writing $2_i$ for the image of $2$ under $i$ iterations of the Frobenius map, and similarly for $3_i$, the simple modules for $H$ can be labelled by $2^m3^n_{r_1,\dots,r_{m+n}}$, where $m,n\geq 0$ are integers, $\{r_1,\dots,r_{m+n}\}\subset \{1,\dots,a\}$ with the $r_i$ distinct, with
\[2^m3^n_{r_1,\dots,r_{m+n}}=\left(\bigotimes_{i=1}^m 2_{r_i}\right) \otimes \left(\bigotimes_{j=m+1}^{m+n} 3_{r_i}\right).\]
Hence for example $12_{2,3,1}=2_2\otimes 2_3\otimes 3_1$ is a simple module for $\PSL_2(3^a)$ for any $a\geq 3$.

We need to understand the restrictions of these simple modules to $L$, in order to understand which ones we can have in the restrictions of minimal modules for $\bG=F_4,E_6,E_7$.

\begin{lemma}\label{lem:sl23restriction}
Let $H=\PSL_2(3^a)$, $a\geq 1$, and let $M$ be a simple module of dimension at most $56$. The restriction of $M$ to $\PSL_2(3)$ is as below.

\begin{center}\begin{tabular}{ccc}
\hline Module & Restriction & Composition factors of restriction
\\\hline $1$ & $1$ & $1$
\\ $3$ & $3$ & $3$
\\ $4=2\otimes 2$ & $3\oplus 1$ & $3,1$
\\ $9=3\otimes 3$ & $3^{\oplus 2}\oplus P(1)$ & $3^2,1^3$
\\ $12=2\otimes 2\otimes 3$ & $3^{\oplus 3}\oplus P(1)$ & $3^3,1^3$
\\ $16=2\otimes 2\otimes 2\otimes 2$ & $3^{\oplus 4}\oplus P(1)\oplus 1$ & $3^4,1^4$
\\ $27=3\otimes 3\otimes 3$ & $3^{\oplus 7}\oplus P(1)^2$ & $3^7,1^6$
\\ $48=2\otimes 2\otimes 2\otimes 2\otimes 3$ & $3^{\oplus 12}\oplus P(1)^{\oplus 4}$ & $3^{12},1^{12}$
\\ \hline
\end{tabular}\end{center}
\end{lemma}

We now move on to extensions. With the labelling above, we have the following easy lemma, which can be found for example in \cite[Corollary 3.9]{andersen1983}.

\begin{lemma}\label{lem:char3sl2cohomology} For any $a>1$, a simple module $M$ has non-trivial $1$-cohomology if and only if $M=4_{i,i+1}$ for some $1\leq i\leq a$, and
\[ \dim(\Ext^1(1,4_{i,i+1}))=\begin{cases} 1& a\geq 3,
\\ 2&a=2.\end{cases}\]
\end{lemma}

We will need more detailed information about extensions between simple modules of low dimension for $H$, and we summarize that which we need now, taken from \cite[Corollary 4.5]{andersen1983}. We restrict to the case when $a\neq 2$, because in this case things are slightly different, with that pesky $2$-dimensional $1$-cohomology group, and second because we describe the full projectives for this group after the lemma anyway.

\begin{lemma}\label{lem:lowdimexts3} Let $H=\PSL_2(3^a)$ for $3\leq a\leq 7$. The following extension groups have dimension $1$, for all $1\leq i,j\leq a$:
\[\begin{gathered} (4_{i,i+1},1),\qquad (1,4_{i,i+1}),\qquad (3_i,4_{i-1,i}),
\\
(4_{i-1,i},3_i),\qquad (4_{i,j},4_{i\pm 1,j}), \qquad (4_{i,j},4_{i,j\pm 1}),
\end{gathered}\]
If $A$ and $B$ are simple modules for $H$ of dimension at most $9$ then $\Ext^1(A,B)=0$ unless $(A,B)$ is on the list above.
\end{lemma}

We now consider certain modules. For $a=2$, the structures of the projective indecomposable modules are as follows:

\[
\begin{array}{c}
1
\\ 4\;\;4
\\ 1\;1\;1\;3_1\;\;3_2
\\ 4\;\;4
\\ 1
\end{array}\qquad
\begin{array}{c}
3_i
\\ 4
\\ 1\;3_{3-i}
\\ 4
\\ 3_i
\end{array}\qquad
\begin{array}{c}
4
\\ 1\;\;1\;\;3_1\;\;3_2
\\ 4\;\;4\;\;4
\\ 1\;\;1\;\;3_1\;\;3_2
\\ 4
\end{array}\]
We see that if a module $M$ has five socle layers then it has a projective summand. More generally, if $M$ has trivial composition factors, then we can use these to prove that $M$ must have more $4$s than pressure arguments suggest.

\begin{lemma}\label{lem:sl29hidetrivials}
Let $H=\PSL_2(9)$. If $M$ is a $kH$-module with no trivial submodules or quotients, and there are $2n-1$ or $2n$ trivial composition factors in $M$, then the number of $4$-dimensional factors in $M$ is at least $2n$.

Furthermore, the only submodules of $P(4)$ consisting of $4$s and $1$s are submodules of a self-dual module $4/1,1/4$. In particular, there is no uniserial module of the form $4/1/4$.
\end{lemma}
\begin{proof}
Let $M$ be a $kH$-module, which we may assume is indecomposable. If $M$ is the $9$-dimensional projective simple then the claim is true. If $M$ has any $3$-dimensional submodules or quotients then we may remove them without affecting the claim, and so we may assume that $M$ is a submodule of copies of $P(4)$.

If $M$ is projective then the result holds, so $M$ is not projective, in which case it has at most four socle layers. Since the fourth socle layer consists solely of copies of $1$ and $3_i$, $M$ must actually have three socle layers. In particular, the trivials are all in the second socle layer, so if there are $2n-1$ or $2n$ of them, there must be at least $n$ copies of the $4$-dimensional module in the socle, and similarly in the top. This completes the proof of the first claim.

The second is easy to see by a computer proof that $4/1,1/4$ is the largest such module. Since it is self-dual, we cannot construct a $4/1/4$ inside it, yielding the second statement.
\end{proof}

\begin{lemma}\label{lem:no414} Let $H=\PSL_2(3^a)$ for some $2\leq a\leq 7$. There does not exist a uniserial module with structure $4_{i,j}/1/4_{m,n}$, where $4_{i,j}$ and $4_{m,n}$ are simple modules of dimension $4$.

As a consequence, if $a\neq 2$ and $M$ is a module with $2i-\alpha$ composition factors of dimension $4$ and $i$ of dimension $1$, for some $i>0$ and $\alpha\geq 0$, then $M$ has a trivial submodule or quotient.
\end{lemma}
\begin{proof}For the first part, we may assume that $a\geq 3$ by Lemma \ref{lem:sl29hidetrivials}. By Lemma \ref{lem:char3sl2cohomology} the only modules with non-trivial $1$-cohomology are $4_{i,i+1}=2_i\otimes 2_{i+1}$, and the module $1/4_{i,i+1}$ is unique for $3\leq a\leq 7$, so by applying a field automorphism we may assume that the socle of our uniserial module is $4_{1,2}$. To prove that there is no module $4_{j,j+1}/1/4_{i,i+1}$ we simply use a computer to compare $\Ext^1(4_{j,j+1},1/4_{1,2})$ and $\Ext^1(4_{j,j+1},4_{1,2})$ for each $j$, and note that they coincide for all $4$-dimensional modules $4_{j,j+1}$. This proves the result because it shows that every module that is an extension of $4_{j,j+1}$ by $1/4_{1,2}$ arises as a module $1,4_{j,j+1}/4_{1,2}$, and hence is not a uniserial module. (This can be checked directly by simply constructing the maximal extension.)

For the second statement, we use induction. If $i=1$ then there are at most two $4$-dimensional factors, and so $M$ has a trivial submodule or quotient as there is no uniserial module $4_{i,j}/1/4_{m,n}$. Let $M$ be a minimal counterexample.

Since $M$ has no trivial submodules, $\soc(M)$ is a sum of $4$-dimensional modules, say $n$ of them. As $\Ext^1(1,4_{i,j})$ has dimension at most $1$, there are at most $n$ trivial submodules of the quotient $M/\soc(M)$. Let $N$ be the preimage in $M$ of the $\{1\}$-radical of the quotient $M/\soc(M)$. Notice that the quotient module $M/N$ has $2i-\alpha-n$ composition factors of dimension $4$ and at least $i-n$ of dimension $1$. Since $M$ has no trivial quotient as it is a quotient of $M$, and has no trivial submodules by construction, $M$ by induction cannot satisfy the hypotheses of the lemma. This is only possible if $M/N$ has no trivial composition factors, so $N$ contains all trivial factors of $M$.

As all trivial composition factors must lie in the second socle layer of $M$ from the above, by Lemma \ref{lem:pressure}, there must be at least $i$ modules of dimension $4$ in $\soc(M)$ and at least $i$ modules of dimension $4$ in $\topp(M)$, but there are at most $2i$ composition factors of dimension $4$ in total. Thus $n=i$, $\alpha=0$, and there are exactly $i$ composition factors in each of the three socle layers. By the first part, we know that $i\geq 2$.

We claim that no such module $M$ can exist, because there is no module $4_{j,j+1}/1/4_{l,l+1}$. To see this, proceed by induction on $\dim(M)$, or equivalently on $i$. The quotient by any $4$-dimensional submodule must have a fixed point by Lemma \ref{lem:pressure}, because the image of $\soc^2(M)$ in the quotient has pressure $-1$. Thus we may quotient out by a submodule $N$ of the form $1/4_{j,j+1}$. Consider the dual $(M/N)^*$ of the quotient module $M/N$. This has $i$ composition factors of dimension $4$ in the socle but only $i-1$ trivial factors in the second socle layer. Thus there exists a particular $4$-dimensional submodule that one may quotient out by and still have no trivial composition factors. This quotient has $2(i-1)$ factors of dimension $4$ and $i-1$ trivial factors, but still has no trivial submodule or quotient. This is a contradiction to our induction hypothesis, so $M$ cannot exist, completing the proof.
\end{proof}

We end with a small lemma, needed at one point in the text.

\begin{lemma}\label{lem:psl227pressure1} Let $p^a=27$. The projective cover of $4_{1,2}$ is
\[ \begin{array}{c}4_{1,2}
\\1\;\;3_2\;\;4_{2,3}\;\;4_{1,3}\;\;12_{2,3,1}
\\1\;\;3_1\;\;3_3\;\;4_{2,3}\;\;4_{1,3}\;\;4_{1,2}\;\;4_{1,2}\;\;4_{1,2}\;\;9_{\
1,3}\;\;12_{1,3,2}
\\1\;\;1\;\;3_2\;\;3_2\;\;4_{2,3}\;\;4_{2,3}\;\;4_{1,3}\;\;4_{1,3}\;\;4_{1,2}\;\;12_{1,2,3}\;\;12_{2,3,1}\;\;12_{2,3,1}
\\1\;\;3_1\;\;3_3\;\;4_{2,3}\;\;4_{1,3}\;\;4_{1,2}\;\;4_{1,2}\;\;4_{1,2}\;\;9_{1,3}\;\;12_{1,3,2}
\\1\;\;3_2\;\;4_{2,3}\;\;4_{1,3}\;\;12_{2,3,1}
\\4_{1,2}
\end{array}\]
Consequently, if $M$ is a self-dual module of pressure $1$ with at least five trivial composition factors then $H$ stabilizes a line or hyperplane of $M$.
\end{lemma}
\begin{proof}
The description of the projective is produced by a Magma calculation.

To see the consequence, let $M$ be a minimal counterexample to the statement. By removing all submodules and quotients not of dimension $4$ from $M$, we may assume that $M$ is a submodule of $P(4_{1,2})$ (up to field automorphism). Since $M$ has pressure $1$ it cannot be the whole of $P(4_{1,2})$, so in particular $M$ is a submodule of $\rad(P(4_{1,2}))$. As the top of $M$ must be $4_{1,2}$, this means that $M$ cannot have $1$, $3_2$, $4_{2,3}$, $4_{1,3}$ or $12_{2,3,1}$ as a quotient either, so $M$ is a submodule of $\soc^5(P(4_{1,2}))$. This has five trivial factors, but one is a quotient, so also needs to be removed, and $M$ has only four trivial factors, which is a contradiction.
\end{proof}

\section{Modules for \texorpdfstring{$\SL_2(p)$}{SL(2,p)}}
\label{sec:sl2p}
Let $p\geq 5$. Since $H=\SL_2(p)$ has a cyclic Sylow $p$-subgroup, there are only finitely many indecomposable modules for it over a field of characteristic $p$. In this section we describe how to construct all indecomposable modules for $H$ in characteristic $p$, using the projective indecomposable modules as a starting point.

There are three blocks for $kH$: one consists of all faithful modules, one of all non-faithful modules other than the Steinberg, and one is simply the Steinberg module (see, for example, \cite[p.47]{andersen1983}). We understand the block containing the Steinberg module, and so we will concentrate on the other blocks.

The Green correspondence \cite[Theorem 11.1]{alperin} shows that the number of non-projective indecomposable modules of dimension congruent to $i$ modulo $p$ for $H$ is the same as that of the normalizer $N_H(P)$ of a Sylow $p$-subgroup $P$ of $H$, a soluble group of order $p(p-1)$ with a centre of order $2$. However, for this group, it is easy to construct the indecomposable modules: the projective modules are all of dimension $p$, and look like truncated polynomial rings $k[X]/(X^p-1)$, hence are uniserial. (There is a fixed $1$-dimensional module $N$ such that the $i$th socle layer of the projective cover of the trivial is $N^{\otimes i-1}$ for $1\leq i\leq p$, where we take $N^{\otimes 0}$ to be the trivial module. All other projectives are found by tensoring this projective by the various $1$-dimensional modules.) Every indecomposable module is a quotient of such a module, and as every simple module for $N_H(P)$ is $1$-dimensional, we see that there are exactly $p-1$ indecomposable modules of dimension $i$ for each $1\leq i\leq p$, with half of these faithful modules for $\SL_2(p)$ and half modules for $\PSL_2(p)$.

In particular, we see that once we have constructed $p(p-1)$ indecomposable modules for $H$ other than the Steinberg module then we must have found them all. Thus we start with the simple and projective modules for $H$, which may be found in \cite[pp. 75--79]{alperin}. Letting $M=L(1)$ be the natural module for $H$, we construct all simple modules using symmetric powers
\[ L(i)=S^i(M)\qquad 0\leq i\leq p-1,\]
with $L(i)$ being of dimension $i+1$. As with the case of $\SL_2(2^a)$, we will normally write the single number $i$ to refer to the simple module of dimension $i$, and so a module $3/5$ for $\SL_2(7)$, for example, is an $8$-dimensional module with $5$-dimensional socle $L(4)$ and $3$-dimensional top $L(2)$. The odd-dimensional simple modules are modules for $\PSL_2(p)$, and the even-dimensional ones are faithful modules for $\SL_2(p)$.

Having defined the simple modules, we consider the projectives: the Steinberg module $L(p-1)$ of dimension $p$ is already projective and is being ignored, and for each simple module $i$ with $1\leq i\leq p-1$, the projective module $P(i)$ has structure
\[ i/((p+1-i)\,,\, (p-1-i))/i,\]
except when $i=1$, in which case $p+1-i$ would have dimension $p$, and we have $1/(p-2)/1$, and when $i=p-1$, so $p-1-i$ would have dimension $0$, and we have $(p-1)/2/(p-1)$.

We represent these in diagrams, with lines linking two composition factors $A$ and $B$ if there is a non-split extension $A/B$ as a subquotient of the module. For example, here are $P(3)$ and $P(5)$ for $\PSL_2(11)$.
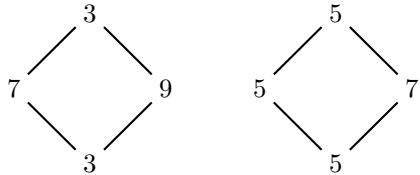
\tikzstyle{every node}=[circle, fill=black!0,
                        inner sep=0pt, minimum width=15pt]

\begin{center}\begin{tikzpicture}[thick,scale=1]
\draw (0,0) -- (1,1);
\draw (0,0) -- (-1,1);
\draw (0,2) -- (1,1);
\draw (0,2) -- (-1,1);
\draw (0,2) node{$3$};
\draw (-1,1) node{$7$};
\draw (1,1) node{$9$};
\draw (0,0) node{$3$};
\end{tikzpicture}\qquad
\begin{tikzpicture}[thick,scale=1]
\draw (0,0) -- (1,1);
\draw (0,0) -- (-1,1);
\draw (0,2) -- (1,1);
\draw (0,2) -- (-1,1);
\draw (0,2) node{$5$};
\draw (-1,1) node{$5$};
\draw (1,1) node{$7$};
\draw (0,0) node{$5$};
\end{tikzpicture}\end{center}

Using these we construct indecomposable modules as follows: we have modules of the form $i/(p+1-i)$ and $i/(p-1-i)$, and also two modules of the form $i/(p-1-i),(p+1-i)$ and $(i+2)/(p-1-i)$. These two indecomposables can be summed together, then quotiented by a diagonal submodule $p-1-i$ to make a new module with four composition factors. We can do that same with the modules $i/(p-1-i),(p+1-i)$ and $(i+2)/(p-3-i),(p-1-i)$ to obtain a module with five composition factors.

It is easier to visualize using diagrams. In the example above, we can remove the socles of the two projectives to obtain modules $3/7,9$ and $5/5,7$, take their direct sum, and then quotient out by a diagonal $7$.
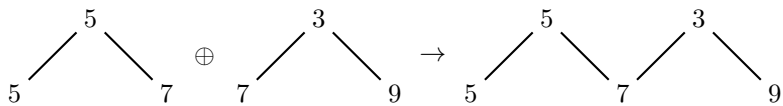
\begin{center}\begin{tikzpicture}[thick,scale=1]
\draw (0,2) -- (1,1);
\draw (0,2) -- (-1,1);
\draw (0,2) node{$5$};
\draw (-1,1) node{$5$};
\draw (1,1) node{$7$};
\draw (1.5,1.5) node{$\oplus$};
\draw (3,2) -- (4,1);
\draw (3,2) -- (2,1);
\draw (3,2) node{$3$};
\draw (2,1) node{$7$};
\draw (4,1) node{$9$};
\draw (4.5,1.5) node{$\rightarrow$};
\draw (6,2) -- (7,1);
\draw (6,2) -- (5,1);
\draw (8,2) -- (7,1);
\draw (8,2) -- (9,1);
\draw (6,2) node{$5$};
\draw (5,1) node{$5$};
\draw (7,1) node{$7$};
\draw (8,2) node{$3$};
\draw (9,1) node{$9$};
\end{tikzpicture}\end{center}
This process certainly produces a module, with quotients both of our original summands, and so this module must be indecomposable. Note that if one tries to do this with say two copies of $3/7,9$ then the fact that $\Ext^1(3,7)$ is $1$-dimensional means that this module splits, so one needs the modules at the top (in this case $3$ and $5$) to be different.

One can continue this process until one constructs an indecomposable module $M$ with all (non-projective) simple modules appearing in the top and the socle of $M$ exactly once. As an example, the diagrams of the two such modules for $p=11$ (one for $\PSL_2(11)$, one for the faithful modules of $\SL_2(11)$) are as follows:
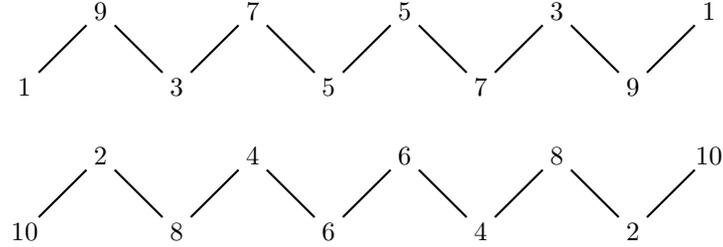
\begin{center}\begin{tikzpicture}[thick,scale=1]
\draw (0,2) -- (1,1);
\draw (0,2) -- (-1,1);
\draw (2,2) -- (1,1);
\draw (2,2) -- (3,1);
\draw (4,2) -- (3,1);
\draw (4,2) -- (5,1);
\draw (6,2) -- (5,1);
\draw (6,2) -- (7,1);
\draw (8,2) -- (7,1);
\draw (0,2) node{$9$};
\draw (2,2) node{$7$};
\draw (4,2) node{$5$};
\draw (6,2) node{$3$};
\draw (8,2) node{$1$};
\draw (-1,1) node{$1$};
\draw (1,1) node{$3$};
\draw (3,1) node{$5$};
\draw (5,1) node{$7$};
\draw (7,1) node{$9$};
\end{tikzpicture}\end{center}

\begin{center}\begin{tikzpicture}[thick,scale=1]
\draw (0,2) -- (1,1);
\draw (0,2) -- (-1,1);
\draw (2,2) -- (1,1);
\draw (2,2) -- (3,1);
\draw (4,2) -- (3,1);
\draw (4,2) -- (5,1);
\draw (6,2) -- (5,1);
\draw (6,2) -- (7,1);
\draw (8,2) -- (7,1);
\draw (0,2) node{$2$};
\draw (2,2) node{$4$};
\draw (4,2) node{$6$};
\draw (6,2) node{$8$};
\draw (8,2) node{$10$};
\draw (-1,1) node{$10$};
\draw (1,1) node{$8$};
\draw (3,1) node{$6$};
\draw (5,1) node{$4$};
\draw (7,1) node{$2$};
\end{tikzpicture}\end{center}
We can take subquotients of these modules and construct new indecomposable modules, and we claim that this constructs all non-projective, indecomposable modules for $\SL_2(p)$, other than the Steinberg module.

First, the non-simple indecomposable subquotients of the module $M$ are in one-to-one correspondence with connected subdiagrams of the diagram with at least one edge, since one notes that no two distinct subdiagrams of the diagram above have the same first and second rows. In the case of simple modules, of course each appears twice as a subdiagram.

The number of connected subdiagrams of each diagram with at least one edge is $(p-1)(p-2)/2$ (i.e., we choose the start and end points), and add in the $p-1$ simple modules (other than the Steinberg), and the $p-1$ non-simple projective modules, to obtain
\[ (p-1)(p-2)+2(p-1)=p(p-1).\]
This is the number of indecomposable modules for the normalizer, and so we must have constructed all indecomposable modules for $\SL_2(p)$.

It is clear from this `zigzag' structure, that for any indecomposable module, if $A$ and $B$ lie in the socle so does any module with dimension between $\dim(A)$ and $\dim(B)$. We have proved the following proposition.

\begin{proposition}\label{prop:simplesl2p}
Let $H=\SL_2(p)$, and let $M$ be an indecomposable module for $H$.
\begin{enumerate}
\item If $M$ has one socle layer then $M$ is simple, and there are $p$ such modules, one of each dimension.
\item If $M$ has three socle layers then $M=P(i)$ for some $1\leq i\leq p-1$.
\item If $M$ has two socle layers then the socle of $M$ consists of modules of dimension $i,i+2,\dots,j$ ($i\leq j$), and the top consists of modules $p-j+\epsilon,p-j+\epsilon+2,\dots,p-i+\delta$, where $\epsilon,\delta=\pm1$. There are $(p-1)(p-2)$ such modules.
\end{enumerate}
\end{proposition}

The indecomposable modules for $\PSL_2(7)$ other than the Steinberg are below, ordered so that the modules in column $i$ have dimension congruent to $i$ modulo $7$.
\begin{center}
\begin{tabular}{ccccccc}
$1$ & $3,5/3,5$ & $3$ & $1,3,5/1,3,5$& $5$ & $3/3$ & $P(1)$
\\ $3/5$&$1,3/5$&$1,3,5/3,5$&$3/3,5$&$1,3/3,5$&$1/5$&$P(3)$
\\ $5/3$&$5/1,3$&$3,5/1,3,5$&$3,5/3$&$3,5/1,3$&$5/1$&$P(5)$
\end{tabular}\end{center}
As another example, the indecomposable modules for $\PSL_2(5)$ are the Steinberg module, together with the following:
\[ 1,\;\;3,\;\;1/3,\;\;3/1,\;\;P(1),\;\;P(3),\;\;3,1/3,\;\;3/1,3,\;\;3/3,\;\;1,3/1,3.\]
Each of these modules $M$ is in Green correspondence with an indecomposable module $V$ of dimension at most $p-1$. Also, Green correspondence means that the restriction of $M$ to $N_H(P)$ is a sum of $V$ and projective modules. These two facts yield the first part of the following result.

\begin{lemma}\label{lem:indJordanblocks} Let $H=\SL_2(p)$ and let $M$ be an indecomposable $kH$-module. Let $u\in H$ have order $p$.
\begin{enumerate}
\item\label{lemi:jorda} $u$ acts on $M$ with at most one block of size different from $p$.
\item\label{lemi:jordb} All blocks have size $p$ if and only if $M$ is projective.
\item\label{lemi:jordc} $u$ acts on $M$ with no blocks of size $p$ if and only if $u$ is simple, but not the Steinberg module, or $u$ has dimension $p-1$.
\end{enumerate}
\end{lemma}
\begin{proof} The first part follows from what was said above. The second part simply states that a module is projective if and only if its restriction to a Sylow $p$-subgroup is projective, which is \cite[Theorem 5.6 and Corollary 9.3]{alperin}. For (\ref{lemi:jordc}), the dimension of any such module must be less than $p$ by (\ref{lemi:jorda}). From the description of the indecomposable modules, we see that if $M$ is not simple then $M$ is the extension of a module of dimension $i$ by one of dimension $p-1-i$, whence the result holds.
\end{proof}

The next lemma is an easy consequence of Lemma \ref{lem:indJordanblocks}.

\begin{lemma}\label{lem:pblockseven}
Let $H=\PSL_2(p)$, and let $M$ be a module for $H$ over a field of characteristic $p$. Let $u\in H$ have order $p$. If $n_i$ denotes the number of Jordan blocks of size $i$ in the action of $u$ on $M$, then
\[ n_p\geq \sum_{1\leq i<(p-1)/2} n_{2i}.\]
\end{lemma}
\begin{proof} If the inequality holds for $M_1$ and $M_2$ then it holds for $M_1\oplus M_2$, so assume that $M$ is indecomposable.

If $M$ has dimension at most $p$ then $M$ is simple or $1/(p-2)/1$, so $u$ acts on $M$ with a single block of odd size, or is $i/(p-1-i)$, so $u$ acts with a single block of size $p-1$. Therefore if there is a block of even size $i<p-1$ then it must come from an indecomposable module of dimension greater than $p$, and so we obtain at least one block of size $p$, as needed.
\end{proof}

For $\bG=E_7$ we must also consider $H=\SL_2(p)$ with $Z(H)=Z(\bG)$. In this case we want a similar result to the above but for faithful modules.

\begin{lemma}\label{lem:faithfulpblockseven}
Let $H=\SL_2(p)$, and let $M$ be a module for $H$ over a field of characteristic $p$ on which the central involution $z$ of $H$ acts as the scalar $-1$. Let $u\in H$ have order $p$. If $n_i$ denotes the number of Jordan blocks of size $i$ in the action of $u$ on $M$, then
\[ n_p\geq \sum_{1\leq i<(p-1)/2} n_{2i-1}.\]
\end{lemma}
\begin{proof}Similar to Lemma \ref{lem:pblockseven}, and omitted.\end{proof}

We often want to understand self-dual modules for $H$, since the minimal module $M(\bG)$ is self-dual for $F_4$ and $E_7$, and the simple adjoint module $L(\bG)^\circ$ is always self-dual. Using the statements above, if $n_i$ is odd, where again $n_i$ is the number of blocks of size $i$ in the action of $u$, there must be a self-dual indecomposable summand of dimension congruent to $i$ modulo $p$.

The next lemma follows from Proposition \ref{prop:simplesl2p} and classifies self-dual indecomposable modules for $\SL_2(p)$. From our zigzag diagrams above, it is clear which the self-dual modules are: choose the same simple module as the start and end points of the subdiagram.

\begin{lemma}\label{lem:selfdualsl2p}
Let $H=\SL_2(p)$, and let $M$ be a self-dual indecomposable module for $H$. If $M$ is not simple or projective, then $M$ has socle (and top) consisting of pairwise non-isomorphic modules $N_1,N_2,\dots,N_r$, where $\dim(N_i)-\dim(N_{i-1})=2$ and $\dim(N_1)+\dim(N_r)=p\pm 1$. In particular, there are exactly $p-1$ non-projective, indecomposable self-dual modules for $\PSL_2(p)$, and exactly $p-1$ non-projective, indecomposable and faithful self-dual modules for $\SL_2(p)$.

Therefore, if $p\equiv 1\bmod 4$, then there is a unique self-dual indecomposable module for $\PSL_2(p)$ of dimension congruent to $2i+1$ for each $0\leq i\leq (p-1)/2$ modulo $p$, and none congruent to $2i$, and there is a unique faithful, self-dual indecomposable module for $\SL_2(p)$ of dimension congruent to $2i$ for each $0\leq i\leq (p-1)/2$ modulo $p$, and none congruent to $2i+1$.

On the other hand, if $p\equiv 3\bmod 4$, then there is a unique self-dual indecomposable module for $\PSL_2(p)$ of dimension congruent to $2i$ for each $0\leq i\leq (p-1)/2$ modulo $p$, and none congruent to $2i+1$, and there is a unique faithful, self-dual indecomposable module for $\SL_2(p)$ of dimension congruent to $2i+1$ for each $0\leq i\leq (p-1)/2$ modulo $p$, and none congruent to $2i$.
\end{lemma}

We can use this to obtain a better handle on which possible Jordan block structures a given unipotent element $u$ can have, given that it lies inside a copy of $\PSL_2(p)$ for $p\equiv 1\bmod 4$. We split the result into two corollaries depending on whether one has modules for $\PSL_2(p)$ or $\SL_2(p)$.

\begin{corollary}\label{cor:blocks1mod4}
Let $H=\PSL_2(p)$ with $p\equiv 1\bmod 4$, and let $M$ be a self-dual module for $H$. Let $u$ be an element of order $p$ in $H$. The action of $u$ has an even number of blocks of a given even size $i$, and there are at least as many blocks of size $p$ as there are blocks of size all even numbers less than $p-1$.
\end{corollary}
\begin{proof} The second statement comes from Lemma \ref{lem:faithfulpblockseven}. That all blocks of even size come in pairs follows from the fact that $u$ acts on any self-dual indecomposable module with only odd blocks by Lemma \ref{lem:selfdualsl2p}. Thus if $u$ has a block of even size it must come from a summand $M_1$ that is not self-dual, and then $M_1^*$ is another summand contributing another block of the same size.
\end{proof}

\begin{corollary}\label{cor:faithfulblocks1mod4}
Let $H=\SL_2(p)$ with $p\equiv 1\bmod 4$, and let $M$ be a self-dual module for $H$ on which the central involution $z$ acts as the scalar $-1$. Let $u$ be an element of order $p$ in $H$. The action of $u$ has an even number of blocks of a given odd size $i$, and there are at least as many blocks of size $p$ as there are blocks of size all odd numbers less than $p$.
\end{corollary}

We now turn to tensor products. By Steinberg's tensor product theorem, simple modules for $\SL_2(p^a)$ are tensor products of Frobenius twists of \emph{$p$-restricted} modules, i.e., $L(i)$ for $i\leq p-1$. These restrict to $\SL_2(p)$ as tensor products of simple modules, so it will come in handy to understand the tensor products of simple modules for $\SL_2(p)$.

The next result gives the tensor product of any two simple modules for $\SL_2(p)$, and will be of great use when computing the restriction of simple $\SL_2(p^a)$-modules to $\SL_2(p)$.

\begin{proposition}\label{prop:tensorproductsl2p}
Let $H=\SL_2(p)$. If $0\leq\mu\leq\lambda\leq p-1$ then $L(\lambda)\otimes L(\mu)$ is given by one of the following:
\begin{enumerate}
\item If $\lambda+\mu<p$ then
\[ L(\lambda)\otimes L(\mu)=L(\lambda-\mu)\oplus L(\lambda-\mu+2)\oplus \cdots \oplus L(\lambda+\mu-2)\oplus L(\lambda+\mu).\]
\item If $\lambda+\mu\geq p$ and $\lambda<p-1$ then
\begin{align*}
 L(\lambda)\otimes L(\mu)=&L(\lambda-\mu)\oplus L(\lambda-\mu+2)\oplus \cdots \oplus L(a)
 \\ &\oplus \begin{cases} P(\lambda+\mu)\oplus P(\lambda+\mu-2)\oplus \cdots \oplus P(p+1)\oplus L(p-1)&\mu\text{ even}
\\ P(\lambda+\mu)\oplus P(\lambda+\mu-2)\oplus \cdots \oplus P(p)&\mu\text{ odd}\end{cases}
\end{align*}
where $a=2p-(\lambda+\mu+4)$.
\item $L(p-1)\otimes L(p-1)=P(1)\oplus P(3)\oplus \cdots \oplus P(p-1)$.
\end{enumerate}
\end{proposition}

This result can be found, for example, in \cite{dotyhenke2005} and explicitly in \cite[Lemma 3.1]{craven2013}.

\section{Modules for \texorpdfstring{$\SL_2(p^a)$}{SL(2,pa)} for \texorpdfstring{$p\geq 5$}{p at least 5} and \texorpdfstring{$a>1$}{a>1}}

As with modules for $\SL_2(3^a)$ we need a notation system for the simple modules, and as in that section, we let $2_1$ denote the natural module $L(1)$, $i_1=S^{i-1}(2_1)=L(i)$ for $2\leq i<p$ be the symmetric powers (the $p$-restricted modules) and let $i_{j+1}$ denote the application of a single Frobenius morphism to $i_j$, i.e., $i_j^{[1]}$. We then write, for a module of dimension $n$ formed as the tensor product of $m$ twisted $p$-restricted modules, $n_{a_1,\dots,a_m}$, in order of increasing dimension of factor; for example, the module $2_1\otimes 3_2$ will be denoted $6_{1,2}$, and $3_1\otimes 3_2\otimes 2_3$ will be denoted $18_{3,1,2}$. If the integer $n$ has a unique decomposition as a product of exactly $m$ integers greater than $1$ such that the module would be for the correct group (i.e., $\PSL_2(p^a)$ or $\SL_2(p^a)$) then we simply write that, so that $6_1$ and $6_{1,2}$ for $\SL_2(49)$ are unambiguous. Sometimes there are modules that could be either for $\SL_2$ or $\PSL_2$, such as 
$12_{1,2}$ for $p\geq 7$, which is either $2_1\otimes 6_2$ or $3_1\otimes 4_2$, but context will tell us which. When there genuinely is ambiguity, for example, $18_{1,2}$ when $p\geq 11$, as it could be $2_1\otimes 9_2$ or $3_1\otimes 6_2$, we label them with subscripts $18_{1,2}^{(1)}$ and $18_{1,2}^{(2)}$ according to the lexicographic ordering on the partitions of $18$, but in these rare cases we remind the reader which is which.

\medskip

We start with some information about extensions of simple $kH$-modules. These are completely determined in \cite[Corollary 4.5]{andersen1983}, where a general formula in terms of the $p$-adic expansion of the highest weights is given. We extract a few special cases which are of use to us.

Of particular interest is which modules have non-trivial $1$-cohomology, since we will often want to prove that we stabilize a line. The next lemma gives this completely.

\begin{lemma}\label{lem:cohomologysimple}
Let $p$ be a prime, $a\geq 1$ be an integer, and let $M$ be a simple module for $H=\SL_2(p^a)$ with non-trivial $1$-cohomology. One of the following holds.
\begin{enumerate}
\item $p^a=2$, $M$ is the trivial module, with $\dim(H^1(H,M))=1$.
\item $p$ is odd and $a=1$, $\dim(M)=p-2$, with $\dim(H^1(H,M))=1$.
\item $p^a=9$, $\dim(M)=4$, with $\dim(H^1(H,M))=2$.
\item $p^a\neq 9$ with $a\geq 2$, $M$ is up to application of a Frobenius map $L(p-2)\otimes L(1)^{[1]}$, where $(-)^{[1]}$ is the Frobenius twist (so that $\dim(M)=2(p-1)$ and $M=2(p-1)_{2,1}$), with $\dim(H^1(H,M))=1$.
\end{enumerate}
\end{lemma}

Just knowing that modules have $1$-cohomology is not going to be enough information. We need more specific information about extensions between simple modules of low dimension, for $p=5,7$ and $a>1$. The next two lemmas of this section furnish us with this information.

\begin{lemma}\label{lem:psl2125smallmodules} Let $H=\PSL_2(5^a)$ for $a=2,3$. The extensions between simple modules of dimension at most $8$ are:
\begin{enumerate}
\item $1$ with $8_{i,i-1}$;
\item $3_i$ with $4_{i,i+1}$, $8_{i,i-1}$;
\item $4_{i,i+1}$ with $3_i$, $8_{i+1,i-1}$ (the latter only for $a=3$);
\item $5_i$ with nothing;
\item for $a=2$, $8_{i,i+1}$ with $3_i$;
\item for $a=3$, $8_{i,i+1}$ with $4_{i-1,i}$;
\item for $a=3$, $8_{i,i-1}$ with $1$, $3_i$.
\end{enumerate}
\end{lemma}

\begin{lemma}\label{lem:psl249smallmodules} Let $H=\PSL_2(49)$. The extensions between simple modules of dimension at most $9$ are:
\begin{enumerate}
\item $1$ with nothing;
\item $3_i$ with $8_{i+1,i}$;
\item $4_{1,2}$ with $5_1$, $5_2$;
\item $5_i$ with $4_{1,2}$;
\item $7_i$ with nothing;
\item $8_{i,i+1}$ with $3_{i+1}$, $9_{1,2}$;
\item $9_{1,2}$ with $8_{1,2}$, $8_{2,1}$.
\end{enumerate}
\end{lemma}

We also require some restrictions of simple $\PSL_2(p^a)$-modules to $\PSL_2(p)$. This is needed because we often understand the action of $\PSL_2(p)$ on the minimal or adjoint modules completely, and want to extend a module for $\PSL_2(p)$ to a module for $\PSL_2(p^a)$. We consider modules of dimension at most $56$ to include the minimal modules of $F_4$, $E_6$ and $E_7$. We use Proposition \ref{prop:tensorproductsl2p} to compute the restrictions of modules for $\SL_2(p^a)$ to $\SL_2(p)$. We list restrictions for $p^a\leq 150$ and dimension up to $56$, as this is all we will need.

\begin{lemma}\label{lem:sl25restriction}
Let $H=\PSL_2(5^a)$ for $a=2,3$, and let $M$ be a simple module of dimension at most $56$. The restriction of $M$ to $L=\PSL_2(5)$ is as in Table \ref{t:tensorproductsp=5}.
\begin{table}\begin{center}\begin{tabular}{ccc}
\hline Module & Restriction & Factors of restriction
\\\hline $1$ & $1$ & $1$
\\ $3$ & $3$ & $3$
\\ $4=2\otimes 2$ & $3\oplus 1$ & $3,1$
\\ $5$ & $5$ & $5$
\\ $8=2\otimes 4$ & $5\oplus 3$ & $5,3$
\\ $9=3\otimes 3$ & $1\oplus 3\oplus 5$ & $5,3,1$
\\ $12=2\otimes 2\otimes 3$ & $5\oplus 3^{\oplus 2}\oplus 1$ & $5,3^2,1$
\\ $15=3\otimes 5$ & $5\oplus P(3)$ & $5,3^3,1$
\\ $16=4\otimes 4$ & $5\oplus P(3)\oplus 1$ & $5,3^3,1^2$
\\ $20=2\otimes 2\otimes 5$ & $5^{\oplus 2}\oplus P(3)$ & $5^2,3^3,1$
\\ $24=2\otimes 4\otimes 3$ & $5^{\oplus 2}\oplus P(3)\oplus 3\oplus 1$ & $5^2,3^4,1^2$
\\ $25=5\otimes 5$ & $5^{\oplus 2}\oplus P(3)\oplus P(1)$ & $5^2,3^4,1^3$
\\ $27=3\otimes 3\otimes 3$ & $5^{\oplus 2}\oplus P(3)\oplus 3^{\oplus 2}\oplus 1$ & $5^2,3^5,1^2$
\\ $40=2\otimes 4\otimes 5$ & $5^{\oplus 3}\oplus P(3)^{\oplus 2}\oplus P(1)$ & $5^3,3^7,1^4$
\\ $45=3\otimes 3\otimes 5$ & $5^{\oplus 4}\oplus P(3)^{\oplus 2}\oplus P(1)$ & $5^4,3^7,1^4$
\\ $48=4\otimes 4\otimes 3$ & $5^{\oplus 4}\oplus P(3)^{\oplus 2}\oplus P(1)\oplus 3$ & $5^4,3^8,1^4$
\\ \hline
\end{tabular}\end{center}\caption{Tensor products of modules for $\PSL_2(5)$}\label{t:tensorproductsp=5}
\end{table}
Consequently, if $V$ is a module for $H$ of dimension at most $56$ such that $V\downarrow_L$ has more trivial than $3$-dimensional composition factors, then $H$ stabilizes a line on $V$.
\end{lemma}
\begin{proof} We prove the last statement: from the table above we see that only the trivial has more $1$s than $3$s in its restriction to $L$. Suppose that the composition factors of $V\downarrow_L$ are $5^i,3^j,1^k$, with $k>j$. Lemma \ref{lem:cohomologysimple} states that the only simple modules with non-trivial $1$-cohomology for $H$ are of dimension $8$. Let $\alpha$ and $\beta$ be the number of trivial and $8$-dimensional composition factors of $V$ respectively. We have that $\alpha\geq k-(j-\beta)>\beta$, so $V$ has negative pressure, and hence $H$ stabilizes a line on $V$ by Lemma \ref{lem:pressure}, as needed.
\end{proof}

For $\PSL_2(25)$, we will need the eigenvalues of an element of order $12$ on the simple modules, so we list them here. These are of course easy to compute.

\begin{lemma}
Let $H=\PSL_2(25)$, and let $x$ be a semisimple element of order $12$ in $H$. Let $\xi$ denote a primitive $12$th root of unity. Choosing $\xi$ so that $x$ acts on the symmetric square of the natural module for $\SL_2(25)$ with eigenvalues $1,\xi^{\pm 1}$, the eigenvalues of $x$ on the various simple modules for $H$ are as follows.
\begin{center}
\begin{tabular}{cc}
\hline Dimension & Eigenvalues
\\\hline $1$&$1$
\\$3$&$1,(\xi,\xi^{11})/(\xi^5,\xi^7)$
\\$4$&$(\xi^2,\xi^{10}),(\xi^3,\xi^9)$
\\$5$&$1,(\xi^2,\xi^{10}),(\xi,\xi^{11})/(\xi^5,\xi^7)$
\\$8$&$(\xi^2,\xi^{10}),(\xi^3,\xi^9),(\xi^4,\xi^8),(\xi,\xi^{11})/(\xi^5,\xi^7)$
\\$9$&$1,(-1)^2,(\xi,\xi^{11}),(\xi^4,\xi^8),(\xi^5,\xi^7)$
\\$15$&$1,(-1)^2,(\xi,\xi^{11}),(\xi^2,\xi^{10}),(\xi^3,\xi^9),(\xi^4,\xi^8),(\xi^5,\xi^7),(\xi,\xi^{11})/(\xi^5,\xi^7)$
\\$16$&$(-1)^2,(\xi,\xi^{11}),(\xi^2,\xi^{10}),(\xi^3,\xi^9)^2,(\xi^4,\xi^8)^2,(\xi^5,\xi^7)$
\\$25$&$1^3,(-1)^2,(\xi,\xi^{11})^2,(\xi^2,\xi^{10})^2,(\xi^3,\xi^9)^2,(\xi^4,\xi^8)^2,(\xi^5,\xi^7)^2$
\\\hline
\end{tabular}
\end{center}
Here, $(\xi,\xi^{11})/(\xi^5,\xi^7)$ means either $(\xi,\xi^{11})$ or $(\xi^5,\xi^7)$, depending on the isomorphism type of the module.
\end{lemma}

We now give the analogue of Lemma \ref{lem:sl25restriction} for $p=7$. Again, we consider the range $p^a\leq 150$, so just $49$ in this case.

\begin{lemma}\label{lem:sl27restriction}
Let $H=\PSL_2(p^2)$, and let $M$ be a simple module for $H$. The restriction of $M$ to $L=\PSL_2(p)$ is as in Table \ref{t:tensorproductsp=7} for $p=7$ and Table \ref{t:tensorproductsp=11} for $p=11$.
\begin{table}
\begin{center}\begin{tabular}{ccc}
\hline Module & Restriction & Factors of restriction
\\\hline $1$ & $1$ & $1$
\\ $3$ & $3$ & $3$
\\ $4=2\otimes 2$ & $3\oplus 1$ & $3,1$
\\ $5$ & $5$ & $5$
\\ $7$ & $7$ & $7$
\\ $8=2\otimes 4$ & $5\oplus 3$ & $5,3$
\\ $9=3\otimes 3$ & $5\oplus 3\oplus 1$ & $5,3,1$
\\ $12=2\otimes 6$ & $7\oplus 5$ & $7,5$
\\ $15=3\otimes 5$ & $7\oplus 5\oplus 3$ & $7,5,3$
\\ $16=4\otimes 4$ & $7\oplus 5\oplus 3\oplus 1$ & $7,5,3,1$
\\ $21=3\otimes 7$ & $7\oplus P(5)$ & $7,5^2,3,1$
\\ $24=4\otimes 6$ & $7\oplus P(5)\oplus 3$ & $7,5^2,3^2,1$
\\ $25=5\otimes 5$ & $7\oplus P(5)\oplus 3\oplus 1$ & $7,5^2,3^2,1^2$
\\ $35=5\otimes 7$ & $7\oplus P(5)\oplus P(3)$ & $7,5^3,3^4,1$
\\ $36=6\otimes 6$ & $7\oplus P(5)\oplus P(3)\oplus 1$ & $7,5^3,3^4,1^2$
\\ $49=7\otimes 7$ & $7^{\oplus 2} \oplus P(5)\oplus P(3)\oplus P(1)$ & $7^2,5^4,3^4,1^3$
\\ \hline
\end{tabular}\end{center}\caption{Tensor products of modules for $\PSL_2(7)$}\label{t:tensorproductsp=7}
\end{table}
\begin{table}
\begin{center}\begin{tabular}{ccc}
\hline Module & Restriction & Factors of restriction
\\\hline $1$ & $1$ & $1$
\\ $3$ & $3$ & $3$
\\ $4=2\otimes 2$ & $3\oplus 1$ & $3,1$
\\ $5$ & $5$ & $5$
\\ $7$ & $7$ & $7$
\\ $8=2\otimes 4$ & $5\oplus 3$ & $5,3$
\\ $9$ & $9$ & $9$
\\ $9=3\otimes 3$ & $5\oplus 3\oplus 1$ & $5,3,1$
\\ $11$ & $11$ & $11$
\\ $12=2\otimes 6$ & $7\oplus 5$ & $7,5$
\\ $15=3\otimes 5$ & $7\oplus 5\oplus 3$ & $7,5,3$
\\ $16=4\otimes 4$ & $7\oplus 5\oplus 3\oplus 1$ & $7,5,3,1$
\\ $16=2\otimes 8$ & $9\oplus 7$ & $9,7$
\\ $21=3\otimes 7$ & $9\oplus 7\oplus 5$ & $9,7,5$
\\ $20=2\otimes 10$ & $11\oplus 9$ & $11,9$
\\ $24=4\otimes 6$ & $9\oplus 7\oplus 5\oplus 3$ & $9,7,5,3$
\\ $25=5\otimes 5$ & $9\oplus 7\oplus 5\oplus 3\oplus 1$ & $9,7,5,3,1$
\\ $27=3\otimes 9$ & $11\oplus 9\oplus 7$ & $11,9,7$
\\ $32=4\otimes 8$ & $11\oplus 9\oplus 7\oplus 5$ & $11,9,7,5$
\\ $33=3\otimes 11$ & $11\oplus P(9)$ & $11,9^2,3,1$
\\ $35=5\otimes 7$ & $11\oplus 9\oplus 7\oplus 5\oplus 3$ & $11,9,7,5,3$
\\ $36=6\otimes 6$ & $11\oplus 9\oplus 7\oplus 5\oplus 3\oplus 1$ & $11,9,7,5,3,1$
\\ $40=4\otimes 10$ & $11\oplus P(9)\oplus 7$ & $11,9^2,7,3,1$
\\ $45=5\otimes 9$ & $11\oplus P(9)\oplus 7\oplus 5$ & $11,9^2,7,5,3,1$
\\ $48=6\otimes 8$ & $11\oplus P(9)\oplus 7\oplus 5\oplus 3$ & $11,9^2,7,5,3^2,1$
\\ $49=7\otimes 7$ & $11\oplus P(9)\oplus 7\oplus 5\oplus 3\oplus 1$ & $11,9^2,7,5,3^2,1^2$
\\ $55=5\otimes 11$ & $11\oplus P(9)\oplus P(7)$ & $11,9^2,7^2,5,3^2,1$
\\ \hline
\end{tabular}\end{center}\caption{Tensor products of modules for $\PSL_2(11)$}\label{t:tensorproductsp=11}
\end{table}
\end{lemma}

\chapter{Some \texorpdfstring{$\PSL_2$s}{PSL2s} inside \texorpdfstring{$E_6$}{E6} in Characteristic 3}
\label{ch:psl2sinchar3}

In this short chapter we lay the groundwork for studying copies of $H=\PSL_2(3^a)$ (for $a\geq 2$) inside $F_4(3^b)$ by embedding $F_4(3^b)$ inside $E_6(3^b)$. Let $p=3$, and let $\bG$ denote $F_4(k)$, with $\hat{\bG}=E_6(k)$. Write $\hat{\mathscr X}$ for the set of positive-dimensional subgroups for $\hat{\bG}$. Suppose that $H$ is contained in $G=\bG^\sigma$. Note that the elements of $\Aut^+(\bG)$ extend to elements of $\Aut^+(\hat{\bG})$. Thus if $H$ is contained in a $\sigma$-stable, $N_{\Aut^+(\hat{\bG})}(H)$-stable subgroup $\bX$ of $\hat{\bG}$, and $\bX\leq \bG$, then $H$ is strongly imprimitive in $\bG$.

Importantly, we may assume that the graph automorphism does not lie in $N_{\Aut^+(\hat{\bG})}(H)$ however, since that does not restrict to an element of $\Aut^+(\bG)$. In particular, this means that we may argue with $M(E_6)$, rather than $M(E_6)\oplus M(E_6)^*$. Write $A$ for the subgroup of $\Aut^+(\hat{\bG})$ generated by inner, diagonal, and $p^i$-power field automorphisms, so not including the graph automorphism.

If $H$ stabilizes a $3$-space on $M(E_6)$, then $H$ lies inside a positive-dimensional subgroup of $\hat{\bG}$, namely the stabilizer $\bY$ of that $3$-space by Lemma \ref{lem:smallsubspaces}. In general, however, the intersection $\bY\cap \bG$ need not be positive dimensional, but we will at least show that $H$ is contained in a $\sigma$-stable, $N_A(H)$-stable, positive-dimensional subgroup of $\bG$, i.e., $H$ is strongly imprimitive in $\bG$.

We start with another useful result. Since the minimal module for $F_4$ is $25$-dimensional rather than $26$-dimensional, the action of $F_4$ on $M(E_6)$ is a uniserial module $1/25/1$. If a subgroup $H$ splits this extension, so acts as $1\oplus 1\oplus 25$, we cannot see this in $F_4$, but we can use the structure of $M(E_6)$ to nevertheless place $H$ inside a member of $\mathscr X$.

\begin{proposition}\label{prop:twotrivials}
Let $H$ be a subgroup of $\bG$. If $H$ centralizes a $2$-space on $M(E_6)$ then $H$ is strongly imprimitive in $\bG$.
\end{proposition}
\begin{proof} Suppose that $H$ centralizes a $2$-space on $M(E_6)$. If $H$ centralizes a line on $M(F_4)$ then we are done by Proposition \ref{prop:fixlineonMG}, so we may assume that $H$ does not. Therefore the $H$-fixed points on $M(E_6)$ are exactly $2$-dimensional. The $2$-space stabilizer $\bY$ is positive dimensional by Lemma \ref{lem:smallsubspaces}, and $N_A(H)$-stable by Proposition \ref{prop:intersectionstabilizers}. Also, since $\bG$ centralizes exactly one line, $\bY\neq \bG$. As it stabilizes the line stabilized by $\bG$, $\bY$ is a proper subgroup of $\bG$, and hence $H$ is strongly imprimitive in $\bG$, as needed.
\end{proof}

\begin{proposition}\label{prop:containedf4e6}
Let $H\cong \PSL_2(3^a)$ be a subgroup of $\bG$ for $a\geq 2$. If $H$ stabilizes a $3$-space on the $25$-dimensional minimal module $M(F_4)$ then $H$ is strongly imprimitive in $\bG$.
\end{proposition}
\begin{proof} Suppose that $H$ stabilizes a $3$-space on $M(E_6)$, and hence $H$, and the stabilizer $\bY$, lie inside some member of $\hat{\mathscr X}$ by Lemma \ref{lem:smallsubspaces}, $\bX$ say. By Propositions \ref{prop:fixlineonLG} and \ref{prop:fixlineonMG}, if $H$ stabilizes a line on either $M(F_4)$ or $L(F_4)$ then the result holds. Since the restriction of $L(E_6)^\circ$ to $F_4$ is the sum of $M(F_4)$ and $L(F_4)$, if $H$ is contained in $F_4$ and stabilizes a line on $L(E_6)^\circ$, it must stabilize a line on either $M(F_4)$ or $L(F_4)$, so we are done. If $H$ is a blueprint for $M(E_6)$ then $H$ is contained in a $\sigma$-stable, $N_A(H)$-stable, positive-dimensional subgroup of $\hat{\bG}$ by Proposition \ref{prop:intersectionstabilizers}, whence $H$ is strongly imprimitive in $\bG$.

If $H\not\leq \bX^0$, then the component group $\bX/\bX^0$ must have a subgroup isomorphic to $\PSL_2(3^a)$ for some $a\geq 2$. By \cite[Corollary 2]{liebeckseitz2004} and Remark (2) shortly before, $H$ must lie in a maximal-rank subgroup of $\hat{\bG}$ (since all others have soluble component group). Examining \cite[Table 10.3]{liebeckseitz2004}, if $\bX$ has maximal rank then $\bX$ must be the normalizer of a torus, $N_{\hat{\bG}}(\bT)$. Furthermore, since $\bX/\bX^0$ is the Weyl group $2\times \PSp_4(3)$, we see that $a=2$.

There is a unique class of subgroups $\PSL_2(9)$ inside $\PSp_4(3)=W(E_6)'$, and since there is a copy of $\PSL_2(9)\cong \Alt(6)$ normalizing a torus in the $A_5$ subgroup, we see that $H$ is conjugate to a subgroup of $A_5$, which of course means that we may replace $\bX$ by a different member of $\hat{\mathscr X}$ and have $H\leq \bX^0$. Thus we assume that $\bX$ is connected.

If $\bX=A_5A_1$ then $\bX$ stabilizes a line on $L(E_6)^\circ$, so we are done by the first paragraph. If $H$ lies in a subgroup $\bX$ that is either a $D_5$-parabolic subgroup or a copy of $F_4$ other than $\bG$, then $H$ stabilizes a line on $M(F_4)$ by Lemma \ref{lem:e6stabs}. Since $H$ already stabilizes a line on $M(F_4)$ since it lies in $\bG$, this means that $H$ stabilizes two different lines, in fact $H$ centralizes them since $H$ is perfect. Thus $H$ centralizes a $2$-space on $M(E_6)$, and therefore we apply Proposition \ref{prop:twotrivials}.

If $\bX$ is a parabolic subgroup, then by Proposition \ref{prop:paraisgood}, either $H$ is contained in a $\sigma$-stable, $N_A(H)$-stable parabolic of $E_6$, whose intersection with $F_4$ is positive dimensional by dimension counting, and still stable under $\sigma$ and $N_A(H)$, or $H$ is $E_6$-completely reducible. If $\bX$ is an $A_5$ parabolic, then this means that $H\leq A_5\leq A_5A_1$, which we have already considered. Similarly, if $\bX$ is an $A_2A_2A_1$-parabolic subgroup then $H\leq A_2A_2A_1\leq A_2A_2A_2$, which we consider below.

If $\bX=C_4$, then $\bX$ acts irreducibly on $M(E_6)$, as the exterior square of the natural module $M(C_4)$ minus a trivial summand. As $F_4$ has a trivial submodule and trivial quotient on $M(E_6)$, if $H\leq C_4\cap F_4$ then the exterior square of an $8$-dimensional module for $H=\PSL_2(3^a)$ or $\SL_2(3^a)$ (in this case the centre must act as a scalar) must have at least one trivial summand, and at least two trivial submodules. If the action of $H$ on the exterior square has two trivial summands, or three trivial submodules, then $H$ must lie inside a line stabilizer of $M(E_6)$ other than $F_4$, and hence stabilize a line on $M(F_4)$, so $H$ has been considered above.

The simple modules of dimension at most $8$ for $H$ have dimensions $1$, $3$ and $4$, and the exterior square of these are $0$, $3^*$ and $3\oplus 3$, none of which has a trivial submodule. Thus we must see repetitions of composition factors in $M(C_4)\downarrow_H$, and either three of one composition factor or two different repeats. If the factors have dimension $3,3,1,1$ then $M(C_4)\downarrow_H$ is semisimple by Lemma \ref{lem:char3sl2cohomology} and has the form $3_1^{\oplus 2}\oplus 1^{\oplus 2}$; $H$ is contained in the algebraic $A_1$ acting as $L(2)^{\oplus 2}\oplus L(0)^{\oplus 2}$ on $M(C_4)$, and this stabilizes the same subspaces of $M(E_6)$ as $H$. Thus $H$ is a blueprint for $M(E_6)$. If the factors have dimension $3,1^5$ or $4,1^4$ then $H$ has at least three trivial submodules on $M(C_4)$ (there is no module $1,1/4/1,1$ for $\PSL_2(9)$), and this results in two trivial submodules on $M(E_6)$, which was handled before. Hence we have an embedding of $\bar H=\SL_2(3^a)$ inside $\bX=C_4$ with $Z(\bar H)=Z(\bX)$.

The faithful simple modules for $\bar H$ of dimension at most $8$ have dimension $2$, $6$ and $8$. Of course, $\Lambda^2(8)$ has a single trivial submodule, so this cannot be the correct action of $\bar H$ on $M(C_4)$. If the composition factors have dimension $6,2$ then $M(C_4)\downarrow_{\bar H}$ is semisimple of the form $6_{i,j}\oplus 2_l$ for some $i,j,l$, and 
\[ \Lambda^2(6_{i,j})=(1/4_{j,j+1}/1)\oplus 9_{i,j},\quad \Lambda^2(2_l)=1.\]
Thus if $H$ stabilizes a $3$-space on $M(E_6)$, it must come from $6_{i,j}\otimes 2_l$, hence $l=i$ or $l=j$. If $l=i$, we have
\[6_{i,j}\otimes 2_i=3_j\oplus 9_{i,j},\]
and $H$ stabilizes a unique $3$-space on $M(F_4)$ and $M(E_6)$, the stabilizer of which is easily seen to be positive dimensional by placing it inside an algebraic $A_1$ inside $C_4$, or by applying Lemma \ref{lem:smallsubspaces}. Since $H$ stabilizes a unique $3$-space, this stabilizer is $\sigma$-stable and $N_A(H)$-stable by Proposition \ref{prop:intersectionstabilizers}, and therefore $H$ is strongly imprimitive. On the other hand, if $l=j$ then
\[6_{i,j}\otimes 2_j=2_i\otimes (2_j/2_{j+1}/2_j),\]
which will only stabilize a $3$-space if $i=j$, which is not allowed, so $H$ cannot stabilize a $3$-space in this case.

Thus we may assume that the composition factors of $M(C_4)\downarrow_H$ all have dimension $2$. There are two orbits of $2$-spaces on $M(C_4)$, with stabilizers a parabolic subgroup and $A_1C_3$. Since we may assume that $H$ does not lie in a parabolic subgroup, all $2$-spaces stabilized by $H$ have $A_1C_3$ stabilizer, and in this case the $2$-space is complemented in the stabilizer, hence in $H$. Thus $M(C_4)$ is a sum of $H$-stable $2$-spaces, and therefore $H$ centralizes at least a $4$-space on the exterior square of $M(C_4)$, a possibility that has already been considered. Thus $\bX=C_4$ is complete.

We are left with the irreducible $G_2$, the $A_2G_2$ subgroup, and the $A_2A_2A_2$ maximal-rank subgroup. In these cases we will show that $H$ is a blueprint for $M(E_6)$.

If $H$ is contained in $\bX=G_2$ then from the list of maximal subgroups in \cite{kleidman1988} either $H$ acts on the natural module (up to Frobenius) as $3_1^{\oplus 2}\oplus 1$, or lies in a diagonal subgroup of $A_1A_1$ acting as $4_{1,i}\oplus 3_1$ for any $i\geq 1$. In both cases $H$ is contained in an algebraic $A_1$ subgroup $\bY$.

The subgroup $\bX$ acts irreducibly as $L(20)$ on $M(E_6)$, and these two copies of $H$ act on $M(E_6)$ as
\[ 3_1^{\oplus 3}\oplus (1/4_{1,2}/1)^{\oplus 3},\qquad 9_{1,i}\oplus (1/4_{1,2}/1)\oplus (4_{1,i}/(2_2\otimes 2_i)/4_{1,i}),\]
where of course $2_2\otimes 2_i$ is $4_{2,i}$ if $i>2$ and $1\oplus 3_2$ if $i=2$. The subgroups $\bY$ containing them act in the same way, and stabilize the same subspaces as $H$, so that $H$ is a blueprint for $M(E_6)$.

If $\bX=A_2G_2$, then $\bX$ acts on $M(E_6)$ as $(10,10)\oplus (02,00)$: if $H$ lies inside the $G_2$ factor then it centralizes a $6$-space on $M(E_6)$, so definitely lies inside a line stabilizer of $F_4$, and hence we may assume that $H$ projects along the $A_2$ factor as $3_1$. Along the $G_2$ factor it can act as either $3_i^{\oplus 2}\oplus 1$ or $4_{i,j}\oplus 3_i$, for any $i,j\geq 1$ with $i\neq j$. In the first case we obtain
\[ 3_1\oplus 9_{1,i}^{\oplus 2}\oplus (1/4_{1,2}/1)\quad \text{and}\quad 3_1^{\oplus 3}\oplus (1/4_{1,2}/1)^{\oplus 3},\]
for $i>1$ and $i=1$ respectively, and in the second case we obtain
\[\begin{gathered} 3_1\oplus (1/4_{1,2}/1)\oplus (4_{1,i}/(2_2\otimes 2_i)/4_{1,i})\oplus (1/4_{1,2}/1),
\\ 9_{1,i}\oplus (4_{1,j}/(2_2\otimes 2_j)/4_{1,j})\oplus (1/4_{1,2}/1),\qquad
9_{1,i}\oplus 12_{i,j,1}\oplus (1/4_{1,2}/1),\end{gathered}\]
for $i=1$, $j=1$ and $i,j\neq 1$ respectively. Again, in all cases, the algebraic subgroup containing $H$ stabilizes the same subspaces as $H$, so again $H$ is a blueprint for $M(E_6)$.

Finally, we have $\bX=A_2A_2A_2$, which acts on $M(E_6)$ as (up to duality) the sum of the three possible configurations of natural times natural times trivial. If we act trivially on one or two of them we obtain, up to field automorphism,
\[ 3_1^{\oplus 6}\oplus 1^{\oplus 9},\qquad 3_1^{\oplus 7}\oplus (1/4_{1,2}/1),\qquad 9_{1,i}\oplus 3_1^{\oplus 3}\oplus 3_i^{\oplus 3},\]
depending on whether we have one non-trivial module $3_1$, two non-trivial and both $3_1$, and two non-trivial and $3_1$ and $3_i$ respectively. In the first two cases, $H$ centralizes a $2$-space on $M(E_6)$, and we already showed that this means that $H$ is strongly imprimitive in $\bG$. In the third case, $H$ does not stabilize a line on $M(E_6)$, so cannot lie in $\bG$.

Thus we may assume that $H$ acts along the first factor as $3_1$, the second as $3_i$ and the third as $3_j$. In this case we have one of 
\[3_1^{\oplus 3}\oplus (1/4_{1,2}/1)^{\oplus 3},\qquad 3_1\oplus (1/4_{1,2}/1)\oplus 9_{1,j}^{\oplus 2},\qquad 27_{1,i,j},\]
depending on whether we have $i=j=1$, $1=i\neq j$, and $1\neq i\neq j\neq 1$ respectively. As with the other cases, each of these is contained in an algebraic $A_1$ stabilizing the same subspaces of $M(E_6)$, so again $H$ is a blueprint for $M(E_6)$.

In the first paragraph we showed that if $H$ is a blueprint for $M(E_6)$ then $H$ is strongly imprimitive, so in all cases we have shown that $H$ is strongly imprimitive.
\end{proof}

\chapter{Proof of the Theorems: Strategy}
\label{ch:strategy}

In this chapter we discuss the techniques that we will use in proving that a given $\PSL_2(p^a)$ subgroup $H$ of $\bG$ is strongly imprimitive in $\bG$.

The first step is usually to use the dimensions of modules and the traces of semisimple elements to produce a list of potential sets of composition factors for the action of $H$ on $M(\bG)$, which in Chapter \ref{ch:notation} we called conspicuous sets of composition factors. For many groups this list is small, but as the sizes of $\bG$ and $H$ grow the number grows larger and we need more efficient methods that cut this number down, for example only considering possible multisets of dimensions that have either no modules of dimension $1$ or more modules of dimension $2(p-1)$ than modules of dimension $1$, at least for $p^a$ odd and not equal to $9$, i.e., modules of positive pressure (see Lemmas \ref{lem:pressure} and \ref{lem:cohomologysimple}).

\medskip

Having done this, we can assume we know the composition factors of $M(\bG)\downarrow_H$, and we have a few ways to proceed.

\begin{enumerate}
\item\label{li:strata} Each semisimple element of $H$ belongs to a semisimple class of $\bG$, and the trace of the element on $M(\bG)$ is often enough to determine this class uniquely, sometimes to a small number of possibilities. In particular, this yields the traces of the semisimple elements on $L(\bG)$ (with potentially several possibilities). Hence the composition factors of $M(\bG)\downarrow_H$ yields a list, often a list with one element, of the possible sets of composition factors for $L(\bG)\downarrow_H$. Sometimes this has no elements, of course only if there is no embedding of $H$ with these composition factors. Other times $L(\bG)\downarrow_H$ has non-positive pressure, so we again show that $H$ is strongly imprimitive. Otherwise, we may analyse both $M(\bG)\downarrow_H$ and $L(\bG)\downarrow_H$ using the techniques below. We will occasionally employ Lemma \ref{lem:relatingvminlg} in this regard.

\item\label{li:stratb} We can easily compute $\Ext^1$ between the composition factors of $M(\bG)\downarrow_H$ using \cite[Corollary 4.5]{andersen1983} and determine if $M(\bG)\downarrow_H$ is semisimple or not. If it is, the action of a unipotent element $u$ must match one of the unipotent classes of $G$, whose actions on $M(\bG)$ and $L(\bG)$ are tabulated in \cite{lawther1995}. If it does not appear, or is generic (see Definition \ref{defn:genericunipotent}), then we are done.

\item\label{li:stratc} Let $V$ be some rational $k\bG$-module. If $V\downarrow_H$ is not semisimple, and $V$ is self-dual (i.e., all cases except when $G=E_6$ and $V=M(E_6)$) then in order for a composition factor to appear in the socle and not be a summand, it must occur with multiplicity at least $2$. This allows us to cut down the possibilities for the socle of $V\downarrow_H$.

\item\label{li:stratd} Let $V$ be as above. If the socle of $V\downarrow_H$ is $W$, then $V$ is a submodule of $P(W)$, where $P(W)$ denotes the projective cover of $W$. In particular, it is a submodule of the $\cf(V)$-radical of $P(W)$, where we recall that $\cf(V)$ is the set of composition factors of $V$. Thus the $\cf(V)$-radical of $P(W)$ needs to contain at least as many copies of each composition factor as there are in $V$, and further analysis of this radical can eliminate more cases.

The Magma commands \texttt{Ext} and \texttt{MaximalExtension} can be used to construct this radical, even without the ambient projective module, so are useful when the projective is too large to compute directly. If \texttt{W} is the putative socle of $V$, and \texttt{I} is a set of simple $kH$-modules, we construct the \texttt{I}-radical of \texttt{V} using the following program.

\begin{verbatim}
function ComputeRadical(W,I)

W0:=Dual(W);
repeat
  n:=Dimension(W0);
  for M in I do
    E,rho:=Ext(W0,Dual(M)); 
    W0:=MaximalExtension(W0,Dual(M),E,rho);
  end for;
until n eq Dimension(W0);
return Dual(W0);

end function;
\end{verbatim}

\item\label{li:strate} We can use Lemma \ref{lem:q/2restricted}: suppose that $H=\PSL_2(p^a)$ embeds in $\bG$, and an algebraic $A_1$-subgroup $\bX$ embeds in $\bG$, such that for some module $V$, the highest weights of the composition factors of both $H$ and $\bX$ on $V$ are the same. Assume furthermore that the composition factors of $V\downarrow_{\bX}$ satisfy the hypotheses of Lemma \ref{lem:q/2restricted}. We wish to conclude that $H$ is a blueprint for $V$. In order to do this, an element $x$ in $H$ of order $(p^a\pm 1)/2$ must be guaranteed to come from a class intersecting $\bX$. If the semisimple class containing $x$ is determined by its eigenvalues on $V$ then this is true, but this is not true for every semisimple class, so we will have to check when we use the lemma.

\item\label{li:stratf} In a similar vein, we can look for elements of $\bG$ that do not lie in $H$ and yet stabilize some eigenspaces of an element of $H$ on a module $V$: if $\zeta_1,\dots,\zeta_r$ are roots of unity and $y$ acts as a scalar on each $\zeta_i$-eigenspace of an element $x\in H$ (i.e., preserves all subspaces of the eigenspace), then $y$ stabilizes any subspace of $V$ on which $x$ acts with eigenvalues some of the $\zeta_i$. In particular, if there is a submodule $W$ of $H$ with this property then $\gen{H,y}$ stabilizes $W$. Of course, it might be that $\gen{H,y}$ is almost simple, say $\PGL_2(p^a)$ for example, so we need to exclude this case by finding other such elements,  proving that the index of $H$ in this group is not $2$, or applying Corollary \ref{cor:sl2nopgl}.

In practice when looking for an element of order $an$ to power to an element of order $n$, we often use Litterick's program \cite[Chapter 7]{litterick} to construct all the eigenvalues of semisimple elements of order $an$ of $V$ first. We then take powers of these eigenvalues to see whether there are elements of order $an$ powering to a fixed element $x$ of order $n$ and stabilizing a given subspace. Of course, this implies that the eigenvalues of $x$ on $V$ determine $x$. If this is not the case, we would have to do it for all possible semisimple classes to which $x$ may belong.

If the order $an$ is too big to construct all semisimple classes of that order, we can employ the preimage trick, described at the end of Section \ref{sec:blueprints}, to obtain the semisimple elements that power to a given one.

A heuristic, but not a formal statement, is that the larger the order of an element, the greater the proportion of elements of that order are blueprints for a given module. For example, $\mathcal X_{M(F_4)}=\{1,\dots,18\}$, so there are semisimple elements of order $17$ in $F_4$ that are not blueprints for $F_4$. However, 228 of the 230 classes of semisimple elements of order $17$ are blueprints for $M(F_4)$. The easiest way to check this is to find elements of order $34$ that square to a given element of order $17$ and have the same eigenspaces. Since elements of order $34$ are blueprints for $M(F_4)$, this shows that a given element of order $17$ is a blueprint.

\item If $\bG=E_6,E_7$ and $p=h-1$ where $h$ is the Coxeter number of $\bG$, then in one case we prove that $H$ stabilizes an $\slf_2$-subalgebra of $L(\bG)$. We can then apply Corollary \ref{cor:sl2coxeter-1} on positive-dimensionality of such a stabilizer.
\end{enumerate}

Some combination of these ideas is usually enough to prove that $H=\PSL_2(p^a)$ lies inside a positive-dimensional subgroup of the algebraic group $\bG$, by showing that $H$ stabilizes some subspace $W$ of some module $V$ for $\bG$ that has positive-dimensional stabilizer. This proves that $H$ is not Lie primitive, but to prove strong imprimitivity we need to know more about $W$. We will sometimes prove directly that the intersection of the stabilizers of the subspaces in the $N_{\Aut^+(\bG)}(H)$-orbit of $W$ is also positive dimensional (often this orbit has length $1$) and then apply Proposition \ref{prop:intersectionstabilizers}. Usually we will prove that $H$ stabilizes a $1$-space on $L(\bG)^\circ$ and apply Proposition \ref{prop:fixlineonLG}, a $1$- or $2$-space on $M(\bG)$ and apply Propositions \ref{prop:fixlineonMG} or \ref{prop:fix2spaceonMG} respectively, or that $H$ is a blueprint for $M(\bG)$ or $L(\bG)$ and apply Proposition \ref{prop:intersectionstabilizers} (we need $M(E_6)\oplus M(E_6)^*$, and we cannot use $M(F_4)$ if $p=2$).

The rest of the chapters will prove these facts for the various $H$ and $\bG$ under consideration.

\chapter{The Proof for \texorpdfstring{$F_4$}{F4}}
\label{ch:f4}
In this chapter, $k$ is an algebraically closed field of characteristic $p\geq 2$ and $\bG=F_4(k)$. Let $H\cong \PSL_2(p^a)$ be a subgroup of $\bG$.

Theorem \ref{thm:goodblueprint} states that if $p=2$ and $a\geq 5$ then $H$ is a blueprint for $L(F_4)$, and if $p$ is odd and $p^a\geq 37$ then $H$ is a blueprint for $M(F_4)$. In both cases, $H$ is strongly imprimitive by Proposition \ref{prop:blueprintissi}. Thus in what follows we may assume that $p^a\leq 31$.

Let $L=\PSL_2(p)\leq H$ and let $u$ denote a unipotent element of $L$ of order $p$. The possibilities for the Jordan block structures of $u$ on $M(F_4)$ and $L(F_4)$ are given in \cite[Tables 3 and 4]{lawther1995}. Recall the definition of a generic unipotent element from Definition \ref{defn:genericunipotent}.

We will prove that $H$ is always strongly imprimitive. Since $H$ cannot act irreducibly on $M(F_4)$ or $L(F_4)$ by \cite{liebeckseitz2004a}, if $H$ is a blueprint for either module then it is strongly imprimitive, by Proposition \ref{prop:intersectionstabilizers}.

\section{Characteristic 2}
\label{sec:f4p=2}
Let $p=2$, so that $a\leq 4$. All semisimple elements of $\bG$ lie inside $D_4$, which centralizes a $2$-space on $M(F_4)$. Furthermore, an element of order $2^a+1$ in $H$ has a fixed point only on the trivial simple module. We therefore see that $M(F_4)\downarrow_H$ has at least two trivial composition factors. In particular, Lemma \ref{lem:char2pressure1} applies in this situation (as an involution cannot act projectively on $M(F_4)$ from \cite[Table 3]{lawther1995}), and so the pressure of $M(F_4)\downarrow_H$ has to be at least $2$ for $H$ not to stabilize a line on $M(F_4)$, and therefore be strongly imprimitive via Proposition \ref{prop:fixlineonMG}.

If $a=1$ then $H$ is soluble, and if $a=2$ then $H$ stabilizes either a $1$- or $2$-space on $M(F_4)$ by \cite[Proposition 5.4]{craven2015un2}, hence is strongly imprimitive by Propositions \ref{prop:fixlineonMG} and \ref{prop:fix2spaceonMG}. 

We first consider $a=3$, then the case $a=4$ afterwards. If one simply wants to show that $H$ must stabilize a $1$- or $2$-space on $M(F_4)$, one can proceed as in the proof of Proposition \ref{prop:e6char2a=3} below, and prove the next result in a few lines. However, we now prove a stronger result.

\begin{proposition}\label{prop:f4char2a=3} Suppose that $p=2$ and $a=3$. If $H$ does not stabilize a line on $L(F_4)$, then up to field automorphism of $H$ the composition factors of $M(F_4)\downarrow_H$ are
\[ 4_{1,3}^2,2_1^4,2_2,2_3^2,1^4,\]
and $H$ stabilizes a $2$-space on $M(F_4)$.
\end{proposition}
\begin{proof} Using a computer, we use the traces of semisimple elements of order at most $17$ (i.e., all of them) to find all conspicuous sets of composition factors for $M(F_4)\downarrow_H$. There are 63 such sets, too many to simply list here, so we need to cut down on this number.

First we exclude those whose traces on semisimple elements imply that there is no corresponding set of composition factors for $L(F_4)\downarrow_H$, or equivalently $L(\lambda_1)\downarrow_H$ since $L(F_4)$ has composition factors $L(\lambda_4)=M(F_4)$ and $L(\lambda_1)$ (see (\ref{li:strata}) from Chapter \ref{ch:strategy}). Thus we are left with 53 conspicuous sets of composition factors for $M(F_4)\downarrow_H$.

These 53 sets fall into orbits, first under the field automorphism of order $3$ of $H$, and second under the graph morphism swapping $L(\lambda_4)$ and $L(\lambda_1)$. The orbits have lengths $2$, $3$ and $6$, and there are ten orbits in total. If $M(F_4)\downarrow_H$ has pressure at most $1$, then as stated at the start of the section, $H$ stabilizes a line on $M(F_4)$ by Lemma \ref{lem:char2pressure1}. Six of these ten orbits contain sets of factors with non-positive pressure, and two more have pressure $1$. Thus for eight of the ten orbits we are guaranteed that $H$ stabilizes a line on either $L(\lambda_4)$ or $L(\lambda_1)$. Since $L(F_4)$ has these modules as two composition factors, this means $H$ stabilizes a line on $L(F_4)$, as claimed.

We are left with two orbits, one of length six and one of length three, with representatives
\[  8,4_{1,3},4_{2,3},2_1^2,2_2,2_3,1^2,\qquad \text{and}\qquad 4_{1,3}^2,2_1^4,2_2,2_3^2,1^4.\]
Thus we may assume that $H$ acts on $M(F_4)$ with one of these two sets of composition factors.

\medskip

\noindent \textbf{Case 1}: Since the only non-trivial simple module appearing more than once is $2_1$, and $M(F_4)$ is self-dual, $M(F_4)\downarrow_H$ is the sum of an indecomposable module with socle $2_1$ and a semisimple module. If the semisimple module has a trivial composition factor then $H$ stabilizes a line on $M(F_4)$, so we may assume that both trivial factors lie in this indecomposable summand. Write $V$ for this indecomposable module, so $\soc(V)=2_1$.

We know that $M(F_4)\downarrow_H$ is self-dual, so if there is a non-split extension between two non-trivial simple modules $A$ and $B$ inside $M(F_4)\downarrow_H$, we must also have one between $B$ and $A$. Thus not both $A$ and $B$ can appear with multiplicity $1$ in the composition factors of $M(F_4)\downarrow_H$.

If $A$ is $4_{2,3}$, then $B$ can only be $2_3$ by Lemma \ref{lem:extforsl22a}. However, both of these appear with multiplicity $1$ in $M(F_4)\downarrow_H$, so there is no extension between them. In particular, $4_{2,3}$ must be a summand of $M(F_4)\downarrow_H$.

The $\{1,2_1,2_2,2_3,4_{1,3}\}$-radical $M_1$ of $P(2_1)$ has three trivial composition factors, and must contain $V$. Since $\topp(V)=2_1$ as well, we may take the $\{2_1\}'$-residual of $M_1$ to obtain another module $M_2$, which also contains $V$. The module $M_2$ has structure
\[ 2_1/1/2_1,2_2/1,4_{1,3}/2_1,\]
and $u$ acts projectively on $M_2$. However, an involution cannot act projectively on $M(F_4)$ (see \cite[Table 3]{lawther1995}) so both trivial factors cannot lie in $V$. Hence $H$ must stabilize a line on $M(F_4)$, as claimed.

\medskip

\noindent \textbf{Case 2}: There are modules with these composition factors that do not stabilize a line, for example
\[ 4_{1,3}^{\oplus 2}\oplus (2_1/1/2_2/1/2_1)\oplus (2_1/1/2_3)\oplus (2_3/1/2_1),\]
and the Jordan blocks of $u$ on this module do appear in \cite[Table 3]{lawther1995}. We claim that $M(F_4)\downarrow_H$ always has a $1$- or $2$-dimensional submodule.

To see the claim, note that otherwise $\soc(M(F_4)\downarrow_H)=4_{1,3}$. Hence $M(F_4)\downarrow_H$ is a submodule of $P(4_{1,3})$, but $P(4_{1,3})$ has structure
\[ 4_{1,3}/2_1/1/2_2/1/2_1/4_{1,3},\]
and there are many reasons why this cannot work: the dimension is $16$, $2_3$ is not involved in it, there are not enough factors that are $1$ or $2_1$, the involution $u$ acts projectively on it, and so on.
\end{proof}

There is a copy of $\SL_2(8)$ inside $F_4$ that does indeed not stabilize a line on $L(F_4)$, inside $\tilde A_2A_2$. The projection of $H$ along  the $\tilde A_2$ factor acts on $M(A_2)$ as $2_1/1$, and the projection along the $A_2$ factor acts on $M(A_2)$ as $1/2_3$. The product of these two modules is an indecomposable module $2_1/1/2_3,4_{1,3}$ with dual $2_3/1,4_{1,3}/2_1$, and the product of $2_1/1$ and its dual is
\[ (2_1/1/2_2/1/2_1)\oplus 1,\]
yielding an embedding into $F_4$ with the required property. (Remember that the trivial in the last decomposition is removed when considering $M(F_4)$, see Appendix \ref{app:actions}.)

\medskip

We now turn to $a=4$. Almost exactly the same result holds in this case.

\begin{proposition}\label{prop:f4char2a=4} Suppose that $p=2$ and $a=4$. If $M(F_4)\downarrow_H$ does not stabilize a line on $L(F_4)$, then up to field automorphism of $H$ and graph automorphism of $\bG$ the composition factors of $M(F_4)\downarrow_H$ are
\[ 4_{1,3}^2,2_1^4,2_2,2_3^2,1^4,\]
and $H$ stabilizes a $2$-space on $M(F_4)$.
\end{proposition}
\begin{proof} We proceed as in Proposition \ref{prop:f4char2a=3}, starting by producing all conspicuous sets of composition factors using the traces of semisimple elements of order up to $17$, this time finding $146$ conspicuous sets of composition factors for $M(F_4)\downarrow_H$. Sixteen of these sets of composition factors have no corresponding set of composition factors on $L(F_4)$, so we reduce to $130$ sets of composition factors. These fall into eighteen orbits under field automorphism of $H$ and graph automorphism of $\bG$: fifteen of length $8$, two of length $4$ and one of length $2$.

If $M(F_4)\downarrow_H$ has pressure at most $1$, then as stated at the start of the section, $H$ stabilizes a line on $M(F_4)$ by Lemma \ref{lem:char2pressure1}. Eleven of these orbits contain a conspicuous set of composition factors with negative pressure, and a further two with factors with pressure $0$, hence $H$ fixes a line on $M(F_4)$, and thus $L(F_4)$. We can exclude factors with pressure $1$ as well, eliminating a further two orbits. There remain three orbits, each of length $8$, with representatives
\[ 8_{1,2,3},4_{1,3},4_{2,3},2_1^2,2_2,2_3,1^2,\qquad 8_{1,2,4},4_{1,3},4_{2,4},2_1^2,2_2,2_3,1^2,\qquad 4_{1,3}^2,2_1^4,2_2,2_3^2,1^4.\]
Thus we assume that $M(F_4)\downarrow_H$ has composition factors one of the above three sets.

\medskip

\noindent \textbf{Cases 1, 2}: We argue as in Proposition \ref{prop:f4char2a=3}, to see that since the only non-trivial simple module appearing more than once is $2_1$, $M(F_4)\downarrow_H$ is a sum of a semisimple module and a self-dual submodule $V$ of $P(2_1)$ that has top $2_1$. For $i=3,4$, the $\{1,2_2,2_3,4_{1,3},4_{2,i},8_{1,2,i}\}$-radical of the quotient module  $P(2_1)/2_1$, lifted back to $P(2_1)$, is
\[ 2_3/1,4_{2,3}/2_2,2_3/1,4_{1,3}/2_1\quad\text{and}\quad 1/2_2,2_3/1,4_{1,3}/2_1,\]
for $i=3$ and $i=4$ respectively.

This module must contain $V$ as a submodule. An involution acts projectively on this module, so if both trivial composition factors of $M(F_4)\downarrow_H$ lie in $V$ then $u$ acts projectively on $M(F_4)\downarrow_H$. However, $u$ cannot act projectively on $M(F_4)$ by \cite[Table 3]{lawther1995}. Thus there is a trivial summand in $M(F_4)\downarrow_H$, as needed.

\medskip

\noindent \textbf{Case 3}: Since we constructed an example of this embedding not stabilizing a line on $M(F_4)$ inside the $A_2\tilde A_2$ subgroup just after Proposition \ref{prop:f4char2a=3}, we will not be able to prove that it stabilizes a line on $M(F_4)$. However, it does stabilize a $2$-space, as it did for $a=3$.

We follow the same proof as for $a=3$ as well. The $\{1,2_1,2_2,2_3,4_{1,3}\}$-radical of $P(4_{1,3})$ is 
\[ 2_2/1/2_1,2_3/4_{1,3},\]
so $\soc(M(F_4)\downarrow_H)$ cannot be just $4_{1,3}$. Hence $H$ stabilizes either a line or a $2$-space on $M(F_4)$, as claimed.
\end{proof}

\section{Characteristic 3}
\label{sec:f4p=3}

For $p=3$, since $p^a\leq 31$, $a\leq 3$. Of course, $\PSL_2(3)$ is soluble, so $a=2,3$. The proof for $a=3$ is significantly easier than for $a=2$. In fact, for $a=2$ we need to embed $H$ inside the $F_4A_1$ subgroup of $E_7$ to prove that $H$ is not Lie primitive.

\begin{proposition}\label{prop:f4char3a=3} If $p=3$ and $a=3$, then $H$ is a blueprint for $M(F_4)$.
\end{proposition}
\begin{proof} Let $x$ be a semisimple element of order $13$ in $H$. Of the $104$ semisimple classes of elements of order $13$ in $F_4$, all but seven contain blueprints for $M(F_4)$. In each case, there are elements of order $26$ that square to them and preserve the eigenspaces on $M(F_4)$, and since all elements of order $26$ are blueprints for $M(F_4)$ by Theorem \ref{thm:blueprintsminimal}, this shows that those classes contain blueprints (see (\ref{li:stratf}) from Chapter \ref{ch:strategy}).

There are 40 conspicuous sets of composition factors for $M(F_4)\downarrow_H$. Removing the conspicuous sets of factors for which $x$ is a blueprint leaves just seven, three up to field automorphism of $H$. Representatives of these orbits are
\[12_{2,3,1},9_{1,2},4_{1,2},\qquad 12_{2,3,1},9_{1,3},4_{1,2},\qquad 4_{1,2}^2,4_{1,3}^2,4_{2,3}^2,1.\]
Let $\zeta$ denote a primitive $13$th root of unity, and let $\theta$ denote a primitive $26$th root of unity with $\theta^2=\zeta$. By choosing $\zeta$ appropriately, $x$ acts on $4_{1,2}$ with eigenvalues $\zeta^{\pm 1}$ and $\zeta^{\pm 2}$.

\medskip

\noindent \textbf{Case 1}: $x$ acts on $M(F_4)$ with eigenvalues
\[ 1,(\zeta^{\pm 1})^2,(\zeta^{\pm 2})^3,(\zeta^{\pm 3})^2,(\zeta^{\pm 4})^2,\zeta^{\pm 5},(\zeta^{\pm 6})^2.\]
There is an element $\hat x$ of order $26$ in $\bG$ that squares to $x$ and has eigenvalues
\[ 1,(\theta^{\pm 1})^2,(\theta^{\pm 2})^3,(\theta^{\pm 3})^2,(\theta^{\pm 4})^2,\theta^{\pm 5},\theta^{\pm 6},(-\theta^{\pm 6}).\]
This does not stabilize all the eigenspaces of $x$, but it only splits the $\zeta^{\pm 6}$-eigenspaces, which are contained inside the $12$-dimensional factor. Hence $\hat x$ stabilizes all subspaces stabilized by $H$. Since $\hat x$ has order $26$, it is a blueprint for $M(F_4)$ by Theorem \ref{thm:blueprintsminimal}.

\medskip

\noindent \textbf{Cases 2, 3}: In the second case, $M(F_4)\downarrow_H$ is semisimple because it is self-dual. In the third case, by applying a field automorphism if necessary, $\soc(M(F_4)\downarrow_H)$ contains $4_{1,2}$.

The trace of $x$ on $M(F_4)$ in both cases is $0$. There is an element $\hat x$ in $\bG$, of order $26$, such that $\hat x^2=x$ and $\hat x$ acts on $M(F_4)$ with eigenvalues
\[ 1,(\theta^{\pm 1})^2,(\theta^{\pm 2})^2,(\theta^{\pm 3})^2,(\theta^{\pm 4})^2,\theta^{\pm 5},(-\theta^{\pm 5}),\theta^{\pm 6},(-\theta^{\pm 6}),\]
so $\hat x$ stabilizes a $4_{1,2}$ in the socle of the second and third cases. Since $\hat x$ is a blueprint for $M(F_4)$, and it stabilizes $4_{1,2}$, this means that the stabilizer $\bY$ of $4_{1,2}$ is positive dimensional. In particular, $H$ lies in a member of $\mathscr X$.

Examining the list of maximal positive-dimensional subgroups of $\bG$ from Appendix \ref{app:actions}, if $H$ acts on $M(F_4)$ with factors $12,9,4$ then the only member of $\mathscr X$ in which $H$ can lie is $A_1C_3$. This subgroup acts with factors of dimension $12$ and $13$, so any positive-dimensional subgroup containing $H$ must also stabilize the $12$. In particular, $\bY$ stabilizes the $12$. Since $\bY$ must act semisimply on $M(F_4)$, $H$ and $\bY$ stabilize the same subspaces of $M(F_4)$, i.e., $H$ is a blueprint for $M(F_4)$. Thus we are left with Case 3.

We now run through the elements of $\mathscr X$, proving that $u$ lies in the generic class $A_2$ (see Definition \ref{defn:genericunipotent}), and thus $H$ is a blueprint for $M(F_4)$ by Lemma \ref{lem:genericunipotent}. See Appendix \ref{app:actions} for the composition factors of the maximal positive-dimensional subgroups of $\bG$ on $M(F_4)$.

We cannot embed $H$ in a maximal parabolic or $A_2\tilde A_2$ as the dimensions of the composition factors of $M(F_4)\downarrow_H$ and $M(F_4)\downarrow_{A_2\tilde A_2}$ are not compatible. We may embed $H$ only in $B_4$, $A_1C_3$ and $A_1G_2$.

Suppose that $H\leq B_4$, which acts as $9\oplus 16$ on $M(F_4)$. Since there is no uniserial module of the form $4_{i,j}/1/4_{m,n}$ by Lemma \ref{lem:no414}, and $H$ acts with composition factors of dimensions $4^6,1$ on $M(F_4)$, the restriction of the $9$ to $H$ must be $4_{i,j}\oplus 4_{m,n}\oplus 1$. The element $u$ acts on $4_{i,j}$ with blocks $3,1$, so $u$ acts with at least three blocks of size $1$ on $M(F_4)$, and at least six blocks of size $3$, one from each factor $4_{i,j}$. Thus $u$ acts on $M(F_4)$ with blocks $3^6,1^3$: from Table \ref{t:unipotentF4} we see that $u$ lies in the generic class $A_2$.

If $H$ embeds in $\bX=A_1C_3$ then we may assume up to field automorphism that $H$ acts on $M(A_1)$ as $2_1$. Since $\bX$ acts on $M(F_4)$ as $(1,100)\oplus (0,010)$, we need $6$-dimensional modules for $H$ whose tensor product with $2_1$ only have $4$-dimensional composition factors, each appearing at most twice, and there are three of these: $2_2^{\oplus 2}\oplus 2_3$, $2_2\oplus 2_3^{\oplus 2}$, and $6_{3,1}$. However, the exterior square of $L(100)$ for $C_3$ is a uniserial module $0/010/0$ (in characteristic 3). Thus the exterior square of one of these modules, minus two trivials, is the other summand of $M(F_4)\downarrow_H$. None of these has the correct exterior square, so $H$ does not embed in $A_1C_3$.

We are left with $A_1G_2$, which acts on $M(F_4)$ with the composition factors $(L(2),L(10))$ of dimension $21$ and $(L(4),L(00))$ of dimension $4$. This embedding of $H$ is impossible: $H$ acts on $M(A_1)$ as, up to field automorphism, $3_1$, and so the action of $H$ on the minimal module for $G_2$ cannot have a trivial or $3$-dimensional composition factor, because the product with $3_1$ would not be correct. But then one cannot make a $7$-dimensional module at all, a contradiction.

Thus if $H$ embeds into $\bG$ with these factors on $M(F_4)$ then it is a blueprint, as needed.
\end{proof}

We now consider $a=2$, which we did not consider in \cite{craven2015un2} because we could not produce a complete answer there. There is no known proof of the result from within $F_4$ itself either, and here we have to move into $E_7$, by embedding $F_4$ in the maximal subgroup $F_4A_1$. The next result establishes as much of the result as we can within $F_4$, and then we embark on a proof of the final part afterwards.

\begin{proposition}\label{prop:f4char3a=2} Suppose that $p=3$ and $a=2$. One of the following holds:
\begin{enumerate}
\item\label{propi:f43a} $H$ stabilizes a line on $M(F_4)$ or $L(F_4)$;
\item\label{propi:f43b} $H$ stabilizes a unique $3$-space on $M(F_4)$;
\item\label{propi:f43c} $H$ centralizes a $2$-space on $M(E_6)$;
\item\label{propi:f43d} up to field automorphism of $H$, the action of $H$ on $M(F_4)$ is
\[ 9\oplus (4\otimes 3_2)\oplus 4.\]
\end{enumerate}
If (\ref{propi:f43a}), (\ref{propi:f43b}) or (\ref{propi:f43c}) hold, then $H$ is strongly imprimitive.\end{proposition}
\begin{proof} First, note that if $H$ stabilizes a line on $M(F_4)$ or $L(F_4)$ then $H$ is strongly imprimitive by Propositions \ref{prop:fixlineonMG} and \ref{prop:fixlineonLG} respectively. If $H$ stabilizes a unique $3$-space on $M(F_4)$ then $H$ is strongly imprimitive by Proposition \ref{prop:containedf4e6}, and if $H$ centralizes a $2$-space on $M(E_6)$ then $H$ is strongly imprimitive by Proposition \ref{prop:twotrivials}. Thus the consequence follows if we can show (\ref{propi:f43a}) to (\ref{propi:f43d}) hold.

\medskip

Using the traces of semisimple elements of orders $2$, $4$ and $5$, one finds, up to field automorphism of $H$, eight conspicuous sets of composition factors, namely
\[ 3_1^6,1^7,\qquad 4^3,3_1^3,1^4,\qquad 4^4,3_1,3_2,1^3,\qquad 9,4^2,3_1,3_2,1^2\]
\[ 9,4^3,3_1,1,\qquad 4,3_1^7,\qquad 4,3_1^6,3_2,\qquad 9^2,4,3_1.\]

\noindent \textbf{Case 1}: Using Lemma \ref{lem:char3sl2cohomology} we can compute the pressure of each set of composition factors, and the first set of factors has pressure $-7$, so in this case $H$ fixes a line on $M(F_4)$ by Lemma \ref{lem:pressure}, so the result holds.

\medskip

\noindent \textbf{Case 2}: The second conspicuous set of composition factors must yield a trivial submodule on $M(F_4)$ by Lemma \ref{lem:sl29hidetrivials}, so again the result holds.

\medskip

\noindent \textbf{Cases 6, 7}: If $M(F_4)\downarrow_H$ is either the sixth or seventh cases, then the trace of an involution in $H$ is $-7$ on $M(F_4)$. The trace of an involution from this class on $L(F_4)$ is $20$ so that $L=\PSL_2(3)$ acts with composition factors $3^8,1^{28}$ on $L(F_4)$. By Lemma \ref{lem:sl23restriction} any non-trivial simple module for $H$ has at most two $3$s for every three $1$s on restriction to $L$, and so $H$ always has at least sixteen trivial composition factors and at most eight non-trivial composition factors. Since $H^1(H,4)$ has dimension $2$ by Lemma \ref{lem:char3sl2cohomology}, $L(F_4)\downarrow_H$ has non-positive pressure and hence has a trivial submodule by Lemma \ref{lem:pressure}. Again, the result holds for these cases.

\medskip

\noindent \textbf{Case 8}: If the factors of $M(F_4)\downarrow_H$ are the eighth case then $M(F_4)\downarrow_H$ must be semisimple and $H$ be stabilize a unique $3$-space on $M(F_4)$, so again the result holds.

\medskip

\noindent \textbf{Case 4}: If $H$ stabilizes a line or a unique $3$-space on $M(F_4)$ then the result holds, so assume that this is not the case. The $9$ must split off as it is the projective Steinberg module. The socle is either $4$ or $3_1\oplus 3_2\oplus 4$, and so it is either $4$, or $3_1$ and $3_2$ are summands of $M(F_4)\downarrow_H$. The structure of $M(F_4)\downarrow_H$ is therefore either
\[ 9\oplus (4/1,1,3_1,3_2/4)\qquad\text{or}\qquad 9\oplus 3_1\oplus 3_2\oplus (4/1,1/4).\]
(In the first possibility, the second summand is not unique determined by the socle structure.)

We claim that, in either case, $M(F_4)\downarrow_H$ has zero $1$-cohomology. If this is true then $M(E_6)\downarrow_H$ must be $1\oplus 1\oplus M(F_4)\downarrow_H$, and therefore $H$ centralizes a $2$-space on $M(E_6)$.

If we have the second possibility, then the module is uniquely determined up to isomorphism, and so this is an easy computer calculation. As we have said, the first module is not determined uniquely up to isomorphism, but the quotient $M_1$ by the socle has structure
\[ 4/1,1,3_1,3_2\]
and is determined uniquely. An easy computer calculation shows that $H^1(H,M_1)=0$. If $M(F_4)\downarrow_H$ has non-trivial $1$-cohomology, then the extension must have the form
\[ 9\oplus (4/1,1,1,3_1,3_2/4),\]
i.e., the $1$ must fall into the second socle layer. But now there is a module $1,1,1/4$, and $H^1(H,4)$ has dimension $2$ by Lemma \ref{lem:char3sl2cohomology}. This proves the claim.

Thus in the second possibility, $H$ centralizes a $2$-space on $M(E_6)$, as the element in $H^1(\bG,M(F_4))$ that yields the module $1/M(F_4)\downarrow_H$ must be zero on restriction to $H$. In particular, the result holds.

\medskip

\noindent \textbf{Case 3}: If the proposition does not hold, then we may assume that $H$ does not stabilize a line or unique $3$-dimensional subspace of $M(F_4)\downarrow_H$. We cannot have $3_1\oplus 3_2$ in the socle, for then $M(F_4)\downarrow_H$ has $3_1\oplus 3_2$ as a summand and the remaining summand $M_1$ has composition factors $4^4,1^3$. The $\{1,4\}$-radical of $P(4)$ is the self-dual module $4/1,1/4$, so $M_1$ must have socle $4^{\oplus 2}$. Furthermore, $M_1$ must have structure
\[ 4,4/1,1,1/4,4,\]
but we see that the $\{1,4\}$-radical of $P(4)^{\oplus 2}$ has four trivial modules in the second socle layer, not three, but still exactly two copies of $4$ in the third layer. Since none of the trivial factors are quotients, one cannot be removed while keeping the two $4$s in the top, so $M_1$ must have a trivial submodule.

Hence $\soc(M(F_4)\downarrow_H)$ is either $4$ or $4^{\oplus 2}$. (If it is $4^{\oplus 3}$, there must be a $4$ as a summand, and we can remove this to still obtain a module with socle $4$ or $4^{\oplus 2}$.) We saw the structure of $P(4)$ just before Lemma \ref{lem:sl29hidetrivials}:
\[ 4/1,1,3_1,3_2/4,4,4/1,1,3_1,3_2/4.\]
If $M(F_4)\downarrow_H$ has five socle layers then it contains $P(4)$, which has dimension $36$, too many dimensions. If it has four socle layers then $1$ or $3_i$ is a quotient, hence a submodule, which is not allowed. Thus it has three socle layers, so has structure
\[ 4,4/1,1,1,3_1,3_2/4,4.\]

We now look at the image of $H$ inside $E_6$, and its action on $M(E_6)$. If $\soc(M(E_6)\downarrow_H)=1$ then, since $P(1)$ has dimension $27$, we have that $M(E_6)\downarrow_H\cong P(1)$. However, the action of $u$ on $M(E_6)$ is clearly now $3^9$, so $u$ acts on $M(F_4)$ with Jordan blocks $3^8,1$, from Table \ref{t:unipotentF4}. But if we remove the top and socle from $P(1)$ we obtain a $25$-dimensional module on which $u$ acts with Jordan blocks $3^7,2^2$, a contradiction.

Thus there exists an $H$-submodule $1\oplus 4$ of the minimal module $M(E_6)$. Notice that $P(4)$ has dimension $36$ and has five socle layers, and $P(1)$ has five socle layers, so since neither of these is contained in the module $M(E_6)\downarrow_H$, we must have that $M(E_6)\downarrow_H$ has at most four socle layers. In particular, since $M(E_6)$ is self-dual, we cannot have a uniserial module $3_i/4/1$ as a subquotient of $M(E_6)\downarrow_H$. To see this, note that $1/4/3_i$ would also have to be a subquotient. As there is a unique $3_i$ in $M(E_6)\downarrow_H$, we would need at least five socle layers in $M(E_6)\downarrow_H$, which is not allowed.

Consider the preimage $W$ of $\soc^2(M(F_4)\downarrow_H)$ in $M(E_6)$, and in particular the $\{1,4\}$-radical of $W$. This is the preimage of a module $(1,1/4)\oplus (1/4)$, and since $\Ext^1(1,1/4,1)=0$, the $\{1,4\}$-radical of $W$ must be a module 
\[(1,1/4)\oplus (1/4/1).\]
(The uniserial module $1/4/1$ is not uniquely determined up to isomorphism.) We need to place both a $3_1$ and a $3_2$ on top of this module, but without constructing a uniserial $3_i/4/1$ as a subquotient. There is only one way to do this:
\[ (1,1,3_1,3_2/4)\oplus (1/4/1).\]
(If $3_i$ were placed diagonally across the two summands, quotienting out by the first summand would yield a uniserial $3_1/4/1$ submodule.)

Since there is no uniserial module $4/1/4$ by Lemma \ref{lem:no414}, no $4$ placed on top of this module $W$ can cover the $1$ in the second summand (again by quotienting out by the first summand we would construct a uniserial module $4/1/4$), and so $M(F_4)\downarrow_H$ has a trivial quotient. This was specifically excluded at the start of the proof, so we obtain a contradiction to the statement that $H^1(H,M(F_4))\neq 0$.

We therefore see that $H$ cannot have non-trivial $1$-cohomology, so $H$ centralizes a $2$-space on $M(E_6)$, as needed for the proposition.

\medskip

\noindent \textbf{Case 5}: The composition factors are $9,4^3,3_1,1$. We may assume that neither $1$ nor $3_1$ lie in the socle of $M(F_4)\downarrow_H$, as else we satisfy the conclusion of the proposition.

The $9$, being projective, splits off as a summand. As in the previous case, we know that $M(F_4)\downarrow_H$ has exactly three socle layers. Since there are three $4$s, two must appear in the socle, and therefore two in the top. Thus one splits off as a summand. 

We therefore see that $M(F_4)\downarrow_H$ must have the form
\[ 9\oplus (4/1,3_1/4)\oplus 4.\]
There is a unique module $4/1,3_1/4$, and it is $4\otimes 3_2$. To see that it is unique, note that above $1,1,3_1/4$ one may place two copies of $4$, which can be understood as the module
\[ (4/1,1/4)\oplus (4/1,3_1/4)\]
quotiented out by a diagonal $4$, where the second summand is $4\otimes 3_2$.

Any submodule of this of codimension $4$ either removes the $4$ from the first or second summand, or diagonally across both. Removing from either the first or second summand is easy to understand. Removing a diagonal $4$ must result in neither trivial factor in the first summand becoming a quotient. Hence the diagonal submodule has a single trivial quotient, so upon removing that it has the form
\[ 4/1,1,3_1/4,\]
with no trivial quotients. This has the wrong form, so we cannot make a module with structure $4/1,3_1/4$ by removing a diagonal quotient.

Thus $4/1,3_1/4$ is the module $4\otimes 3_2$, and this is case (\ref{propi:f43d}) in the proposition.
\end{proof}

Before we attack this last case above, we prove a small lemma about it, reducing us to the case where $H$ is Lie primitive.

\begin{lemma}\label{lem:f4p=3a=2Lieprim} Suppose that $p=3$ and $a=2$. Either $H$ is strongly imprimitive, or $H$ is Lie primitive in $\bG$.
\end{lemma}
\begin{proof} By Proposition \ref{prop:f4char3a=2}, either $H$ is strongly imprimitive or $M(F_4)\downarrow_H$ is
\[ 9\oplus (4\otimes 3_2)\oplus 4.\]
If $H$ is Lie primitive then we are done, so assume that $H$ is contained in a member $\bX$ of $\mathscr X$. We will show that $\bX$ always stabilizes the $9$, and hence $H$ is strongly imprimitive by Proposition \ref{prop:intersectionstabilizers}.

We first consider $\bX$ a parabolic subgroup. The composition factors of $M(F_4)\downarrow_{\bX}$ are given in Appendix \ref{app:actions}. As the $B_3$-parabolic has three trivial factors on $M(F_4)$, the $C_3$-parabolic has two factors $M(C_3)$, which are not compatible with the dimensions, and the two $A_2A_1$-parabolics have composition factors either of dimension $2$ or two of dimension $3$, $H$ cannot lie in any parabolic. Thus $H$ is $\bG$-irreducible, and so $\bX$ is a maximal-rank subgroup or $A_1G_2$.

Suppose that $\bX=B_4$. As $\bX$ has summands of dimension $9$ and $16$, the summand $9$ of $M(F_4)\downarrow_H$ is stabilized by $\bX$, as claimed. (In fact, this cannot work: $H$ must act irreducibly on $M(B_4)$. In this case, the composition factors of $H$ on the spin module are $4^2,3_1,3_2,1$, which are not correct. But we do not need this to prove the result.)

The summands of $A_2A_2$ on $M(F_4)$ are of dimension $9$, $9$ and $7$, so $H$ cannot embed in this subgroup.

The composition factors of $\bX=A_1G_2$ on $M(F_4)$ are of dimension $21$ and $4$. The $21$-dimensional module is $(L(A_1),M(G_2))$, and this must restrict to $H$ as $(3_1\otimes 3_2)\oplus (4\otimes 3_2)$, up to field automorphism. We see that $H$ must act on $M(G_2)$ as $4\oplus 3_2$, and on $M(A_1)$ as $2_1$ (so that $H$ acts on $L(A_1)$ as $3_1$).

Let $\bY$ denote a diagonal $A_1$ subgroup of $A_1G_2$ acting on $M(A_1)$ as $L(1)$ and on $M(G_2)$ as $L(4)\oplus L(6)$. The action of $\bY$ on $M(F_4)$ has a summand $L(2)\otimes L(6)=L(8)$, which means that $\bY$ stabilizes the irreducible $9$-space on $M(F_4)$, as needed.

Finally, if $H$ is contained in $\bX=C_3A_1$, then $H$ acts on the $13$-dimensional module $L(010)$ for the $C_3$ factor, which is a summand of $M(F_4)\downarrow_{\bX}$, as $9\oplus 4$. Since $C_3$ is a classical group and $H$ is $\bG$-irreducible, hence $C_3$-irreducible, $H$ is contained in an $A_1$-subgroup $\bY$ of $C_3$ (see Proposition \ref{prop:morphismextensiontrueforclassical}). There are no $13$-dimensional simple modules for $\bY$, and $H$ acts as $4\oplus 9$, so $\bY$ must also act as a sum of two modules, of dimensions $4$ and $9$. In particular, $\bY$ stabilizes the irreducible $9$-space of $M(F_4)\downarrow_H$, as claimed.
\end{proof}

\begin{proposition}\label{prop:f4p3a2end} Suppose that $p=3$ and $a=2$. The subgroup $H$ is always strongly imprimitive in $\bG$.
\end{proposition}
\begin{proof} From Proposition \ref{prop:f4char3a=2} we may assume that $M(F_4)\downarrow_H$ is the module
\[ 9\oplus (4\otimes 3_1)\oplus 4,\]
up to field automorphism. From Lemma \ref{lem:f4p=3a=2Lieprim} we may assume that $H$ is Lie primitive in $\bG$.

Consider the $F_4A_1$ maximal subgroup $\bX$ of $E_7$. We may construct a subgroup $J\cong \SL_2(9)$ in $\bX$ by projecting $J$ along $F_4$ as the subgroup $H$, and along $A_1$ irreducibly, acting on $M(A_1)$ as $2_1$. Since $H$ is Lie primitive in $\bG$, the only positive-dimensional subgroups of $\bX$ containing $J$ are contained in $H\cdot A_1$, and have dimension at most $3$.

\medskip

\noindent \textbf{Step 1}: $\bX$ stabilizes a unique $2$-space $V$ on $M(E_7)$.

\noindent From \cite[Table 10.2]{liebeckseitz2004}, the action of $\bX$ on $M(E_7)$ has structure
\[ (0000,1)/(0001,1),(0000,3)/(0000,1),\]
and so the claim is clear.

\medskip

\noindent \textbf{Step 2}: $J$ stabilizes more than one $2$-space on $M(E_7)$. It acts as $2_2$ on $V$ and $2_1$ on some other subspace $W$.

\noindent This is a simple calculation. The action of $J$ on $(M(F_4),M(A_1))$ is
\[ (9\oplus (4\otimes 3_1)\oplus 4)\otimes 2_1=P(6_{1,2})\oplus P(6_{2,1})\oplus (2_2/2_1/2_2)\oplus 6_{2,1}\oplus 2_2.\]  
Thus the action of $J$ on the second socle layer of $M(E_7)\downarrow_{\bX}$ has three copies of $2_2$ as submodules. Since $\Ext^1(2_2,2_1)\cong\Ext^1(4,1)$ has dimension $2$ by Lemma \ref{lem:char3sl2cohomology}, there is no module $2_2,2_2,2_2/2_1$, this means that $M(E_7)\downarrow_J$ has a submodule $2_2$, as well as a submodule $2_1$.

\medskip

The subspace $W$ is stabilized by $J$ but not by $\bX$. Since $\dim(E_7)=133$, the dimension of the stabilizer of $W$ in $E_7$ is at least $133-56-55=22$. Thus $J$ is contained in a positive-dimensional subgroup $\bY_0$ of $E_7$, not equal to $\bX$ by construction, and of dimension at least $22$. It is not contained in $\bX$ either, since the largest dimension of a positive-dimensional subgroup of $\bX$ containing $J$ is $3$. Thus $\bY_0\not\leq \bX$.

Hence we replace $\bY_0$ by some maximal positive-dimensional subgroup $\bY$ of $E_7$, which has dimension at least $22$ and is not equal to $\bX$. 

\medskip

\noindent \textbf{Step 3}: $J$ is $E_7$-irreducible, and contained in either $\bY=D_6A_1$ or $\bY=F_4A_1$.

\noindent Suppose that $\bY$ is some other maximal positive-dimensional subgroup. If $J$ is not contained in $\bY^0$ then $J$ must be a subgroup of $\bY/\bY^0$. \cite[Table 10.3]{liebeckseitz2004} shows that if $\bY$ has maximal rank then $\bY$ is the normalizer of a torus, but this has dimension $7$, which is too small. \cite[Theorem 1 and Remark (2)]{liebeckseitz2004} shows that there are no other candidates. Thus $J\leq \bY^0$ so we may assume that $\bY$ is a maximal \emph{connected} positive-dimensional subgroup of $E_7$.

Thus $\bY$ is a maximal parabolic, a maximal-rank subgroup $A_7$, $D_6A_1$ or $A_2A_5$ from \cite[Table 10.4]{liebeckseitz2004}, or a subgroup $G_2C_3$ or $F_4A_1$ ($A_1G_2$ has dimension $17$, which is too small).

\medskip

If $\dim(\bY)$ has dimension at least $82$ then $\dim(\bY\cap \bX)\geq 82+55-133=4$. Since $J\leq \bY\cap\bX$ this contradicts the fact that $J$ is contained in no subgroup of $\bX$ of dimension greater than $3$.

We remind the reader of the dimensions of the maximal parabolic subgroups of $E_7$.

\begin{center}\begin{tabular}{cc}
\hline Parabolic & Dimension
\\ \hline $E_6$ & $106$
\\ $D_6$ & $100$
\\ $A_6$ & $90$
\\ $A_5A_1$ & $86$
\\ $D_5A_1$ & $91$
\\ $A_4A_2$ & $83$
\\ $A_3A_2A_1$ & $80$
\\ \hline
\end{tabular}\end{center}

Thus we need only consider the $A_3A_2A_1$-parabolic subgroup. Notice that $Z(J)=Z(E_7)$. In order for $Z(J)=Z(E_7)$, $J$ must project faithfully along the $A_1$ and the $A_3$ factors of the parabolic, but there are no faithful $4$-dimensional simple modules for $J$ in characteristic $3$. Hence $J$ acts reducibly on $M(A_4)$, and thus lies in a $A_2A_1A_1A_1$-parabolic subgroup. This lies in another parabolic subgroup, so we obtain a contradiction.

Thus we now have that $J$ is $E_7$-irreducible, establishing one of our claims.

Suppose that $H$ lies in one of $A_7$, $A_2A_5$ and $G_2C_3$. Then Proposition \ref{prop:morphismextensiontrueforclassical} implies that $H$ is contained in an $E_7$-irreducible $A_1$ subgroup. Consulting \cite[Table 7]{thomas2016}, we find that all $E_7$-irreducible $A_1$ subgroups lie in $D_6A_1$, so $\bY$ may be replaced by $D_6A_1$ in these cases.

Thus $H$ lies inside $D_6A_1$, or $F_4A_1$, as needed.

\medskip

\noindent \textbf{Step 4}: $H$ is contained in $\bY=D_6A_1$.

\noindent Suppose that $\bY$ is of type $F_4A_1$, with $\bY\neq \bX$. By making another choice for $W$, we may assume that $\bY$ stabilizes $W$. Since $\bY$ acts on $W$ as $(0,1)$, this means that $J$ embeds in $\bY$ acting on $M(A_1)$ as $2_2$. In order for $Z(J)$ and $Z(E_7)$ to coincide, $Z(J)$ must lie in the kernel of the projection along $F_4$.

The possible composition factors of $J$ on $M(F_4)$ are given in Proposition \ref{prop:f4char3a=2}, and are up to field automorphism
\[ 3_1^6,1^7,\quad 4^3,3_1^3,1^4,\quad 4^4,3_1,3_2,1^3,\quad 9,4^2,3_1,3_2,1^2,\]
\[ 9,4^3,3_1,1,\quad 4,3_1^7,\quad 4,3_1^6,3_2,\quad 9^2,4,3_1.\]
We tensor these factors (and their images under the field automorphism) by $2_2$, then add $2_2^2,2_1$ to them to give the factors of $M(E_7)\downarrow_J$, which we computed above to be
\[ 6_{1,2}^3,6_{2,1}^3,2_1^4,2_2^6.\]
This yields a unique set of composition factors for $M(F_4)\downarrow_J$, which are $4^3,3_1^3,1^4$. Thus $J$ fixes a line on $M(F_4)$ by the proof of Proposition \ref{prop:f4char3a=2}. Hence the projection of $J$ along the $F_4$ factor of $\bY$ lies in a positive-dimensional subgroup of $F_4$ of dimension at least $52-26=26$, but not a parabolic since $J$ is $\bY$-irreducible. Thus this projection lies inside $B_4$. Thus $J$ lies inside a $B_4A_1$ subgroup of $\bY$, and we may apply Proposition \ref{prop:morphismextensiontrueforclassical} to see that $J$ is contained in an $E_7$-irreducible $A_1$ subgroup. As before, this means that $J\leq D_6A_1$, as claimed.

\noindent \textbf{Step 5}: Conclusion.

\noindent Thus we may assume that $\bY=D_6A_1$. If $J\leq D_6$ then it is not $E_7$-irreducible, so $J$ acts faithfully on $M(A_1)$. This means that $Z(J)$ acts trivially on $M(D_6)$.

Since $J$ is $\bY$-irreducible, by Proposition \ref{prop:morphismextensiontrueforclassical}, the projection of $J$ along $D_6$ is contained in a $D_6$-irreducible $A_1$ subgroup. We see from \cite[Lemma 3.4]{thomas2016} that the action of $J$ on $M(D_6)$ must be a sum of inequivalent simple modules. This is only possible if $J$ acts on $M(D_6)$ as $3_i\oplus 9$ for some $i=1,2$, and of course $J$ acts along $M(A_1)$ as $2_j$ for some $j=1,2$.

We consult \cite[Table 7]{thomas2016}, and see that the $A_1$ subgroup of $E_7$ must be subgroup 10, with $r=j$, $s=1$, $t=2$, $u=i$.

Working up to field automorphism, we may fix $j=1$, and let $u=i=1,2$. This yields two sets of composition factors for $M(E_7)\downarrow_J$, easily computable from \cite[Table 12]{thomas2016}. They are
\[6_{1,2}^3,6_{2,1}^2,2_1^7,2_2^6,\qquad 6_{1,2}^5,6_{2,1},2_1^4,2_2^6.\]
On the other hand, we see from the composition factors of $M(E_7)\downarrow_J$ above that they must be $6_{1,2}^3,6_{2,1}^3,2_1^4,2_2^6$. Up to field automorphism, these do not coincide, and this yields a contradiction.

The only assumption was that the positive-dimensional subgroups of $\bX$ containing $J$ had dimension at most $3$, so this is not the case. In particular, they have positive-dimensional projection along $F_4$, so the projection $H$ of $J$ along $F_4$ is contained in a positive-dimensional subgroup of $F_4$, as desired.
\end{proof}

\section{Characteristic At Least 5}
\label{sec:f4atleast5}

Let $p\geq 5$, and recall that $H=\PSL_2(p^a)$ for some $a\geq 1$, with $p^a\leq 36$, with $u\in H$ of order $p$, as detailed at the start of this chapter. The possible actions of $u$ on $M(F_4)$ are given in \cite[Table 3]{lawther1995}; by Lemma \ref{lem:genericunipotent} we may assume that our unipotent class is not generic. This leaves us with the following three unipotent classes:
\begin{enumerate}
\item $C_3$, $p=7$, acting as $7^2,6^2$;
\item $F_4(a_2)$, $p=7$, acting as $7^3,5$;
\item $F_4$, $p=13$, acting as $13^2$.
\end{enumerate}

This proves the following result immediately.

\begin{proposition}\label{prop:f4char>3firstreduction} If $p^a\neq 7,13$ then $H$ is a blueprint for $M(F_4)$.
\end{proposition}

For $p=7$ we have the following result.

\begin{proposition}\label{prop:f4char7} If $p^a=7$ then $H$ stabilizes a line on either $M(F_4)$ or $L(F_4)$.
\end{proposition}
\begin{proof}
We use the traces of elements of orders $2$, $3$ and $4$ to produce the possible composition factors of $M(F_4)\downarrow_H$. These are
\[ 3^6,1^8,\qquad 5,3^7,\qquad 5^3,3^3,1^2,\qquad 7,5^3,3,1,\qquad 7^3,1^5.\]
We saw in Section \ref{sec:sl2p} that the only indecomposable module with a trivial composition factor but no trivial submodule or quotient is $P(3)=5/1,3/5$. This immediately tells us that in the first, third and fifth cases, $H$ stabilizes a line on $M(F_4)$. (Indeed, in the first and fifth cases all trivial factors are summands.)

The case $7,5^3,3,1$ yields traces of elements of orders $2$, $3$ and $4$ of $2$, $-1$ and $-2$ respectively. This corresponds to traces on $L(F_4)$ of $-4$ for an involution, $-2$ or $7$ for an element of order $3$, and finally $4$ for an element of order $4$. There is no set of composition factors that are compatible with this, so this case cannot occur.

If the composition factors are $5,3^7$, then the traces of the elements of orders $2$ and $3$ yield a unique conspicuous set of composition factors on $L(F_4)$, which is $5^7,3,1^{14}$. This has negative pressure and so $H$ stabilizes a line on $L(F_4)$ by Lemma \ref{lem:pressure}.
\end{proof}

For $p=13$ we are left with one open possibility, which we will prove yields a Serre embedding (see Definition \ref{defn:serreembedding}).

\begin{proposition}\label{prop:f4char13}
Suppose that $p^a=13$. Either $H$ is a blueprint for $M(F_4)$, or $u$ is a regular unipotent element and $M(F_4)\downarrow_H$ and $L(F_4)\downarrow_H$ are given by
\[ P(9)=9/3,5/9 \qquad \text{and}\qquad P(3)\oplus P(11)=(3/9,11/3)\oplus (11/1,3/11)\]
respectively. In particular, $H$ is a Serre embedding.
\end{proposition}
\begin{proof} From the list above, the regular unipotent class is the only non-generic one for $p=13$, so if $H$ is not a blueprint for $M(F_4)$ then $u$ is regular and in particular acts projectively on $M(F_4)$ and $L(F_4)$ by \cite[Tables 3 and 4]{lawther1995}, hence both modules must restrict to $H$ as projectives. The projective indecomposable modules for $H$ are
\[ 1/11/1,\quad 3/9,11/3,\quad 5/7,9/5,\quad 7/5,7/7,\quad 9/3,5/9,\quad 11/1,3/11,\quad 13.\]
Thus there are eight possible projective modules of dimension $26$, two of which are conspicuous for $M(F_4)$: $P(5)$ and $P(9)$. The first of these does not have corresponding factors on $L(F_4)$ (see (\ref{li:strata}) from Chapter \ref{ch:strategy}), and the second has factors $11^3,9,3^3,1$, which yield the projective module $P(3)\oplus P(11)$, as claimed.
\end{proof}

We now summarize the results of this chapter.

If $p=2$ we showed that $H$ always stabilizes a line or $2$-space on $M(F_4)=L(\lambda_4)$ or on $M(F_4)^\tau=L(\lambda_1)$, and hence is strongly imprimitive by Propositions \ref{prop:fixlineonMG} and \ref{prop:fix2spaceonMG}. If $p=3$ then $H$ is strongly imprimitive by Proposition \ref{prop:f4p3a2end} for $a=2$, and for $a=3$ we showed in Proposition \ref{prop:f4char3a=3} that $H$ is a blueprint for $M(F_4)$, hence strongly imprimitive by Proposition \ref{prop:blueprintissi}.

If $p\geq 5$ then we showed that $H$ is a blueprint for $M(F_4)$ unless $p=7$ -- in which case $H$ stabilizes a line on either $M(F_4)$ or $L(F_4)$ -- or $p=13$, in which case $H$ is a Serre embedding. By Propositions \ref{prop:fixlineonMG}, \ref{prop:fixlineonLG} and \ref{prop:blueprintissi}, $H$ is strongly imprimitive or a Serre embedding.

We now complete the proof of Theorem \ref{thm:f4}. Any maximal subgroup of the finite group that is almost simple with socle a copy of $\PSL_2(p^a)$ is either strongly imprimitive or it is not. If it is not, then it is a Serre embedding by the above reasoning, so we assume that it is strongly imprimitive. Therefore it is the fixed points of a maximal subgroup from \cite[Corollary 5]{liebeckseitz2004}. It cannot be a parabolic, and no $A_1$-type subgroup appears in \cite[Table 5.1]{liebecksaxlseitz1992}, so it cannot be the fixed points of a maximal-rank subgroup. This leaves only the $A_1$ for $p\geq 13$, as claimed.

\chapter{The Proof for \texorpdfstring{$E_6$}{E6}}
\label{ch:e6}
In this chapter, $k$ is an algebraically closed field of characteristic $p\geq 2$ and $\bG=E_6(k)$, by which we mean the simply connected form, i.e., $|Z(\bG)|=3$ if $p\neq 3$ and $\bG'=\bG$.  Let $H\cong \PSL_2(p^a)$ be a subgroup of $\bG$.

Theorem \ref{thm:goodblueprint} states that if $p=2$ and $a\geq 5$ then $H$ is a blueprint for $M(E_6)\oplus M(E_6)^*$, and if $p$ is odd and $p^a\geq 37$ then $H$ is a blueprint for $M(E_6)\oplus M(E_6)^*$. In both cases, $H$ is strongly imprimitive by Proposition \ref{prop:blueprintissi}. Thus in what follows we may assume that $p^a\leq 31$.

Let $L=\PSL_2(p)\leq H$ and let $u$ denote a unipotent element of $L$ of order $p$. The possibilities for the Jordan block structures of $u$ on $M(E_6)$ and $L(E_6)$ are given in \cite[Tables 5 and 6]{lawther1995}, and for $L(E_6)^\circ$ in Lemma \ref{lem:unipotentE6}. Recall the definition of a generic unipotent element from Definition \ref{defn:genericunipotent}.

We will prove that $H$ is always strongly imprimitive. Since $H$ cannot act irreducibly on $M(E_6)$ or $L(E_6)^\circ$ by \cite{liebeckseitz2004a}, if $H$ is a blueprint for either module then it is strongly imprimitive, by Proposition \ref{prop:intersectionstabilizers}.

\section{Characteristic 2}

Let $p=2$. This case is easy, as we do not aim to produce the same depth of result. If $a=1$ then $H$ is soluble, and if $a=2$ then \cite[Proposition 5.4]{craven2015un2} shows that $H$ stabilizes a $1$- or $2$-space on $M(E_6)$ or its dual. We show that $\SL_2(8)$ stabilizes a line or $2$-space on $M(E_6)$ or its dual, and we show that $\SL_2(16)$ is a blueprint for $M(E_6)\oplus M(E_6)^*$.

For $a=4$ we use the idea that, while not every semisimple element of order $17$ in $F_4$ is a blueprint for $M(F_4)$, since $17$ is very close to the maximum of the integers in $\mathcal X_{M(F_4)}$ given in Proposition \ref{prop:numbersforF4}, almost all elements of order $17$ are (see (\ref{li:stratf}) from Chapter \ref{ch:strategy}). This statement passes through to $M(E_6)$ and $M(E_6)\oplus M(E_6)^*$, since our real semisimple elements lie in $F_4$, via Lemma \ref{lem:f4blueprintisblueprint}.

We start with $a=3$.

\begin{proposition}\label{prop:e6char2a=3} Suppose that $p=2$ and $a=3$. Then $H$ stabilizes a line or $2$-space on $M(E_6)$ or $M(E_6)^*$.
\end{proposition}
\begin{proof} Suppose that $\soc(M(E_6)\downarrow_H)$ and $\soc(M(E_6)^*\downarrow_H)$ have neither $1$s nor $2$s, so $M(E_6)\downarrow_H$ is a submodule of $P(4)$s and $8$s. The projective cover of $4_{i,i+1}$ is
\[ 4_{i,i+1}/2_{i+1}/1/2_{i-1}/1/2_{i+1}/4_{i,i+1},\]
and thus $M(E_6)\downarrow_H$ is a sum of projectives $P(4_{i,i+1})$ and $8$s, but this has even dimension, which is not correct.
\end{proof}

Now we move on to $a=4$, where we use semisimple elements of order $17$ that are blueprints for $M(E_6)$, as suggested earlier.

\begin{proposition}\label{prop:e6char2a=4} If $p=2$ and $a=4$ then $H$ is a blueprint for $M(E_6)\oplus M(E_6)^*$.
\end{proposition}
\begin{proof} Every real semisimple element of $E_6$ is (conjugate to an element) in $F_4$ by Proposition \ref{prop:containedf4e6}, and by Lemma \ref{lem:f4blueprintisblueprint} if an element of $F_4$ is a blueprint for $M(F_4)$ then it is for $M(E_6)\oplus M(E_6)^*$ as well.

Of the $230$ semisimple classes in $F_4$ of elements of order $17$, all but two are blueprints for $M(F_4)$ (see (\ref{li:stratf}) from Chapter \ref{ch:strategy}). These two classes have representatives $x$ and $x^3$, where $x$ has eigenvalues
\[ 1^3,(\zeta_{17}^{\pm 1})^2,(\zeta_{17}^{\pm 2})^2,(\zeta_{17}^{\pm 3}),(\zeta_{17}^{\pm 4})^2,(\zeta_{17}^{\pm 5}),(\zeta_{17}^{\pm 6}),(\zeta_{17}^{\pm 7}),(\zeta_{17}^{\pm 8})^2\]
on $M(E_6)$. We thus may assume that every element of $H$ of order $17$ is conjugate to either $x$ or $x^3$.

There are $107766$ possible sets of composition factors for a $kH$-module of dimension $27$, but for none of them does an element of order $17$ have the eigenvalues above (up to algebraic conjugacy). Thus a semisimple element of $H$ of order $17$ is always a blueprint for $M(E_6)\oplus M(E_6)^*$, and so the result holds.
\end{proof}

\section{Characteristic 3}

Let $p=3$. From the remarks at the start of the chapter, we may assume that $a\leq 3$. Of course, $\PSL_2(3)$ is soluble, and $\PSL_2(3^2)$ was proved to stabilize a line on either $M(E_6)$ or $L(E_6)^\circ$ in \cite[Proposition 6.2]{craven2015un2}. Thus we may assume that $a=3$.

In the previous section we exploited the fact that most semisimple elements of order $17$ are blueprints for $M(E_6)$. We will do the same here with order $13$ elements. There are $104$ classes of semisimple elements of order $13$ in $F_4$. All but seven of these are blueprints for $M(F_4)$. We can easily see this by finding elements of order $26$ that have the same number of distinct eigenvalues on $M(F_4)$ as their square, and noting that elements of order $26$ are blueprints for $M(F_4)$ by Proposition \ref{prop:numbersforF4}. By Proposition \ref{prop:f4real} real semisimple elements of $E_6$ conjugate into $F_4$, and by Lemma \ref{lem:f4blueprintisblueprint} blueprints for $M(F_4)$ are blueprints for $M(E_6)\oplus M(E_6)^*$.

\begin{proposition}\label{prop:e6char3} Suppose that $p=3$ and $a=3$. Either $H$ is a blueprint for $M(E_6)\oplus M(E_6)^*$ or $H$ stabilizes a line on $M(E_6)$ or $M(E_6)^*$.
\end{proposition}
\begin{proof} This is easier than the case of $F_4$, but will start in exactly the same way. There are fifty conspicuous sets of composition factors for $M(E_6)\downarrow_H$, but for only seven of these do the elements of order $13$ come from semisimple classes that are not blueprints for $M(E_6)\oplus M(E_6)^*$, three up to field automorphism of $H$. These are
\[12_{2,3,1},9_{1,2},4_{1,2},1^2,\qquad 12_{2,3,1},9_{1,3},4_{1,2},1^2,\qquad 4_{1,2}^2,4_{1,3}^2,4_{2,3}^2,1^3.\]
The first two have pressure $-1$ and so $H$ must stabilize a line on $M(E_6)$ by Lemma \ref{lem:pressure}. The third must stabilize a line on $M(E_6)$ or $M(E_6)^*$ by Lemma \ref{lem:no414} with $i=3$, $\alpha=0$. Thus $H$ stabilizes a line on $M(E_6)$ or $M(E_6)^*$, as needed.
\end{proof}

\section{Characteristic At Least 5}
\label{sec:e6p>3}

Let $p\geq 5$, and recall that $H=\PSL_2(p^a)$ for some $a\geq 1$, with $p^a\leq 36$, with $u\in L\leq H$ of order $p$, where $L=\PSL_2(p)$. The possible actions of $u$ on $M(E_6)$ are given in \cite[Table 5]{lawther1995}; by Lemma \ref{lem:genericunipotent} we may assume that our unipotent class is not generic, leaving us with the following seven unipotent classes:
\begin{enumerate}
\item $A_4$, $p=5$, acting as $5^5,1^2$;
\item $A_4+A_1$, $p=5$, acting as $5^5,2$;
\item $A_5$, $p=7$, acting as $7^2,6^2,1$;
\item $D_5(a_1)$, $p=7$, acting as $7^3,3,2,1$;
\item $E_6(a_3)$, $p=7$, acting as $7^3,5,1$;
\item $E_6(a_1)$, $p=11$, acting as $11^2,5$;
\item $E_6$, $p=13$, acting as $13^2,1$.
\end{enumerate}

We now go prime by prime, starting with $p=5$.

\begin{proposition}\label{prop:e6char5} Suppose that $p=5$.
\begin{enumerate}
\item If $a=1$ then $H$ stabilizes a line on either $M(E_6)$ or $L(E_6)$. 
\item If $a=2$ then $H$ stabilizes a line on $M(E_6)$, $M(E_6)^*$ or $L(E_6)$.
\end{enumerate}
\end{proposition}
\begin{proof}
Suppose that $a=1$. The conspicuous sets of composition factors of $M(E_6)\downarrow_H$ are
\[ 3^6,1^9,\qquad 5,3^7,1,\qquad 5^3,3^3,1^3.\]
The first set of composition factors has pressure $-3$, so in this case $H$ stabilizes a line on $M(E_6)$ by Lemma \ref{lem:pressure}. In the second case we switch to $L(E_6)$ (see (\ref{li:strata}) from Chapter \ref{ch:strategy}), on which $H$ acts with composition factors
\[ 5^8,3^8,1^{14}\qquad\text{or}\qquad 5^{11},3^5,1^8.\]
In either case, we see that $H$ stabilizes a line on $L(E_6)$, as needed. The third set of composition factors has pressure $0$, so might only stabilize a hyperplane on $M(E_6)$. However, the only indecomposable modules with a trivial composition factor but no trivial submodule are submodules of $P(3)=3/1,3/3$, so in order not to stabilize a line, $M(E_6)\downarrow_H$ must be
\[ 5^{\oplus 3}\oplus (1/3)^{\oplus 3},\]
on which $u$ acts as $5^3,4^3$, but this does not appear on \cite[Table 5]{lawther1995}, so $H$ does indeed stabilize a line (and hyperplane) on $M(E_6)$.

\medskip

Now suppose that $a=2$. By Lemma \ref{lem:sl25restriction}, if $M(E_6)\downarrow_L$ has more trivial than $3$-dimensional factors then $H$ stabilizes a line on $M(E_6)$. Thus if $M(E_6)\downarrow_L$ is the first set of composition factors then $H$ stabilizes a line on $M(E_6)$, and if $M(E_6)\downarrow_L$ is the second set of composition factors then $H$ stabilizes a line on $L(E_6)$.

We therefore assume that $M(E_6)\downarrow_L$ has factors $5^3,3^3,1^3$. At this point it seems easiest to use the traces of semisimple elements of order at most $13$, finding eighteen conspicuous sets of composition factors, each with at least one trivial factor and with non-positive pressure, so in all cases $H$ stabilizes either a line or a hyperplane on $M(E_6)$.
\end{proof}

For $p=7$ we do not need to go past $a=1$, which makes this easier than the previous case.

\begin{proposition}\label{prop:e6char7} If $p=7$ and $a=1$, then $H$ stabilizes a line or hyperplane on $M(E_6)$.
\end{proposition}
\begin{proof} The conspicuous sets of composition factors are, as for $p=5$, the same as for $F_4$ but with an extra trivial factor, namely
\[ 3^6,1^9,\qquad 5,3^7,1,\qquad 5^3,3^3,1^3,\qquad 7,5^3,3,1^2,\qquad 7^3,1^6.\]
The only indecomposable module that has a trivial composition factor but no trivial submodule or quotient is $P(5)=5/1,3/5$, thus all of these sets of composition factors stabilize either a line or hyperplane on $M(E_6)$.
\end{proof}

For $p=11$, we see the first use of the idea of fixing an $\slf_2$-subalgebra (see Section \ref{sec:sl2subalgebra}).

\begin{proposition}\label{prop:e6char11} Suppose that $p=11$. Either $H$ is a blueprint for both $M(E_6)$ and $L(E_6)$, or $H$ has a trivial summand on $M(E_6)$, or $H$ acts on $M(E_6)$ and $L(E_6)$ as
\[ P(9)\oplus 5\qquad\text{and}\qquad 11^{\oplus 2}\oplus P(7)\oplus P(5)\oplus 9\oplus 3\]
respectively. Furthermore, $H$ stabilizes a unique $3$-space on $L(E_6)$, which is an $\slf_2$-subalgebra of $L(E_6)$.
\end{proposition}
\begin{proof} Examining \cite[Tables 5 and 6]{lawther1995}, we see that there are only two unipotent classes of elements of order $11$ that are not generic for both $M(E_6)$ and $L(E_6)$, namely $D_5$ (generic for $M(E_6)$) and $E_6(a_1)$ (not generic for either), and generic unipotent elements are blueprints by Lemma \ref{lem:genericunipotent}. If $u$ belongs to class $D_5$, then it acts on $M(E_6)$ with Jordan blocks $11,9,5,1^2$. Since there are two Jordan blocks of size $1$ and only one of size $11$, $M(E_6)\downarrow_H$ must have a trivial summand. To see this, each non-trivial indecomposable summand of dimension congruent to $1$ modulo $11$ has dimension $12$ and uses up a block of size $11$, by Lemma \ref{lem:indJordanblocks}.

We therefore assume that $u$ belongs to class $E_6(a_1)$, so acts as $11^2,5$ on $M(E_6)$ and as $11^6,9,3$ on $L(E_6)$. There are five indecomposable modules of dimension congruent to $5$ modulo $11$, which up to duality are
\[ 5,\qquad 7,5,3/5,7,\qquad 9,7,5,3/1,3,5,7,9,\]
with the last one having dimension $49$, which is not allowed, and the second one having dimension $27$, with trace of an involution $-1$, which is also not allowed (see Appendix \ref{app:traces}). Thus $M(E_6)\downarrow_H$ is the sum of $5$ and a $22$-dimensional projective module.

We now use traces of semisimple elements of orders at most $6$ to see which sums of projectives and a $5$ are conspicuous, finding two, namely
\[ 11\oplus P(1)\oplus 5\qquad\text{and}\qquad P(9)\oplus 5.\]
The first stabilizes a line on $M(E_6)$ but does not have a trivial summand, hence lies inside a $D_5$-parabolic subgroup, acting on $M(E_6)$ uniserially as $10/16/1$, and the image of $H$ inside the $D_5$-Levi subgroup must act as $1/9$ on $M(D_5)$, which is not possible as $M(D_5)$ is a self-dual module. Thus the first case does not exist, and $M(E_6)\downarrow_H$ must be the second.

The corresponding sets of composition factors on $L(E_6)$ (see (\ref{li:strata}) from Chapter \ref{ch:strategy}) are
\[ 11,9^3,7^4,3^3,1^3\qquad\text{and}\qquad 11^2,9,7^3,5^4,3^2.\]
Since $L(E_6)$ is self-dual and there is a unique self-dual module congruent to each dimension modulo $11$ by Lemma \ref{lem:selfdualsl2p}, the fact that $u$ acts on $L(E_6)$ with blocks $11^6,9,3$ means that $9$ and $3$ must be summands of $L(E_6)\downarrow_H$. The first set of factors cannot form a projective and these summands, but the second case can, yielding
\[ 11^{\oplus 2}\oplus P(7)\oplus P(5)\oplus 9\oplus 3.\]
By Proposition \ref{prop:sl2ifsplitoff}, the $3$-dimensional summand is an $\slf_2$-subalgebra, as claimed.
\end{proof}

When $p=13$, the only non-generic class is the regular unipotent class. We will show more generally that if $H$ contains a regular unipotent element then $H$ either lies in $F_4$, or $p=13$ and $H$ is a non-$\bG$-completely reducible subgroup in a $D_5$-parabolic subgroup of $\bG$.

\begin{proposition}\label{prop:e6char>11} Suppose that $p\geq 13$. If $H$ contains a regular unipotent element then $H$ is contained in a conjugate of $F_4$, or $p=13$ and $H$ is a non-$\bG$-completely reducible subgroup of the $D_5$-parabolic subgroup acting on $M(E_6)$ as
\[ (1/11/1)\oplus (9/5).\]
or its dual. In either case, $H$ stabilizes a line on $M(E_6)$.

If $H$ does not contain a regular unipotent element, then $H$ is a blueprint for $M(E_6)\oplus M(E_6)^*$.
\end{proposition}
\begin{proof} Suppose first that $u\in H$ is not regular. From \cite[Table 5]{lawther1995}, we see that $u$ is generic for $M(E_6)$, whence $u$ and therefore $H$ are blueprints for $M(E_6)\oplus M(E_6)^*$ by Lemma \ref{lem:genericunipotent}. Hence the result holds. Thus in what follows we may assume that $u$ is a regular unipotent element.

Suppose that $p\geq 17$: the action of a regular unipotent element on $M(E_6)$ has Jordan blocks $17,9,1$. Thus if $p\geq 19$, we must have that 
\[M(E_6)\downarrow_H=17\oplus 9\oplus 1\]
by Lemma \ref{lem:indJordanblocks}, and so $H$ lies inside either $F_4$, as desired, or a $D_5$-parabolic subgroup. However, a $D_5$-parabolic has composition factors of dimensions $10$, $16$ and $1$. These are incompatible with the composition factors of $M(E_6)\downarrow_H$ above, so $H\leq F_4$. For $p=17$, the Jordan blocks $17,1$ from the action of $u$ could come from an $18$-dimensional indecomposable module for $H$. However, the $9$ is definitely a summand, so in particular $H$ has three composition factors on $M(E_6)$. However, $u$ is contained in the regular class, which is generic for $p=17$, hence $H$ is a blueprint for $M(E_6)$ by Lemma \ref{lem:genericunipotent}. In particular $H$ is contained in a member $\bX$ of $\mathscr X$, which are listed in Appendix \ref{app:actions}. Since $\bX$ contains a regular unipotent element (eliminating all reductive maximal subgroups except for $F_4$ using \cite{lawther2009}) and must have at most three composition factors on $M(E_6)$, and if it does have three then one has dimension $9$ (eliminating all parabolic subgroups), we must have $H\leq F_4$, as claimed.

\medskip

We therefore have that $p=13$, and $u$ acts on $M(E_6)$ with factors $13^2,1$. Suppose that the block of size $1$ in the action of $u$ arises from a trivial summand in $M(E_6)\downarrow_H$. From the proof of Proposition \ref{prop:f4char13} we see that the conspicuous sets of composition factors are
\[ (5/7,9/5)\oplus 1\qquad\text{and}\qquad (9/3,5/9)\oplus 1.\]
Since there is no $10$-dimensional quotient not involving the trivial summand, these structures are incompatible with coming from a $D_5$-parabolic subgroup, and so $H\leq F_4$, as needed.

We thus assume that $M(E_6)\downarrow_H$ has no trivial summand. We therefore have a projective of dimension $13$ (either $P(1)$ or $13$, both with a trace of $1$ for the involution) and a module $i/(p+1-i)$, with a trace of $\pm 2$. As the trace of an involution $t\in H$ on $M(E_6)$ is either $3$ or $-5$ from Appendix \ref{app:traces}, we see that $t$ has a trace of $+2$ on $i/(p+1-i)$, and hence $i=5,9$. This means that, up to duality, $M(E_6)\downarrow_H$ is either
\[ 13\oplus (5/9)\qquad \text{or}\qquad (1/11/1)\oplus (9/5).\]
The second case has the structure claimed in the proposition. Using Lemma \ref{lem:e6stabs}, we see that $H$ cannot lie in $F_4$ as $F_4$ has a trivial summand on $M(E_6)$, and thus $H$ lies in the $D_5$-parabolic subgroup. However, it again does not lie in the $D_5$ Levi subgroup because that has a trivial summand on $M(E_6)$, hence $H$ is a non-$\bG$-completely reducible subgroup.

Thus we are left to eliminate the first case. Here we take the Borel subgroup $B$ of $H$: the exact module structure of $B=\gen u\rtimes \gen x$ on the $27$-dimensional module $M(E_6)$ is up to duality as follows, where $\zeta$ is a cube root of unity. The action of $u$ on $M(E_6)$ has blocks $13^2,1$, and the element $x$ acts on each trivial composition factor of $M(E_6)\downarrow_{\gen u}$ as a $6$th root of unity. The socle structure of the action of $B$ on $M(E_6)$ is given below.
\[ \begin{array}{c}
1\\-\zeta\\\zeta^2\\-1\\\zeta\\-\zeta^2\\1\\-\zeta\\\zeta^2\\-1\\\zeta\\-\zeta^2\\1
\end{array}\quad
\begin{array}{c}
\zeta\\-\zeta^2\\1\\-\zeta\\\zeta^2\\-1\\\zeta\\-\zeta^2\\1\\-\zeta\\\zeta^2\\-1\\\zeta
\end{array}\quad \zeta^2\]
Since $F_4$ acts on $M(E_6)$ as $M(F_4)\oplus 1$, the point that $B$ stabilizes cannot be a $F_4$ point. Thus it is either a $D_5$-parabolic point or a $B_4$ point, but either way $H$ lies inside a $D_5$-parabolic subgroup, either one stabilizing a line or one stabilizing a hyperplane.

Let $v$ be a unipotent element of $D_5$ contained in the image of $L$ inside the $D_5$-Levi subgroup. Thus $v$ acts on the $10$ and $16$ as subquotients of the action of $u$ on $M(E_6)$, namely $13^2,1$. Therefore $v$ acts on both the $10$ and the $16$ with at most three Jordan blocks, and if it has three then one is of size $1$. 

We can read off the unipotent classes of $D_5$ from the table for $D_6$, \cite[Table 6]{lawther2009}, which shows that there are only three unipotent classes, $A_4$, $D_5(a_1)$ and $D_5$, that act with at most three blocks on the $10$. From the embedding of the $D_5$-Levi subgroup into $E_6$ we can easily deduce the actions of these on the $16$, as we just consult \cite[Table 5]{lawther1995} which lists the block sizes for the classes for $E_6$, and look for the unipotent classes with these names. This gives us the list below.
\begin{center}
\begin{tabular}{cccc}
\hline Class & $A_4$ & $D_5(a_1)$ & $D_5$
\\\hline Action on $10$ & $5^2$ & $7,3$ & $9,1$
\\ Action on $16$ & $7,5,3,1$ & $7^2,2$ & $9,7$
\\ \hline
\end{tabular}
\end{center}
We therefore see that $v$ comes from the regular class $D_5$, and so the image $\bar B$ of $B$ in $D_5$, which contains $v$, must act on the self-dual module $10$ as
\[ 1\oplus (\zeta^2/-1/\zeta/-\zeta^2/1/-\zeta/\zeta^2/-1/\zeta).\]
This is a submodule of the action of $B$ above, and we therefore see that $\bar B$ acts on the $16$ with eigenvalues
\[ (1,-1)^2,(\zeta,\zeta^2,-\zeta,-\zeta^2)^3;\]
these cannot form modules of dimension $9$ and $7$, since a module of dimension $9$ needs exactly three $\pm 1$ eigenvalues, and a module of dimension $7$ needs at least two $\pm 1$s.

This proves that $H$ cannot embed with these composition factors, and completes the proof of the proposition.
\end{proof}

We will construct this non-$\bG$-completely reducible subgroup of the $D_5$-parabolic subgroup when $p=13$; the same construction works for the $E_6$-parabolic subgroup of $E_7$ and $p=19$.

Let $(\bG,p,\bX,\bY)$ be either $(E_6,13,D_5,B_4)$ or $(E_7,19,E_6,F_4)$. One of the stabilizers of a point on $M(\bG)$ is a subgroup that is the extension of a unipotent group by $\bY$, so let $H$ be a copy of $\PSL_2(p)$ inside $\bY$ that contains a regular unipotent element, the fixed points of a principal $\PSL_2$ subgroup of $\bY$. This copy of $H$ embeds in $\bX$, of course, and the action of $\bX$ on the unipotent radical of the $\bX$-parabolic subgroup is as a single simple module, so that the $1$-cohomology is easy to compute. We see that the restriction of this simple module to $H$ contains a summand of dimension $p-2$, hence the $1$-cohomology of $H$ on the unipotent radical is $1$-dimensional. There is an action of the torus of the $\bX$-parabolic subgroup outside of $\bX$ on this cohomology group, and this yields two conjugacy classes of subgroups $H$ in the $\bX$-parabolic subgroup, one inside $\bX$ and another class of complements. Given the composition factors of $H$ on $M(\bG)$, together with the table from \cite{lawther1995}, there is a unique possible module structure for $M(\bG)\downarrow_H$ if $H$ does not lie inside $\bX$ but merely the $\bX$-parabolic subgroup of $\bG$, and the action of a non-trivial unipotent element of $H$ on this module has Jordan blocks $13^2,1$ and $19^2,18$ (for $p=13,19$ respectively), consistent with coming from the regular unipotent class of $\bG$.

\medskip

We now summarize the results of this chapter and prove Theorem \ref{thm:e6}. Let $H=\PSL_2(p^a)$ inside $\bG=E_6$.

If $p=2$ then $H$ is a blueprint for $M(E_6)\oplus M(E_6)^*$ when $a\geq 5$ by Theorem \ref{thm:goodblueprint}, and when $a=4$ by Proposition \ref{prop:e6char2a=4}. Thus $H$ is strongly imprimitive by Proposition \ref{prop:blueprintissi}. If $a=3$ or $a=2$ then $H$ stabilizes a line or a $2$-space on $M(E_6)$ or its dual by Proposition \ref{prop:e6char2a=3} and \cite[Proposition 5.4]{craven2015un2} respectively. Therefore $H$ is strongly imprimitive by Propositions \ref{prop:fixlineonMG} and \ref{prop:fix2spaceonMG}.

If $p=3$ then $H$ is a blueprint for $M(E_6)\oplus M(E_6)^*$ when $a\geq 4$ by Theorem \ref{thm:goodblueprint}, and for $a=3$ it is either a blueprint for that module, or stabilizes a line on it by Proposition \ref{prop:e6char3}. If $a=2$ then $H$ stabilizes a line on $M(E_6)\oplus M(E_6)^*$ by \cite[Proposition 6.2]{craven2015un2}. Thus $H$ is strongly imprimitive by Propositions \ref{prop:blueprintissi}, \ref{prop:fixlineonMG} and \ref{prop:fix2spaceonMG}.

If $p\geq 5$ and $p^a\geq 36$ then we apply Theorem \ref{thm:goodblueprint} again to obtain that $H$ is a blueprint for $M(E_6)\oplus M(E_6)^*$, so strongly imprimitive as above. If $p^a=5,7,25$ then $H$ stabilizes a line on one of $M(E_6)$, $M(E_6)^*$ and $L(E_6)$ by Propositions \ref{prop:e6char5} and \ref{prop:e6char7}, so $H$ is strongly imprimitive by Proposition \ref{prop:fixlineonMG} as above, and also Proposition \ref{prop:fixlineonLG} for $L(E_6)$.

For $p\geq 13$, $H$ is either a blueprint for $M(E_6)\oplus M(E_6)^*$ or stabilizes a line on $M(E_6)$ by Proposition \ref{prop:e6char>11}, so is strongly imprimitive. If $p=11$, if $H$ is a blueprint for $M(E_6)\oplus M(E_6)^*$ or stabilizes a line on $M(E_6)$ then $H$ is strongly imprimitive. Proposition \ref{prop:e6char11} states that if these do not hold then $H$ stabilizes a unique $3$-space on $L(E_6)$ that is an $\slf_2$-subalgebra of $L(E_6)$, and so $H$ is again strongly imprimitive by Corollary \ref{cor:sl2coxeter-1}. Thus $H$ is strongly imprimitive in all cases.

\medskip

We now complete the proof of Theorem \ref{thm:e6}. Any maximal subgroup of the finite group that is almost simple with socle a copy of $\PSL_2(p^a)$ is strongly imprimitive by the above arguments. Therefore it is the fixed points of a maximal subgroup from \cite[Corollary 5]{liebeckseitz2004}. It cannot be a parabolic, and no $A_1$-type subgroup appears in \cite[Table 5.1]{liebecksaxlseitz1992}, so it cannot be the fixed points of a maximal-rank subgroup. There are no other $A_1$ subgroups in \cite[Corollary 5]{liebeckseitz2004}, so there are no such maximal subgroups.

\chapter{The Proof for \texorpdfstring{$E_7$}{E7} in Characteristic 2}

In this chapter, $k$ is an algebraically closed field of characteristic $2$ and $\bG=E_7(k)$. Let $H=\SL_2(2^a)$, and let $u$ be an element of order $2$ in $H$. The case of $p=2$ is very different from odd characteristic because if $p$ is odd then a copy of $\PSL_2(p^a)$ inside the simple group of type $E_7$ can lift in the simply connected group to either $\PSL_2(p^a)\times 2$ or $\SL_2(p^a)$, and the two possibilities require very different strategies. In characteristic $2$ there is no such bifurcation.

Theorem \ref{thm:goodblueprint} states that if $a\geq 7$ then $H$ is a blueprint for $M(E_7)$, and therefore $H$ is strongly imprimitive by Proposition \ref{prop:blueprintissi}. Thus in what follows we may assume that $2^a\leq 64$.

The case $a=2$ is done in \cite[Proposition 5.4]{craven2015un2} so we may assume that $3\leq a\leq 6$. Furthermore, if $a=5,6$ and $M(E_7)\downarrow_H$ has at least six trivial composition factors, then $H$ is a blueprint for $M(E_7)$ by Proposition \ref{prop:f4a4}.

We can use a computer to find which semisimple elements are blueprints for $M(E_7)$ even when they have order smaller than $77$, or $30$ when they centralize a $6$-space. For example, of the $2430$ classes of elements of order $17$, $1892$ of them are blueprints for $M(E_7)$, which helps reduce the number of conspicuous sets of composition factors that need to be considered when $a=4$.

\medskip

We will proceed by splitting the (many) possibilities according to the number of trivial composition factors. Since $H$ is strongly imprimitive if $H$ stabilizes either a $1$- or $2$-space on $M(E_7)$ by Propositions \ref{prop:fixlineonMG} and \ref{prop:fix2spaceonMG}, one particularly troublesome case might be when $M(E_7)\downarrow_H$ contains no composition factors at all of dimensions $1$ or $2$.

In this case if we cannot prove that $H$ is a blueprint for $M(E_7)$ then we have to switch to the Lie algebra $L(E_7)^\circ$, which we recall has dimension $132$, not $133$ in the case $p=2$. We address this situation now: luckily there are very few such sets of composition factors that are conspicuous.

\begin{proposition}\label{prop:e7char2not1or2} Suppose that $3\leq a\leq 6$. Suppose that there are no $1$- or $2$-dimensional composition factors in $M(E_7)\downarrow_H$.
\begin{enumerate}
\item We cannot have $a=3,4$.
\item If $a=5,6$ then $H$ is a blueprint for $M(E_7)$.
\end{enumerate}

\end{proposition}
\begin{proof}
The trace of an element of order $3$ on $M(E_7)$ is one of $-25,-7,2,20$ (see Appendix \ref{app:traces}), and so if $H$ acts on $M(E_7)$ with no composition factors of dimension $1$ or $2$, then the dimensions of the composition factors are one of seven possibilities:
\begin{equation} 32,16,4^2,\qquad 32,8,4^4,\qquad 16^3,8,\qquad 16^2,8^2,4^2,\qquad 16,8^3,4^4,\qquad 8^7,\qquad 8^4,4^6. \label{eq:dimsfore7no12}\end{equation}
For these, no general arguments about stabilizing $1$- and $2$-spaces can work, and $u$ must act projectively on $M(E_7)$, so cannot be generic. We will just have to deal with them case by case, switching to the Lie algebra on one occasion. From now on we assume that $M(E_7)\downarrow_H$ has composition factors of dimensions from (\ref{eq:dimsfore7no12}). Fortunately there are very few of these.

\medskip

\noindent \textbf{Case $a=3$}: The only conspicuous set of composition factors for $M(E_7)\downarrow_H$ (subject to our dimension requirements) is
\[ 8^4,4_{1,2}^2,4_{1,3}^2,4_{2,3}^2,\]
which does not have a corresponding set of composition factors on $L(E_7)^\circ$ (see (\ref{li:strata}) from Chapter \ref{ch:strategy}).
\medskip

\noindent \textbf{Case $a=4$}: We find no conspicuous sets of composition factors for $M(E_7)\downarrow_H$ at all.

\medskip

\noindent \textbf{Case $a=5$}: Up to field automorphism of $H$, there are two conspicuous sets of composition factors for $M(E_7)\downarrow_H$, namely
\[ 16_{1234},8_{135},8_{235},8_{345},4_{13}^2,4_{23},4_{34}\qquad\text{and}\qquad 8_{123}^2,8_{124},8_{135},4_{12},4_{13},4_{14}^2,4_{15}^2.\]
Let $x\in H$ have order $31$. In both of the cases above, the semisimple class of $\bG$ to which $x$ belongs is determined by the eigenvalues of $x$ on $M(E_7)$ (as we have a list of all 53503 semisimple classes). It is easy to check using the preimage trick from the end of Section \ref{sec:blueprints}, that there is an element $\hat x$ in $\bG$ of order $93$ and such that $\hat x^3=x$, with $\hat x$ and $x$ having the same number of distinct eigenvalues on $M(E_7)$. Since $\hat x$ is a blueprint for $M(E_7)$ by Theorem \ref{thm:blueprintsminimal}, $x$, and therefore $H$, are blueprints for $M(E_7)$.

\medskip

\noindent \textbf{Case $a=6$}: We do not have lists of the traces of elements of orders $63$ and $65$, but we can check whether a given matrix possesses the eigenvalues of a semisimple element of order $63$ on $M(E_7)$ by using the preimage trick from the end of Section \ref{sec:blueprints}. Doing this to the seven possible sets of dimensions in (\ref{eq:dimsfore7no12}) yields the following table. In this, the number of sets of composition factors up to field automorphism is given in the second column, and those that are conspicuous using elements of order up to $21$ and $63$ are given in the third and fourth columns respectively.
\begin{center}\begin{tabular}{cccc}
\hline Dimensions & Number of modules & Conspicuous up to $21$ & Conspicuous for $63$
\\\hline $32,16,4^2$ & $1800$ & $1$ & $0$
\\ $32,8,4^4$ & $61200$ & $5$ & $0$
\\ $16^3,8$ & $2270$ & $1$ & $0$
\\ $16^2,8^2,4^2$ & $504240$ & $32$ & $0$
\\ $16,8^3,4^4$ & $11781000$ & $159$ & $2$
\\ $8^7$ & $109660$ & $1$ & $0$
\\ $8^4,4^6$ & $57206136$ & $934$ & $9$
\\ \hline
\end{tabular}\end{center}
We thus simply need to check whether for a given conspicuous set of composition factors, that any conjugacy class of elements of order $63$ with the correct eigenvalues on $M(E_7)$ is a blueprint for $M(E_7)$. This can easily be done with a computer, finding elements of order $315=5\cdot 63$ that have the same number of eigenvalues on $M(E_7)$ and power to our element of order $63$, and so we prove the result.
\end{proof}

We have now dealt with the case where $M(E_7)\downarrow_H$ has no $1$- or $2$-dimensional composition factors. We generally cannot prove that $H$ stabilizes a line on $M(E_7)$, and often want to prove that $H$ stabilizes a $2$-space on $M(E_7)$. This for example could happen if $M(E_7)\downarrow_H$ possesses $2$-dimensional factors but no trivial factors. This is the next proposition.

\begin{proposition}\label{prop:e7char2no1s} Suppose that $M(E_7)\downarrow_H$ has at least one $2$-dimensional composition factor and no trivial composition factors.
\begin{enumerate}
\item If $a=3$ then $H$ stabilizes a $2$-space on $M(E_7)$.
\item If $a=4,5,6$ then $H$ stabilizes a $2$-space on $M(E_7)$ or is a blueprint for $M(E_7)$.
\end{enumerate}
\end{proposition}
\begin{proof} \textbf{Case $a=3$}: Any $8$s split off, so we just consider the $4$s and $2$s. The projective cover of $4_{i,i+1}$ is
\[ P(4_{i,i+1})=4_{i,i+1}/2_{i+1}/1/2_{i-1}/1/2_{i+1}/4_{i,i+1},\]
and from this we see that no module can have a $2$-dimensional composition factor, no trivial composition factor, and not have a $2$-dimensional submodule or quotient. This proves the case $a=3$..

\medskip

\noindent \textbf{Case $a=4$}: We first compute the conspicuous sets of composition factors, finding 81 sets up to field automorphism of $H$. As stated at the start of this chapter, using the preimage trick from the end of Section \ref{sec:blueprints}, we can determine if an element of order $17$ is a blueprint for $M(E_7)$, and 1892 classes out of 2430 classes are. All but fifteen of the 81 conspicuous sets of composition factors are blueprints for $M(E_7)$ via an element of order $17$.

We can also compute which have positive $2_i$-pressure (or no $2_i$) for each $i$, and find that only eighteen of the 81 sets of factors do. (Those of non-positive $2_i$-pressure must have a submodule $2_i$ by Lemma \ref{lem:pressure}.) The intersection of these two short lists has just two sets of composition factors on it (up to field automorphism), and so we consider these two:
\[ 8_{1,2,3},4_{1,2}^2,4_{1,3}^3,4_{2,3}^3,4_{2,4},2_1^2,2_2^2,2_3^2,\qquad 
8_{1,2,4},8_{2,3,4},4_{1,2}^2,4_{1,3},4_{1,4}^2,4_{2,3}^3,2_1^2,2_3^2.\]
For the first of these, $H$ must stabilize a $2$-space on $M(E_7)$. To see this, notice that otherwise the socle can only consist of summands of $M(E_7)\downarrow_H$ and a submodule of $4_{1,2}\oplus 4_{1,3}\oplus 4_{2,3}$. The $\cf(M(E_7)\downarrow_H)$-radicals of $P(4_{1,2})$, $P(4_{1,3})$ and $P(4_{2,3})$ are
\[ 4_{2,4}/2_2/4_{1,2},\quad 4_{1,3}/2_3/4_{2,3}/2_1,2_3/4_{1,3},\quad 2_3/4_{1,3},4_{2,3}/2_3,8_{1,2,3}/4_{2,3}.\]
Thus at most a single $2_1$ can lie in $M(E_7)\downarrow_H$, which is a contradiction.

The second case is even easier, given that the corresponding submodules are
\[ 4_{1,2}/8_{1,2,4}/4_{1,2},\quad 4_{1,3}/2_1/4_{1,4},\quad 4_{1,3}/2_3/4_{2,3}.\]
This completes the proof for $a=4$.

\medskip

\noindent \textbf{Case $a=5$}: There are 30 possible multisets of dimensions for the composition factors of $M(E_7)\downarrow_H$ that have at least one $2$, no $1$s, and have the correct trace of an element of order $3$ (see Appendix \ref{app:traces}). If $H$ stabilizes a $2$-space on $M(E_7)$ then the result holds, so we suppose that this is not the case. If $H$ does not stabilize a $2$-space, then we need two $4$s in the dimensions, removing ten multisets of dimensions from the list. We can also apply Lemma \ref{lem:largest2and4}, which shows that if there are no $8$s in $M(E_7)\downarrow_H$ then we need at least as many $4$s as $2$s, removing another three. Since any $4$-dimensional factor has $\{2_i\}$-pressure at most $2$, there needs to be more than half as many $4$s as $2$s in all cases; this brings us down to ten. These are
\begin{equation}
\begin{gathered} 4^{10},2^8,\quad 8,4^9,2^6,\quad 16,4^7,2^6,\quad 8^2,4^8,2^4,\quad  16,8,4^6,2^4,\\
16^2,4^4,2^4,\quad 8^3,4^7,2^2,\quad 16,8^2,4^5,2^2,\quad 32,4^5,2^2,\quad  16^2,8,4^3,2^2.
\end{gathered}
\label{eq:setsofdimensionse7char2}\end{equation}
In these cases we switch to proving that $H$ is a blueprint for $M(E_7)$. (This could be done for the other cases but the amount of extra work is significant and so this has not been done.)

We give a table listing: the total number of possible sets of composition factors for $M(E_7)\downarrow_H$ with dimensions each of the options from (\ref{eq:setsofdimensionse7char2}); then those that are conspicuous; and finally those for which an element of order $31$ in $H$ is a blueprint for $M(E_7)$. These numbers are all up to a field automorphism of $M(E_7)$.

\begin{center}
\begin{tabular}{cccc}
\hline Case & Number & Conspicuous & 31 is blueprint
\\ \hline
$4^{10},2^8$ & 9145422 & 23 & 23
\\$8,4^9,2^6$ & 20420400 & 32 & 32
\\$16,4^7,2^6$ & 2402400 & 3 & 2
\\$8^2,4^8,2^4$ & 18718700 & 52 & 51
\\$16,8,4^6,2^4$ & 3503500 & 12 & 8
\\$16^2,4^4,2^4$ & 150150 & 2 & 2
\\ $8^3,4^7,2^2$ & 7550400 & 22 & 21
\\ $16,8^2,4^5,2^2$ & 1651650 & 20 & 19
\\$32,4^5,2^2$ & 6006 & 0 & 0
\\$16^2,8,4^3,2^2$ & 99000 & 4 & 4
\\\hline\end{tabular}
\end{center}
Thus, up to field automorphism, there are eight conspicuous sets of composition factors for $M(E_7)\downarrow_H$ that have dimensions from (\ref{eq:setsofdimensionse7char2}). They are as follows:
\[ 16_{1,2,3,4},4_{1,3},4_{1,4},4_{1,5},4_{2,3},4_{2,4}^2,4_{2,5},2_1^2,2_2^2,2_3^2,\]
\[ 8_{1,2,5},8_{2,3,5},4_{1,3},4_{1,4},4_{1,5}^2,4_{2,3}^2,4_{3,5}^2,2_1^2,2_3^2,\]
\[ 16_{1,3,4,5},8_{2,4,5},4_{1,5},4_{2,3},4_{2,5},4_{3,5},4_{4,5}^2,2_1^3,2_3,\] \[16_{1,3,4,5},8_{2,4,5},4_{1,3},4_{1,4},4_{2,3},4_{3,5},4_{4,5}^2,2_1^2,2_3^2,\]
\[
16_{1,2,3,4},8_{1,4,5},4_{1,5},4_{2,3},4_{2,4},4_{2,5},4_{3,4},4_{3,5},2_1^2,2_3^2,\] \[16_{1,2,3,5},8_{2,3,4},4_{1,3}^2,4_{1,4},4_{2,3},4_{3,5},4_{4,5},2_1^2,2_3,2_4,\]
\[ 8_{1,2,5},8_{1,3,4},8_{2,4,5},4_{1,2},4_{1,5},4_{2,3},4_{2,4},4_{3,4},4_{4,5}^2,2_1,2_3,\] \[16_{1,2,4,5},8_{1,2,4},8_{1,2,5},4_{1,2},4_{1,5},4_{2,3},4_{2,4},4_{3,4},2_1,2_2.\]
Recall from Lemma \ref{lem:extforsl22a} that $2_i$ has extensions only with $4_{i,j}$ for $j\neq i,i+1$.

\medskip

\noindent \textbf{Case 1}: This has $2_3$-pressure $0$, so $M(E_7)\downarrow_H$ has a submodule $2_3$.

\medskip

\noindent \textbf{Cases 3, 4, 5, 7, 8}: These have non-positive $2_1$-pressure, so $M(E_7)\downarrow_H$ has a submodule $2_1$.

\medskip

\noindent \textbf{Case 6}: This has $2_4$-pressure $0$, so $M(E_7)\downarrow_H$ has a submodule $2_4$.

\medskip

\noindent \textbf{Case 2}: Suppose that $H$ does not stabilize a $2$-space on $M(E_7)$. The only composition factors appearing with multiplicity greater than $1$ are $4_{1,5}$, $4_{2,3}$ and $4_{3,5}$, together with $2_1$ and $2_3$. Thus $M(E_7)\downarrow_H$ is a sum of simple summands and submodules of $P(4_{1,5})$, $P(4_{2,3})$ and $P(4_{3,5})$. The $\cf(M(E_7)\downarrow_H)$-radicals of these three projectives are
\[ 4_{1,4}/2_1/4_{1,5},\quad 2_3/4_{2,3},4_{3,5}/2_3,8_{2,3,5}/4_{2,3},\quad \text{and}\quad 2_3/4_{2,3},4_{3,5}/2_3,8_{2,3,5}/4_{3,5}\]
respectively. There is only one $2_1$ in the sum of these three modules, so there must be a $2$-dimensional factor in $M(E_7)\downarrow_H$, as claimed.

\medskip

\noindent \textbf{Case $a=6$}: If $H$ stabilizes a $2$-space on $M(E_7)$ then the result holds, so assume otherwise. We have the same ten multisets of dimensions of composition factors for $M(E_7)\downarrow_H$ from (\ref{eq:setsofdimensionse7char2}), and we perform the same analysis as before.

\begin{center}
\begin{tabular}{ccccc}
\hline Case & Number & Consp.\ to 21 & Consp.\ to 63 & 63 is blueprint
\\ \hline
$4^{10},2^8$ & 420696342 & 68 & 41 & 41
\\$8,4^9,2^6$ & 1258472670 & 369 & 76 & 76
\\$16,4^7,2^6$ & 134306100 & 121 & 9 & 9
\\$8^2,4^8,2^4$ & 1410195600 & 1068 & 104 & 104
\\$16,8,4^6,2^4$ & 244188000 & 750 & 38 & 38
\\$16^2,4^4,2^4$ & 7712064 & 89 & 12 & 12
\\ $8^3,4^7,2^2$ & 626749200 & 983 & 90 & 90
\\ $16,8^2,4^5,2^2$ & 128200860 & 1097 & 80 & 80
\\$32,4^5,2^2$ & 244188 & 15 & 1 & 1
\\$16^2,8,4^3,2^2$ & 5712000 & 208 & 24 & 24
\\\hline\end{tabular}
\end{center}
As every conspicuous set of composition factors for $M(E_7)\downarrow_H$ has an element of order $63$ that is a blueprint for $M(E_7)$, $H$ is always a blueprint for $M(E_7)$, as needed.
\end{proof}

From now on we therefore assume that $M(E_7)\downarrow_H$ possesses at least one (hence two as the dimension is even) trivial composition factor. We split this investigation up into three propositions, depending on whether $H$ has two, four, or at least six trivial composition factors on $M(E_7)$.

\begin{proposition}\label{prop:e7char22trivials} Suppose that $M(E_7)\downarrow_H$ has exactly two trivial composition factors.
\begin{enumerate}
\item If $a=3$ then $H$ stabilizes a subspace of dimension at most $2$ of $M(E_7)$.
\item If $a=4,5$ then $H$ is strongly imprimitive.
\item If $a=6$ then $H$ is either a blueprint for $M(E_7)$ or stabilizes a subspace of dimension at most $2$ of $M(E_7)$.
\end{enumerate}
\end{proposition}
\begin{proof} \textbf{Case $a=3$}: As we have seen in Proposition \ref{prop:e7char2no1s}, the projective cover of $4_{i,i+1}$ is 
\[ P(4_{i,i+1})=4_{i,i+1}/2_{i+1}/1/2_{i-1}/1/2_{i+1}/4_{i,i+1},\]
whence if $M(E_7)\downarrow_H$ has no $1$- or $2$-dimensional submodules or quotients, it is a sum of copies of $8$ and $P(4_{i,i+1})$ for various $i$. In particular, since $\dim(P(4_{i,i+1}))=16$, we must have $P(4),8^5$, as there are exactly two trivial factors. Thus the composition factors of $M(E_7)\downarrow_H$ are $8^5,4^2,2^3,1^2$, on which an element of order $3$ acts with trace $-4$. From Appendix \ref{app:traces} we see that this is not a trace of an element of order $3$ on $M(E_7)$. This completes the proof for $a=3$.

\medskip

\noindent \textbf{Case $a=4$}: If $H$ stabilizes a subspace of dimension $1$ or $2$ on $M(E_7)$, then $H$ is strongly imprimitive by Propositions \ref{prop:fixlineonMG} and \ref{prop:fix2spaceonMG} respectively. Thus we assume that this is not the case. We may also assume that $H$ is not a blueprint for $M(E_7)$ by Proposition \ref{prop:blueprintissi}. In addition $H$ cannot stabilize a line on $L(E_7)^\circ$, else it is strongly imprimitive by Proposition \ref{prop:fixlineonLG}.

Using all semisimple elements, there are (up to field automorphism) $113$ conspicuous sets of composition factors for $M(E_7)\downarrow_H$ with exactly two trivial composition factors. Only $80$ of these have corresponding factors on $L(E_7)^\circ$, and of these only $39$ have either no $2_i$ or positive $2_i$-pressure for every $i$. One can eliminate three more as they have no $4$-dimensional factors appearing with multiplicity greater than $1$, so must stabilize either a line or $2$-space as $M(E_7)$ is self-dual. We also exclude those whose corresponding factors on $L(E_7)^\circ$ have pressure less than $6$ (and therefore stabilize a line on $L(E_7)^\circ$ using Lemma  \ref{lem:char2pressuregeneral} and the fact that $u$ must act with at least six blocks of size $1$ on $L(E_7)^\circ$ from \cite[Table 8]{lawther1995}). This leaves $25$ conspicuous sets of composition factors for $M(E_7)\downarrow_H$.

This is, however, still too many to list. Next, we remove all blueprints and `near blueprints'. Let $x$ be an element of order $17$ in $H$, and note that the eigenvalues of $x$ on $M(E_7)$ determine its conjugacy class in $\bG$, by a check of the list of all 2430 classes. Using the preimage trick from the end of Section \ref{sec:blueprints}, we look for elements $\hat x$ of order $85$ in $\bG$ whose fifth power is $x$ and have at most one more eigenvalue on $M(E_7)$ than $x$. If $\hat x$ has the same number of eigenvalues on $M(E_7)$ as $x$ then of course $\hat x$ stabilizes all subspaces of $M(E_7)$ that $H$ stabilizes. If $\hat x$ has one more eigenvalue, then since $\hat x$ must be real (as all semisimple elements of $E_7$ are real) it must be the $(+1)$-eigenspace of $x$ on $M(E_7)$ that is split in two. As $x$ only has a $(+1)$-eigenspace on the trivial simple module, this means that $\hat x$ stabilizes every non-trivial simple submodule of $M(E_7)\downarrow_H$.

Since $o(\hat x)=85$, $\hat x$ is a blueprint for $M(E_7)$ by Theorem \ref{thm:blueprintsminimal}, and therefore there is a positive-dimensional subgroup of $\bG$ stabilizing every non-trivial simple submodule of $M(E_7)\downarrow_H$, certainly enough to guarantee that $H$ is strongly imprimitive by Proposition \ref{prop:intersectionstabilizers}.

Nine sets of factors have such elements, so we are down to sixteen sets of composition factors. We can remove a few more by introducing a general module, which we will use again in similar proofs.

Let $W$ be the subquotient obtained from $M(E_7)\downarrow_H$ by quotienting out by the $\{8_{i,j,l},16\}$-radical and taking the $\{8_{i,j,l},16\}$-residual, and removing any $4$-dimensional simple summands. Since $H$ can be assumed not to stabilize a line or $2$-space on $M(E_7)$, the socle and top of $W$ consists of $4$-dimensional modules, and the factors of $\soc(W)$ (and $\topp(W)$) consist of $4$-dimensional simple modules that occur with multiplicity at least $2$ in $M(E_7)\downarrow_H$, and hence in $W$. Let $S_1,\dots,S_r$ be the $4$-dimensional simple modules that appear in $M(E_7)\downarrow_H$ with multiplicity at least $2$. (If a module appears more than twice, we take the floor of half of its multiplicity, since this is the maximum number of times it may appear in the socle.)

We construct the largest submodule of $P(S_1\oplus \cdots\oplus S_r)$ that consists solely of composition factors from $M(E_7)\downarrow_H$, and then remove all quotients not isomorphic to one of the $S_i$ to form a module $W'$; certainly $W\leq W'$. Thus $W'$ must have at least two trivial factors, and all the requisite $2$-dimensional factors.

Of the sixteen sets of composition factors we had before this test, only six remain after checking that $W'$ contains enough $1$- and $2$-dimensional composition factors. These are as follows:
\[8_{1,3,4},4_{1,3}^2,4_{1,4}^3,4_{2,3}^2,4_{2,4},2_1^3,2_2,2_3^2,2_4,1^2,\] \[8_{1,3,4},4_{1,2},4_{1,3},4_{1,4},4_{2,3},4_{2,4}^2,4_{3,4}^2,2_1^3,2_2,2_3,2_4^2,1^2,\]
\[ 8_{1,2,4},4_{1,2}^2,4_{1,3}^2,4_{1,4},4_{2,3},4_{2,4},4_{3,4},2_1^2,2_2^2,2_3^2,2_4,1^2,\] \[16,4_{1,2},4_{1,3},4_{1,4},4_{2,4}^2,4_{3,4},2_1^2,2_2^2,2_3^2,2_4,1^2,\]
\[8_{1,3,4}^2,4_{1,2},4_{1,3}^2,4_{1,4}^3,4_{2,3},2_1^2,2_2,2_3^2,1^2,\] \[8_{1,2,4},8_{1,3,4},4_{1,2},4_{1,3},4_{1,4}^2,4_{2,3}^2,4_{2,4},2_1^2,2_2,2_3,2_4,1^2\]

We can eliminate some more using module structures.

\noindent \textbf{Case 1}: Suppose that $4_{1,4}$ lies in the socle of $M(E_7)\downarrow_H$. If it is a summand, we quotient it out and ignore it, so suppose it is a submodule but not a summand, and let $U$ denote the $\{1,2_i\}$-radical of the quotient module $M(E_7)\downarrow_H/4_{1,4}$, lifted to $M(E_7)\downarrow_H$. This module $U$ is a submodule of $P(4_{1,4})$ and so we use Lemma \ref{lem:largestsubmodofP4}, seeing that $U$ is a submodule of
\[ 2_1/1/2_2/1/2_1/4_{1,4};\]
if both trivials are in $U$ then the quotient $M(E_7)\downarrow_H/U$ has $2_3$-pressure $0$, so has $2_3$ as a submodule, a contradiction from the definition of $U$. If there is a single trivial in $U$ then first replace $U$ by the $7$-dimensional submodule $1/2_1/4_{1,4}$ of $U$, and since $M(E_7)$ is self-dual, there is a (unique) corresponding submodule $U'$ such that the quotient $M(E_7)\downarrow_H/U'$ is isomorphic to $U^*$. If $U\leq U'$ then the quotient $U'/U$ has no trivials and again it has $2_3$-pressure $0$, so we obtain a contradiction. Thus $U'$ does not contain $U$, and we claim that in this case an involution $u$ must act with exactly two trivial Jordan blocks, which contravenes the possible actions given in \cite[Table 7]{lawther1995}. To see this, first let $M$ denote the $\{1\}'$-residual modulo the $\{1\}'$-radical of $M(E_7)\downarrow_H$, so it is a submodule of $P(1)$, as otherwise it is simply $1^{\oplus 2}$, with this impossible by \cite[Table 7]{lawther1995}. The submodule $U$ of $M(E_7)\downarrow_H$ has image inside $M$ which is just $\soc(M)$, and the image of $U'$ has image inside $M$ which is simply $\rad(M)$. It is therefore clear that the image of $U'$ contains the image of $U$ and, since $U$ is uniserial, $U'$ contains $U$.

Hence $U$ is a submodule of $2_1/4_{1,4}$. Thus we can remove any $4_{1,4}$ in the socle and top, perhaps remove two copies of $2_1$ that are now in the socle and top, and assume that the resulting module $V'$ is a self-dual submodule of $P(4_{1,3})\oplus P(4_{2,3})$.

We now give the three modules obtained from the following procedure, given a socle $S$ that is a submodule of $4_{1,3}\oplus 4_{2,3}$:
\begin{enumerate}
\item Take the preimage $S_1$ in $P(S)$ of the radical of the quotient module $P(S)/S$ corresponding to all composition factors of $M(E_7)\downarrow_H$ other than those in $S$;
\item Take the preimage $S_2$ in $P(S)$ of the $\cf(S)$-radical of the quotient $P(S)/S_1$;
\item Take the $\cf(S)'$-residual $S_3$ of $S_2$.
\end{enumerate}
This must contain the module $V'$, so we examine the composition factors of the modules $S_3$ for the choices of $S$, which are
\[ 4_{1,3},4_{1,3}/2_1,2_3/4_{1,4},4_{2,3}/2_1,2_3/4_{1,3},\qquad 4_{2,3}/2_3/1/2_4/1/2_3/4_{2,3},\] 
\[4_{2,3}/2_3/1,4_{1,3}/2_1,2_4/1,4_{1,3},4_{1,4},4_{2,3}/2_1,2_3,2_3/4_{1,3},4_{2,3}.\]
None of these has $2_2$ as a composition factor, and this yields a contradiction.

\medskip

\noindent \textbf{Case 2}: $W'$ might have enough $2$-dimensional factors, but in order to have three copies of $2_1$ in $W'$ we need both $4_{2,4}$ and $4_{3,4}$ in the socle, whence they cannot appear elsewhere in the module (which they can do in our construction of $W'$). With this restriction, that $4_{2,4}$ and $4_{3,4}$ can only appear in the socle and top of $M(E_7)\downarrow_H$, the module $W'$ becomes
\[ 4_{3,4}/2_3,2_4/1,1,4_{1,3},4_{2,4}/2_1,2_1,2_2,2_3/1,1,4_{1,2},4_{2,4},4_{3,4}/2_2,2_4,2_4/4_{2,4},4_{3,4},\]
which does not have three copies of $2_1$ in it, a contradiction.

\medskip

\noindent \textbf{Case 3}: The socle of $W$ cannot simply be $4_{1,2}$ as there is no $2_1$ in its contribution to $W'$. If it is $4_{1,2}\oplus 4_{1,3}$ then, arguing as in the previous case, we obtain
\[4_{1,2}/2_2/1,4_{1,3},4_{1,3}/2_1,2_3,2_3/1,4_{1,2},4_{1,4},4_{2,3}/2_1,2_2,2_3,8_{1,2,4}/4_{1,2},4_{1,3},\]
and if it is just $4_{1,3}$ then we obtain
\[ 4_{1,3},4_{1,3}/2_1,2_3/4_{1,4},4_{2,3}/2_1,2_3/4_{1,3}.\]
In neither case is $2_4$ a composition factor of this module, so it does not contain all of $W$, a contradiction.

\medskip

\noindent \textbf{Case 4}: The socle of $W$ must be $4_{2,4}$, and if so then no $4_{2,4}$ can appear outside of the socle and top of $W$. Taking the radical of the quotient $P(4_{2,4})/4_{2,4}$ with factors all other composition factors of $M(E_7)\downarrow_H$, then adding on as many copies of $4_{2,4}$ on top of that, then taking the $\{4_{2,4}\}'$-residual of this (since the socle of $W$ must be $4_{2,4}$), we end up with
\[ 4_{2,4},4_{2,4}/2_2,2_4/4_{1,2},4_{3,4}/2_2,2_4/4_{2,4},\]
which is clearly wrong.

\medskip
\noindent \textbf{Case 5}: Choose $\zeta$ a primitive $17$th root of unity so that $x$ acts on $2_1$ with eigenvalues $\zeta^{\pm 1}$, then $x$ acts on $M(E_7)$ with eigenvalues
\[ 1^2, (\zeta^{\pm 1})^3,(\zeta^{\pm 2})^2,(\zeta^{\pm 3})^5,(\zeta^{\pm 4})^4,(\zeta^{\pm 5})^4,(\zeta^{\pm 6})^3,(\zeta^{\pm 7})^3,(\zeta^{\pm 8})^3,\]
and there is an element $\hat x$ of order $85$ in $\bG$ that powers to $x$ and has almost the same eigenspaces, except that it splits the $\zeta^{\pm 1}$ and $1$-eigenspaces, so has twenty distinct eigenvalues on $M(E_7)$. An easy calculation shows that the only composition factors of $M(E_7)$ on which $x$ has $1$ or $\zeta^{\pm 1}$ as an eigenvalue are $1$, $2_1$ and $4_{1,2}$: $1$ and $2_1$ are not submodules of $M(E_7)\downarrow_H$ by assumption. If $4_{1,2}$ is a submodule then it is a summand, so there must be another factor in the socle, and one that is not a summand. Thus the simultaneous stabilizer of all simple submodules that are not summands is positive dimensional and $N_{\Aut^+(\bG)}(H)$-stable, hence $H$ is strongly imprimitive by Proposition \ref{prop:intersectionstabilizers}.

\medskip

\noindent \textbf{Case 6}: The possible factors of $\soc(W)$ are $4_{1,4}$ and $4_{2,3}$, with both required for all of the $2$-dimensional factors to be present, as an examination of $W'$ proves. In this case, we do as in Cases 3 and 4 to find that $W$ is a submodule of
\[ 4_{1,4}/2_1/1,4_{1,4}/2_1,2_2/1,1,4_{1,3},4_{1,4}/2_1,2_3,8_{1,3,4}/4_{1,4},4_{2,3},\]
which does not have a copy of $2_4$ in it, a contradiction.

\medskip

\noindent \textbf{Case $a=5$}: If $H$ stabilizes a line or $2$-space on $M(E_7)$, then $H$ is strongly imprimitive by Propositions \ref{prop:fixlineonMG} and \ref{prop:fix2spaceonMG} respectively. Thus we may assume that this is not the case. As there are exactly two trivial composition factors in $M(E_7)\downarrow_H$, we need at least three $2$-dimensional composition factors in order not to stabilize a line by Lemma \ref{lem:pressure}, and at least two composition factors of dimension $4$, to avoid fixing a line or $2$-space on $M(E_7)$. There are seventeen possible sets of dimensions of composition factors with these properties that also have the correct trace of an element of order $3$ (given in Appendix \ref{app:traces}). If there are exactly two $4$s in $M(E_7)\downarrow_H$ then we can use Lemma \ref{lem:largestsubmodofP4} to see that we can have exactly three $2$s, thus eliminating two of these cases, and if there are three $4$s we can have at most eight $2$s, eliminating two more.

We now give a table listing the possible sets of dimensions, together with the number of sets of composition factors (up to field automorphism) with those dimensions, and those that are conspicuous. Furthermore, we list those sets of composition factors for which an element of order $31$ is a blueprint for $M(E_7)$, which can be checked from an element $\hat x$ of order $93$ such that $\hat x^3=x$ using the preimage trick from the end of Section \ref{sec:blueprints}. (Elements of order $93$ are blueprints for $M(E_7)$ by Theorem \ref{thm:blueprintsminimal}.) In addition, we list those sets of factors for which there exists an element $\hat x$ of order $93$ in $\bG$, cubing to $x$, and such that $\hat x$ has one more distinct eigenvalue on $M(E_7)$ than $x$ (but there is not one with the same number of distinct eigenvalues). This last condition does not ensure that $H$ is a blueprint for $M(E_7)$, but does show that $H$ lies inside a positive-dimensional subgroup of $\bG$ stabilizing every simple submodule of $\soc(M(E_7)\downarrow_H)$ not of dimension $1$ or $32$.

To see this, if $\hat x$ has one more eigenvalue than $x$ then, since $\hat x$ must be real as it lies in $E_7$, all eigenspaces are preserved except for the $1$-eigenspace. Only the trivial and $32$-dimensional have $1$ as an eigenvalue for $x$. Also, $\hat x$ is a blueprint for $M(E_7)$, so let $\bX$ be a positive-dimensional subgroup of $\bG$ stabilizing the same subspaces of $M(E_7)$ as $\hat x$. We see that $\bX$ stabilizes any simple submodule of $M(E_7)\downarrow_H$ not of dimension $1$ or $32$. If there is such a submodule, we then apply Proposition \ref{prop:intersectionstabilizers} to see that $H$ is strongly imprimitive. If there is not, then $H$ either stabilizes a line on $M(E_7)$, contrary to our assumption, or it stabilizes only a $32$-space, but this is impossible as $\dim(M(E_7))=56$ and all $32$-dimensional factors -- Steinberg modules -- are summands.
\begin{center}
\begin{tabular}{ccccc}
\hline Case & Number & Consp. & 31 is blueprint & One more eigenvalue
\\ \hline
$4^6,2^{15},1^2$ & 3879876 & 5 & 5 & 0
\\$8,4^5,2^{13},1^2$ & 9529520 & 2 & 2 &0
\\$8^2,4^4,2^{11},1^2$ & 10735725 & 13 & 12 & 1
\\ $4^9,2^9,1^2$ & 6952660 & 16 & 16 & 0
\\ $8,4^8,2^7,1^2$ & 16044600 & 30 & 23 & 0
\\ $16,4^6,2^7,1^2$ & 1651650 & 9 & 3 & 1
\\ $8^2,4^7,2^5,1^2$ & 15855840 & 54 & 29 & 12
\\ $16,8,4^5,2^5,1^2$ & 2522520 & 24 & 10 & 3
\\ $16^2,4^3,2^5,1^2$ & 83160 & 6 & 6 & 0
\\ $8^3,4^6,2^3,1^2$ & 7707700 & 22 & 5 & 9
\\ $16,8^2,4^4,2^3,1^2$ & 1376375 & 19 & 14 & 3
\\ $32,4^4,2^3,1^2$ & 5005 & 1 & 0 & 0
\\ $16^2,8,4^2,2^3,1^2$ & 57750 & 3 & 1 & 1
\\\hline\end{tabular}
\end{center}
Excluding both those that are blueprints and where there is an element with one more eigenvalue on $M(E_7)$, we are left with 48 conspicuous sets of composition factors. 22 of these 48 have no corresponding set of composition factors on $L(E_7)^\circ$, or one with pressure at most $5$. In the first case $H$ cannot embed in $\bG$, and in the second $H$ stabilizes a line on $L(E_7)^\circ$ by Lemma \ref{lem:char2pressuregeneral}, so is strongly imprimitive by Proposition \ref{prop:fixlineonLG}.

We are left with 26 conspicuous sets of composition factors for $M(E_7)\downarrow_H$, still too many to list. Just as in the $a=4$ case, let $W$ be the subquotient obtained from $M(E_7)\downarrow_H$ by quotienting out by the $\{8_{i,j,l},16_{i,j,l,m},32\}$-radical and taking the $\{8_{i,j,l},16_{i,j,l,m},32\}$-residual, and remove any $4$-dimensional simple summands. Since $H$ can be assumed not to stabilize a line or $2$-space on $M(E_7)$, the socle of $W$ consists of $4$-dimensional modules, and the factors of $\soc(W)$ consist of $4$-dimensional simple modules that occur with multiplicity at least $2$ in $M(E_7)\downarrow_H$, and hence $W$. Let $S_1,\dots,S_r$ be the $4$-dimensional simple modules that appear in $M(E_7)\downarrow_H$ with multiplicity at least $2$. (Note that no composition factor of $M(E_7)\downarrow_H$, in the twenty-eight remaining sets of factors, appears with multiplicity greater than $3$, so we need only one copy of each $S_i$.)

We construct the largest submodule $W'$ of $P(S_1\oplus \cdots\oplus S_r)$ that consists solely of composition factors from $M(E_7)\downarrow_H$; certainly $W\leq W'$. Thus $W'$ must have at least two trivial factors, and all the requisite $2$-dimensional factors. In fact, only ten out of the 28 cases yield modules $W'$ with any trivial factors, with two even being the zero module (as there are no such $S_i$). Another seven can be removed for not having the correct $2$-dimensional factors, leaving the following three sets of factors:
\[ 8_{1,3,5},8_{1,4,5},4_{1,2},4_{1,3}^2,4_{1,4},4_{1,5}^2,4_{2,3},2_1^2,2_2,2_3^2,1^2,\]
\[8_{1,2,4},8_{1,3,5},4_{1,2},4_{1,3},4_{1,5}^2,4_{2,4}^2,4_{3,4},2_1^2,2_2,2_4^2,1^2,\]
\[ 16_{1,3,4,5},8_{1,3,4},8_{1,4,5},4_{1,2},4_{1,4},4_{1,5}^2,2_1^2,2_2,1^2.\]
In these final three cases we need to consider a preimage $\hat x$ that does not stabilize all eigenspaces on $M(E_7)$, but does stabilize those that make up some submodule of $M(E_7)\downarrow_H$. Let $\zeta$ be a primitive $31$st root of unity, chosen so that $x$ acts with eigenvalues $\zeta^{\pm 1}$ on $2_1$. In all three cases, $x$ has $31$ eigenvalues on $M(E_7)$.

In the first case, the fewest number of eigenvalues for a preimage $\hat x$ of order $93$ is $35$, with the four eigenvalues of $x$ not being stabilized being $\zeta^{\pm 14}$, $\zeta^{\pm 15}$. In the second case, $\hat x$ can take $34$ eigenvalues, with the three eigenvalues of $x$ not being stabilized being $1$ and $\zeta^{\pm 11}$ (there are four options for $\hat x$, two with this property). In the third case, the fewest number of eigenvalues for $\hat x$ is $34$, with the three eigenvalues of $x$ not stabilized by $\hat x$ being $1,\zeta^{\pm 2}$ (there are four options for $\hat x$, two with this property).

The eigenvalues of $x$ on the $4$-dimensional modules in the sets above are as follows:
\begin{center}
\begin{tabular}{cc}
\hline Module & Eigenvalues
\\ \hline $4_{1,2}$ & $\zeta^{\pm 1},\zeta^{\pm 3}$
\\ $4_{1,3}$ & $\zeta^{\pm 3},\zeta^{\pm 5}$
\\ $4_{1,4}$ & $\zeta^{\pm 7},\zeta^{\pm 9}$
\\ $4_{1,5}$ & $\zeta^{\pm 14},\zeta^{\pm 15}$
\\ $4_{2,3}$ & $\zeta^{\pm 2},\zeta^{\pm 6}$
\\ $4_{2,4}$ & $\zeta^{\pm 6},\zeta^{\pm 10}$
\\ \hline
\end{tabular}
\end{center}
In the second and third cases all simple $4$-dimensional submodules are stabilized by $\hat x$ (having chosen the correct one), hence a positive-dimensional subgroup, so we just have to show that one exists, and then apply Proposition \ref{prop:intersectionstabilizers}. However, there are no $1$- and $2$-dimensional submodules by assumption. The $8$- and $16$-dimensional factors appear with multiplicity $1$, so if they lie in the socle then they are summands. Thus there must be a $4$-dimensional submodule, and we are done.

In the first case, all simple $4$-dimensional submodules other than a copy of $4_{1,5}$ are preserved by $\hat x$, so we need to find an $N_{\Aut^+(\bG)}(H)$-stable collection of $4$-dimensional submodules that avoids $4_{1,5}$. There are no extensions between $4_{1,2}$ and any of $4_{1,3}$, $4_{1,5}$, $2_1$ or $2_3$, and so since all other composition factors are multiplicity free, $4_{1,2}$ must split off as a summand. The modules $4_{1,2}$ and $4_{1,5}$ do not lie in the same $\Aut(H)$-orbit of simple modules, so the $N_{\Aut^+(\bG)}(H)$-orbit of this summand must consist entirely of $4$-dimensional submodules of $M(E_7)\downarrow_H$ stabilized by $\hat x$. Thus we may apply Proposition \ref{prop:intersectionstabilizers} again, and this completes the proof.

\medskip

\noindent \textbf{Case $a=6$}: We have exactly the same possible dimensions for composition factors for $M(E_7)\downarrow_H$ as for $a=5$. The traces of semisimple elements of $\bG$ of order up to $21$ are known, but not $63$ or $65$, so we can check if a set of composition factors are conspicuous for elements of order up to $21$. Letting $x$ be an element of order $63$ in $H$, we use the preimage trick from the end of Section \ref{sec:blueprints} first to see if the composition factors are conspicuous up to $63$, and then use the preimage trick again to see if there exists an element $\hat x$ of order $63\cdot 5=195$ with the same eigenspaces as $x$ and with $\hat x^5=x$. If this is the case, then since $\hat x$ is a blueprint for $M(E_7)$ by Theorem \ref{thm:blueprintsminimal}, $H$ is a blueprint for $M(E_7)$, as needed.
\begin{center}
\begin{tabular}{ccccc}
\hline Case & Number & Consp.\ to 21 & Consp.\ to 63 & 63 is blueprint
\\ \hline
$4^6,2^{15},1^2$ & 100155870 & 6 & 6 & 6
\\$8,4^5,2^{13},1^2$ & 332095680 & 22 & 3 & 3
\\$8^2,4^4,2^{11},1^2$ & 467812800 & 60 & 18 & 18
\\ $4^9,2^9,1^2$ & 272669110 & 164 & 21 & 21
\\ $8,4^8,2^7,1^2$ & 844192800 & 1201 & 40 & 40
\\ $16,4^6,2^7,1^2$ & 76744800 & 254 & 16 & 16
\\ $8^2,4^7,2^5,1^2$ & 1025589600 & 3079 & 93 & 93
\\ $16,8,4^5,2^5,1^2$ & 146512800 & 1203 & 59 & 59
\\ $16^2,4^3,2^5,1^2$ & 3427200 & 53 & 20 & 20
\\ $8^3,4^6,2^3,1^2$ & 557110500 & 2665 & 54 & 54
\\ $16,8^2,4^4,2^3,1^2$ & 89964000 & 996 & 63 & 63
\\ $32,4^4,2^3,1^2$ & 171360 & 14 & 5 & 5
\\ $16^2,8,4^2,2^3,1^2$ & 2688000 & 58 & 18 & 18
\\\hline\end{tabular}
\end{center}
In every case, we find that the element of order $63$ is a blueprint for $M(E_7)$. This completes the proof for $a=6$.\end{proof}

\begin{proposition}\label{prop:e7char24trivials} Suppose that $M(E_7)\downarrow_H$ has exactly four trivial composition factors.
\begin{enumerate}
\item If $a=3$ then $M(E_7)\downarrow_H$ has a submodule of dimension at most $2$.
\item If $a=4$ then $H$ is strongly imprimitive.
\item If $a=5,6$ then $H$ is either a blueprint for $M(E_7)$ or stabilizes a subspace of dimension at most $2$ of $M(E_7)$.
\end{enumerate}
\end{proposition}
\begin{proof} \textbf{Case $a=3$}: The proof is the same as for Proposition \ref{prop:e7char22trivials}. In this case, the only option is $P(4)^2,8^3$, as there are four trivial factors. Thus the composition factors of $M(E_7)\downarrow_H$ have dimensions $8^3,4^4,2^6,1^4$, on which an element of order $3$ acts with trace $-1$. From Appendix \ref{app:traces} we see that this is not allowed, completing the proof.

\medskip

\noindent \textbf{Case $a=4$}: If $H$ stabilizes a subspace of dimension $1$ or $2$ on $M(E_7)$, then $H$ is strongly imprimitive by Propositions \ref{prop:fixlineonMG} and \ref{prop:fix2spaceonMG} respectively. Thus we assume that this is not the case. We may also assume that $H$ is not a blueprint for $M(E_7)$ by Proposition \ref{prop:blueprintissi}. In addition $H$ cannot stabilize a line on $L(E_7)^\circ$, else it is strongly imprimitive by Proposition \ref{prop:fixlineonLG}.

Using all semisimple elements, there are (up to field automorphism) $114$ conspicuous sets of composition factors for $M(E_7)\downarrow_H$ with exactly four trivial composition factors. Only $94$ of these have corresponding factors on $L(E_7)^\circ$, and of these only $81$ have either no $2_i$ or positive $2_i$-pressure for every $i$. One can eliminate nine more as they have no $4$-dimensional factors appearing with multiplicity greater than $1$, so must stabilize either a line or $2$-space as $M(E_7)$ is self-dual. We also exclude those whose corresponding factors on $L(E_7)^\circ$ have pressure less than $6$ (and therefore stabilize a line on $L(E_7)^\circ$ using Lemma  \ref{lem:char2pressuregeneral} and the fact that $u$ must act with at least six blocks of size $1$ on $L(E_7)^\circ$ from \cite[Table 8]{lawther1995}). This leaves $50$ conspicuous sets of composition factors for $M(E_7)\downarrow_H$.

Let $x$ be an element of order $17$ in $H$. As with the proof of Proposition \ref{prop:e7char22trivials}, we may exclude those sets of composition factors that are blueprints or near blueprints for $M(E_7)$, as in these cases $H$ is strongly imprimitive. This time there are eleven such sets of factors, bringing us down to $39$ that still need to be checked.

In the same proof, we introduced the subquotient $W$ of $M(E_7)\downarrow_H$, which is the $\{8_{i,j,l},16\}$-residual of $M(E_7)\downarrow_H$ modulo its $\{8_{i,j,l},16\}$-radical, with all $4$-dimensional summands removed. We also constructed a module $W'$, which is the $\cf(M(E_7)\downarrow_H)$-radical of $P(S_1\oplus \cdots S_r)$, where $S_1,\dots,S_r$ are all $4$-dimensional composition factors of $M(E_7)\downarrow_H$ appearing with multiplicity the floor of half of their multiplicity in $M(E_7)\downarrow_H$. We also removed all quotients from $W'$ that are not isomorphic to one of the $S_i$. By construction $W\leq W'$, so in particular the number of copies of $1$ and $2_i$ in $W'$ are at least those in $W$ (and hence $M(E_7)\downarrow_H$).

Computing the module $W'$ for all $39$ sets of composition factors yields eighteen sets where $W'$ does not contain enough factors $1$ and $2_i$. But still $21$ remain, which is too many to list.

As we saw when considering the case with two trivial factors in Proposition \ref{prop:e7char22trivials}, construction of the module $W'$ does not take into account that if a $4$-dimensional factor lies in the socle of $W$ and has multiplicity exactly $2$ in $M(E_7)\downarrow_H$ then it cannot appear anywhere other than the socle or the top of $W$. Including this, and ranging over all possible socles rather than just the largest one, yields a collection of modules for each case, all smaller than the original $W'$, and another thirteen that no longer have enough $1$- or $2$-dimensional factors, bringing us down to eight. The last eight cases are as follows:
\[ 4_{1,3}^3,4_{1,4},4_{2,3}^3,4_{2,4},2_1^4,2_2^2,2_3^3,2_4,1^4,\quad 4_{1,2},4_{1,3}^3,4_{1,4}^2,4_{2,3},4_{2,4},2_1^4,2_2^2,2_3^3,2_4,1^4,\]
\[8_{1,3,4},4_{1,3}^2,4_{1,4}^2,4_{2,3},4_{3,4}^2,2_1^3,2_2,2_3^2,2_4^2,1^4,\quad 8_{1,3,4},4_{1,2},4_{1,3}^2,4_{1,4}^3,4_{2,3},2_1^4,2_2^2,2_3^2,1^4,\]
\[8_{1,3,4},8_{2,3,4},4_{1,4}^3,4_{2,4},4_{3,4}^2,2_1^3,2_2,2_4^2,1^4,\quad 8_{1,3,4}^2,4_{1,4}^3,4_{2,4},4_{3,4}^2,2_1^3,2_2,2_4^2,1^4,\]
\[8_{1,2,3},8_{1,3,4},4_{1,2},4_{1,3},4_{1,4}^2,4_{2,3}^2,2_1^2,2_2,2_3^2,2_4,1^4,\] \[8_{1,3,4}^2,4_{1,3},4_{1,4}^3,4_{2,3},4_{3,4},2_1^3,2_2^2,2_3,1^4.\]

\noindent \textbf{Case 1}: The socle of $W$ can be either $4_{1,3}$ or $4_{1,3}\oplus 4_{2,3}$. If the socle of $W$ is $4_{1,3}$ then the module $W'$ in which $W$ can be found is
\[\begin{array}{c} 4_{1,3}
\\2_1\;2_3
\\1\;4_{1,4}\;4_{2,3}
\\2_1\;2_2\;2_3\;2_4
\\1\;1\;4_{1,3}\;4_{1,3}\;4_{2,4}
\\2_1\;2_2\;2_3\;2_4
\\1\;4_{1,4}\;4_{2,3}
\\2_1\;2_3
\\4_{1,3}\end{array}
\]
This is self-dual, so has a simple top, and since it is $64$-dimensional, $W$ must be contained in $\rad(W')$, and indeed in the $\{4_{1,3},4_{2,3}\}'$-residual of this, which is
\[ 4_{2,3}/2_3/1,4_{1,3},4_{1,3}/2_1,2_3,2_4/1,4_{1,4},4_{2,3}/2_1,2_3/4_{1,3},\]
which has no $2_2$, so $4_{1,3}$ cannot be the socle. If $4_{1,3}\oplus 4_{2,3}$ is the socle, then the module $W'$ is the sum of the one above and
\[ 4_{1,3},4_{2,3}/2_1,2_3/1,4_{1,4}/2_1,2_4/1,4_{1,3}/2_3/4_{2,3},\]
which also has no $2_2$. The same statement about the top $4_{1,3}$ not appearing in $W$ remains true, and so we take the same residual (this is why we took the $\{4_{1,3},4_{2,3}\}'$-residual rather than the $\{4_{1,3}\}'$-residual above) and see no $2_2$ again. Thus $H$ must stabilize a $1$- or $2$-space on $M(E_7)$.

\medskip

\noindent \textbf{Case 2}: The socle of $W'$ must be $4_{1,3}$, and indeed $W'$ is the same module as in the previous case, so the same method works there.

\medskip

\noindent \textbf{Cases 4, 8}: The module $W'$ is the self-dual module
\[\begin{array}{c}
4_{1,4}
\\2_1\;8_{1,3,4}
\\1\;4_{1,3}\;4_{1,4}
\\2_1\;2_2\;2_3
\\1\;1\;4_{2,3}
\\2_1\;2_2\;2_3
\\1\;4_{1,3}\;4_{1,4}
\\2_1\;8_{1,3,4}
\\4_{1,4}
\end{array}\]
which has two copies of $8_{1,3,4}$, so as in the first two cases we can take the $\{4_{1,4}\}'$-residual of $\rad(W')$ to obtain a module
\[ 4_{1,4}/2_1/1/2_2/1,4_{1,4}/2_1,8_{1,3,4}/4_{1,4},\]
which cannot work for several reasons, so that $H$ stabilizes a line or $2$-space on $M(E_7)$. The exact same module appears as $W'$ in the eighth case as well, so this method works there.

\medskip

We have therefore eliminated the first, second, fourth and eighth cases, and will look at semisimple elements in the third, fifth, sixth and seventh cases.

\medskip

\noindent \textbf{Case 3}: The element $x$ acts on $M(E_7)$ with eigenvalues
\[ 1^4,(\zeta^{\pm 1})^3,(\zeta^{\pm 2})^2,(\zeta^{\pm 3})^3,(\zeta^{\pm 4})^5,(\zeta^{\pm 5})^5,(\zeta^{\pm 6})^2,(\zeta^{\pm 7})^2,(\zeta^{\pm 8})^4,\]
and there exists an element $\hat x$ of order $85$ in $\bG$ that powers to $x$ and has nineteen distinct eigenvalues on $M(E_7)$, only splitting the $\zeta^{\pm 2}$-eigenspaces. In $M(E_7)\downarrow_H$, these lie in $2_2$ and $4_{2,3}$, the latter of which can only lie in the socle if it is a summand. Hence every other simple submodule is preserved by an element of order $85$, and therefore a positive-dimensional subgroup of $\bG$ by Theorem \ref{thm:blueprintsminimal}. Since there must be a submodule that is not a summand, there is an $N_{\Aut^+(\bG)}(H)$-orbit of simple submodules whose stabilizer is positive dimensional. Thus $H$ is strongly imprimitive by Proposition \ref{prop:intersectionstabilizers}.

\medskip

\noindent \textbf{Case 5}: This time, we find an element $\hat x$ of order $85$ that powers to $x$ and only disturbs the $1$- and $\zeta^{\pm 1}$-eigenspaces. Since these only lie in the trivial and $2_1$, the simultaneous stabilizer of every simple submodule of $M(E_7)\downarrow_H$ contains $\hat x$, and hence is a positive-dimensional subgroup of $\bG$ by Theorem \ref{thm:blueprintsminimal}. In particular, $H$ is strongly imprimitive by Proposition \ref{prop:intersectionstabilizers}.

\medskip

\noindent \textbf{Case 6}: There are eight elements of order $85$ in $\bG$ that power to $x$ and have nineteen eigenvalues on $M(E_7)$: four split the $\zeta^{\pm 6}$-eigenspace and the other four split the $\zeta^{\pm 8}$-eigenspace. The $\zeta^{\pm 6}$-eigenspace is contributed to by $4_{2,4}$ and $8_{1,3,4}$ from $M(E_7)\downarrow_H$, and the $\zeta^{\pm 8}$-eigenspace is contributed to by $2_4$ and $4_{1,4}$. Thus by Theorem \ref{thm:blueprintsminimal}, the simultaneous stabilizer of all simple submodules of $M(E_7)\downarrow_H$ not isomorphic to those on the first list is positive dimensional, and the same holds for the second list.

If there are $8$-dimensional simple modules in the socle, then the simultaneous stabilizer of all of these is positive dimensional, and hence $H$ is strongly imprimitive by Proposition \ref{prop:intersectionstabilizers}. Thus the socle of $M(E_7)\downarrow_H$ consists entirely of $4$-dimensional modules. Since $4_{2,4}$ must be a summand if it is a submodule, the simultaneous stabilizer of the collection of all simple submodules of $M(E_7)\downarrow_H$ that are not summands is positive dimensional, and of course is $N_{\Aut^+(\bG)}(H)$-stable. Thus we apply Proposition \ref{prop:intersectionstabilizers} again to see that $H$ is strongly imprimitive.

\medskip

\noindent \textbf{Case 7}: By examining the possibilities for $W'$, we see that $W'$ must have $4_{1,4}\oplus 4_{2,3}$ in the socle, as the other cases cannot yield the enough copies of $1$ and $2_i$. Thus $\soc(M(E_7)\downarrow_H)$ must contain $4_{1,4}\oplus 4_{2,3}$, as the $8$-dimensional factors must be summands if they are submodules. This also means that $4_{1,4}$ cannot be a summand of $M(E_7)\downarrow_H$, and therefore there is a unique submodule $V$ isomorphic to $4_{1,4}$.

We claim that the submodule $4_{1,4}$ is $N_{\Aut^+(\bG)}(H)$-invariant, so suppose the contrary. The $\Aut(H)$-orbit of $4_{1,4}$ contains $4_{1,2}$ and $4_{2,3}$. If $4_{1,2}$ is a submodule it is a summand, so $W$ cannot be sent to such a submodule by an element of $N_{\Aut^+(\bG)}(H)$. Thus we may assume that there is an element $\phi$ mapping $W$ to a submodule $4_{2,3}$. Any element of $\Aut(H)$ swapping $4_{1,4}$ and $4_{2,3}$ maps $4_{1,2}$ to $4_{3,4}$, but $4_{3,4}$ is not a composition factor of $M(E_7)\downarrow_H$, so $\phi$ cannot extend to a map on all of $M(E_7)$. This contradiction means that $V$ is $N_{\Aut^+(\bG)}(H)$-stable.

The smallest number of eigenvalues that an element $\hat x$ of order $85$ in $\bG$ powering to $x$ has on $M(E_7)$ is 23, and there is one that splits the $\zeta^{\pm 1}$, $\zeta^{\pm 2}$ and $\zeta^{\pm 5}$-eigenspaces. The element $x$ has eigenvalues $\zeta^{\pm 7},\zeta^{\pm 8}$ on $4_{1,4}$, so $\hat x$ stabilizes $V$. Thus by Theorem \ref{thm:blueprintsminimal}, the stabilizer of $V$ is positive dimensional, and is $N_{\Aut^+(\bG)}(H)$-stable by the above argument, so $H$ is strongly imprimitive by Proposition \ref{prop:intersectionstabilizers}.

\medskip

\noindent \textbf{Case $a=5$}: There are 28 possible sets of dimensions for the factors of $M(E_7)\downarrow_H$ that have a trace of an element of order $3$ that is one of $-25$, $-7$, $2$ and $25$ (see Appendix \ref{app:traces}). We exclude those of non-positive pressure, using Lemma \ref{lem:pressure} -- bringing us down to sixteen sets -- and those that do not have three $4$s as needed by Lemma \ref{lem:largestsubmodofP4}. We apply this lemma again to see that we cannot have more than four $2$-dimensional factors per $4$-dimensional factor (minus one $4$-dimensional factor), and this brings us down to six possible sets of dimensions, given in the table below.
\begin{center}
\begin{tabular}{cccc}
\hline Case & Number & Conspicuous & 31 is blueprint
\\ \hline $4^5,2^{16},1^4$ & 1939938 & 3 & 3
\\ $4^8,2^{10},1^4$ & 4866862 & 14 & 14
\\ $8,4^7,2^8,1^4$ & 11325600 & 30 & 29
\\ $16,4^5,2^8,1^4$ & 990990 & 7 & 3
\\ $8^2,4^6,2^6,1^4$ & 11561550 & 45 & 45
\\ $16,8,4^4,2^6,1^4$ & 1501500 & 19 & 18
\\\hline\end{tabular}
\end{center}

This leaves just six sets of composition factors (up to field automorphism) that are not guaranteed to be blueprints for $M(E_7)\downarrow_H$. These are
\[ 8_{1,4,5},4_{1,3}^2,4_{1,5}^2,4_{2,3},4_{3,5}^2,2_1^2,2_2^2,2_3^2,2_4^2,1^4,\]
\[16_{1,2,3,5},4_{1,3}^2,4_{2,3},4_{2,5}^2,2_1^2,2_2^2,2_3^2,2_4^2,1^4,\]
\[16_{1,2,3,4},4_{1,3},4_{1,4},4_{1,5},4_{3,4}^2,2_1^2,2_2^2,2_3,2_4^2,2_5,1^4,\]
\[16_{1,3,4,5},4_{1,4},4_{1,5},4_{2,3},4_{4,5}^2,2_1^3,2_2,2_3^2,2_5^2,1^4,\]
\[16_{1,2,4,5},4_{2,4}^2,4_{2,5},4_{3,4},4_{3,5},2_1^3,2_2,2_4^2,2_5^2,1^4,\]
\[16_{1,2,3,4},8_{1,3,4},4_{1,4},4_{2,3},4_{2,4},4_{3,4},2_1^2,2_2^2,2_3,2_5,1^4.\]
The easiest way to eliminate these is to consider the modules $W$ and $W'$ from the proof of Proposition \ref{prop:e7char22trivials}: in each of the six cases, we have at most two trivial factors in $W'$, and so we cannot have $W\leq W'$. Hence in these six cases $H$ must always stabilize a line or $2$-space on $M(E_7)$, as claimed in the proposition.

\medskip

\noindent \textbf{Case $a=6$}: We proceed in the same way as in Proposition \ref{prop:e7char22trivials} for $a=6$. We have exactly the same possible dimensions for composition factors for $M(E_7)\downarrow_H$ as for $a=5$. The traces of semisimple elements of $\bG$ of order up to $21$ are known, but not $63$ or $65$, so we can check if a set of composition factors are conspicuous for elements of order up to $21$. Letting $x$ be an element of order $63$ in $H$, we use the preimage trick from the end of Section \ref{sec:blueprints} first to see if the composition factors are conspicuous up to $63$, and then use the preimage trick again to see if there exists an element $\hat x$ of order $63\cdot 5=195$ with the same eigenspaces as $x$ and with $\hat x^5=x$. If this is the case, then since $\hat x$ is a blueprint for $M(E_7)$ by Theorem \ref{thm:blueprintsminimal}, $H$ is a blueprint for $M(E_7)$, as needed.
\begin{center}
\begin{tabular}{ccccc}
\hline Case & Number & Consp.\ to 21 & Consp.\ to 63 & 63 is blueprint
\\ \hline $4^5,2^{16},1^4$ & 39437442 & 4 & 4 & 4
\\ $4^8,2^{10},1^4$ & 160048350 & 170 & 19 & 19
\\ $8,4^7,2^8,1^4$ & 498841200 & 792 & 47 & 47
\\ $16,4^5,2^8,1^4$ & 37414170 & 61 & 12 & 12
\\ $8^2,4^6,2^6,1^4$ & 626754246 & 1484 & 85 & 85
\\ $16,8,4^4,2^6,1^4$ & 70686000 & 146 & 37 & 37
\\\hline\end{tabular}
\end{center}
This completes the proof for $a=6$.\end{proof}

We are left with $H$ having at least six trivial composition factors, where by the remarks at the start of this chapter we noted that if $a=5,6$ then $H$ is always a blueprint for $M(E_7)$.

\begin{proposition}\label{prop:e7char2>5trivials} Suppose that $a\geq 3$ and $M(E_7)\downarrow_H$ has at least six trivial composition factors. 
\begin{enumerate}
\item If $a=3$ then $M(E_7)\downarrow_H$ either has a $1$- or $2$-dimensional submodule or is
\[ 8\oplus P(4_{1,2})\oplus P(4_{2,3})\oplus P(4_{1,3}).\]
\item If $a=4$ then $H$ is a blueprint for $M(E_7)$ or $H$ stabilizes a subspace of dimension at most $2$ on $M(E_7)$.
\item If $a\geq 5$ then $H$ is a blueprint for $M(E_7)$.
\end{enumerate}
\end{proposition}
\begin{proof} \textbf{Case $a=3$}: We use the proof of the previous proposition to note that the only possibility is that $M(E_7)\downarrow_H$ is the sum of three projectives $P(4_{i,j})$ and a single summand $8$. We therefore consider the ten possible such modules, and note that only one has a conspicuous set of composition factors for $M(E_7)\downarrow_H$, the one mentioned. This completes the proof for $a=3$.

\medskip

\noindent \textbf{Case $a=4$}: As in previous chapters and earlier this chapter, we note that most classes of elements of order $17$ are blueprints for $M(E_7)$. In order to restrict the number of conspicuous sets of composition factors, we assume that an element of order $17$ is not a blueprint for $M(E_7)$. If $H$ stabilizes a $1$-space or $2$-space on $M(E_7)$ then we are done, so we assume that this is not the case either. By Lemma \ref{lem:pressure}, this means that $H$ has positive pressure on $M(E_7)$, and also there are two $4$-dimensional composition factors, else $H$ would stabilize a $2$-space on $M(E_7)$, by Lemma \ref{lem:extforsl22a}.

Remove any $8$s and $16$s in the top and socle of $M(E_7)\downarrow_H$, together with any simple summands of dimension $4$, leaving a self-dual module $W$ whose top and socle consist of $4$-dimensional modules, with $W$ having all trivial factors in $M(E_7)\downarrow_H$.

The projectives $P(4_{1,2})$ and $P(4_{1,3})$ both have exactly four trivial composition factors, and have dimension $64$. Therefore we cannot have the whole projective as a submodule of $W$. Thus we remove the simple top, then any $1$-, $2$- and $8$-dimensional modules from the top of each projective module (i.e., take the $\{1,2_i,8_{i,j,l}\}$-radical of $\rad(P(4_{\alpha,\beta}))$) to find the following modules:
\begin{equation}
\begin{gathered}4_{1,2},4_{2,4}/2_2,2_4/1,4_{3,4}/2_3,2_4/1,4_{1,2},4_{2,4}/2_2,8_{1,2,4}/4_{1,2};\\
4_{1,4},4_{2,3}/2_1,2_3/1,1,4_{1,3},4_{1,3},4_{2,4}/2_1,2_2,2_3,2_4/1,4_{1,4},4_{2,3}/2_1,2_3/4_{1,3}.
\end{gathered}
\label{eq:e7char2a=46trivsprojsubmodules}
\end{equation}
Thus $W$ is a submodule of a sum of these modules and their images under field automorphisms.

Since $M(E_7)\downarrow_H$ has at least six trivial composition factors, the socle of $W$ cannot be simple, and if it has only two factors they must both be $4_{1,3}$ or $4_{2,4}$. This means that we need either three $4$-dimensional factors appearing in both $\soc(W)$ and $\topp(W)$, or one of $4_{1,3}^4$, $4_{2,4}^4$ or $4_{1,3}^2,4_{2,4}^2$ as composition factors of $W$.

Using the traces of non-blueprint semisimple elements of order $17$, and traces of all elements of order $3$, $5$ and $15$, we end up with ten conspicuous sets of composition factors with at least six trivials, positive pressure, and at least two $4$s, up to field automorphism. These are
\[ 4_{1,3}^2,4_{1,4}^2,4_{2,3}^2,2_1^4,2_2^2,2_3^4,2_4^2,1^8,\quad
4_{1,3}^2,4_{1,4}^2,4_{2,3},4_{3,4}^2,2_1^4,2_2^3,2_3^2,2_4^2,1^6,\]
\[8_{2,3,4},4_{1,3}^2,4_{1,4}^2,4_{3,4}^2,2_1^4,2_2^3,2_4^2,1^6,\quad
4_{1,3}^3,4_{1,4},4_{2,3}^3,2_1^4,2_2^2,2_3^4,2_4,1^6,\]
\[8_{1,3,4},4_{1,3}^3,4_{1,4}^2,4_{2,4},2_1^4,2_2,2_3^2,2_4^2,1^6,\quad
8_{1,2,4},8_{2,3,4},4_{1,3}^3,4_{2,4}^2,2_1^4,2_2,2_3^2,1^6,\]
\[16,4_{1,3},4_{1,4},4_{3,4}^2,2_1^3,2_2^3,2_3,2_4^2,1^6,\quad
8_{2,3,4},4_{1,2},4_{1,3},4_{2,3}^2,4_{3,4}^2,2_1^3,2_2^2,2_3^2,2_4^2,1^6,\]
\[8_{1,3,4}^2,4_{1,3}^2,4_{1,4},4_{2,3},4_{2,4},2_1^3,2_3^2,2_4^2,1^6,\quad
8_{1,2,3}^2,4_{1,2},4_{1,3},4_{1,4},4_{2,3}^2,2_1^2,2_2^2,2_3^2,2_4,1^6. \]
By our restrictions on the multiplicities of $4$-dimensional factors of $W$ (and hence $M(E_7)\downarrow_H$), in all but the first, second, third and sixth cases $H$ must stabilize either a $1$-space or a $2$-space on $M(E_7)$, as needed. In the sixth case $H$ has pressure $1$, and $4_{1,3}\oplus 4_{2,4}$ is the socle of $W$. However, all six trivial factors in the sum of the modules from (\ref{eq:e7char2a=46trivsprojsubmodules}) must be present in $W$, so the third socle layer of $W$ has a submodule $1^{\oplus 2}$, contradicting Lemma \ref{lem:pressure}.

In the first, second and third cases, all $4$-dimensional factors that appear with multiplicity greater than $1$ must appear in the socle of $W$ and hence $M(E_7)\downarrow_H$. In the first case, $W$ possesses eight trivial composition factors, but $\soc(W)$ contains at most $4_{1,3}\oplus 4_{1,4}\oplus 4_{2,3}$, which can support only seven trivial factors, as we can see from the submodules of $P(4_{i,j})$ in (\ref{eq:e7char2a=46trivsprojsubmodules}). This yields a contradiction in the first case.

In the second and third cases, in order to obtain six trivial factors in $W$, the socle of $W$ must be $4_{1,3}\oplus 4_{1,4}\oplus 4_{3,4}$ and $4_{1,3}\oplus 4_{1,4}\oplus 4_{3,4}$ respectively. In particular, all of these composition factors can only appear in the socle and top of $W$. In the second case we take the preimages of the $\{1,2_i,4_{2,3}\}$-radicals of the quotient modules $P(4_{1,3})/4_{1,3}$, $P(4_{1,4})/4_{1,4}$ and $P(4_{3,4})/4_{3,4}$ to produce three modules in whose direct sum $\rad(W)$ is a submodule. These submodules are
\[1/2_2,2_3,2_4/1,4_{2,3}/2_1,2_3/4_{1,3},\quad
2_1/1/2_2/1/2_1/4_{1,4},\quad 
2_4/1/2_1/1/2_4/4_{3,4}.\]
These have six trivial composition factors, as does $W$, so all trivial factors in the above modules must occur in $W$. However, we obtain $W$ by adding on top only modules $4_{i,j}$, which have trivial $1$-cohomology by Lemma \ref{lem:extforsl22a}. Thus the trivial quotient of the module above must yield a trivial quotient of $W$, which is a contradiction. (The socle and top of $W$ consist entirely of $4$-dimensional modules, as stated before.)

In the third case we do the same thing, but with the $\{1,2_1,2_2,2_4,8_{2,3,4}\}$-radicals, to obtain
\[2_1/4_{1,3},\quad 2_1/1/2_2/1/2_1/4_{1,4},\quad
2_4/1/2_1/1/2_4,8_{2,3,4}/4_{3,4},\]
and clearly we have a contradiction here as there are not enough trivial factors in the sum of these modules.

\medskip

\noindent \textbf{Case $a\geq 5$}: This proof is easy, and was mentioned at the start of the chapter. By Proposition \ref{prop:f4a4}, if a semisimple element $x$ has order at least $31$ and at least a $6$-dimensional eigenspace on $M(E_7)$, then $x$ is a blueprint for $M(E_7)$. This clearly holds for $a\geq 5$, and so $H$ is a blueprint for $M(E_7)$ in these cases.
\end{proof}

We now give a summary of what we have proved. We are proving that unless $a=3$ and $M(E_7)\downarrow_H$ is
\[ 8\oplus P(4_{1,2})\oplus P(4_{1,3})\oplus P(4_{2,3}),\]
$H$ is always strongly imprimitive.

If $H$ stabilizes a line on $M(E_7)$ then $H$ is strongly imprimitive by Proposition \ref{prop:fixlineonMG}, if $H$ stabilizes a $2$-space on $M(E_7)$ then $H$ is again strongly imprimitive, this time by Proposition \ref{prop:fix2spaceonMG}. If $H$ is a blueprint for $M(E_7)$ then the same statement holds, by Proposition \ref{prop:blueprintissi}.

If $a\geq 7$ then $H$ is a blueprint for $M(E_7)$ by Theorem \ref{thm:goodblueprint}, so we may assume that $a\leq 6$. If $a=1$ then $H$ is soluble. If $a=2$ then by \cite[Proposition 5.4]{craven2015un2}, $H$ stabilizes either a line or $2$-space on $M(E_7)$, so $3\leq a\leq 5$.

The results depend on the composition factors of $M(E_7)\downarrow_H$.
\begin{itemize}
\item If $H$ has neither $1$- nor $2$-dimensional composition factors on $M(E_7)$, then $a=5,6$ and $H$ is a blueprint for $M(E_7)$, by Proposition \ref{prop:e7char2not1or2}.
\item If $H$ has $2$-dimensional factors but no $1$-dimensional factors on $M(E_7)$, then either $H$ is a blueprint for $M(E_7)$ or $H$ stabilizes a $2$-space on $M(E_7)$, by Proposition \ref{prop:e7char2no1s}.
\item If $H$ has exactly two trivial factors on $M(E_7)$, then in Proposition \ref{prop:e7char22trivials} we show directly that $H$ is strongly imprimitive for $a=4,5$, and for $a=3,6$ we show that either $H$ is a blueprint for $M(E_7)$ or $H$ stabilizes a $1$- or $2$-space on $M(E_7)$.
\item If $H$ has exactly four trivial factors on $M(E_7)$, then in Proposition \ref{prop:e7char24trivials} we show directly that $H$ is strongly imprimitive for $a=4$, and for $a=3,5,6$ we show that either $H$ is a blueprint for $M(E_7)$ or $H$ stabilizes a $1$- or $2$-space on $M(E_7)$.
\item If $H$ has at least six trivial factors on $M(E_7)$, then Proposition \ref{prop:e7char2>5trivials} shows that either $H$ is a blueprint for $M(E_7)$ or $H$ stabilizes a $1$- or $2$-space on $M(E_7)$, with the one exception given above for $a=3$. 
\end{itemize}
This proves that $H$ is always strongly imprimitive unless we are in the exceptional case above for $p^a=8$.

\chapter{The Proof for \texorpdfstring{$E_7$}{E7} in Odd Characteristic: \texorpdfstring{$\PSL_2$}{PSL2} Embedding}
\label{ch:e7oddpsl}

In this chapter, $k$ is an algebraically closed field of characteristic $p\geq 3$ and $\bG=E_7(k)$, by which we mean the simply connected form, i.e., $|Z(\bG)|=2$ and $\bG'=\bG$. Let $H\cong \PSL_2(p^a)$ be a subgroup of $\bG$.

Theorem \ref{thm:goodblueprint} states that if $p^a\geq 150$ then $H$ is a blueprint for $M(E_7)$. In this case, $H$ is strongly imprimitive by Proposition \ref{prop:blueprintissi}. Thus in what follows we may assume that $p^a\leq 150$.

Let $L=\PSL_2(p)\leq H$ and let $u$ denote a unipotent element of $L$ of order $p$. The possibilities for the Jordan block structures of $u$ on $M(E_6)$ and $L(E_6)$ are given in \cite[Tables 7 and 8]{lawther1995}. Recall the definition of a generic unipotent element from Definition \ref{defn:genericunipotent}.

By Proposition \ref{prop:f4a4}, if a semisimple element $x$ has order at least $31$ in $\bG$ and centralizes a $6$-space on $M(E_7)$, then $x$ is a blueprint for $M(E_7)$. Also, any semisimple element in $H$ has a $1$-dimensional $1$-eigenspace on every odd-dimensional simple module. Hence, if $H$ has at least six odd-dimensional composition factors on $M(E_7)$ and $p^a\geq 60$ (so that $H$ possesses an element of order at least $31$), then $H$ is a blueprint for $M(E_7)$. This normally ends up being the case for $p^a\geq 60$.

\section{Characteristic 3}

Let $p=3$, so that $H=\PSL_2(3^a)$ for some $a=2,3,4$. The case $a=2$ was considered in \cite[Proposition 6.2]{craven2015un2}, so we exclude this. If $a=4$ then we may assume that there are fewer than six odd-dimensional composition factors in $M(E_7)\downarrow_H$, by the discussion at the start of this chapter.

We begin by computing the composition factors of $M(E_7)\downarrow_L$, which depends only on the trace of an involution on $M(E_7)$, which is $\pm 8$ (its centralizer is $A_1D_6$, which acts with composition factors of dimension $24$ and $32$, which must be the $(+1)$- and $(-1)$-eigenspaces). This means that there are eight more of one factor than the other, so $3^{12},1^{20}$ and $3^{16},1^8$. From Lemma \ref{lem:sl23restriction} we can see the possible dimensions of composition factors for $M(E_7)\downarrow_H$: if $M(E_7)\downarrow_L$ has factors $3^{16},1^8$ then we must have at least eight $3$-dimensional factors in $M(E_7)\downarrow_H$, and if the factors are $3^{12},1^{20}$ then as only $9$ and $1$ for $H$ have more $1$s than $3$s on restriction to $L$, we need at least eight of these in $M(E_7)\downarrow_H$, and again have at least eight odd-dimensional composition factors in $M(E_7)\downarrow_H$. This gives us the first proposition.

\begin{proposition}\label{prop:e7char3a=4} Let $p=3$ and $a=4$. A semisimple element of order $41$ in $H$ is always a blueprint for $M(E_7)$, and hence $H$ is always a blueprint for $M(E_7)$.
\end{proposition}

We turn to $a=3$, where we cannot quite get the same result, but we come close.

\begin{proposition}\label{prop:e7char3a=3} Let $p=3$ and $a=3$. Either $H$ is a blueprint for $M(E_7)$ or $H$ stabilizes a line on either $M(E_7)$ or $L(E_7)$.
\end{proposition}
\begin{proof} As with $F_4$ and $E_6$, we want to discount conspicuous sets of composition factors where a semisimple element is a blueprint for $M(E_7)$. We have already seen in Proposition \ref{prop:f4char3a=3} that there are $97$ classes of semisimple elements of order $13$ that are blueprints for the minimal module for $F_4$, and there are exactly $188$ classes of semisimple elements of order $13$ in $E_7$ whose $1$-eigenspace is at least $8$-dimensional (which it must be by the discussion at the start of this section), leaving $91$ classes to which an element of order $13$ in $H$ can belong.

Using this, we find up to field automorphism eight conspicuous sets of composition factors, two of which have negative pressure so will not be displayed. The other six are 
\[ 9_{1,3},4_{1,3}^9,3_1,1^8,\qquad 4_{1,2}^4,4_{1,3}^4,4_{2,3}^4,1^8,\qquad 9_{2,3}^3,4_{1,2}^5,4_{1,3},1^5,\]\[ 4_{1,2}^6,3_1^9,3_2,1^2,\qquad 4_{1,2},4_{1,3}^5,4_{2,3},3_1^9,1,\qquad 9_{2,3},4_{1,2}^5,3_1^5,3_2^4.\]

\noindent \textbf{Cases 1, 2, 3}: The first and third cases have pressure $1$, so $H$ stabilizes a line on $M(E_7)$ by Lemma \ref{lem:psl227pressure1}. In the second case $H$ stabilizes a line on $M(E_7)$ by Lemma \ref{lem:no414}.

\medskip

\noindent \textbf{Case 4}: If $H$ does not stabilize a line on $M(E_7)$ then we by quotienting out any $3$s in the socle may assume that the socle consists of copies of $4_{1,2}$, and the $\{1,3_1,3_2,4_{1,2}\}$-radical of $P(4_{1,2})$ is
\[ 4_{1,2}/1,3_2/4_{1,2},\]
but since there is only one $3_2$ in $M(E_7)\downarrow_H$ we cannot cover both trivials in this way, thus $H$ stabilizes a line on $M(E_7)$. (Alternatively, the factors of $H$ on $L(E_7)$ are $4_{1,2}^{16},3_1^{10},3_2^6,1^{21}$, so $H$ stabilizes a line on $L(E_7)$.)

\medskip

\noindent \textbf{Case 5}: The corresponding set of composition factors on $L(E_7)$ (see (\ref{li:strata}) from Chapter \ref{ch:strategy}) is
\[4_{1,2}^5,4_{2,3}^5,4_{1,3}^{11},3_1^{10},3_3,1^{16}.\]
This module has pressure $5$, and so we cannot simply use Lemma \ref{lem:pressure} to find a trivial submodule. However, the largest submodule of $P(4_{i,j})$ with composition factors from among those of $L(E_7)\downarrow_H$ has three trivial composition factors for all pairs $i,j$, and so we need at least six $4$s in the socle of $L(E_7)\downarrow_H$ (once we remove all $3$s), contradicting the fact that the module has pressure $5$. Thus $H$ stabilizes a line on $L(E_7)$.

\medskip

\noindent \textbf{Case 6}: The $9$ splits off and we may quotient out by the $\{3_1,3_2\}$-radical to obtain a module with copies of $4_{1,2}$ in the socle. On this we can only place copies of $3_2$, and so $u$ would act on $M(E_7)$ as $3^{17},1^5$: $u$ acts on $9_{1,2}$ with blocks $3^3$ and $3_i$ with a single block $3$, and on $4_{1,2}$ with blocks $3,1$, so since there must be a subquotient $4_{1,2}^{\oplus 5}$ from the above radical, the action must be $3^{17},1^5$ on the whole module. However, this is not a valid unipotent action in \cite[Table 7]{lawther1995}, so $H$ cannot embed with these factors.
\end{proof}

\section{Characteristic At Least 5}

We now let $p\geq 5$, let $H=\PSL_2(p^a)$ with $a\geq 1$, let $L=\PSL_2(p)\leq H$ and let $u\in L$ have order $p$. We begin by producing a list of all unipotent classes to which $u$ can belong, excluding those that come from generic classes (see Lemma \ref{lem:genericunipotent}) and those that fail Lemma \ref{lem:pblockseven}. Moreover, we make a few remarks now about indecomposable modules for $L$, which can cut down our list.

When $p=5,13,17$, we use Corollary \ref{cor:blocks1mod4}, so for these primes the number of blocks of each even size is even. For $p=7,11,19,23$, there exists a unique self-dual indecomposable module for $L$ of dimension congruent to a given even number modulo $p$ by Lemma \ref{lem:selfdualsl2p}.

For $p=11$, the self-dual indecomposable module of dimension congruent to $6$ modulo $p$ has socle structure
\[ 1,3,5,7,9/1,3,5,7,9\]
and has dimension $50$. The trace of an involution on the module is $0$, and since involutions have trace $\pm 8$ on $M(E_7)$ (see Appendix \ref{app:traces}), we would need a trace of $\pm 8$ from the remaining factors of $M(E_7)\downarrow_L$, a module of dimension $6$, so not possible. Thus this is not a summand of $M(E_7)\downarrow_L$.

For $p=19$, the indecomposable module of dimension congruent to $6$ modulo $p$ has socle structure
\[ 5,7,9,11,13,15/5,7,9,11,13,15\]
and has dimension $120$, so cannot be a summand of $M(E_7)\downarrow_L$.

For $p=23$ the indecomposable module of dimension congruent to $10$ modulo $p$ has socle structure
\[ 3,5,7,9,11,13,15,17,19,21/3,5,7,9,11,13,15,17,19,21\]
and has dimension $240$, so cannot be a summand of $M(E_7)\downarrow_L$.

We now list the unipotent classes of interest using \cite[Table 7]{lawther1995}.

\begin{enumerate}
\item\label{li:pslfirst5} $A_3+A_2$, $p=5$, acting as $5^6,4^2,3^4,2^2,1^2$;
\item $A_4$, $p=5$, acting as $5^{10},1^6$;
\item $A_4+A_1$, $p=5$, acting as $5^{10},2^2,1^2$;
\item\label{li:psllast5} $A_4+A_2$, $p=5$, acting as $5^{10},3^2$;
\item\label{li:pslfirst7} $(A_5)''$, $p=7$, acting as $7^2,6^7$;
\item $D_4+A_1$, $p=7$, acting as $7^6,2^5,1^4$; 
\item $D_5(a_1)$, $p=7$, acting as $7^6,3^2,2^2,1^4$;
\item $(A_5)'$, $p=7$, acting as $7^4,6^4,1^4$;
\item $A_5+A_1$, $p=7$, acting as $7^4,6^3,5^2$;
\item $D_5(a_1)+A_1$, $p=7$, acting as $7^6,4,2^5$;
\item $D_6(a_2)$, $p=7$, acting as $7^6,5^2,4$;
\item $E_6(a_3)$, $p=7$, acting as $7^6,5^2,1^4$;
\item\label{li:psllastbutone7} $E_7(a_5)$, $p=7$, acting as $7^6,6,4^2$;
\item\label{li:psllast7} $A_6$, $p=7$, acting as $7^8$;
\item\label{li:pslfirst11} $D_6$, $p=11$, acting as $11^4,10,1^2$;
\item\label{li:pslsecond11} $E_6(a_1)$, $p=11$, acting as $11^4,5^2,1^2$;
\item\label{li:psllast11} $E_7(a_3)$, $p=11$, acting as $11^4,10,2$;
\item\label{li:pslfirst13} $E_6$, $p=13$, acting as $13^4,1^4$;
\item\label{li:pslfirst19} $E_7$, $p=19$, acting as $19^2,18$.
\end{enumerate}

We start with $p=5$. Because of the small number of possible sets of factors for $M(E_7)\downarrow_L$, we do not need to use the list of unipotent classes above for this prime.

\begin{proposition}\label{prop:e7char5} Let $p=5$ and $a\geq 1$.
\begin{enumerate}
\item\label{propi:e75a} If $a=1,2$ then $H$ stabilizes a line on either $M(E_7)$ or $L(E_7)$.
\item If $a=3$ then $H$ is a blueprint for $M(E_7)$.
\end{enumerate}
\end{proposition}
\begin{proof} \textbf{Case $a=1$}: The traces of elements of orders $2$ and $3$ yield conspicuous sets of composition factors for $H$ of
\[ 3^{12},1^{20},\qquad 5^2,3^{14},1^4,\qquad 5^9,3^3,1^2,\qquad 5^6,3^6,1^8.\]
As we saw in the proof of Proposition \ref{prop:e6char5}, only $P(3)=3/1,3/3$ has a trivial composition factor and no trivial submodule or quotient, and so in the first, third and fourth cases $H$ stabilizes a line on $M(E_7)$. However, in the third case this means that $H$ cannot embed with these factors at all: as it stabilizes a line on $M(E_7)$, from Lemma \ref{lem:e7stabs} we see that $H$ lies in either an $E_6$-parabolic subgroup, with factors $1,1,27,27^*$ or a $B_5$-subgroup, factors $1,1,11^2,32$, neither of which is compatible with $5^9,3^3,1^2$, so this case cannot occur. (In particular, it cannot occur for the subgroup $L$ when $a\geq 2$.)

For the remaining case of $5^2,3^{14},1^4$, we switch to the Lie algebra (see (\ref{li:strata}) from Chapter \ref{ch:strategy}). There are two possibilities for the corresponding composition factors of $L(E_7)$ (since an element of order $3$ with trace $2$ on $M(E_7)$ can have trace either $-2$ or $7$ on $L(E_7)$): $5^{10},3^{22},1^{17}$ and $5^{13},3^{19},1^{11}$. Both of these must have trivial submodules as again we can only cover a $1$ by $3/1,3/3$. This proves (\ref{propi:e75a}) for $a=1$.

\medskip

\noindent \textbf{Case $a=2$}: Recall from Lemma \ref{lem:cohomologysimple} that the only simple modules for $H$ with non-trivial $1$-cohomology when $a=2$ have dimension $8$ and restrict to $L\cong \PSL_2(5)$ as $5\oplus 3$ by Lemma \ref{lem:sl25restriction}. If $M(E_7)\downarrow_L$ has composition factors $3^{12},1^{20}$ then $M(E_7)\downarrow_H$ has at least eight trivial composition factors and can have no factors of dimension $8$, so $M(E_7)\downarrow_H$ has eight trivial summands proving (\ref{propi:e75a}) for this set of composition factors for $M(E_7)\downarrow_L$. Similarly, if the composition factors of $M(E_7)\downarrow_L$ are $5^6,3^6,1^8$ then $M(E_7)\downarrow_H$ must have at least two trivial composition factors by Lemma \ref{lem:sl25restriction}, and for every composition factor of dimension $8$ we must have another trivial factor, so $M(E_7)\downarrow_H$ always has pressure at most $-2$, so (\ref{propi:e75a}) holds when $M(E_7)\downarrow_L$ has this set of composition factors.

To finish the proof of (\ref{propi:e75a}) for $a=2$, we thus may assume that $M(E_7)\downarrow_L$ has factors $5^2,3^{14},1^4$. We have at most two $8$s in $M(E_7)\downarrow_H$ since there are only two $5$s in $M(E_7)\downarrow_L$, and hence there can be at most a single trivial composition factor in $M(E_7)\downarrow_H$, else $H$ stabilizes a line on $M(E_7)$. We thus get two cases: there is a trivial composition factor in $M(E_7)\downarrow_H$ and there is not.

If there is a trivial factor then we have $8^2,1$ in $M(E_7)\downarrow_H$, and the remaining factors of $M(E_7)\downarrow_H$ restrict to $L$ as $3^{12},1^3$, so we need factors $8^2,4^3,3^9,1$. For $a=2$, there is no such set of composition factors, so we cannot have a trivial composition factor in $M(E_7)\downarrow_H$.

There are up to field automorphism five conspicuous sets of composition factors for $M(E_7)\downarrow_H$ with the correct restriction to $L$ and no trivial composition factors:
\[ 5_1^2,4^4,3_1^{10},\quad 8_{2,1},5_1,4^4,3_1^9,\quad 8_{2,1}^2,4^4,3_1^7,3_2,\quad 8_{2,1}^2,4^4,3_1^6,3_2^2,\quad 
15_{2,1}^2,4^2,3_1^4,3_2^2.\]
The fact that $4$ has an extension with $3_1$ and $3_2$ (see Lemma \ref{lem:psl2125smallmodules}) makes deducing the module structure difficult, and so we turn to the Lie algebra $L(E_7)$ in all cases (see (\ref{li:strata}) from Chapter \ref{ch:strategy}). These are
\[ 8_{2,1}^3,5_1^{10},4^5,3_1^{11},1^6,\qquad 8_{2,1}^5,5_1^5,4^6,3_1^{10},3_2,1^{11},\qquad 9,8_{2,1}^7,5_1^5,4^3,3_1^5,3_2^3,1^7,\]
\[ 9^3,8_{2,1}^4,5_1^3,4^8,3_1^6,3_2,1^6,\qquad 15_{2,1}^2,9_{1,2}^4,8_{1,2}^2,8_{2,1},5_1^3,5_2,4^3,3_1^2,3_2,1^2:\]
each of these except the last has non-positive pressure, as needed. For the final case, which has pressure $1$, we take the $\cf(L(E_7)\downarrow_H)$-radical of $P(8_{2,1})$, which has structure
\[9/4,8_{1,2},8_{2,1}/1,3_1,9/8_{1,2}.\]
Since this can support only a single trivial factor, $H$ must again stabilize a line on $L(E_7)$, as claimed. This completes the proof of (\ref{propi:e75a}).

\medskip

\noindent \textbf{Case $a=3$}: From Lemma \ref{lem:sl25restriction} we see the following facts: first, in any even-dimensional composition factor of $M(E_7)\downarrow_H$ there are the same number of $3$s as $5$s and $1$s combined on restriction to $L$, and second, in any odd-dimensional factor of $M(E_7)\downarrow_H$, the number of $5$s and $1$s combined is at most one more than the number of $3$s, on restriction to $L$. This means that if $M(E_7)\downarrow_L$ has factors $5^6,3^6,1^8$, there must be at least six odd-dimensional composition factors. Lemma \ref{lem:sl25restriction} easily shows that if the factors of $M(E_7)\downarrow_L$ are $3^{12},1^{20}$ then there must be at least eight trivial factors in $M(E_7)\downarrow_H$, and if we have $5^2,3^{14},1^4$ then we have at least six $3$s in $M(E_7)\downarrow_H$. In all cases we have at least six odd-dimensional composition factors, so an element of order $63$ in $H$ has a $1$-eigenspace of dimension at least $6$. This means that $H$ is a blueprint for $M(E_7)$ by Proposition \ref{prop:f4a4}, as needed for the result.
\end{proof}

Having completed $p=5$, we now move on to $p=7$. This time $p^a=7,49$ will need to be considered, but $p^a=343$ is above $2\cdot v(E_7)=150$.

\begin{proposition}\label{prop:e7char7} Suppose that $p=7$.
\begin{enumerate}
\item If $a=1$ then either $H$ stabilizes a line on $M(E_7)$ or $L(E_7)$, or the actions of $H$ on $M(E_7)$ and $L(E_7)$ are
\[ 7^{\oplus 4}\oplus P(3)^{\oplus 2}\qquad\text{and}\qquad 7^{\oplus 5}\oplus P(5)^{\oplus 6}\oplus P(3)\]
respectively.
\item If $a=2$ then either $H$ is a blueprint for $M(E_7)$ or $H$ stabilizes a line on $M(E_7)$.
\end{enumerate}
\end{proposition}
\begin{proof} We first compute the possible sets of composition factors for $M(E_7)\downarrow_H$ when $a=1$, using the traces of elements of orders $2$, $3$ and $4$. There are seven of these, given by
\[ 3^{12},1^{20},\quad 5^2,3^{14},1^4,\quad 5^6,3^6,1^8,\quad 7,5^9,3,1,\quad 7^2,5^6,3^2,1^6,\quad 7^4,5^2,3^6,\quad 7^6,1^{14}.\]

\noindent \textbf{Cases 1, 2, 3, 5}: As in the case of $p=5$, the only indecomposable module with a trivial composition factor but no trivial submodule or quotient is $P(5)=5/1,3/5$, so either $H$ stabilizes a line on $M(E_7)$ or we have twice as many $5$s as $1$s and as many $3$s as $1$s. Thus all but the fourth and sixth cases must stabilize lines on $M(E_7)$.

\medskip

\noindent\textbf{Case 4}: If the factors are $7,5^9,3,1$, then $H$ cannot stabilize a line or hyperplane on $M(E_7)$, since by Lemma \ref{lem:e7stabs} the line stabilizers for $M(E_7)$ are contained in either an $E_6$-parabolic subgroup -- composition factors $27,27^*,1^2$ -- or a subgroup $q^{1+32}B_5(q)\cdot (q-1)$ -- composition factors $32,11^2,1^2$ -- neither of which can work. As there are no self-extensions of the $5$, $M(E_7)\downarrow_H$ is
\[ 7\oplus P(5)\oplus 5^{\oplus 7},\]
with $u$ acting with Jordan blocks $7^3,5^7$, which is not in \cite[Table 7]{lawther1995}, so there does not exist an embedding of $H$ into $\bG$ with these factors.

\medskip

\noindent\textbf{Case 6}: We are left with $7^4,5^2,3^6$. Here, the lack of trivials means $M(E_7)\downarrow_H$ is a sum of modules of the form
\[ (3/5)\oplus (5/3),\qquad 3/3\qquad P(3)=3/3,5/3,\qquad (3/3,5)\oplus (3,5/3),\qquad 3,5/3,5,\]
as we saw in the example after Proposition \ref{prop:simplesl2p}; we also saw that $u$ acts on each of these with at most one Jordan block of size not equal to $7$. Therefore we need an even number of $1$s, $4$s and $7$s in the Jordan block structure of $u$, with at least four more $7$s than the $1$s, $2$s and $4$s combined. Examining the list above, we see only two examples of this, namely (\ref{li:psllastbutone7}) and (\ref{li:psllast7}). This yields the two possible embeddings $M(E_7)\downarrow_H$ to be
\[ 7^{\oplus 4}\oplus (3/3)\oplus (3/3,5)\oplus (3,5/3),\qquad\text{and}\qquad 7^{\oplus 4}\oplus P(3)^{\oplus 2}.\]
The traces of semisimple elements of orders $3$ and $4$ on $M(E_7)$ yield two possibilities each for the semisimple class, and so we get four possible sets of composition factors for $L(E_7)\downarrow_H$, namely
\[ 7^6,5^{10},3^{10},1^{11},\qquad 7^8,5^{10},3^6,1^9,\qquad 7^3,5^{13},3^{13},1^8,\qquad
7^5,5^{13},3^9,1^6.\]
Apart from the last one, each of these has enough trivials and not enough $5$s to ensure that $H$ stabilizes a line on $L(E_7)$, since $P(5)=5/1,3/5$ is the only indecomposable module for $H$ with a trivial factor but no trivial submodule or quotient (see the example after Proposition \ref{prop:simplesl2p}). In this case, we obtain the action given in the statement of the proposition, and we know that $u$ acts on $L(E_7)$ with blocks $7^{19}$, so lies in class $A_6$, which acts on $M(E_7)$ with blocks $7^8$.

We now remove the first possible action on $M(E_7)$, using the simple fact that for $p\geq 5$, the symmetric square of $M(E_7)$ is the sum of $L(E_7)$ and the $1463$-dimensional module $L(2\lambda_1)$, Lemma \ref{lem:relatingvminlg}. The symmetric square of the first module is a sum of projectives and
\[ (3,5/1,3,5)^{\oplus 2}\oplus (1,3,5/3,5)^{\oplus 2}\oplus 5\oplus 3\oplus 1.\]
Since $u$ comes from class $E_7(a_5)$, and this acts on $L(E_7)$ as $7^{17},5,3^3$, we must have two of the summands of dimension $17$ in $L(E_7)\downarrow_H$, hence $H$ stabilizes a line on $L(E_7)$. This completes the proof for $a=1$.

\medskip

Now let $a=2$, so that $H=\PSL_2(49)$, and recall that $L\leq H$ is a copy of $\PSL_2(7)$. At the start of this proof we gave the conspicuous sets of composition factors for $M(E_7)\downarrow_L$, and from Lemmas \ref{lem:cohomologysimple} and \ref{lem:sl27restriction} we see that the only simple modules for $H$ with non-trivial $1$-cohomology have dimension $12$ and restrict to $L$ as $7\oplus 5$, and only the trivial module for $H$ restricts to $L$ with more $1$s than $3$s. These two facts mean that if $M(E_7)\downarrow_L$ has factors the first, third, fifth and seventh cases then $M(E_7)\downarrow_H$ has trivial composition factors, and these are summands except for the fifth case, and there we have at least four trivials and at most two $12$s, so pressure at most $-2$. We proved for $a=1$ that the fourth case for $M(E_7)\downarrow_L$ does not occur, so $M(E_7)\downarrow_L$ has factors either $5^2,3^{14},1^4$ or $7^4,5^2,3^6$.

In the case of $5^2,3^{14},1^4$, from Lemma \ref{lem:sl27restriction}, apart from $3$, there are no simple modules for $H$ whose restriction to $L$ has more $3$s than other factors, and the composition factors of $M(E_7)\downarrow_H$ have dimensions $1$, $3$, $4$, $5$, $8$ and $9$. In particular, this means that $M(E_7)\downarrow_H$ has at least eight $3$-dimensional composition factors. By Lemma \ref{lem:psl249smallmodules}, of these modules only $8$s can have an extension with $3$s, with there being at most two of those, so the $3$-pressure is at least $6$. This means that $M(E_7)\downarrow_H$ has at least four $3$-dimensional summands, so the action of the unipotent element $u$ on $M(E_7)$ has at least four Jordan blocks of size $3$. There are no non-generic unipotent classes with this property, as we saw in the list at the start of this section, and so $H$ is a blueprint for $M(E_7)$.

Thus we end with $M(E_7)\downarrow_L$ being $7^4,5^2,3^6$, and Lemma \ref{lem:sl27restriction} implies that $H$ has at least two $7$s and four $3$s on $M(E_7)$. The remaining composition factors have dimension $3$, $5$, $7$, $8$, $12$ or $15$: using the traces of semisimple elements, we find exactly four conspicuous sets of composition factors for $M(E_7)\downarrow_H$ with the correct restriction to $L$, and for each of these the eigenvalues of an element of order $24$ on $M(E_7)$ determine its conjugacy class, and this lies inside $F_4$, hence the element is a blueprint for $M(E_7)$ by Lemma \ref{lem:f4blueprintisblueprint}. Thus $H$ is a blueprint for $M(E_7)$, as needed.
\end{proof}

\begin{proposition}\label{prop:e7char11} Suppose that $p=11$. If $a=1,2$ then either $H$ is a blueprint for $M(E_7)$ or $H$ stabilizes a line on $M(E_7)$.
\end{proposition}
\begin{proof} \textbf{Case $a=1$}: If the unipotent class of $\bG$ to which $u$ belongs is generic for $M(E_7)$ then $H$ is a blueprint for $M(E_7)$ by Lemma \ref{lem:genericunipotent}. Thus we may assume that $u$ is not generic, so the class to which $u$ belongs is given in cases (\ref{li:pslfirst11}) to (\ref{li:psllast11}) from the start of this section.

In each case we either have $5^2$ or $10$ in the action of $u$. A single block of size $10$ (as $M(E_7)$ is self-dual) must come from a self-dual indecomposable module, which must be $5/5$ by Lemma \ref{lem:selfdualsl2p}. For the two blocks $5^2$ in the action of $u$, these come from two modules of dimension congruent to $5$ modulo $11$ by Lemma \ref{lem:indJordanblocks}. From Proposition \ref{prop:simplesl2p}, we see that these modules are $5$ itself, $5,7/3,5,7$ and its dual, of dimension $27$, and $3,5,7,9/1,3,5,7,9$ and its dual, of dimension $49$. Thus if we are in case (\ref{li:pslsecond11}), so $u$ belongs to class $E_6(a_1)$ acting as $11^4,5^2,1^2$, we have a summand $5^{\oplus 2}$ of $M(E_7)\downarrow_H$, or 
\[ M(E_7)\downarrow_H=(5,7/3,5,7)\oplus (3,5,7/5,7)\oplus 1^{\oplus 2}.\]
An involution $x\in H$ acts with trace $0$ on this module, but involutions act on $M(E_7)$ with trace $\pm 8$ (see Appendix \ref{app:traces}), so it is not allowed.

Therefore, the Jordan blocks $10$ or $5^2$ always correspond to a module $5/5$ or $5^{\oplus 2}$.

If $u$ comes from class $E_7(a_3)$, so case (\ref{li:psllast11}), acting as $11^4,10,2$, the self-dual module that can contribute a block of size $2$ to the action of $u$ is $5,7/5,7$, so that $M(E_7)\downarrow_H$ has composition factors including $7^2,5^4$. The only conspicuous sets of composition factors with this many $7$s and $5$s are
\[ 9^2,7^2,5^4,1^4\qquad\text{and}\qquad 7^2,5^6,3^2,1^6,\]
with the latter being incompatible with the unipotent action and the former implying that $M(E_7)\downarrow_H$ is 
\begin{equation}
P(1)^{\oplus 2}\oplus (5/5)\oplus (5,7/5,7).\label{eq:case1forp=11}
\end{equation} Clearly therefore $H$ stabilizes a line on $M(E_7)$.

For $u$ belonging to one of the classes $D_6$ or $E_6(a_1)$, so cases (\ref{li:pslfirst11}) and (\ref{li:pslsecond11}), the blocks $5^2$ must come from a summand $5^{\oplus 2}$. In both cases, $u$ has two blocks of size $1$ on $M(E_7)$. If these come from two trivial summands, then $M(E_7)\downarrow_H$ is the sum of $(5/5)\oplus 1^{\oplus 2}$ or $5^{\oplus 2}\oplus 1^{\oplus 2}$ and a projective module. The only such modules with conspicuous sets of composition factors of $M(E_7)\downarrow_H$ are
\begin{equation} 11^{\oplus 2}\oplus P(1)^{\oplus 2}\oplus 1^{\oplus 2}\oplus M\quad\text{and}\quad P(9)^{\oplus 2}\oplus 1^{\oplus 2}\oplus M,\label{eq:case2forp=11}\end{equation}
with $M$ either $5^{\oplus 2}$ or $5/5$.

Thus suppose that $1^{\oplus 2}$ is not a summand of $M(E_7)\downarrow_H$. The indecomposable modules of dimension congruent to $1$ modulo $11$ are $1$, $3/9$ and its dual, and $5/7$ and its dual. Thus we may assume that $M(E_7)\downarrow_H$ has one of $(3/9)\oplus (9/3)$ or $(5/7)\oplus (7/5)$ as a summand. Again, the remainder of $M(E_7)\downarrow_H$ must be projective, and the only possibilities with conspicuous sets of composition factors are
\begin{equation} P(1)^{\oplus 2}\oplus (5/7)\oplus (7/5)\oplus M\quad \text{and}\quad P(1)^{\oplus 2}\oplus (3/9)\oplus (9/3)\oplus M,\label{eq:case3forp=11}\end{equation}
where $M$ is above.

In each of the cases in (\ref{eq:case1forp=11}), (\ref{eq:case2forp=11}) and (\ref{eq:case3forp=11}), $H$ stabilizes a line on $M(E_7)$, as claimed in the proposition.

\medskip

\noindent \textbf{Case $a=2$}: If $u$ is not from cases (\ref{li:pslfirst11}) to (\ref{li:psllast11}) then $u$ is generic and therefore $H$ is a blueprint by Lemma \ref{lem:genericunipotent}, so again we assume $u$ comes from one of these classes, and $L$ acts on $M(E_7)$ as one of the modules in (\ref{eq:case1forp=11}), (\ref{eq:case2forp=11}) and (\ref{eq:case3forp=11}). 

The only simple $kH$-module that restricts to $L$ with more $1$-dimensional factors than $3$-dimensional factors is the trivial module, by Lemma \ref{lem:sl27restriction}. Thus in every case $M(E_7)\downarrow_H$ has at least two trivial composition factors, and in the first option from (\ref{eq:case2forp=11}) $M(E_7)\downarrow_H$ has at six trivial composition factors.

In all cases, either $M(E_7)\downarrow_L$ has no $11$-dimensional factors at all, or has four more trivial factors than $11$-dimensional factors. Since the simple $kH$-modules with non-trivial $1$-cohomology are $20_{i,j}$ by Lemma \ref{lem:cohomologysimple}, and each restricts to $L$ as $11\oplus 9$, we see that $H$ has negative pressure in all cases, in fact has trivial summands in all cases. This completes the proof for $a=2$.
\end{proof}
It is possible to show in the above case that $M(E_7)\downarrow_H$ has at least six odd-dimensional factors for $p^a=121$, and therefore $H$ is a blueprint for $M(E_7)$, but this is not needed for our result.

\medskip

For $p=13$, since $169>150=2\cdot v(E_7)$, we need only consider $a=1$.

\begin{lemma}\label{prop:e7char13} If $p=13$ and $a=1$ then either $H$ is a blueprint for $M(E_7)$ or $H$ stabilizes a line on either $M(E_7)$ or $L(E_7)$.
\end{lemma}
\begin{proof} If $u$ lies in a generic unipotent class then $H$ is a blueprint for $M(E_7)$ by Lemma \ref{lem:genericunipotent}, so assume that the unipotent class is not generic. This means $u$ acts with Jordan blocks as in case (\ref{li:pslfirst13}) of the list of possible unipotent actions at the start of this section, so $13^4,1^4$.

Suppose that $H$ has no trivial summand on $M(E_7)$, so that the $1^4$ in the action of $u$ all comes from indecomposables of dimension congruent to $1$ modulo $13$: these are of the form $i/(14-i)$ for some odd $i$ by Proposition \ref{prop:simplesl2p}. Thus $M(E_7)\downarrow_H$ is the sum of two dual pairs of these modules.

Using the traces of elements of orders $2$, $3$ and $4$ from Appendix \ref{app:traces}, there is a unique conspicuous set of composition factors for $M(E_7)\downarrow_H$ of the form just described:
\[ (11/3)\oplus (3/11)\oplus (7/7)\oplus (7/7).\]
Using the corresponding traces on $L(E_7)$, we can determine two possibilities for the composition factors of $H$ on $L(E_7)$ (see (\ref{li:strata}) from Chapter \ref{ch:strategy}): these are
\[ 13^3,11,9^3,7^5,5^3,3,1^3\qquad\text{and}\qquad 13^2,11^3,9^4,7,5^4,3^3,1^2.\]
The only module for $\SL_2(p)$ that has a trivial factor but no trivial submodule or quotient is $P(11)=11/1,3/11$, so for both of these sets of composition factors, $H$ stabilizes a line on $L(E_7)$. Thus $H$ is either a blueprint for $M(E_7)$, or $H$ has a trivial summand on $M(E_7)$, or $H$ stabilizes a line on $L(E_7)$. Thus the proof is complete.
\end{proof}

The last case is $p=19$, where again we only have $a=1$.
\begin{proposition}\label{prop:e7char19} Suppose that $p=19$ and $a=1$. If $H$ is not a blueprint for $M(E_7)$ then $H$ centralizes a $2$-space on $M(E_7)$ and is a non-$\bG$-completely reducible subgroup of the $E_6$-parabolic subgroup acting on $M(E_7)$ as
\[ P(1)^{\oplus 2}\oplus (9/9).\]
\end{proposition}
\begin{proof} If $H$ is not a blueprint for $M(E_7)$ then in particular $u$ is non-generic, and so we are in case (\ref{li:pslfirst19}) from the list at the start of this section, i.e., $u$ is regular and acts with Jordan blocks $19^2,18$. We need a self-dual indecomposable module of dimension congruent to $18$ modulo $19$, and there is only one of these by Lemma \ref{lem:selfdualsl2p}, namely $9/9$, and the remainder of the module is projective. If $x$ denotes an involution in $H$ then $x$ has trace $\pm 8$ on $M(E_7)$ (see Appendix \ref{app:traces}), and has trace $2$ on $9/9$, leaving a trace of $6$ or $-10$ on the remaining projective summand. The trace of $x$ on $P(i)$ for $3\leq i\leq 17$ is $\pm 2$, the trace of $x$ on $19$ is $-1$, and on $P(1)$ it is $3$. Thus $M(E_7)\downarrow_H$ is $P(1)^{\oplus 2}\oplus (9/9)$, as needed.
\end{proof}

This non-$\bG$-completely reducible subgroup was constructed at the end of Chapter \ref{ch:e6}.

\medskip

We now check that in all cases $H$ is strongly imprimitive, or $p^a=7$. In all cases, $H$ is either a blueprint for $M(E_7)$, stabilizes a line on $M(E_7)$, or stabilizes a line on $L(E_7)$, and hence $H$ is strongly imprimitive by Propositions \ref{prop:blueprintissi}, \ref{prop:fixlineonMG} or \ref{prop:fixlineonLG} respectively.

Unless $p$ is one of $3,5,7,11,13,19$, then $u\in H$ is generic for $M(E_7)$, hence $H$ is a blueprint for $M(E_7)$ by Lemma \ref{lem:genericunipotent}. Thus we assume $p$ is one of these primes.

If $p=3$ then \cite[Proposition 6.2]{craven2015un2} (which proves that if $p^a=9$ then $H$ stabilizes a line on either $M(E_7)$ or $L(E_7)$) Propositions \ref{prop:e7char3a=4} and \ref{prop:e7char3a=3} prove the conclusion, and if $p=5$ then Proposition \ref{prop:e7char5} proves the conclusion. If $p=7$ then Proposition \ref{prop:e7char7} gives us the result, with the exception of $p^a=7$ and $H$ having actions
\[ 7^{\oplus 4}\oplus P(3)^{\oplus 2}\qquad\text{and}\qquad 7^{\oplus 5}\oplus P(5)^{\oplus 6}\oplus P(3)\]
on $M(E_7)$ and $L(E_7)$ respectively.

If $p=11$ then Proposition \ref{prop:e7char11} proves the result. If $p=13$ then we apply Proposition \ref{prop:e7char13}, and finally Proposition \ref{prop:e7char19} deals with $p=19$. Thus either $H$ is strongly imprimitive or $p^a=7$ and we have that single possibility, as needed for this part of the proof of Theorem \ref{thm:e7}.

\chapter{The Proof for \texorpdfstring{$E_7$}{E7} in Odd Characteristic: \texorpdfstring{$\SL_2$}{SL2} Embedding}
\label{ch:e7oddsl}

In this chapter, $k$ is an algebraically closed field of characteristic $p\geq 3$ and $\bG=E_7(k)$, by which we mean the simply connected form, i.e., $|Z(\bG)|=2$ and $\bG'=\bG$. Let $H\cong \SL_2(p^a)$ be a subgroup of $\bG$ such that $Z(H)=Z(\bG)$.

Theorem \ref{thm:goodblueprint} states that if $p^a\geq 150$ then $H$ is a blueprint for $M(E_7)$. In this case, $H$ is strongly imprimitive by Proposition \ref{prop:blueprintissi}. Thus in what follows we may assume that $p^a\leq 150$.

Let $L=\SL_2(p)\leq H$ and let $u$ denote a unipotent element of $L$ of order $p$. The possibilities for the Jordan block structures of $u$ on $M(E_6)$ and $L(E_6)$ are given in \cite[Tables 7 and 8]{lawther1995}. Recall the definition of a generic unipotent element from Definition \ref{defn:genericunipotent}.

On $M(E_7)$, since we consider $\SL_2(p^a)$ rather than $\PSL_2(p^a)$, there can be no trivial composition factors in $M(E_7)\downarrow_H$, but rather $2$-dimensional factors. We will thus normally aim to show that $H$ is a blueprint for $M(E_7)$, that $H$ stabilizes a line on $L(E_7)$, or that $H$ stabilizes a $2$-space on $M(E_7)$ (and then apply Propositions \ref{prop:blueprintissi}, \ref{prop:fixlineonLG} and \ref{prop:fix2spaceonMG} respectively to show that $H$ is strongly imprimitive). This will not always be possible, and we show directly that $H$ is strongly imprimitive.

\section{Characteristic 3}

Let $p=3$, so that $H=\SL_2(3^a)$ for some $a=2,3,4$. The case $a=2$ was considered in \cite[Proposition 6.2]{craven2015un2}, so we exclude this.

We first attack the case of $a=3$. We will prove that one of a variety of conditions holds, each of which is sufficient to prove that $H$ is strongly imprimitive.

\begin{proposition}\label{prop:e7char3a=3sl} If $p=3$ and $a=3$ then $H$ is strongly imprimitive.
\end{proposition}
\begin{proof} If $H$ stabilizes a $2$-space on $M(E_7)$ then $H$ is strongly imprimitive by Proposition \ref{prop:fix2spaceonMG}. Similarly, if $H$ stabilizes a $1$-space on $L(E_7)$ then $H$ is strongly imprimitive by Proposition \ref{prop:fixlineonLG}. Thus for the rest of this proof we may assume that $H$ does not stabilize a $2$-space on $M(E_7)$, or a line on $L(E_7)$.

There are $284$ conspicuous sets of composition factors for $M(E_7)\downarrow_H$, but only $137$ of these have corresponding sets of factors on $L(E_7)$, each of these being unique. Of these, $77$ have either no $2_i$ or positive $2_i$-pressure for each $i=1,2,3$, and of these only $67$ have either no trivial or positive pressure on $L(E_7)$, so we may eliminate those via Lemma \ref{lem:pressure}. Exactly one of these remaining sets of factors is invariant under the field automorphism, so we are left with $23$ sets of composition factors up to field automorphism.

The module $2_i$ has non-split extensions only with $2_{i\pm 1}$, $6_{i-1,i}$ and $8$, so if there are $2_i$s in $M(E_7)\downarrow_H$ but no $6_{i,i-1}$ or $8$ appearing with multiplicity $2$ or above, then $H$ stabilizes a $2$-space on $M(E_7)$: two sets of composition factors (up to field automorphism) satisfy this, so we are down to $21$ sets of composition factors.

There are, up to field automorphism, six sets of composition factors with no $8$s. Since $2_i$ only has extensions with the $8$, $2_{i\pm 1}$ and $6_{i-1,i}$, the socle must consist of copies of $6_{i-1,i}$ for various $i$, plus modules we can quotient out by without a $2_i$ appearing in the socle. In each case there is a unique $i$ such that $6_{i-1,i}$ appears with multiplicity at least $2$, so this must be the socle and all $2_i$ must be stacked on top of it in some way. In each case we cannot place enough copies of $2_i$ on top of each $6_{i-1,i}$, and so $H$ must stabilize a $2$-space on $M(E_7)$ in all these cases. This reduces us to fifteen sets of composition factors.

The remaining fifteen sets have at least one $8$ in $M(E_7)\downarrow_H$. There are five conspicuous sets of composition factors with a single $8$, up to field automorphism, namely
\[ 8,6_{1,2}^2,6_{3,2}^2,2_1^2,2_2^6,2_3^4,\qquad 8,6_{1,2}^2,6_{2,1}^2,6_{3,1},2_1^2,2_2^5,2_3^2,\qquad 8,6_{1,2}^5,2_1^2,2_2^5,2_3^2,\]
\[8,6_{1,2}^2,6_{2,1}^2,6_{3,1},6_{3,2},2_1^2,2_2^2,2_3^2,\qquad  18_{3,1,2},8,6_{1,2}^2,6_{3,2},2_1^2,2_2^3,2_3.\]
In each case, the only $6_{i-1,i}$ that appears more than once in $M(E_7)\downarrow_H$ is $6_{1,2}$. Since $8$ appears exactly once, if it is a submodule (or quotient) of $M(E_7)\downarrow_H$ then it is a summand, so we may assume that it is not. Let $W$ denote the $\{6_{i,i+1},8\}$-radical of $M(E_7)\downarrow_H$ modulo its $\{6_{i,i+1},8\}$-residual, so that $W$ has socle and top $6_{1,2}^{\oplus i}$ for some $i$.

We take the $\{2_i,6_{i,j},8\}$-radical of $P(6_{1,2})$, then take the $\{6_{1,2}\}'$-residual of this radical, and obtain the module
\[ 6_{1,2}/2_2/2_3/2_2/6_{1,2}.\]
The module $W$ must be a submodule of a sum of these, but $2_1$ does not appear in it. Thus $H$ must stabilize a $2$-space on $M(E_7)$, except possibly in the fifth case, as there is a module $18_{3,1,2}$. In this case, $W$ must be a submodule of $P(6_{1,2})$. We take the $\{2_i,6_{3,2},8,18_{3,1,2}\}$-radical of the quotient module $P(6_{1,2})/6_{1,2}$, and lift to $P(6_{1,2})$, to obtain the module
\[ 2_1,2_2/2_3,8/18_{3,1,2},2_2,6_{3,2}/6_{1,2}.\]
This has only four $2$-dimensional composition factors, and $M(E_7)\downarrow_H$ should possess six, and so we obtain a contradiction.

Thus we are down to ten sets of composition factors for $M(E_7)\downarrow_H$, which are below.
\[ 8^2,6_{2,1}^2,6_{1,2},6_{3,1},2_1^4,2_2^3,2_3,\qquad 8^2,6_{2,1},6_{1,3},6_{3,1}^2,2_1^4,2_2^2,2_3^2,\] \[8^2,6_{3,2}^2,6_{1,2},6_{2,1},6_{3,1},6_{1,3},2_1,2_2,\qquad 8^2,6_{3,2}^2,6_{2,1},6_{3,1}^2,2_1^2,2_2^2,2_3,\] \[8^2,6_{1,3},6_{3,1},6_{2,3}^2,6_{1,2},2_1^2,2_2^2,2_3,\qquad 8^2,6_{3,2},6_{2,1},6_{2,3},6_{1,3}^2,2_1^2,2_2^2,2_3,\]
\[ 18_{1,2,3},8^2,6_{2,1},6_{3,1},2_1^2,2_2^2,2_3,\qquad 18_{1,2,3},8^3,6_{1,3},2_1^3,2_2,\] \[18_{1,2,3},8^2,6_{2,1},6_{1,3},2_1^3,2_2^2.\qquad 8^4,(2_1,2_2,2_3)^4.\]
Recall that $2_i$ has extensions with $2_{i\pm 1}$, $6_{i-1,i}$ and $8$, and no other simple modules.

\medskip

\noindent\textbf{Cases 1, 6, 7, 9}: In the first, sixth, seventh and ninth cases, $6_{i-1,i}$ occurs with multiplicity at most $1$ for all $i$, and so if such a module occurs in the socle of $M(E_7)\downarrow_H$ then it is a summand, and may be ignored. Therefore, for these four sets of composition factors, we can remove all quotients and submodules from $M(E_7)\downarrow_H$ other than $8$, and yield a submodule $W$ of $P(8)$, which contains all composition factors of dimension $2$ in $M(E_7)\downarrow_H$.

The $\{2_i,6_{i,i+1},6_{i+1,i},18_{1,2,3}\}$-radical of the quotient module $P(8)/8$, lifted to $P(8)$, is
\[ 2_1,2_2,2_3,6_{1,2},6_{2,3},6_{3,1}/2_1,2_2,2_3,6_{2,1},6_{3,2},6_{1,3}/8,\]
and so if $W$ has two $8$s and three factors $2_i$ for some $i$, then $H$ must stabilize a $2$-space on $M(E_7)$. This eliminates the first and ninth cases. We also see that in the sixth and seventh cases, $W$ has at most four socle layers. Since $W$ has at most four socle layers and is self-dual, none of the $6_{i,i+1}$ can occur in $W$. Also, any $2_i$ that occurs with multiplicity $1$ in $M(E_7)\downarrow_H$ must occur in the second socle layer, and cannot have any extensions with a $2_{i\pm 1}$ in the third socle layer. However, in the module above, an explicit check shows that both $2_{i\pm 1}$ in the third layer have an extension with the $2_i$ in the second layer, so cannot exist in $W$. In other words, we cannot have two $2_i$ in $W$, eliminating the sixth and seventh cases.

\medskip

\noindent\textbf{Case 8}: The only factors appearing with multiplicity at least $2$ are $8$ and $2$-dimensional factors. Since $H$ does not stabilize a $2$-space on $M(E_7)$ by assumption, $M(E_7)\downarrow_H$ is a submodule of $P(8)$ and possibly some summands $6_{i,j}$ and $18_{1,2,3}$. However, the $\{2_1,2_2,6_{1,3},8,18_{1,2,3}\}$-radical of $P(8)$ is
\[ 8/2_1,2_2,6_{1,3}/8,\]
which doesn't have enough copies of $2_1$.

\medskip

\noindent\textbf{Case 5}: Remove all summands $6_{1,3}$, $6_{3,1}$ and $6_{1,2}$ from $M(E_7)\downarrow_H$ to yield a summand $W$ of $M(E_7)\downarrow_H$ with socle a submodule of $8\oplus 6_{2,3}$ and all $2_i$ in it.

The preimages of the  $\{2_1,2_2,2_3,6_{1,2},6_{1,3},6_{2,3}\}$-radical of the quotient module $P(8)/8$ and the $\{2_1,2_2,2_3,6_{1,2},6_{1,3},8\}$-radical of $P(6_{2,3})/6_{2,3}$ in their respective projectives are
\[ 2_1,2_2,2_3,6_{2,3}/2_1,2_2,2_3,6_{1,3}/8
\qquad\text{and}\qquad
2_2,2_3/2_1,8/2_3,6_{1,3}/6_{2,3},\]
and the fact that there is a single $2_3$ means it must lie in the second socle layer, and has no extensions with the other $2_i$ composition factors. But then all other $2_i$ lie in the second socle layer as well, and that means we cannot fit enough $2$s in, so $H$ stabilizes a $2$-space on $M(E_7)$.

\medskip

\noindent\textbf{Cases 2, 3, 4}: Here we will apply Corollary \ref{cor:sl2nopgl}. Let $y$ be the diagonal matrix with entries $\zeta^2,\zeta^{-2}$ for $\zeta$ a primitive $26$th root of unity, so that $y$ has order $13$ and acts with eigenvalues $\zeta^{\pm 2}$ on $2_1$. The eigenvalues of $y^2$ on $6_{3,1}$ and $8$ are
\[ \zeta^{\pm 4},\zeta^{\pm 8},\zeta^{\pm 12},\quad\text{and}\quad 1^2,\zeta^{\pm 4},\zeta^{\pm 10},\zeta^{\pm 12}\]
respectively. We first show that in the second and fourth cases, $M(E_7)\downarrow_H$ possesses a unique submodule $8$ (under our standing assumption that $H$ does not stabilize a $2$-space on $M(E_7)$), and in the third $M(E_7)\downarrow_H$ possesses a summand $6_{3,1}$.

In the third case, $6_{3,1}$ has no extension with $8$ or $6_{3,2}$, the only factors to appear with multiplicity at least $2$, so $6_{3,1}$ must split off as a summand. Thus we consider the second and fourth cases. If there is a copy of $8$ in the socle that is not a summand, then there must be another copy in the top and we are done, so assume this is not the case. Thus in the two cases, $M(E_7)\downarrow_H$ is the sum of some $8$s and a submodule of $P(6_{3,1})$ and $P(6_{3,2}\oplus 6_{3,1})$ respectively. The $\cf(M(E_7)\downarrow_H)$-radical of $P(6_{3,1})$ in Case 2 is
\[6_{3,1}/2_1,2_3/2_2,8/2_1,6_{2,1}/6_{3,1},\]
which does not contain enough $2$-dimensional factors. Similarly, the $\cf(M(E_7)\downarrow_H)$-radicals of $P(6_{3,1})$ and $P(6_{3,2})$ in Case 4 are
\[6_{2,1}/6_{3,1},8/2_1,2_3,6_{3,2}/2_2,8/2_1,6_{2,1}/6_{3,1}\qquad \text{and}\qquad 6_{3,1}/2_1,6_{2,1}/8/6_{3,2}\]
respectively. There are not enough copies of $2_2$ in this case. Thus in these two cases $H$ stabilizes a unique, irreducible $8$-space $W$ on $M(E_7)$, which is necessarily $N_{\Aut^+(G)}(H)$-stable.

We also claim that the $6_{3,1}$ in the third case is $N_{\Aut^+(G)}(H)$-stable-stable. If it were not, it could only be sent to a summand $6_{1,2}$ (as all other $6$-dimensional factors lie in the other $\Aut(H)$-orbit of simple $kH$-modules), but then $6_{1,2}$ must be sent to $6_{3,1}$. There is no automorphism of $H$ that swaps $6_{1,2}$ and $6_{3,1}$, so this subspace must be stable under elements of $N_{\Aut^+(G)}(H)$-stable. In this case, let $W$ denote this $6$-space.

Suppose that we can find an element $\hat y$ of order $26$ in $\bG\setminus H$ that has no $(-1)$-eigenspace on $L(E_7)$ and that stabilizes $W$. Suppose that $\hat y$ normalizes $H$. The group $\gen{H,\hat y}$ cannot be $\PGL_2(p^a)$ (modulo $Z(\bG)$) by Corollary \ref{cor:sl2nopgl}, and since the composition factors of $M(E_7)\downarrow_H$ are not stable under a field automorphism of $H$, $\hat y$ cannot induce a field automorphism on $H$. Thus $\hat y$ centralizes $H$, but then $\hat y^2=y$ centralizes $H$, which is wrong. Thus $\gen{H,\hat y}\not\leq N_{\bG}(H)$. Since $\gen{H,\hat y}$ stabilizes $W$, $\gen{H,\hat y}$ does not have the same type as $\bG$ either. We apply Proposition \ref{prop:maximalnotinP} to see that either $H$ is strongly imprimitive or $\gen{H,\hat y}$ is contained in a member of $\mathscr P$ (up to taking normalizers). By Proposition \ref{prop:maximalnotinP}, this must be ${}^2\!G_2(27)$, but Ree groups do not have irreducible $6$- or $8$-dimensional modules, and so it must be the case that $H$ is strongly imprimitive.

We therefore find, in each case, an element $\hat y$ in $\bG\setminus H$ of order $26$, squaring to $y^2$, and stabilizing the eigenspaces of the particular stabilized submodule. On $M(E_7)$ these have eigenvalues
\[1^4,(-\zeta^{\pm 1})^6,(\zeta^{\pm 2})^4,(\zeta^{\pm 3}),(-\zeta^{\pm 3})^2,(\zeta^{\pm 4})^2,(-\zeta^{\pm 4})^3,(\zeta^{\pm 5})^3,(-\zeta^{\pm 6})^5,\]
\[ 1^4,(\zeta^{\pm 1})^2,(-\zeta^{\pm 1})^2,(-\zeta^{\pm 2})^5,(\zeta^{\pm 3})^3,(-\zeta^{\pm 3}),(\zeta^{\pm 4})^4,(\zeta^{\pm 5}),(-\zeta^{\pm 5})^3,(-\zeta^{\pm 6})^5,\;\; \text{and}\]
\[ 1^4,(\zeta^{\pm 1})^2,(-\zeta^{\pm 1}),(-\zeta^{\pm 2})^6,(\zeta^{\pm 3})^4,(-\zeta^{\pm 3}),(\zeta^{\pm 4})^5,(-\zeta^{\pm 5})^3,(-\zeta^{\pm 6})^4,\]
respectively. Thus the result holds.

\medskip

\noindent \textbf{Case 10}: The $\{2_i,8\}$-radical of $P(8)$ is
\[ 8/2_1,2_2,2_3/2_1,2_2,2_3/8,\]
so $\soc(M(E_7)\downarrow_H)$ is $8^{\oplus 2}$, as we are assuming that $H$ does not stabilize a $2$-space on $M(E_7)$. We keep the notation of $\zeta$, $y$ and $\hat y$ from Cases 2, 3, 4, above. In this case $y^2$ has trace $4$ on $M(E_7)$, and we may choose $\hat y$ to have eigenvalues
\[ 1^8,(-\zeta^{\pm 1})^4,(\zeta^{\pm 2})^4,(\zeta^{\pm 3})^2,(-\zeta^{\pm 3})^2,(\zeta^{\pm 4})^2,(-\zeta^{\pm 4})^2,(\zeta^{\pm 5})^4,(-\zeta^{\pm 6})^4\]
on $M(E_7)$. The same proof therefore holds. (Of course, $\gen{H,\hat y}$ stabilizes all $8$-dimensional submodules, so it is $N_{\Aut^+(\bG)}(H)$-stable.)
\end{proof}

\begin{proposition}\label{prop:e7char3a=4sl} Suppose that $p=3$ and $a=4$. Either $H$ is a blueprint for $M(E_7)$ or the composition factors of $M(E_7)\downarrow_H$ are
\[ 18_{4,2,3},8_{1,2,3}^2,8_{1,2,4},6_{1,3},6_{4,1},2_1,\]
the $18$-dimensional composition factor is a summand of $M(E_7)\downarrow_H$, and the stabilizer of this unique irreducible $18$-space is positive dimensional.
\end{proposition}
\begin{proof}
Using semisimple elements of order up to $41$, one whittles down the 55 million or so possible sets of composition factors for a module of dimension $56$ to just $190$ up to field automorphism. Using the preimage trick from the end of Section \ref{sec:blueprints}, we can also check the traces of elements of order $80$, and two of these sets of composition factors are not conspicuous for these elements. This leaves 188 conspicuous sets of composition factors for $M(E_7)\downarrow_H$.

Of these, we consider an element $y$ of order $41$ in $H$, and whether there exists an element of order $123$ in $\bG$ cubing to $y$ and stabilizing the same subspaces of $M(E_7)$ as $y$. If this is true then $y$ (and therefore $H$) is a blueprint for $M(E_7)$ by Theorem \ref{thm:blueprintsminimal}. 

As we have a list of the semisimple classes of elements of order $41$, we know that the class of $y$ is determined by its trace on $M(E_7)$. Thus we simply identify $y$ with an element of a maximal torus $\bT$ of $\bG$ with these eigenvalues, and we may use the preimage trick to check whether one of the $3^7$ elements of order $123$ that cube to $y$ has the same number of distinct eigenvalues on $M(E_7)$. (To find such an element, simply choose at random from the elements of order $41$ until an element is found with the correct eigenvalues on $M(E_7)$.)

Indeed, a computer check shows that this is true for $187$ of the $188$ semisimple elements involved. The remaining one comes from the conspicuous set of composition factors in the statement of the proposition,
\[ 18_{4,2,3},8_{1,2,3}^2,8_{1,2,4},6_{1,3},6_{4,1},2_1.\]
For these composition factors, $y$ has $38$ distinct eigenvalues on $M(E_7)$, and there exist elements of order $123$ cubing to $y$ and with $40$ distinct eigenvalues, but none with $38$. If $\zeta$ is a primitive root of unity then $y$ can be chosen to have the following eigenvalues.
\begin{center}\begin{tabular}{cc}
\hline Module & Eigenvalues
\\\hline $2_1$& $\zeta^{\pm 1}$ 
\\ $6_{1,3}$&$\zeta^{\pm 1},\zeta^{\pm 17},\zeta^{\pm 19}$
\\$6_{4,1}$&$\zeta^{\pm 12},\zeta^{\pm 14},\zeta^{\pm 16}$
\\ $8_{1,2,3}$&$\zeta^{\pm 5},\zeta^{\pm 7},\zeta^{\pm 11},\zeta^{\pm 13}$
\\ $8_{1,2,4}$&$\zeta^{\pm 10},\zeta^{\pm 12},\zeta^{\pm 16},\zeta^{\pm 18}$
\\ $18_{4,2,3}$&$\zeta^{\pm 2},\zeta^{\pm 3},\zeta^{\pm 4},\zeta^{\pm 8},\zeta^{\pm 9},\zeta^{\pm 10},\zeta^{\pm 14},\zeta^{\pm 15},\zeta^{\pm 20}$
\\ \hline \end{tabular}\end{center}
There exists an element $\hat y$ of order $123$ cubing to $y$ and stabilizing all eigenspaces except for the $\zeta^{\pm 1}$-eigenspaces.

Since $8_{1,2,3}$ is the only composition factor to occur with multiplicity greater than $1$, any other factor in the socle must be a summand. The module $8_{1,2,3}$ only has extensions with $2_1$, $6_{1,3}$ and $8_{1,2,4}$ from the composition factors of $M(E_7)\downarrow_H$, and so the structure of $M(E_7)\downarrow_H$ must be $W\oplus 6_{4,1}\oplus 18_{4,2,3}$, where $W$ consists of the remaining factors.

Since $\hat y$ stabilizes all but the $\zeta^{\pm 1}$-eigenspaces, it stabilizes the $32\oplus 6\oplus 18$ decomposition above. If $\bY$ denotes an infinite subgroup of $\bG$ containing $\hat y$ and stabilizing the same subspaces of $M(E_7)$ as $\hat y$ (which exists by Theorem \ref{thm:blueprintsminimal}), then $\bX=\gen{\bY,H}$ certainly stabilizes the $6$- and $18$-dimensional summands of $M(E_7)\downarrow_H$ and is positive dimensional.\end{proof}

\section{Characteristic At Least 5}

We now let $p\geq 5$: we still have $p^a\leq 150$, as we stated at the start of this chapter. As all unipotent classes are generic for $M(E_7)$ for all $p\geq 29$, we only need consider
\[ p^a=5,7,11,13,17,19,23,25,49,121,125.\]

As for $\PSL_2(p^a)$, there are some restrictions we can place on the possible actions of a unipotent element $u$, beyond that of appearing on \cite[Table 7]{lawther1995}, given by Lemma \ref{lem:faithfulpblockseven}. This yields twenty-nine possible non-generic classes for various primes, as given below.

\begin{enumerate}
\item\label{li:sl2first5} $(A_3+A_1)''$, $p=5$, acting as $5^2,4^8,2^7$;
\item $D_4(a_1)+A_1$, $p=5$, acting as $5^6,4,3^4,2^5$;
\item $A_3+A_2$, $p=5$, acting as $5^6,4^2,3^4,2^2,1^2$;
\item $A_4$, $p=5$, acting as $5^{10},1^6$;
\item $A_3+A_2+A_1$, $p=5$, acting as $5^6,4^4,2^5$;
\item $A_4+A_1$, $p=5$, acting as $5^{10},2^2,1^2$;
\item\label{li:sl2last5} $A_4+A_2$, $p=5$, acting as $5^{10},3^2$;
\item\label{li:sl2first7} $(A_5)''$, $p=7$, acting as $7^2,6^7$;
\item $D_4+A_1$, $p=7$, acting as $7^6,2^5,1^4$; 
\item $D_5(a_1)$, $p=7$, acting as $7^6,3^2,2^2,1^4$;
\item $(A_5)'$, $p=7$, acting as $7^4,6^4,1^4$;
\item $A_5+A_1$, $p=7$, acting as $7^4,6^3,5^2$;
\item $D_5(a_1)+A_1$, $p=7$, acting as $7^6,4,2^5$;
\item $D_6(a_2)$, $p=7$, acting as $7^6,5^2,4$;
\item $E_6(a_3)$, $p=7$, acting as $7^6,5^2,1^4$;
\item $E_7(a_5)$, $p=7$, acting as $7^6,6,4^2$;
\item\label{li:sl2last7} $A_6$, $p=7$, acting as $7^8$;
\item\label{li:sl2first11} $E_7(a_4)$, $p=11$, acting as $11^2,10,8,6,4^2,2$;
\item $D_6$, $p=11$, acting as $11^4,10,1^2$;
\item $E_6(a_1)$, $p=11$, acting as $11^4,5^2,1^2$;
\item\label{li:sl2last11} $E_7(a_3)$, $p=11$, acting as $11^4,10,2$;
\item\label{li:sl2first13} $E_6$, $p=13$, acting as $13^4,1^4$;
\item\label{li:sl2second13} $E_7(a_3)$, $p=13$, acting as $13^2,12,10,6,2$;
\item\label{li:sl2last13} $E_7(a_2)$, $p=13$, acting as $13^4,4$;
\item\label{li:sl2first17} $E_7(a_2)$, $p=17$, acting as $17^2,10,8,4$;
\item\label{li:sl2last17} $E_7(a_1)$, $p=17$, acting as $17^2,16,6$;
\item\label{li:sl2first19} $E_7(a_1)$, $p=19$, acting as $19^2,12,6$;
\item\label{li:sl2last19} $E_7$, $p=19$, acting as $19^2,18$;
\item\label{li:sl2first23} $E_7$, $p=23$, acting as $23^2,10$.
\end{enumerate}

Thus we need to consider $p=5,7,11,13,17,19,23$, and we will examine each in turn.

\begin{proposition}\label{prop:e7char5sl} Suppose that $p=5$.
\begin{enumerate}
\item If $a=1$ then $H$ stabilizes a $2$-space on $M(E_7)$.
\item If $a=2$ then either $H$ is strongly imprimitive or the actions of $H$ on $M(E_7)$ and $L(E_7)$ are
\[ 10_{2,1}\oplus (12_{1,2}/2_1,4_1,6_{1,2},10_{1,2}/12_{1,2})\]
and
\[(15_{2,1}/8_{2,1}/1,3_2/8_{2,1}/15_{2,1})\oplus (15_{1,2}/8_{1,2}/1,3_1/8_{1,2}/15_{1,2})\oplus (3_2/8_{2,1},16/3_2)\oplus 3_1\]
respectively. Furthermore, in the latter case $H$ stabilizes an $\slf_2$ subalgebra of $L(E_7)$.
\item If $a=3$ then an element of order $63$ in $H$ is a blueprint for $M(E_7)$, and hence $H$ is a blueprint for $M(E_7)$.
\end{enumerate}
\end{proposition}
\begin{proof} First let $a=1$. Only the element of order $3$ is important here, and it has trace one of $-25,-7,2,20$ (see Appendix \ref{app:traces}), with the last case not possible, and so the conspicuous sets of composition factors for $M(E_7)\downarrow_H$ are as follows:
\[ 4,2^{26},\qquad 4^7,2^{14},\qquad 4^{10},2^8.\]
As $P(4)=4/2/4$, each of these must stabilize a $2$-space on $M(E_7)$, as claimed.

When $a=2$, we will prove that, with the single exception in the proposition, $H$ must satisfy one of a number of conditions, each of which implies strong imprimitivity.

There are $106$ conspicuous sets of composition factors for $M(E_7)\downarrow_H$, but fifty of these have no corresponding set of composition factors for $L(E_7)$ (see (\ref{li:strata}) from Chapter \ref{ch:strategy}), so can be ignored. If $H$ stabilizes a line on $L(E_7)$ or a $2$-space on $M(E_7)$ then $H$ is strongly imprimitive by Propositions \ref{prop:fixlineonLG} or \ref{prop:fix2spaceonMG} respectively. Thus we may exclude those sets of factors with non-positive $2_i$-pressure (and at least one $2_i$) by Lemma \ref{lem:pressure}. This leaves eighteen sets of factors, nine up to field automorphism of $H$.

Exactly one of these has two possible sets of composition factors on $L(E_7)$, one of which has a single trivial and a single $8$-dimensional, so that second option will be ignored as having pressure $0$ on $L(E_7)$. Of the other eight, which all have a unique corresponding set of composition factors on $L(E_7)$, three have non-positive pressure and trivial factors, so $H$ stabilizes a line on $L(E_7)$.

This leaves six conspicuous sets of composition factors on $M(E_7)$ up to field automorphism. These are
\[
12_{1,2}^2,10_{1,2},4_2^2,4_1^2,2_1^3,\quad 
12_{1,2},10_{1,2},6_{2,1},6_{1,2}^2,4_2^3,2_1^2,\]
\[12_{1,2},10_{1,2}^2,6_{2,1},6_{1,2},4_2,4_1,2_2,2_1,\quad
12_{1,2},12_{2,1},10_{2,1},6_{2,1},6_{1,2},4_2,4_1,2_1,\]
\[12_{1,2}^2,10_{2,1},10_{1,2},6_{1,2},4_1,2_1,\quad 
12_{2,1}^2,10_{1,2},6_{1,2}^3,4_1.
\]
The simple modules that have extensions with $2_1$ are $4_2$, $6_{2,1}$ and $12_{1,2}$, by \cite[Corollary 4.5]{andersen1983} (see also Lemma \ref{lem:psl2125smallmodules} for those of dimension at most $8$).

\medskip

\noindent\textbf{Case 1}: This has pressure $1$, and we may assume that we have a submodule of $P(12_{1,2})$ or $P(4_2)$ with three copies of $2_1$. The $\{2_1,4_1,4_2,10_{1,2},12_{1,2}\}$-radicals of these two modules are $10_{1,2}/12_{1,2}/2_1,4_1,10_{1,2}/12_{1,2}$ and $4_2/2_1/4_2$, so $H$ stabilizes a $2$-space on $M(E_7)$. Thus $H$ is strongly imprimitive by Proposition \ref{prop:fix2spaceonMG}.

\medskip

\noindent\textbf{Case 2}: The corresponding submodule of $P(4_2)$ in the second case is also $4_2/2_1/4_2$, so again this stabilizes a $2$-space on $M(E_7)$. Thus $H$ is again strongly imprimitive by Proposition \ref{prop:fix2spaceonMG}.

\medskip

\noindent\textbf{Cases 3, 4}: The only module appearing with multiplicity greater than $1$ in the third case is $10_{1,2}$, so unless a module has a non-trivial extension with $10_{1,2}$, it must split off as a summand. Thus $H$ stabilizes a $2$-space on $M(E_7)$. In the fourth case $M(E_7)\downarrow_H$ must be semisimple, so again $H$ stabilizes a $2$-space on $M(E_7)$.

\medskip

\noindent\textbf{Case 6}: There are no extensions between the simple modules involved so $M(E_7)\downarrow_H$ is semisimple. The corresponding composition factors for $L(E_7)\downarrow_H$ have no extensions between them either and so the restriction is also semisimple, acting as
\[ 15_{1,2}^{\oplus 3}\oplus 15_{2,1}^{\oplus 2}\oplus 9^{\oplus 3}\oplus 5_1\oplus 5_2\oplus 3_1^{\oplus 4}\oplus 3_2^{\oplus 3}.\]
Write $x$ for an element of order $13$ in $H$. Choosing $\zeta$ a primitive $13$th root of unity appropriately (so that $x$ acts on $2_1$ with eigenvalues $\zeta^{\pm 1}$) the eigenvalues of $x$ on $M(E_7)$ are
\[ 1^4,(\zeta^{\pm 1})^7,(\zeta^{\pm 2})^6,(\zeta^{\pm 3})^3,(\zeta^{\pm 4})^6,\zeta^{\pm 5},(\zeta^{\pm 6})^3.\]

Looking through the elements of order $26$ in $E_7$, we find one $\hat x$ that squares to $x$, and if $\theta$ is a primitive $26$th root of $1$ with $\theta^2=\zeta$, we have that the eigenvalues of $\hat x$ are
\[ 1^4,(\theta^{\pm 1})^7,(\theta^{\pm 2})^6,(\theta^{\pm 3})^3,(\theta^{\pm 4})^6,\theta^{\pm 5},(\theta^{\pm 6})^2,(-\theta^{\pm 6}).\]
This stabilizes the $4$-space of $M(E_7)$ stabilized by $H$ (as well as the $6$-spaces and the sum of the $12$- and $10$-spaces). Write $Y$ for the stabilizer of this $4$-space: since there is a unique irreducible $4$-space of $M(E_7)$ stabilized by $H$, $Y$ is $N_{\Aut^+(\bG)}(H)$-stable by Proposition \ref{prop:intersectionstabilizers}.

We apply Propositions \ref{prop:maximalnotinP} and \ref{prop:possiblemaximals} to $H$ and $Y$. Since $Y$ stabilizes a $4$-space on $M(E_7)$, certainly $Y$ does not contain the Rudvalis simple group. Thus while $H$ is not a maximal member of $\mathscr P$, the results still apply and either $Y\leq N_{\bG}(H)$ or $H$ is strongly imprimitive. As the composition factors of $M(E_7)\downarrow_H$ are not stable under a field automorphism of $H$, if $Y\leq N_{\bG}(H)$ then $Y\cong \PGL_2(25)$ (as $\Out(H)$ has order $4$, with generators a diagonal automorphism and a field automorphism).

The eigenvalues of $\hat x$ on $L(E_7)$ are
\[ 1^{17},(\theta^{\pm 1})^{10},(\theta^{\pm 2})^{13},(\theta^{\pm 3})^{12},(\theta^{\pm 4})^6,(\theta^{\pm 5})^8,(-\theta^{\pm 5})^3,(\theta^{\pm 6})^4,(-\theta^{\pm 6})^2,\]
so $\gen{H,\hat x}$ is not $\PGL_2(25)$ modulo $Z(\bG)$ by Corollary \ref{cor:sl2nopgl}. Hence $Y\not\leq N_{\bG}(H)$ and so $H$ is strongly imprimitive.

\medskip

\noindent\textbf{Case 5}: The $10_{2,1}$ must split off as it has no extensions with $12_{1,2}$, but the rest of the composition factors of $M(E_7)\downarrow_H$ can lie above $12_{1,2}$, and there is a unique module
\[ 10_{2,1}\oplus (12_{1,2}/2_1,4_1,6_{1,2},10_{1,2}/12_{1,2}),\]
with $u$ acting with Jordan blocks $5^{10},3^2$, so unipotent class $A_4+A_2$. The action of $u$ on the direct sum of the composition factors of $M(E_7)\downarrow_H$ has block structure $5^8,4^2,2^4$, and examining \cite[Table 7]{lawther1995}, we see that there are only two possible actions for $u$ that have at least eight blocks of size $5$ and at most fourteen blocks: $5^{10},3^2$ and $5^{10},2^2,1^2$. Thus the $4_1$ and $6_{1,2}$ cannot be summands (as they both have a $4$ in the action of $u$), and we assume that the $2_1$ is not a submodule, so if $M(E_7)\downarrow_H$ is not the module above then only the $10_{1,2}$ can be removed from the non-simple summand. However, the $\{2_1,4_1,6_{1,2},12_{1,2}\}$-radical of $P(12_{1,2})$ is
\[ 12_{1,2}/2_1,4_1,6_{1,2}/12_{1,2},\]
but with a $6_{1,2}$ quotient, which is not allowed. Thus $M(E_7)\downarrow_H$ is as above, $u$ lies in class $A_4+A_2$, and in particular the symmetric square of this has $L(E_7)\downarrow_H$ as a summand (since $S^2(M(E_7))=L(E_7)\oplus L(2\lambda_1)$ by Lemma \ref{lem:relatingvminlg}).

The composition factors of $L(E_7)\downarrow_H$ are
\[ 16,15_{1,2}^2,15_{2,1}^2,8_{1,2}^2,8_{2,1}^3,3_1^2,3_2^3,1^2,\]
and $u$ must act on $L(E_7)$ with Jordan blocks $5^{26},3$ from \cite[Table 8]{lawther1995}. There are only six isomorphism types of indecomposable module appearing as a summand of $S^2(M(E_7)\downarrow_H)$ whose composition factors appear on the list above, and these have structures
\[ 3_1,\quad 15_{2,1},\quad 8_{2,1}/1,3_2/8_{2,1},\quad 15_{2,1}/8_{2,1}/1,3_2/8_{2,1}/15_{2,1},\]
\[ 15_{1,2}/8_{1,2}/1,3_1/8_{1,2}/15_{1,2},\quad 3_2/8_{2,1},16/3_2,\]
with the second module appearing only once. There is only one way to assemble these summands into a module with the right unipotent action and composition factors, and this is
\[(15_{2,1}/8_{2,1}/1,3_2/8_{2,1}/15_{2,1})\oplus (15_{1,2}/8_{1,2}/1,3_1/8_{1,2}/15_{1,2})\oplus (3_2/8_{2,1},16/3_2)\oplus 3_1.\]
In particular, this has $3_1$ as a summand and so this is an $\slf_2$-subalgebra by Proposition \ref{prop:sl2ifsplitoff}, and also $3_2$ as a submodule and subalgebra, but not necessarily a copy of $\slf_2$. This is as described in the proposition, and so completes the proof for $a=2$.

\medskip

Finally, let $a=3$. Using the traces of semisimple elements of order up to $31$ there are 434 conspicuous sets of composition factors, 146 up to field automorphism. We now use the preimage trick from the end of Section \ref{sec:blueprints}: of these 146 sets of composition factors, for 145 of them the eigenvalues of an element of order $21$ determines its class in $E_7$, but for the last one there are two possibilities. (We have a list of all such classes and can check this manually.)

Checking the traces of elements of order $63$, we find that twelve of these sets of factors are not conspicuous for elements of order $63$. In addition, for the set of composition factors whose element of order $21$ does not belong to a single semisimple class in $E_7$, one of those two classes is incompatible with the trace of an element of order $63$. 

The remaining $134$ conspicuous sets of composition factors all have preimages of order $5\cdot 63=315$ that have the same number of eigenspaces on $M(E_7)$. Elements of order $315$ in $\bG$ are blueprints for $M(E_7)$ by Theorem \ref{thm:blueprintsminimal}, and therefore the elements of order $63$ in $H$ are always blueprints for $M(E_7)$, and hence $H$ is as well.
\end{proof}

Theorem \ref{thm:sl2subalgebra} requires $p\geq 11$ for $E_7$, so just because $H$ stabilizes an $\slf_2$-subalgebra of $L(E_7)$ does not mean that $H$ is contained inside a positive-dimensional subgroup of $\bG$, and this could yield a Lie primitive subgroup of $G$.

In the introduction we claimed that this potential $\SL_2(25)$ is Lie primitive if it exists, so we must show that it cannot be contained in any positive-dimensional subgroup. By consideration of composition factors and summand dimensions, it can only lie inside a $D_6$-parabolic subgroup; then one can proceed either by showing that the $12$-dimensional factor cannot support a symmetric bilinear form, or by noting that if the subgroup lies inside the $D_6$-parabolic subgroup then there is another subgroup with the same composition factors on $M(E_7)$ inside the $D_6$-Levi subgroup, hence acting semisimply, but the action of the unipotent element would be $5^8,4^2,2^4$, which does not appear in \cite[Table 7]{lawther1995}.

The next case is $p=7$, where we again cannot prove that there are no maximal $\SL_2(7)$s in all cases.

\begin{proposition}\label{prop:e7char7sl} Suppose that $p=7$.
\begin{enumerate}
\item Let $a=1$. Either $H$ is a blueprint for $M(E_7)$, or $H$ stabilizes a $2$-space on $M(E_7)$, or stabilizes a $1$-space on $L(E_7)$, or the actions of $H$ on $M(E_7)$ and $L(E_7)$ are
\[ P(6)^{\oplus 2}\oplus P(4)\oplus 6\oplus 4^{\oplus 2} \quad\text{and}\quad 7^{\oplus 5}\oplus P(5)^{\oplus 3}\oplus P(3)^{\oplus 3}\oplus 5\oplus 3^{\oplus 3}\]
respectively.
\item Let $a=2$. Either $H$ is a blueprint for $M(E_7)$, or $H$ stabilizes a $2$-space on $M(E_7)$, or stabilizes a $1$-space on $L(E_7)$.
\end{enumerate}
\end{proposition}
\begin{proof} We start with $a=1$. The conspicuous sets of composition factors are
\[ 4,2^{26},\qquad  4^7,2^{14}, \qquad 6,4^9,2^7,\qquad 6^3,4^7,2^5,\] \[6^4,4^3,2^{10},\qquad 6^5,4^5,2^3,\qquad 6^6,4,2^8,\qquad 6^7,4^3,2.\]
As the projective indecomposable modules are
\[ P(2)=2/4,6/2,\qquad P(4)=4/2,4/4,\qquad P(6)=6/2/6,\]
we look through the list above, checking to see whether we have enough $4$s and $6$s (three of the first or two of the second) to cover all $2$s; this leaves the sixth and eighth cases of $6^5,4^5,2^3$ and $6^7,4^3,2$ to deal with.

\medskip

\noindent\textbf{Case 8}: We switch to $L(E_7)$, and there is only one corresponding set of composition factors on $L(E_7)$, namely $7,5^{15},3^{10},1^{21}$, which means that $H$ stabilizes a line on $L(E_7)$.

\medskip

\noindent\textbf{Case 6}: The only possible structure that does not stabilize a $2$-space on $M(E_7)$ and also yields a unipotent action from the list at the start of this section is $P(6)^{\oplus 2}\oplus P(4)\oplus 6\oplus 4^{\oplus 2}$, with $u$ lying in class $E_7(a_5)$. The factors of $L(E_7)\downarrow_H$ aren't uniquely determined, and can be any one of
\[ 7^5,5^{15},3^2,1^{17},\qquad 7^2,5^{18},3^5,1^{14},\qquad 7^8,5^7,3^{12},1^6,\qquad 7^5,5^{10},3^{15},1^3.\]
The first three of these must stabilize a line on $L(E_7)$, but the last one could in theory not, with module action
\[ 7^{\oplus 5}\oplus P(5)^{\oplus 3}\oplus P(3)^{\oplus 3}\oplus 5\oplus 3^{\oplus 3},\]
this action compatible with the action of $u$ on $M(E_7)$.

\medskip

Now let $a=2$, and recall that $L=\SL_2(7)$. The eigenvalues of an element $y$ of order $25$ on $M(E_7)$ are enough to determine the semisimple class of $E_7$ to which $y$ belongs. This allows us to apply Lemma \ref{lem:q/2restricted} to see that if there is an $A_1$ subgroup with $24$-restricted composition factors, then any subgroup $H$ of $\bG$ whose composition factors on $M(E_7)$ match the restriction of this $A_1$ to $\SL_2(49)$ is a blueprint for $M(E_7)$.

There are $150$ conspicuous sets of composition factors for $M(E_7)\downarrow_H$, and this is too many to analyse one at a time, but we can eliminate those of negative $2_i$-pressure using Lemma \ref{lem:pressure}. Of the $150$, only $92$ of them have a corresponding set of composition factors for $L(E_7)$, and only $42$ of these have either no trivial factors or positive pressure. Of these $42$, only $26$ have either no $2_i$ or positive $2_i$-pressure for $i=1,2$, so we have thirteen conspicuous sets of composition factors left to deal with, up to field automorphism. Six of these have a $2_1$ composition factor:
\[ 6_2^5,6_{1,2},4_2^4,2_1^2,\quad 18_{1,2}^2,10_{2,1},4_1,4_2,2_1,\quad 28_{1,2},10_{2,1},6_2^2,4_1,2_1,\]\[ 14_{1,2},10_{1,2},6_2^2,6_{1,2},4_2^3,2_1,\quad 18_{2,1},10_{2,1}^2,6_1,6_{1,2},4_2,2_1,\quad 14_{2,1},10_{2,1}^3,6_{1,2},4_2,2_1.
\]
\medskip

\noindent\textbf{Case 1}: Consider the diagonal $A_1$ inside the $A_1A_1$ maximal subgroup of $\bG$, acting along each factor as $L(1)$: this acts on $M(E_7)$ as
\[ (L(6)\otimes L(3))\oplus (L(2)\otimes L(5))\oplus (L(4)\otimes L(1)),\]
which has factors $L(5)^5,L(3)^4,L(7)^2,L(9)$, up to field automorphism the same as the first case above. 

\medskip

\noindent\textbf{Case 4}: Inside $C_3G_2$, consider an $A_1$ subgroup $\bX$ acting along the first factor as $L(5)$ and along the second as $L(2)\oplus L(8)$. The action of $\bX$ on $M(E_7)$ is
\[ (L(5)/L(7)/L(5))\oplus L(3)\oplus L(11)\oplus L(13)\oplus (L(3)/L(9)/L(3)),\]
so the restriction to $\SL_2(49)$ has composition factors $14_{2,1},10_{2,1},6_1^2,6_{2,1},4_1^3,2_2$, up to field automorphism a match for the fourth case.

\medskip

\noindent\textbf{Case 5}: We note that the composition factors are the same as
\[2_2\otimes ((6_1\otimes 2_2)\oplus (5_1/4_{1,2}/5_1))\oplus 4_2:\]
inside $A_1F_4$, let $\bX$ be a copy of $A_1$ acting along the first factor as $L(7)$ and along the second factor as $L(12)\oplus (L(4)/L(8)/L(4))$. This subgroup of $F_4$ exists inside the $A_1C_3$ subgroup, acting irreducibly on the minimal modules of both subgroups as $L(7)$ and $L(5)$. The composition factors of $M(E_7)\downarrow_{\bX}$ match the fifth case.

\medskip

\noindent\textbf{Case 6}: As in the fifth case, we note that the composition factors are the same as
\[2_2\otimes (7_1\oplus (5_1/4_{1,2}/5_1)\oplus 5_1)\oplus 4_2:\]
inside $A_1F_4$, let $\bX$ be a copy of $A_1$ acting along the first factor as $L(7)$ and along the second factor as $L(6)\oplus (L(4)/L(8)/L(4))\oplus L(4)$. This subgroup of $F_4$ exists inside the $A_1G_2$ subgroup, acting irreducibly on the minimal modules of both subgroups as $L(1)$ and $L(6)$. The composition factors of $M(E_7)\downarrow_{\bX}$ match the sixth case.

Since these embeddings have factors up to $L(21)$, in all cases $H$ is a blueprint for $M(E_7)$ by Lemma \ref{lem:q/2restricted}.

\medskip

\noindent\textbf{Case 3}: Inside $C_3G_2$, let $\bX$ be an $A_1$ subgroup acting along the $C_3$ as $L(1)\oplus L(21)$ and acting along $G_2$ as $L(6)$. The action of $\bX$ on $M(E_7)$ is
\[ L(27)\oplus (L(5)/L(7)/L(5))\oplus L(3)\oplus L(29).\]
Up to field automorphism, the composition factors match the third case. While these are not $24$-restricted, they are close: checking the weight spaces against the eigenvalues of the $\SL_2(49)$ contained within it, all weight spaces that have the same eigenvalues when restricted to $\SL_2(49)$ are contained within the $L(27)\oplus L(29)$. Thus if $H$ is not a blueprint then $H$ is contained in a positive-dimensional subgroup $\bY$ with composition factors of dimension $38,6,6,4,2$. Notice that therefore $\bY$, and hence $H$, must lie in either $C_3G_2$ itself -- and we know from above that in that case $H$ is a blueprint for $M(E_7)$ -- or inside $A_1F_4$, but it is easily seen to not be possible to place $H$ inside this subgroup by the action of $H$ on $M(E_7)$. Thus $H$ is indeed a blueprint for $M(E_7)$.

\medskip

\noindent\textbf{Case 2}: Here if $H$ is not semisimple -- and hence stabilizes a $2$-space of $M(E_7)$ -- then the action of $H$ on $M(E_7)$ is
\[ (18_{1,2}/2_1,4_1/18_{1,2})\oplus 10_{2,1}\oplus 4_2;\]
we claim that such a subgroup $H$ must be Lie primitive, so we check the members of $\mathscr X$ (see Appendix \ref{app:actions}). To see this, first the dimensions of the composition factors are not compatible with coming from any maximal parabolic, so that $H$ must be contained in a reductive maximal subgroup, where the dimensions and multiplicities exclude $A_7$ and $A_2$, and it is easy to see that it doesn't lie in the $A_1A_1$.

For $H$ to lie in $A_1D_6$, $18_{1,2}$ would have to lie in the product of a module for $\SL_2$ of dimension $2$ and a module for $\PSL_2$ of dimension $12$: this is possible, but only with $2_1$ being tensored by $2_1\otimes 6_2$, and this yields $6_2$ as well, which is not in $M(E_7)\downarrow_H$.

If $H$ lies in $A_2A_5$ then we must have that $H$ acts on the natural modules along each factor as $3_1$ and $6_2$, whence the module $L(00)\otimes L(\lambda_3)$ is $(4_2/6_{1,2}/4_2)\oplus 6_2$, which is obviously not correct.

If $H$ lies in $G_2C_3$, then the tensor product of the two minimal modules for these groups must have composition factors $18_{1,2}^2,4_i$ for some $i$, and this is obviously impossible.

For $H$ to lie in $A_1G_2$, $18_{1,2}$ would have to be a composition factor of the tensor product of a module for $\SL_2(49)$ of dimension $4$ and a module for $\PSL_2(49)$ of dimension $7$, but this is not possible.

We finally have $H$ in $A_1F_4$, where $H$ must act on the natural module as $2_i$ for some $i$, yielding $4_i$ as a composition factor of $M(E_7)\downarrow_H$, and the rest of the module must be $2_i\otimes M$ for some $26$-dimensional module $M$ for $\PSL_2(49)$, and this cannot yield $18_{1,2}^2$, so $H$ cannot embed in this subgroup either.

\medskip

Having proved this, we now show that the stabilizer of the $4_2$ submodule is positive dimensional, which yields a contradiction. There are eight elements of order $50$ in a maximal torus of $\bG$ squaring to $y$ of order $25$ and preserving the eigenspaces making up the $4_2$: these eight elements generate a subgroup $Z_{50}\times Z_2\times Z_2$ of order $200$, and so the stabilizer of the $4$-space contains $H$ as a subgroup of index at least $4$, ruling out the possibility that it is almost simple with socle $H$. Thus we now apply Proposition \ref{prop:possiblemaximals}, which states that $\PSL_2(49)$ lies in $\mathscr P$. Hence by Proposition \ref{prop:maximalnotinP}, the stabilizer of the $4$-space is contained in a member of $\mathscr X$, but then $H$ is not Lie primitive, a contradiction.

Thus $H$ acts semisimply on $M(E_7)$, hence stabilizes a $2$-space on $M(E_7)$.

\bigskip

We now give the seven conspicuous sets of composition factors with no $2_i$ in them.
\[18_{2,1}^2,6_1,6_{2,1},4_1^2,\quad 18_{2,1},14_{2,1},10_{2,1},6_{2,1},4_1^2,\quad 42_{1,2},6_{2,1},4_1^2,\quad 28_{1,2},14_{1,2},6_1,4_1^2\]
\[ 12_{2,1}^2,10_{1,2},6_{1,2}^3,4_1,\quad 28_{2,1},12_{2,1},6_{2,1}^2,4_1,\quad 28_{1,2},18_{2,1},10_{1,2}.\]
\noindent \textbf{Cases 1, 2}: Inside $A_2A_5$, let a subgroup $\bX$ of type $A_1$ act along the two factors as $L(14)$ and $L(5)$ respectively. The composition factors of $M(E_7)\downarrow_{\bX}$ are $L(21)^2,L(5),L(3)^2,L(9)$. This is the first case.

Inside $C_3G_2$, consider an $A_1$ subgroup $\bX$ acting along the first factor as $L(5)$ and along the second as $L(14)\oplus L(8)$. The action of $\bX$ on $M(E_7)$ is
\[ L(19)\oplus L(13)\oplus L(11)\oplus (L(3)/L(9)/L(3)).\]
The composition factors of $\bX$ on $M(E_7)$ match the second case.

\noindent \textbf{Case 3}: Inside $C_3G_2$, consider an $A_1$ subgroup $\bX$ acting along the first factor as $L(5)$ and along the second as $L(42)$. The action of $\bX$ on $M(E_7)$ is
\[ L(47)\oplus (L(3)/L(9)/L(3)),\]
a match for the third case, but of course these are not $24$-restricted, but do satisfy the second condition of Lemma \ref{lem:q/2restricted}, so that $H$ is a blueprint for $M(E_7)$.

\medskip

\noindent \textbf{Case 4}: Inside $F_4A_1$, let $\bX$ denote an $A_1$ subgroup acting along the second factor as $L(1)$, and along the first factor as $L(4)\oplus L(44)$, which exists inside the $A_1G_2$ subgroup of $F_4$. The action of $\bX$ on $M(E_7)$ is
\[  L(43)\oplus L(45)\oplus L(5)\oplus L(3)^{\oplus 2},\]
so this is the fourth case. Of course, these are not $24$-restricted, so we proceed as in the second case of the previous set of composition factors, looking for elements of order $50$ in $\bG$ squaring to $y$ and stabilizing the eigenspaces the comprise the $4_1$ in the socle. Again, we find a subgroup $Z_{50}\times Z_2\times Z_2$, and we conclude as before that the stabilizer of the $4_1$ is a positive-dimensional subgroup of $\bG$.

We claim that $H$ is a blueprint for $M(E_7)$. With the dimensions and multiplicities of the composition factors, the only maximal positive-dimensional subgroups it can lie in are $D_6A_1$, $C_3G_2$, $A_1G_2$, $A_1F_4$ and $A_1A_1$, with the last one clearly impossible.

If $H\leq D_6A_1$ then $14_{1,2}\oplus 6_1\oplus 4_1$ is a tensor product of a $12$-dimensional and a $2$-dimensional module, so must be $2_1\otimes (7_2\oplus 5_1)$. Thus $H$ lies inside the product of the $A_1$ and a product of two orthogonal groups, $\Spin_7\times \Spin_5$, and there is a unique action of an $A_1$ subgroup inside these of acting as $L(42)\oplus L(0)$ and $L(42)$ on the two relevant modules of $\Spin_7$, and as $L(3)$ and $L(4)$ on the two modules of $\Spin_5$. This $A_1$ stabilizes the same subspaces of $M(E_7)$ as $H$, so $H$ is a blueprint for $M(E_7)$. The same statement holds from above for $A_1F_4$.

If $H\leq C_3G_2$ then a similar analysis shows that $H$ acts on the minimal modules of the two factors as $2_1\oplus 4_1$ and $7_2$ respectively, and the $A_1$ acting along each factor as $L(1)\oplus L(3)$ and $L(42)$ again stabilizes the same subspaces of $M(E_7)$ as $H$, so $H$ is a blueprint for $M(E_7)$.

We are left with $A_1G_2$, where in order to find the $28_{1,2}$ we must have that $H$ acts on the two natural modules as $2_1$ and $7_2$ respectively, so that the factor $(L(3),L(10))$ in $M(E_7)\downarrow_{A_1G_2}$ yields $28_{1,2}$, but then the other factor of $(L(1),L(01))$ yields two copies of $6_{1,2}$, which is not correct. (Indeed, this is how we obtain the sixth case above.)

Thus, whenever $H$ is a subgroup of a positive-dimensional subgroup of $\bG$, it is a blueprint for $M(E_7)$, as needed.

\medskip

\noindent \textbf{Case 5}: Inside the maximal subgroup $A_1G_2$, let $\bX$ be an $A_1$ subgroup acting along the first factor as $L(7)$ and along the second as $L(2)^{\oplus 2}\oplus L(0)$: the composition factors of $\bX$ on $M(E_7)$ are $L(9)^3,L(11),L(21),L(23)^2$, matching up with the fifth case, up to field automorphism. Since they are all $24$-restricted, $H$ is a blueprint for $M(E_7)$ by Lemma \ref{lem:q/2restricted}.

\medskip

\noindent \textbf{Case 6}: Inside the same subgroup $A_1G_2$, let $\bX$ instead be an $A_1$ subgroup acting along the first factor as $L(7)$ and along the second as $L(6)$: the composition factors of $\bX$ on $M(E_7)$ are $L(3),L(9)^2,L(17),L(27)$, matching up with the sixth case, but no longer $24$-restricted, but $H$ is still a blueprint for $M(E_7)$ by Lemma \ref{lem:q/2restricted}(\ref{lemi:q2resb}).

\medskip

\noindent \textbf{Case 7}: Inside the maximal $A_1A_1$ subgroup, take a diagonal $A_1$ as we have done before, but this time acting as $L(1)$ and $L(7)$ along the two factors. The composition factors of this on $M(E_7)$ are $L(27),L(13),L(37)$, of course not $24$-restricted, and even Lemma \ref{lem:q/2restricted}(\ref{lemi:q2resb}) doesn't work in this case. If $H$ is contained in a positive-dimensional subgroup other than $A_1A_1$ then the dimensions of the composition factors and multiplicities show that it can only come from $A_1G_2$, and in order to get $28_{1,2}$ appearing, $H$ must act along $A_1$ as $2_1$ and $G_2$ as $7_2$. However, the other factor of dimension $28$ must have $6_{1,2}$ as a composition factor, which is not allowed. Thus $H\leq A_1A_1$, and so $H$ is a blueprint for $M(E_7)$.

It remains to show that $H$ is always contained inside a member of $\mathscr X$. Of course, since $M(E_7)\downarrow_H$ is multiplicity free it is semisimple, and so the $10_{1,2}$ is a submodule. As with previous cases, we find more than one element of order $50$ in $\bG$ squaring to $y$ and stabilizing the eigenspaces that comprise $10_{1,2}$. The subgroup generated by these is $Z_{50}\times Z_2$, and we wish to apply Corollary \ref{cor:sl2nopgl}, so we need to find an element $\hat y$ of order $50$ in this subgroup whose action on $L(E_7)$ has no $(-1)$-eigenspace, but this is easy: its eigenspaces are
\[ 1^7,(\theta^{\pm 2})^6,(\theta^{\pm 4})^4,(\theta^{\pm 5}),(\theta^{\pm 6})^2,(\theta^{\pm 7})^3,(\theta^{\pm 8}),(\theta^{\pm 9})^5,(\theta^{\pm 11})^6,(\theta^{\pm 13})^6,\]
\[(\theta^{\pm 15})^5,(\theta^{\pm 16}),(\theta^{\pm 17})^3,(\theta^{\pm 18})^2,(\theta^{\pm 19})^2,(\theta^{\pm 20})^4,(\theta^{\pm 21}),(\theta^{\pm 22})^5,(\theta^{\pm 24})^6,\]
where $\theta$ is a primitive $50$th root of unity and $y=\hat y^2$ acts on $2_1$ with eigenvalues $\theta^{\pm 2}$.

\bigskip

We have thus shown that all cases are blueprints for $M(E_7)$, stabilize a $2$-space on $M(E_7)$, or stabilize a line on $L(E_7)$, as needed.
\end{proof}

\begin{proposition}\label{prop:e7char11sl} Suppose that $p=11$.
\begin{enumerate}
\item Let $a=1$. If $H$ is not a blueprint for $M(E_7)$, then $H$ stabilizes a unique $2$-space on $M(E_7)$ or a line on $L(E_7)$.
\item Let $a=2$. If $H$ is not a blueprint for $M(E_7)$, then $H$ stabilizes a unique $2$-space on $M(E_7)$.
\end{enumerate}
\end{proposition}
\begin{proof} Let $a=1$. As $p=11$ the action of $u$ is one of cases (\ref{li:sl2first11}) to (\ref{li:sl2last11}) in the list above. In the first unipotent class, $E_7(a_4)$ acting as $11^2,10,8,6,4^2,2$, there are single blocks of size $10$, $8$, $6$ and $2$, which must come from simple summands of those dimensions. Since we cannot have a faithful indecomposable module of dimension $11+4=15$, the $4$s must also come from simple summands, and so $M(E_7)\downarrow_H$ is a single projective plus a semisimple module. The conspicuous such sets of composition factors yield
\[ P(10)\oplus 10\oplus 8\oplus 6\oplus 4^{\oplus 2}\oplus 2\qquad\text{and}\qquad P(4)\oplus 10\oplus 8\oplus 6\oplus 4^{\oplus 2}\oplus 2.\]
The second of these has corresponding set of factors on $L(E_7)$ given by $9^8,7^2,5^7,1^{12}$, and the action of $u$ on $L(E_7)$ is $11^8,9,7^2,5^2,3^4$. Blocks of size $3$ come from, up to duality,
\[ 3,\qquad 7,9/1,3,5,\qquad 3,5,7/5,7,9,\]
and so $H$ cannot embed with these factors and this action of $u$. (The other set of composition factors yields
\[ 11^{\oplus 4}\oplus P(9)\oplus P(7)\oplus 9\oplus 7^{\oplus 2}\oplus 5^{\oplus 2}\oplus 3^{\oplus 4}.)\]

\medskip

If $u$ comes from class $D_6$, acting as $11^4,10,1^2$, then the $10$ must come from a simple summand, and the $1^2$ comes from (up to duality) $6/6$, $4/8$ or $2/10$. The sum of one of these plus its dual, a single projective indecomposable module, and the $10$, must be conspicuous. There are five such conspicuous sets of composition factors, all of which have corresponding sets of factors on $L(E_7)$, three of them having two different sets. The action of $u$ on $L(E_7)$ is $11^{10},10^2,1^3$, so certainly $H$ has a trivial summand on $L(E_7)$, with the $10^2$ coming from $(5/5)^{\oplus 2}$, $(3/7)\oplus (7/3)$ or $(1/9)\oplus (9/1)$, the other $1^2$ being either semisimple, $(9/3)\oplus (3/9)$ or $(5/7)\oplus (7/5)$, with the rest of the module being projective. There are $937$ such sets of composition factors, and when taking the intersection of that list with those of the corresponding sets of composition factors to our list for $M(E_7)\downarrow_H$, we find two members:
\[ P(4)\oplus (10/2)\oplus (2/10)\oplus 10\qquad\text{and}\qquad P(6)\oplus (6/6)^{\oplus 2}\oplus 10,\]
with corresponding embeddings
\[ 11^{\oplus 4}\oplus P(5)\oplus P(3)\oplus (3/9)\oplus (9/3)\oplus (3/7)\oplus (7/3)\oplus 1\]
and
\[ 11^{\oplus 6}\oplus P(5)\oplus P(3)\oplus (3/7)\oplus (7/3)\oplus 1^{\oplus 3}.\]
Of course, these stabilize a line on $L(E_7)$, as needed.

\medskip

If $u$ comes from class $E_6(a_1)$ acting as $11^4,5^2,1^2$, then the two blocks of size $5$ must come from summands of dimension $16$, so $(4,6/6)\oplus (6/4,6)$, and the two $1$s must come from summands of dimension $12$. The conspicuous such sets of composition factors yield
\[ (4,6/6)\oplus (6/4,6)\oplus (10/2)\oplus (2/10),\qquad (4,6/6)\oplus (6/4,6) \oplus (4/8)\oplus (8/4),\]
\[(4,6/6)\oplus (6/4,6)\oplus (6/6)^{\oplus 2}.\]
Of these, only the last has a corresponding set of composition factors on $L(E_7)$, and this is $11^2,9^8,5^5,1^{14}$; however, $u$ acts as $11^{10},9,5^2,3,1$ on $L(E_7)$, so we need a $3$ as a summand of $L(E_7)\downarrow_H$, which is not possible. Thus $H$ does not embed with $u$ from this class.

\medskip

Finally, if $u$ acts as $11^4,10,2$, coming from $E_7(a_3)$, then the $2$ and $10$ must come from simple summands, so we need two projectives plus $10\oplus 2$. There are three such conspicuous sets of composition factors, yielding
\[ P(4)^{\oplus 2}\oplus 10\oplus 2,\qquad P(4)\oplus P(10)\oplus 10\oplus 2,\qquad P(6)\oplus P(10)\oplus 10\oplus 2.\]
Each of these has corresponding sets of composition factors on $L(E_7)$, with the third having two. However, $u$ acts on $L(E_7)$ with blocks $11^{11},9,3$, so $L(E_7)\downarrow_H$ is the sum of $9\oplus 3$ and projectives, and since we have $11^{11}$ in the action of $u$, the number of summands of $L(E_7)\downarrow_H$ that are either $11$s or $P(1)$s must be odd. In particular, any trivial composition factors lie either in $P(1)$s or $P(9)$s. The four corresponding sets of factors are
\[ 9^8,7^3,5^7,1^{12},\quad 11^4,9^3,7^4,5^3,3^6,1,\quad 11^3,9^2,7^6,5^5,3^5,\quad 11^2,9^4,7^7,5,3^6,1^3;\]
the first and second cases need an odd number of $P(1)$s, and the second case can have no $P(9)$s as it has no $3$s, leading to a contradiction. The fourth case cannot work with the unipotent class either, but the third case yields
\[ 11^{\oplus 3}\oplus P(7)^{\oplus 2}\oplus P(5)\oplus P(3)\oplus 9\oplus 3.\]
Thus $P(6)\oplus P(10)\oplus 10\oplus 2$ is the only acceptable embedding of $H$ into $M(E_7)$ with this unipotent action.

\bigskip

Now let $a=2$, and recall that $L$ is a copy of $\SL_2(11)$ inside $H$. We have traces of semisimple elements of order up to $40$, and can use the preimage trick from Section \ref{sec:blueprints} to find traces of elements of orders $60$ and $120$ as well. We would like that the eigenvalues of an element $y$ in $H$ of order $120$ on $M(E_7)$ uniquely determine the semisimple class of $E_7$ to which $y$ belongs. This is not true in general for all classes, but will be true for the particular classes that arise from conspicuous sets of composition factors for $M(E_7)\downarrow_H$, by a check using the preimage trick.

Suppose that $u$ comes from class $E_7(a_4)$, so that the composition factors of $M(E_7)\downarrow_L$ are $10^3,8,6,4^2,2$. There are, up to field automorphism, twenty conspicuous sets of composition factors for $M(E_7)\downarrow_H$, using semisimple elements of order up to $40$; using the preimage trick, we can eliminate seven of these from contention, as they fail the trace of an element of order either $60$ or $120$, leaving thirteen. Note that, for these thirteen remaining sets of composition factors, the eigenvalues of an element $y$ of order $120$ on $M(E_7)$ determine the semisimple class to which $y$ belongs.

Nine of the thirteen have a $2$-dimensional composition factor, and are
\[ 14_{1,2},10_2^3,4_2^2,2_1^2,\quad 10_1^3,8_1,6_1,4_1^2,2_1,2_2,\quad 22_{2,1},18_{2,1}^{(1)},6_1,4_1^2,2_1,\] \[22_{2,1},14_{2,1},10_1,4_1^2,2_1,\quad 10_1^3,8_1,6_1,6_{2,1},4_1,2_1,\quad
30_{1,2}^{(1)},10_2,10_{1,2},4_1,2_1,\] \[22_{1,2},18_{1,2}^{(1)},10_{1,2},4_1,2_1,\quad 18_{1,2}^{(1)},10_2^2,6_2,6_{1,2},4_2,2_1,\quad 22_{1,2},18_{1,2}^{(2)},10_2,4_2,2_1.\]
(Recall that $18_{i,j}^{(1)}=2_i\otimes 9_j$, $18_{i,j}^{(2)}=3_i\otimes 6_j$ and $30_{i,j}^{(1)}=3_i\otimes 10_j$.)
The simple modules with non-trivial extensions with $2_1$ are $10_2$, $18_{2,1}^{(1)}$ and $30_{1,2}^{(1)}$, and in order for $H$ not to stabilize a $2$-space, one of these must occur with multiplicity $2$: thus all cases must stabilize a $2$-space on $M(E_7)$ except for the first and eighth. However, even in the first case the $\{2_1,4_2,10_2,14_{1,2}\}$-radical of $P(10_2)$ is simply $10_2/2_1/10_2$, so there must be a $2_1$ submodule of $M(E_7)\downarrow_H$, so that case also stabilizes a $2$-space.

For the eighth case, up to field automorphism we find this inside $D_6A_1$, by taking an $A_1$ subgroup acting on the natural modules for the two factors as $9_1\oplus 3_1=L(8)\oplus L(2)$ and $2_2=L(11)$ respectively. This yields $18_{2,1}^{(1)}\oplus 6_{2,1}=L(19)\oplus L(13)$, and for the spin module for the $D_6$ term, we need a module with unipotent action $11^2,6,4$ (one sees this from the entry for $D_6(a_1)$ in \cite[Table 7]{lawther1995}) and composition factors $10_1^2,6_1,4_1,2_2=L(9)^2,L(5),L(3),L(11)$ (obtained from the traces of semisimple elements). Thus the restriction of the spin module to this $A_1$ subgroup must be
\[ (L(9)/L(11)/L(9))\oplus L(5)\oplus L(3),\]
so we apply Lemma \ref{lem:q/2restricted} to see that $H$ is a blueprint for $M(E_7)$, as needed.

The remaining four conspicuous sets of composition factors, which have no $2$-dimensional composition factor, are
\[ 36_{2,1},10_1,6_{1,2},4_1,\quad 42_{1,2},10_1,4_1,\quad 22_{2,1},10_1,10_{2,1},8_1,6_{1,2},\quad 28_{2,1},22_{2,1},6_{2,1}.\]
For the last of these, consider a diagonal $A_1$ inside $A_1G_2$, acting as $2_2=L(11)$ along $A_1$ and as $7_1=L(6)$ along the $G_2$ factor. The composition factors on $M(E_7)$ are $L(39),L(21),L(13)$, matching up with the fourth case above, so we satisfy the conditions of Lemma \ref{lem:q/2restricted}. Consider an $A_1$ subgroup of $G_2C_3$, acting as $7_2=L(66)$ along the $G_2$ factor and as $6_1=L(5)$ along the $C_3$ factor: the composition factors on $M(E_7)$ are $L(71),L(9),L(3)$, matching up with the second case above. Again we apply Lemma \ref{lem:q/2restricted}, this time the second statement, and therefore $H$ is a blueprint for $M(E_7)$, as claimed.

For the third case, we find an $A_1$ subgroup $\bY$ inside $D_6A_1$ that works: consider $\bY$ acting as $2_1=L(1)$ along the $A_1$ factor, and as $9_1\oplus 3_2=L(8)\oplus L(22)$ along the second factor. This second $A_1$ is contained diagonally as an irreducible subgroup inside the product of orthogonal groups $\Spin_9\times \Spin_3$, i.e., $B_4A_1$, so its action on the $32$-dimensional half-spin module is as the tensor product of the spins. The action on $\Spin_3$ must be as $2_2=L(11)$, and the action on the $\Spin_9$ has unipotent factors $11,5$ and can be seen to be $11_1\oplus 5_1=L(10)\oplus L(4)$. Thus the subgroup $\bY$ has composition factors
\[ L(21),L(15),L(9),L(7),L(23),\]
and hence satisfies Lemma \ref{lem:q/2restricted}, with the $\SL_2(121)$ inside $\bY$ having the same factors on $M(E_7)$ as the third case.

We are left with $36_{2,1},10_1,6_{1,2},4_1$. We first note that if there is such a subgroup $H$ then $H$ is not contained in a positive-dimensional subgroup of $\bG$: to see this, notice that since $M(E_7)\downarrow_H$ is multiplicity-free and contains a $36$-dimensional composition factor, the only positive-dimensional subgroups that could contain it are $G_2C_3$ and $A_1F_4$. If $H\leq G_2C_3$, the module $(L(10),L(100))$ has dimension $42$, so must restrict to $H$ as $36_{2,1}\oplus 6_{1,2}$. However, $36_{2,1}=9_1\otimes 4_2$ is not a composition factor of any tensor product of a $6$-dimensional module and a $7$-dimensional module, which is a contradiction. Since $A_1F_4$ acts on $M(E_7)$ with factors $(L(1),L(0001))$ and $(L(3),L(0000))$, we see that if $H\leq A_1F_4$ then the projection of $H$ along $A_1$ must act on the natural module as $2_1$, so that $L(3)$ restricts to $H$ as $4_1$. However, again $36_{2,1}$ is not a composition factor of any tensor product of $2_1$ and a module of dimension at most $26$, so $H$ cannot lie in $A_1F_4$ either. Thus $H$ cannot lie in a positive-dimensional subgroup of $\bG$.

We therefore must have that $N_{\bG}(H)$ contains $\SL_2(121)$ with index at most $2$ by Proposition \ref{prop:possiblemaximals} and the fact that the composition factors of $M(E_7)\downarrow_H$ are not stable under a field automorphism of $H$, so $\Aut_{\bar G}(H)$ can only induce a diagonal automorphism on $H$. Recalling that $y$ is an element of order $120$ in $H$, chosen so that the eigenvalues of $y$ on $2_1$ are $\zeta^{\pm 1}$ for $\zeta$ a primitive $120$th root of unity, we consider the elements of order $120$ in a maximal torus of $\bG$ squaring to $y^2$, noting that the $\zeta^2$- and $\zeta^6$-eigenspaces of $y^2$ on $M(E_7)$ are both $3$-dimensional and coincide with the $\zeta$- and $\zeta^3$-eigenspaces of $y$ respectively. We thus look for elements that square to $y^2$ and preserve the $\zeta^2$- and $\zeta^6$-eigenspaces of $y^2$; of course, $y$ is one of these elements, and we find four elements with $3$-dimensional $\zeta$- and $\zeta^3$-eigenspaces (and therefore four with $3$-dimensional $(-\zeta)$- and $(-\zeta^3)$-eigenspaces), which together generate a subgroup $Z_{120}\times Z_2\times Z_2$ of the torus. Thus the stabilizer in $\bG$ of the $4$-dimensional submodule $4_1$ of $M(E_7)$ contains $H$ with index at least $4$, a contradiction, and so $H$ doesn't exist. This completes the proof of the proposition when $u$ comes from class $E_7(a_4)$.

\medskip

Suppose that $u$ comes from class $D_6$, and that the composition factors of $M(E_7)\downarrow_L$ are $10,6^7,4$: the trace of an element of order $5$ is $6$, and this is enough to seriously restrict the possibilities. Using other semisimple elements of order up to $20$, we find up to field automorphism a single conspicuous set of composition factors, namely $10_1,6_1^7,4_1$, and this must be semisimple as there are no non-trivial extensions between the factors. However, this is incompatible with the action of $u$, so $H$ cannot embed with this restriction to $L$. The other set of composition factors are the same as for $E_7(a_4)$, which we have already considered above.

Suppose that $u$ comes from $E_7(a_3)$, so that $M(E_7)\downarrow_L$ has composition factors $10^3,6^3,4,2^2$. Checking traces of elements of order up to $40$ yields, up to field automorphism, only five conspicuous sets of composition factors, which are
\[ 10_1^3,6_1^3,4_1,2_1^2,\qquad 10_2^3,6_2^3,4_2,2_1^2,\qquad 10_1^3,10_{2,1},6_1^2,2_1,2_2,\] \[10_1,10_2^2,6_1,6_2^2,6_{1,2},2_1,\qquad 22_{1,2},10_2,10_{1,2},6_2^2,2_1.\]
The last two of these fail the traces of elements of order $60$, so do not exist. The $2_1$-pressures of the remaining three are $-2$, $1$ and $-1$ respectively, as only $10_2$ from these simple modules has an extension with $2_1$, so only the second need not stabilize a $2$-space. In this case the $\{2_1,4_2,6_2,10_2\}$-radical of $P(10_2)$ is $10_2/2_1/10_2$, so $M(E_7)\downarrow_H$ has a $2$-dimensional submodule.
\end{proof}

We now have that $p\geq 13$, for which we only need consider $a=1$. For $p=17,19,23$ we will get an $\slf_2$-subalgebra being a possible outcome, and for $p=19$ we will also get a Serre embedding (see Definition \ref{defn:serreembedding}).

We begin with $p=13$.

\begin{proposition}\label{prop:e7char13sl} Suppose that $p=13$ and $a=1$. Then $H$ is strongly imprimitive.
\end{proposition}
\begin{proof} In all but one case, we will show that $H$ is either a blueprint for $M(E_7)$, or $H$ stabilizes a $2$-space on $M(E_7)$, or stabilizes a line on $L(E_7)$. In these cases, $H$ is strongly imprimitive by Propositions \ref{prop:blueprintissi}, \ref{prop:fix2spaceonMG} and \ref{prop:fixlineonLG} respectively. Thus we will show one of these three properties, except the last case, where we argue directly.

There are three possibilities for the action of $u$ on $M(E_7)$, namely cases (\ref{li:sl2first13}), (\ref{li:sl2second13}) and (\ref{li:sl2last13}) from the list at the start of the section, with the last of these being semiregular.

In the first case, $u$ acts as $13^4,1^4$, and so $M(E_7)\downarrow_H$ is the sum of four modules of dimension $14$, which are
\[ 2/12,\quad 4/10,\quad 6/8,\quad 8/6,\quad 10/4,\quad 12/2.\]
There are only two conspicuous sets of composition factors for $M(E_7)\downarrow_H$ consisting of dual pairs of these, and they are
\[ (2/12)\oplus (12/2)\oplus (4/10)\oplus (10/4)\qquad\text{and}\qquad (4/10)\oplus (10/4)\oplus (6/8)\oplus (8/6).\]
Neither of these has a corresponding set of composition factors on $L(E_7)$, so $H$ does not embed in $\bG$ with $u$ coming from class $E_6$.

\medskip

If $u$ comes from class $E_7(a_3)$ then it acts as $13^2,12,10,6,2$. The single block of size $2$ must come from a self-dual indecomposable module of dimension congruent to $2$ modulo $13$, and the two of these are $2$ itself -- so $H$ stabilizes a $2$-space on $M(E_7)$ -- or a $28$-dimensional module $6,8/6,8$; from here, the blocks of sizes $6$ and $10$ must come from simple summands and the $12$ comes either from a $12$ or a $6/6$, yielding two possible sets of composition factors, neither of which is conspicuous, having trace $-1$ for an element of order $3$. If $H$ stabilizes more than a single $2$, then we must also have $P(2)$, and the block of size $12$ comes from either a simple $12$ or a $6/6$. In either case, again the trace of an element of order $3$ is $-1$, so $H$ stabilizes a unique $2$-space on $M(E_7)$. This completes the proof for the second action of $u$.

\medskip

The final unipotent class to consider is the semiregular $E_7(a_2)$, acting with Jordan blocks $13^4,4$ on the minimal module. The single block of size $4$ comes either from a summand $4$ or from the indecomposable module $4,6,8,10/4,6,8,10$, which is conspicuous, but we saw above that it has no corresponding set of factors for $L(E_7)$. Hence $M(E_7)\downarrow_H$ has a $4$ as a summand, with two projective indecomposable summands. The conspicuous such sets of composition factors yield
\[ P(2)\oplus P(4)\oplus 4,\qquad P(12)\oplus P(10)\oplus 4,\qquad P(10)\oplus P(8)\oplus 4,\qquad P(6)\oplus P(4)\oplus 4.\]
The first and second of these cannot occur because they do not have corresponding factors on $L(E_7)$.

In the fourth case we again switch to $L(E_7)$, and find two corresponding sets of composition factors, namely
\[ 13,11,9^2,7^9,5^3,3^4,1,\qquad 11^3,9^3,7^5,5^4,3^6.\]
The action of $u$ on $L(E_7)$ must be $13^{10},3$, and the single $3$ in this action comes from a summand isomorphic to either $3$ or $5,7,9/5,7,9$. The first case cannot occur as the single $1$ must lie in a $P(11)$, but this cannot occur. In the second case, the lack of trivial factors means there can be no $P(11)$s, so we must have $P(3)^{\oplus 3}$ in $L(E_7)\downarrow_H$. In this case, there are then no $3$s or $9$s remaining, so the summand contributing the $3$ to the action of $u$ cannot occur, which is a contradiction. Thus $H$ cannot embed with these factors either.

In the third case, there are again two corresponding sets of composition factors on $L(E_7)$, namely
\[  13,11^3,9^2,7^7,5,3^4,1^3,\qquad  11^5,9^3,7^3,5^2,3^6,1^2.\]
In the first of these, the single $5$ means one has no $P(5)$ and at most one $P(7)$, but these are the only two projectives containing $7$, so we cannot use up the seven $7$s. The second case does have a unique possibility, however, of
\[ P(11)^{\oplus 2}\oplus P(9)\oplus P(7)\oplus P(3)\oplus 3.\]
Although it has $3$ as a summand, the presence of a $P(3)$ means that we cannot guarantee that it is an $\slf_2$-subalgebra of $L(E_7)$ using Proposition \ref{prop:sl2ifsplitoff}, although the $3\oplus 3$ in the socle of $L(E_7)\downarrow_H$ does form a subalgebra.

Let $x$ denote an element of order $14$ in $H$ and let $\zeta$ be a primitive $14$th root of unity, arranged so that the eigenvalues of $x$ on $4$ are $\zeta^{\pm 1},\zeta^{\pm 3}$. The eigenvalues of $x$ on $M(E_7)$ are
\[ (-1)^8,(\zeta^{\pm 1})^9,(\zeta^{\pm 3})^8,(\zeta^{\pm 5})^7.\]
Let $\theta$ denote a primitive $28$th root of unity with $\theta^2=\zeta$. Looking through the eigenvalues of elements of order $28$ in $\bG$, we find $\hat x\in \bG$ such that $\hat x^2=x$ and $\hat x$ has eigenvalues
\[ (\pm \I)^4,(\theta^{\pm 1})^9,(\theta^{\pm 3})^8,(\theta^{\pm 5})^5,(-\theta^{\pm 5})^2\]
on $M(E_7)$. This stabilizes the eigenspaces intersecting the $4$, and so $\hat x$ stabilizes the $4$-space stabilized by $H$.

Let $K=\gen{H,\hat x}$. Since $H$ is a maximal member of $\mathscr P$ by Proposition \ref{prop:possiblemaximals}, we wish to apply Proposition \ref{prop:maximalnotinP}. Thus either $K\leq N_{\bG}(H)$ or $H$ is strongly imprimitive, since $H$ and $K$ stabilize a unique $4$-space on $M(E_7)$ and so is $N_{\Aut^+(\bG)}(H)$-stable. Since $\Aut(H)\cong \PGL_2(13)$, if we show that $K\neq \PGL_2(13)$ then we are done.

The action of $\hat x$ on $L(E_7)$ has eigenvalues
\[ (1)^{21},(\theta^{\pm 2})^{18},(-\theta^{\pm 2}),(\theta^{\pm 4})^{14},(-\theta^{\pm 4})^4,(\theta^{\pm 6})^{13},(-\theta^{\pm 6})^6.\]
In particular, $\hat x$ has no eigenvalue $-1$ on $L(E_7)$, so $K\not\cong \PGL_2(13)$ by Corollary \ref{cor:sl2nopgl}. Hence $K\not\leq N_{\bG}(H)$, so $H$ is strongly imprimitive.
\end{proof}

\begin{proposition}\label{prop:e7char17sl} If $p=17$ and $a=1$, then $H$ is strongly imprimitive.
\end{proposition}
\begin{proof} Note that if $H$ is a blueprint for $M(E_7)$ then $H$ is strongly imprimitive by Proposition \ref{prop:blueprintissi}. We will show that this is the case, or will argue directly that $H$ is strongly imprimitive.

There are two non-generic unipotent classes, cases (\ref{li:sl2first17}) and (\ref{li:sl2last17}) above, where there are single Jordan blocks of sizes $4,6,8,10,16$. Apart from the simple modules of dimension congruent to $4,6,8,10,16$ modulo $17$, the self-dual indecomposable modules congruent to those dimensions have dimensions $72$, $108$, $144$, $114$ and $16$ respectively, so only the $16$ might not come from a simple summand.

For $u$ belonging to class $E_7(a_2)$, so acting as $17^2,10,8,4$, we therefore have a single projective plus $10\oplus 8\oplus 4$. Applying the traces of semisimple elements yields two possibilities:
\[ P(16)\oplus 10\oplus 8\oplus 4\qquad\text{and}\qquad P(10)\oplus 10\oplus 8\oplus 4.\]
The second of these does not yield an action on $L(E_7)$ as the traces do not match up, but the first of these has a unique set of composition factors on $L(E_7)$ which yields
\[ 17\oplus P(15)\oplus P(11)\oplus 15\oplus 11\oplus 9\oplus 7\oplus 3^{\oplus 2}.\]
The subspace $3^{\oplus 2}$ is an $H$-invariant Lie subalgebra of $L(E_7)$, but we proceed as in the case of $p=13$, finding an element of order $36$ that preserves the $4$-dimensional submodule.

Let $x$ be an element of $H$ of order $18$ and $\zeta$ be a primitive $18$th root of unity, arranging our choices so that the eigenvalues of $x$ on $4$ are $\zeta^{\pm 1}$ and $\zeta^{\pm 3}$. The eigenvalues of $x$ on $M(E_7)$ are
\[ (-1)^6,(\zeta^{\pm 1})^6,(\zeta^{\pm 3})^7,(\zeta^{\pm 5})^6,(\zeta^{\pm 7})^6.\]
Letting $\theta$ denote a primitive $36$th root of unity squaring to $\zeta$, we find an element $\hat x$ of order $36$ in $\bG$ with $\hat x^2=x$ and with eigenvalues on $M(E_7)$ given by
\[ (\pm \I)^3,(\theta^{\pm 1})^6,(\theta^{\pm 3})^7,(\theta^{\pm 5})^5,-\theta^{\pm 5},(\theta^{\pm 7})^4,(-\theta^{\pm 7})^2.\]
We see immediately that $\hat x$ preserves the $4$-space stabilized by $H$, and we proceed exactly as in the proof of Proposition \ref{prop:e7char13sl}.

The eigenvalues of $\hat x$ on $L(E_7)$ are
\[ (1)^{15},(\theta^{\pm 2})^{14},(-\theta^{\pm 2}),(\theta^{\pm 4})^{12},(-\theta^{\pm 4})^3,(\theta^{\pm 6})^{11},(-\theta^{\pm 2})^3,(\theta^{\pm 8})^8,(-\theta^{\pm 8})^7.\]
Again, this has no eigenvalue $-1$, so $K=\gen{H,\hat x}$ is not $\PGL_2(17)$ modulo $Z(\bG)$. Thus $K$ is strongly imprimitive since $H$ is a maximal member of $\mathscr P$, by Proposition \ref{prop:possiblemaximals}, and hence $H$ is strongly imprimitive by Proposition \ref{prop:maximalnotinP}.

\medskip

For $u$ belonging to class $E_7(a_1)$, so acting on $M(E_7)$ with Jordan blocks $17^2,16,6$, the $6$ must come from a simple summand, but the $16$ comes from either a simple summand or $8/8$. Thus our embedding of $H$ is either a single projective plus $(8/8)\oplus 6$ or a single projective plus $16\oplus 6$. Using traces, the two options are
\[ P(12)\oplus 16\oplus 6\qquad\text{and}\qquad P(4)\oplus 16\oplus 6.\]
The second of these has no corresponding set of composition factors on $L(E_7)$, but the first has a unique set, which implies that $L(E_7)\downarrow_H$ is
\[ 17\oplus P(15)\oplus P(11)\oplus P(7)\oplus 11\oplus 3,\]
and so the $3$ is an $\slf_2$-subalgebra by Proposition \ref{prop:sl2ifsplitoff}. We complete the proof that $H$ is strongly imprimitive using Corollary \ref{cor:sl2coxeter-1}.
\end{proof}

\begin{proposition}\label{prop:e7char19sl}
Suppose that $p=19$ and $a=1$. If $H$ is not a blueprint for $M(E_7)$ then $H$ stabilizes a unique $3$-space on $L(E_7)$ that is an $\slf_2$-subalgebra of $L(E_7)$, or $H$ is a Serre embedding.
\end{proposition}
\begin{proof}
When $p=19$, there are two non-generic unipotent classes, $E_7(a_1)$ and $E_7$, cases (\ref{li:sl2first19}) and (\ref{li:sl2last19}) above. As $p\equiv 3\bmod 4$ there is a unique self-dual indecomposable module congruent to any given integer modulo $p$. For $E_7(a_1)$ the $12$ and $6$ in the action of $u$ must therefore come from simple summands, leaving a single projective module of dimension $38$. Only two possibilities yield conspicuous sets of composition factors, namely $P(16)\oplus 12\oplus 6$ and $P(4)\oplus 12\oplus 6$. The second of these has no corresponding set of composition factors on $L(E_7)$, with the first of these yielding the unique action
\[ P(15)\oplus P(11)\oplus 19\oplus 17\oplus 11\oplus 7\oplus 3,\]
with the $3$ being an $\slf_2$-subalgebra of $L(E_7)$ by Proposition \ref{prop:sl2ifsplitoff}.

The remaining case is $u$ coming from the regular class, where as with the $E_7(a_1)$ case the $10$ from the action of $u$ must yield a simple summand, with the rest projective. There are again two conspicuous sets of composition factors for $M(E_7)\downarrow_H$, coming from $P(4)\oplus 18$ and $P(10)\oplus 18$. The first has no corresponding set of composition factors on $L(E_7)$, and the second has a single set, which since $u$ is projective on $L(E_7)$, must be arranged so that $L(E_7)\downarrow_H$ is
\[ 19\oplus P(15)\oplus P(11)\oplus P(3).\]
While the $3$-dimensional submodule is a subalgebra of $L(E_7)$, it is not obviously an $\slf_2$-subalgebra because we cannot apply Proposition \ref{prop:sl2ifsplitoff}. This is a Serre embedding as defined in Definition \ref{defn:serreembedding}, as needed.
\end{proof}

The last case is $p=23$ and the regular unipotent class, to conclude this section, chapter and article.

\begin{proposition}\label{prop:e7char23sl}
Suppose that $p=23$ and $a=1$. If $H$ is not a blueprint for $M(E_7)$ then $H$ stabilizes a unique $3$-space on $L(E_7)$ that is an $\slf_2$-subalgebra of $L(E_7)$.
\end{proposition}
\begin{proof} The only non-generic unipotent class for $p=23$ and $M(E_7)$ is the regular class, with Jordan blocks $23^2,10$, case (\ref{li:sl2first23}) above. The $10$ in the action of $u$ must come from a simple summand, leaving a single projective module of dimension $46$. Only two possibilities yield conspicuous sets of composition factors, namely $20,18,10,4^2$ and $18^2,10,6,4$. The first of these has no corresponding set of composition factors on $L(E_7)$, and the second of these yields the unique action
\[ P(19)\oplus P(11)\oplus 23\oplus 15 \oplus 3,\]
with the $3$ being an $\slf_2$-subalgebra of $L(E_7)$ by Proposition \ref{prop:sl2ifsplitoff}.
\end{proof}

We now conclude the proof of Theorem \ref{thm:e7}. We start by showing that if $p$ is an odd prime with $p^a\neq 7,25$ then $H=\SL_2(p^a)$ with $Z(\bG)=Z(H)$ then $H$ is strongly imprimitive.

If $p^a>150$ then $H$ is a blueprint for $M(E_7)$, whence we are done by Proposition \ref{prop:blueprintissi}. For $p=3$, we have $a=2,3,4$. If $a=2$ then $H$ stabilizes a $2$-space on $M(E_7)$ by \cite[Proposition 6.2]{craven2015un2}, so we are done by Proposition \ref{prop:fix2spaceonMG}. If $a=3$ then $H$ is strongly imprimitive by Proposition \ref{prop:e7char3a=3sl}, and if $a=4$ then $H$ is strongly imprimitive by Proposition \ref{prop:e7char3a=4sl}. This needs Proposition \ref{prop:blueprintissi} if $H$ is a blueprint for $M(E_7)$, and Proposition \ref{prop:intersectionstabilizers} if $H$ stabilizes a unique $18$-space.

Thus $p\geq 5$. If $p=5$ then Proposition \ref{prop:e7char5sl} shows that $H$ stabilizes a $2$-space on $M(E_7)$ for $a=1$, is a blueprint for $M(E_7)$ if $a=3$, and is therefore strongly imprimitive by Propositions \ref{prop:fix2spaceonMG} and \ref{prop:blueprintissi}. If $a=2$ then $H$ is strongly imprimitive except for one specific action on $M(E_7)$ and $L(E_7)$.

If $p=7$ then by Proposition \ref{prop:e7char7sl}, if $a=2$ then $H$ is a blueprint for $M(E_7)$, or $H$ stabilizes a $2$-space on $M(E_7)$ or a line on $L(E_7)$. The first two possibilities yield strong imprimitivity as we have seen above, and stabilizing a line yields strong imprimitivity by Proposition \ref{prop:fixlineonLG}. For $a=1$, one of these three conditions hold -- and hence $H$ is strongly imprimitive -- or $M(E_7)\downarrow_H$ and $L(E_7)\downarrow_H$ have a specific action.

For $p=11$, Proposition \ref{prop:e7char11sl} shows that again, $H$ is either a blueprint for $M(E_7)$, stabilizes a $2$-space on $M(E_7)$ or a line on $L(E_7)$, and hence is strongly imprimitive. For $p=13,17$, Propositions \ref{prop:e7char13sl} and \ref{prop:e7char17sl} states directly that $H$ is strongly imprimitive.

For $p=19$, if $H$ is a blueprint for $M(E_7)$ then $H$ is strongly imprimitive as we have seen above, and if $H$ stabilizes a unique $3$-space that is an $\slf_2$-subalgebra of $L(E_7)$ then $H$ is strongly imprimitive by Corollary \ref{cor:sl2coxeter-1}. The other alternative is that $H$ is a Serre embedding, as stated in Theorem \ref{thm:e7}. Finally, if $p=23$ then $H$ is either a blueprint for $M(E_7)$ or stabilizes a unique $3$-space that is an $\slf_2$-subalgebra of $L(E_7)$, so again $H$ is strongly imprimitive.

\medskip

We now complete the proof of Theorem \ref{thm:e7}. Suppose that $H\cong \PSL_2(p^a)$ is a subgroup of the simple group $G=E_7(p^b)$ such that $N_{\bar G}(H)$ is maximal in an almost simple group $\bar G$ with socle $G$. The corresponding subgroup ($\PSL_2(2^a)$ if $p=2$, and $\PSL_2(p^a)\times 2$ or $\SL_2(p^a)$ for $p$ odd) of the (simply connected) algebraic group $\bG$ is either strongly imprimitive or not strongly imprimitive. If it is not strongly imprimitive then $p^a$ is one of $7$, $8$ and $25$, as we have seen in the summaries of this and the previous two chapters. Thus we may assume that $H$ is strongly imprimitive. Since it is maximal, this means that $H$ is the fixed points $\bH^\sigma$ of a Frobenius endomorphism of a positive-dimensional subgroup $\bH$ of $\bG$. By \cite[Corollary 2]{liebeckseitz2004}, $\bH$ is either maximal rank, a maximal parabolic (which is obviously impossible), $(2^2\times D_4)\cdot \Sym(3)$ (again, obviously impossible) or appears in \cite[Table 1]{liebeckseitz2004}.

If $\bH$ is maximal-rank then we examine \cite[Table 5.1]{liebecksaxlseitz1992}, and the result holds. If $\bH$ appears in \cite[Table 1]{liebeckseitz2004} then $\bH$ must be a product of type $A_1$ subgroups. There are two subgroups $A_1$, which appear in Theorem \ref{thm:e7}, and one subgroup $A_1A_1$. However, as we see from \cite[Table 10.2]{liebeckseitz2004}, the two $A_1$ factors are not interchangeable (they have different actions on $M(E_7)$) and so a Frobenius endomorphism cannot have fixed points $\PSL_2(p^a)$ on this subgroup, only a product of $\PSL_2$ subgroups.

This completes the proof of Theorem \ref{thm:e7}, and therefore concludes the whole proof.

\appendix
%    Include appendix "chapters" here.
\chapter[Actions of Maximal Positive-Dimensional Subgroups]{Actions of Maximal Positive-Dimensional Subgroups on Minimal and Adjoint Modules}
\label{app:actions}

In this appendix we collate information on the actions of the reductive and parabolic maximal subgroups of positive dimension on the minimal and adjoint modules for the algebraic groups $F_4$, $E_6$ and $E_7$ that we have used in the text, other than those in Lemmas \ref{lem:e6stabs} and \ref{lem:e7stabs}. These have been documented in many places, but we give them here as well for ease of reference.

We need information for $F_4$ and $E_6$ in characteristic $3$, and for $E_6$ in characteristics $7$ and $11$. We list the composition factors of every maximal closed, connected subgroup of positive dimension in these characteristics, taken from \cite{liebeckseitz2004}, on $M(\bG)$, and $L(\bG)$ for $\bG=F_4$. We list the reductive subgroups first, and then the parabolics. Write $M^\pm$ to mean both a module $M$ and its dual $M^*$ are composition factors.

We begin with the table for $F_4$ in characteristic $3$.

\begin{center}
\begin{tabular}{ccc}
\hline Subgroup & Factors on $M(F_4)$ & Factors on $L(F_4)$
\\\hline $B_4$ & $1000,0001$ & $0100,0001$
\\ $\tilde A_1C_3$ & $(1,100),(0,010)$ & $(2,000),(0,200),(1,001)$
\\ $A_2\tilde A_2$ & $(10,10),(01,01),(00,11)$ & $(11,00),(00,11),(10,02),(01,20),(00,00)^2$
\\ $A_1G_2$ & $(2,10),(4,00)$ & $(2,00),(0,01),(0,10)^2,(4,10)$
\\ \hline $B_3$ & $100,001^2,000^2$ & $100^2,010,001^2,000$
\\ $C_3$ & $100^2,010$ & $200,001^2,000^3$
\\\multirow{2}{*}{$A_2\tilde A_1$} & $(10,1)^\pm,(10,0)^\pm,$
 & $(11,0),(10,2)^\pm,(10,1)^\pm,(10,0)^\pm,$
 \\ & $(00,2),(00,1)^2$ & $(00,2),(00,1)^2,(00,0)^2$
\\\multirow{2}{*}{$\tilde A_2A_1$} & \multirow{2}{*}{$(10,1)^\pm,(10,0)^\pm,(11,0)$} & $(20,1)^\pm,(20,0)^\pm,(11,0),$
\\ &&$(00,2),(00,1)^2,(00,0)^2$
\\ \hline
\end{tabular}
\end{center}

Next, the subgroups of $E_6$ in characteristic $3$.

\begin{center}
\begin{tabular}{cc}
\hline Subgroup & Factors on $M(E_6)$
\\ \hline $A_5A_1$ & $(\lambda_4,0),(\lambda_1,1)$
\\ $A_2A_2A_2$ & $(10,01,00),(00,10,01),(01,00,10)$
\\ $F_4$ & $0001,0000^2$
\\ $C_4$ & $0100$
\\ $G_2A_2$ & $(10,10),(00,02)$
\\ $G_2$ (2 classes) & $20$
\\\hline  $D_5$ & $\lambda_1,\lambda_4,0$
\\ $A_5$ & $\lambda_1^2,\lambda_4$
\\ $A_4A_1$ & $(1000,1),(0001,0),(0010,0),(0000,1)$
\\ $A_2A_2A_1$ & $(10,01,0),(01,00,1),(00,10,1),(01,00,0),(00,10,0)$
\\ \hline
\end{tabular}
\end{center}

Finally, the subgroups of $E_7$ in characteristics $7$ and $11$.

\begin{center}
\begin{tabular}{cc}
\hline Subgroup & Factors on $M(E_7)$
\\ \hline $D_6A_1$ & $(\lambda_1,1),(\lambda_5,0)$
\\$A_7$ & $\lambda_2^\pm$
\\$A_5A_2$ & $(\lambda_1,10)^\pm,(\lambda_3,00)$
\\$C_3G_2$ & $(001,00),(100,10)$
\\$G_2A_1$ & $(01,1),(10,3)$
\\$F_4A_1$ & $(0001,1),(3,0000)$
\\$A_2$ & $60^\pm$
\\$A_1A_1$ & $(6,3),(4,1),(2,5)$
\\ \hline $E_6$ & $\lambda_1^\pm,0^2$
\\ $D_6$ & $\lambda_1^2,\lambda_5$
\\$A_6$ & $\lambda_1^\pm,\lambda_2^\pm$
\\$A_5A_1$ & $(\lambda_1,1)^\pm,(\lambda_1,0)^\pm,(\lambda_3,0)$
\\ $A_4A_2$ & $(10,1000)^\pm,(10,0000)^\pm,(00,0100)^\pm$
\\ $A_3A_2A_1$ & $(000,10,1)^\pm,(010,00,1),(100,10,0)^\pm,(100,00,0)^\pm$
\\ \hline
\end{tabular}
\end{center}

\chapter{Traces of Small-Order Semisimple Elements}
\label{app:traces}
We use the traces of semisimple elements on $M(\bG)$ and $L(\bG)$, of fairly large order compared with many similar papers in the literature. In this chapter we give a few tables for the real elements of orders at most $5$ for $\bG=F_4,E_6,E_7$ and $M(\bG)$ and $L(\bG)$. For elements of order $4$, we write the trace of the element followed by that of its square. Write $\omega$ for the sum of a $5$th root of unity and its inverse. We only list traces of elements of order $5$ up to algebraic conjugacy.

For $E_6$, by Proposition \ref{prop:f4real}, all real semisimple elements lie in $F_4$, so we list the class in $F_4$, its trace on $M(F_4)$, $L(F_4)$ and $L(E_6)$ in the first table. The trace on $M(E_6)$ is that on $M(F_4)$ plus $1$. For $p=3$, subtract $1$ from the trace on $M(F_4)$ and $L(E_6)$. For $E_7$, we simply list the class on $M(E_7)$ and $L(E_7)$ in the second table.

\begin{center}
\begin{tabular}{cccc}
\hline Order & Trace on $M(F_4)$ & Trace on $L(F_4)$ & Trace on $L(E_6)$
\\ \hline $2$ & $2$ & $-4$ & $-2$
\\ & $-6$ & $20$ & $14$
\\ \hline $3$ & $8$ & $7$ & $15$
\\ & $-1$ & $7$ & $6$
\\ & $-1$ & $-2$ & $-3$
\\ \hline $4$ & $14,2$ & $20,-4$ & $34,-2$
\\ & $6,-6$ & $8,20$ & $14,14$
\\ & $2,2$ & $0,-4$ & $2,-2$
\\ & $-2,2$ & $4,-4$ & $2,-2$
\\ & $-2,-6$ & $0,20$ & $-2,14$
\\ \hline $5$ & $1$ & $2$ & $3$
\\ & $\omega-1$ & $-2\omega+1$ & $-\omega$
\\ & $7\omega+7$ & $\omega+15$ & $8\omega+22$
\\ & $6\omega+14$ & $13\omega+21$ & $19\omega+35$
\\ & $3\omega+5$ & $4\omega+4$ & $7\omega+9$
\\ \hline
\end{tabular}
\end{center}

\begin{center}
\begin{tabular}{ccc}
\hline Order & Trace on $M(E_7)$ & Trace on $L(E_7)$
\\ \hline $2$ & $8$ & $5$
\\ & $-8$ & $5$
\\ & $-56$ & $133$
\\\hline $3$ & $20$ & $34$
\\ & $2$ & $7$
\\ & $2$ & $-2$
\\ & $-7$ & $7$
\\ & $-25$ & $52$
\\ \hline $4$ & $32,8$ & $65,5$
\\ & $16,-8$ & $29,5$
\\ & $8,8$ & $9,5$
\\ & $0,8$ & $1,5$
\\ & $0,-8$ & $5,5$
\\ & $0,-8$ & $-3,5$
\\ & $0,-56$ & $25,133$
\\ & $0,-56$ & $-7,133$
\\ & $-8,8$ & $9,5$
\\ & $-16,-8$ & $29,5$
\\ & $-32,8$ & $65,5$
\\ \hline $5$& $6$ & $8$
\\ & $26\omega-1$ & $-27\omega+52$
\\ & $14\omega+18$ & $22\omega+39$
\\ & $14\omega-7$ & $-28\omega+14$
\\ & $12\omega+32$ & $31\omega+66$
\\ & $8\omega+5$ & $4\omega+15$
\\ & $7\omega-8$ & $-9\omega+21$
\\ & $6\omega+14$ & $13\omega+22$
\\ & $5\omega-4$ & $-10\omega+3$
\\ & $2\omega+2$ & $\omega+1$
\\ & $\omega-1$ & $-2\omega+2$
\\ \hline
\end{tabular}
\end{center}

\backmatter
\newcommand{\etalchar}[1]{$^{#1}$}
\providecommand{\bysame}{\leavevmode\hbox to3em{\hrulefill}\thinspace}
\providecommand{\MR}{\relax\ifhmode\unskip\space\fi MR }
\providecommand{\MRhref}[2]{%
  \href{http://www.ams.org/mathscinet-getitem?mr=#1}{#2}
}
\providecommand{\href}[2]{#2}

\end{document}